\documentclass[a4paper,12pt]{article}

\usepackage[top=3.5cm, bottom=3.5cm, left=2.5cm, right=2.5cm]{geometry} 
\usepackage{amsmath}
\usepackage{amsthm}
\usepackage{amssymb}
  \usepackage{paralist}
  \usepackage{graphics} 
  \usepackage{epsfig} 
\usepackage{graphicx}  \usepackage{epstopdf}
\usepackage[colorlinks=false]{hyperref}
\usepackage[nottoc, notlof, notlot]{tocbibind}

\usepackage{float}
\usepackage{tikz}
\usetikzlibrary{arrows}
\usetikzlibrary[patterns]

\usepackage[boxed]{algorithm2e}
\usepackage{setspace}
\usepackage{etoolbox}
\AtBeginEnvironment{algorithm}{\setstretch{1.2}}

\usepackage{array}
\usepackage{tabularx}
\usepackage{ltablex} 
\usepackage{caption}

\usepackage{authblk}
\usepackage{mathrsfs}
\usepackage{thmtools}
\usepackage{bbm}

\usepackage{stmaryrd}

\makeatletter
\newsavebox{\@brx}
\newcommand{\llangle}[1][]{\savebox{\@brx}{\(\m@th{#1\langle}\)}%
  \mathopen{\copy\@brx\kern-0.5\wd\@brx\usebox{\@brx}}}
\newcommand{\rrangle}[1][]{\savebox{\@brx}{\(\m@th{#1\rangle}\)}%
  \mathclose{\copy\@brx\kern-0.5\wd\@brx\usebox{\@brx}}}
\makeatother

\usepackage{cleveref}
\crefname{I}{}{}
\creflabelformat{I}{#2I-#1#3}
\crefname{Ieq}{}{}
\creflabelformat{Ieq}{#2(I-#1)#3}
\crefname{II}{}{}
\creflabelformat{II}{#2#1#3}
\crefname{IIeq}{}{}
\creflabelformat{IIeq}{#2(#1)#3}

\usepackage{comment}

\usepackage{nameref,zref-xr}
\zxrsetup{toltxlabel}
\zexternaldocument*{part1}

\newcommand{\R}{\mathbb{R}}
\newcommand{\N}{\mathbb{N}}
\newcommand{\1}{\mathbbm{1}}
\newcommand{\E}{\mathbb{E}}
\newcommand{\PP}{\mathbb{P}}
\newcommand{\e}{\mathrm{e}}
\newcommand{\dd}{\mathrm{d}}
\newcommand{\pb}{\mathcal{P}}
\DeclareMathOperator{\Tr}{Tr}
\DeclareMathOperator{\Dom}{Dom}
\DeclareMathOperator{\argmin}{argmin}
\DeclareMathOperator{\diag}{diag}

\newtheorem{theorem}{Theorem}[section]
\newtheorem{corollary}{Corollary}

\newtheorem{lemma}[theorem]{Lemma}
\newtheorem{proposition}{Proposition}

\newtheorem{assumption}{Assumption}
\theoremstyle{definition}
\newtheorem{definition}[theorem]{Definition}
\newtheorem{remark}{Remark}
\newtheorem{example}{Example}

\title{\bf \LARGE Propagation of chaos: a review of models, methods and applications. \\ II. Applications}

\author[1]{Louis-Pierre \textsc{Chaintron}}
\author[2]{Antoine \textsc{Diez}}
\affil[1]{\small
DMA, \'Ecole Normale Sup\'erieure

45 rue d'Ulm

75005 Paris, France

\url{lchaintron@clipper.ens.fr}
\bigskip
}
\affil[2]{\small
Department of Mathematics, 

Imperial College London, South Kensington Campus,

London, SW7 2AZ, UK

\url{antoine.diez18@imperial.ac.uk}
\bigskip
}

\begin{document}
\maketitle

\begin{abstract}
The notion of propagation of chaos for large systems of interacting particles originates in statistical physics and has recently become a central notion in many areas of applied mathematics. The present review describes old and new methods as well as several important results in the field. The models considered include the McKean-Vlasov diffusion, the mean-field jump models and the Boltzmann models. The first part of this review is an introduction to modelling aspects of stochastic particle systems and to the notion of propagation of chaos. The second part presents concrete applications and a more detailed study of some of the important models in the field.
\end{abstract}



\newpage

\tableofcontents

\section{Introduction}

The second part of this review is devoted to many classical and recent modelling problems which are based on the simulation of large systems of interacting particles. This approach was initiated in the second half of the XIXth century by Boltzmann who proposed to model a gas as a myriad of elementary particles undergoing a simple Newtonian dynamics. When the number of particles grows to infinity, Boltzmann's \emph{kinetic theory of gases} is able to explain complex thermodynamics phenomena which previously had only a phenomenological interpretation. 

Beyond the contributions in Physics, the ideas of Boltzmann have had a profound influence on the development of mathematical concepts both in Probability and in Analysis. In the middle of the XXth century, Kac and later McKean introduced a proper mathematical formalisation of the concepts introduced by Boltzmann. The notions of \emph{Kac's chaos} and \emph{propagation of chaos} together with the probabilistic models of Kac and McKean are the foundations of the \emph{mathematical kinetic theory}. The derivation of the Boltzmann equation for rarefied gas dynamics as well as the other classical equations of statistical physics have long motivated the development of the theory. Since the last two decades, there is an ever growing number of applications of these ideas in wider range of domains, from the modelling of large animal societies, to socio-economic models or numerical methods in data sciences. 

The first part of this review introduced the tools, concepts and some of the main abstract models for the study of many-particle systems. Throughout this second part, references to the first part are indicated by ``I-'' (for instance Section \cref{sec:models} refers to the second section of the first part). In the second part of this review, the analysis is specialized on the one hand to the classical models introduced by Kac and McKean and their recent developments and on the other hand to a gallery of recent applications in applied mathematics and beyond. 

\subsubsection*{Outline} The outline of the article is as follows. 

Section \cref{sec:summary} summarises the content of the first part of this review. 

Section \cref{sec:mckeanreview} and Section \cref{sec:boltzmannreview} are devoted to the review of the main results in the literature respectively for McKean-Vlasov models and Boltzmann-Kac models. We emphasize that although none of the results presented are new, we include some proofs that we did not find or hardly found in the literature in this form, in particular: the proofs of McKean's and Kac's theorems (Section \cref{sec:mckeancoupling} and Section~\cref{sec:kactheorem}), the functional law of large numbers by martingale arguments (Section \cref{sec:martingalecompactness}) and the proof of propagation of chaos for Boltzmann models via coupling methods (Section~\cref{sec:couplingBoltzmann}). 

Section \cref{sec:applications} is an introductory section to various recent modelling problems and practical applications of the concept of propagation of chaos. A selection of examples which motivate and often extend the results of the previous sections is presented, including some open problems and current research trends.   

Several appendices complete this work. Generalised high-order expansions of the particle generators against monomial test functions are shown in Appendix \cref{appendix:generatorestimates} and a technical lemma in Appendix \cref{appendix:combinatoriallemma}. Finally, for the reader's convenience, we collect in Appendix \cref{appendix:tightness} useful tightness criteria.

\subsection*{Notations and conventions}

\subsubsection*{Sets}

\noindent\begin{tabularx}{\linewidth}{@{}lX}
    $C( I , E)$ & The set of continuous functions from a time interval $I=[0,T]$ to a set $E$, endowed with the uniform topology.   \\
    $C_b(E)$, $C_b^k(E)$ & Respectively the set of real-valued bounded continuous functions and the set of functions with $k\geq1$ bounded continuous derivatives on a set $E$. \\
    $C_c(E)$ & The set of real-valued continuous functions with compact support on a locally compact space $E$. \\
    $C_0(E)$ & The set of real-valued continuous functions vanishing at infinity on a locally compact space $E$, i.e. $\varphi\in C_0(E)$ when for all $\varepsilon>0$, there exists a compact set $K_\varepsilon\subset E$ such that $|\varphi(x)|<\varepsilon$ for all $x\in E$ outside $K_\varepsilon$.\\ 
    $D( I , E)$ & The space of functions which are right continuous and have left limit everywhere from a time interval $I=[0,T]$ to a set $E$, endowed with the Skorokhod $J1$ topology. This is the space of \emph{c\`adl\`ag} functions. This space is also called the \emph{Skorokhod space} or the \emph{path space}. \\
    $L^p(E)$ or $L^p_\mu(E)$ & The set of measurable functions $\varphi$ defined almost everywhere on a measured space $(E,\mu)$ such that the $|\varphi|^p$ is integrable for $p\geq1$. When $p=+\infty$, this is the set of functions with a bounded essential supremum. We do not specify the dependency in $\mu$ when no confusion is possible. \\
    $\mathcal{M}_d(\R)$ & The set of $d$-dimensional square real matrices. \\
    $\mathcal{M}(E)$ & The set of signed measures on a measurable space $E$. \\
    $\mathcal{M}^+(E)$ & The set of positive measures on a measurable space $E$. \\
    $\pb(E)$ & The set of probability measures on a space $E$.\\
    $\pb_p(E)$ & The set of probability measures with bounded moment of order $p\geq1$ on a space $E$.\\
    $\widehat{\pb}_N(E)$ & The set of empirical measures of size $N$ over  a set $E$, that is measures of the form $\mu = \frac{1}{N}\sum_{i=1}^N \delta_{x^i}$, where $x^i\in E$. \\
    $\R_+$ & The set $[0,+\infty)$.\\
    $\mathfrak{S}_N$ & The permutation group of the set $\{1,\ldots,N\}$.\\
    $\mathbb{S}^{d-1}$ & The sphere of dimension $d-1$. 
\end{tabularx}

\subsubsection*{Generic elements and operations}

\noindent\begin{tabularx}{\linewidth}{@{}lX}

$C$ & A generic nonnegative constant, the value of which may change from line to line. \\
$C(a_1,\ldots a_n)$ & A generic nonnegative constant which depends on some fixed parameters denoted by $a_1,\ldots,a_n$. Its value may change from line to line.\\
$\diag(x)$ & The $d$-dimensional diagonal matrix whose diagonal coefficients $x_1,\ldots,x_d$ are the components of the $d$-dimensional vector $x$. \\ 
$\nabla\cdot V$ & The divergence of a vector field $V:\R^d\to\R^d$ or of a matrix field $V:\R^d\to\mathcal{M}_d(\R)$, respectively defined by $\nabla\cdot V = \sum_{i=1}^d \partial_{x_i} V_i$ or componentwise by $(\nabla\cdot V)_i = \sum_{j=1}^d \partial_{x_j}V_{ij}$. \\
$A:B$ and $\|A\|$ & The Frobenius inner product of two matrices $A,B\in\mathcal{M}_d(\R)$ defined by $A:B:=\sum_{i=1}^d\sum_{j=1}^d A_{ij}B_{ij}$ and the associated norm $\|A\|:=\sqrt{A:A}$. \\
$\nabla^2 V$ & The Hessian matrix of a scalar field $V:\R^d\to \R$ defined componentwise by $(\nabla^2 V)_{ij} = \partial^2_{x_i, x_j} V$. \\
$I_d$ & The $d$-dimensional identity matrix.\\
$\mathrm{Id}$ & The identity operator on a vector space.\\
$\langle x, y\rangle$ or $x\cdot y$ & The Euclidean inner product of two vectors $x,y\in\R^d$ defined by $\langle x,y\rangle \equiv x\cdot y := \sum_{i=1}^d x^i y^i$. One notation or the other may be preferred for typographical reasons in certain cases. \\ 
$M_{ij}$ & The $(i,j)$ (respectively row and column indexes) component of a matrix $M$. \\ 
$\mathsf{P}(u)$ & The projection matrix $\mathsf{P}(u):=I_d-\frac{u\otimes u}{|u|^2}$ on the plane orthogonal to a vector $u\in\R^d$.\\
$\varphi\in C_b(E)$ & A generic bounded continuous test function on $E$.\\
$\varphi_N\in C_b(E^N)$ & A generic bounded continuous test function on the product space $E^N$. \\
$\Phi\in C_b(\pb(E))$ & A generic bounded continuous test function on the set of probability measures on $E$. \\
$u\otimes v$, $\mu\otimes\nu$ or $\varphi\otimes\psi$ & Respectively, the matrix tensor product of two vectors $u,v\in\R^d$ defined componentwise by $(u\otimes v)_{ij} = u_iv_j$; the product measure on $E\times F$ of two measures $\mu,\nu$ respectively on $E$ and $F$; the product function on $E\times F$ defined by $(\varphi\otimes\psi)(x,y)=\varphi(x)\psi(y)$ for two real-valued function $\varphi,\psi$ respectively on $E$ and $F$. \\
$\Tr M$ & The trace of the matrix $M$. \\
$M^\mathrm{T}$ & The transpose of the matrix $M$.\\
$\mathbf{x}^N=(x^1,\ldots,x^N)$ & A generic element of a product space $E^N$. The components are indexed with a superscript. \\
$\mathbf{x}^{M,N}=(x^1,\ldots,x^M)$ & The $M$-dimensional vector in $E^M$ constructed by taking the $M$ first components of $\mathbf{x}^N$. \\
$x=(x_1,\ldots,x_d)^\mathrm{T}$ and $|x|$ & A generic element of a $d$-dimensional space and its norm. The coordinates are indexed with a subscript. The norm of $x$ denoted by $|x|$ is the Euclidean norm. 
\end{tabularx}

\subsubsection*{Probability and measures}

\noindent\begin{tabularx}{\linewidth}{@{}lX}
$K\star\mu$ & The convolution of a function $K:E\times F\to G$ with a measure $\mu$ on $F$ defined as the function $K\star\mu:x\in E \mapsto \int_F K(x,y)\mu(\dd y)\in G$. When $E=F=G=\R^d$ and $K:\R^d\to\R^d$, we write $K\star\mu(x)=\int_{\R^d}K(x-y)\mu(\dd y)$. \\
$\delta_x$ & The Dirac measure at the point $x$. \\
$\mu_{\mathbf{x}^N}$ & The empirical measure defined by $\mu_{\mathbf{x}^N}=\frac{1}{N}\sum_{i=1}^N\delta_{x^i}$ where $\mathbf{x}^N=(x^1,\ldots,x^N)$. \\
$\E_\mu[\varphi]$ & Alternative expression for $\langle \mu,\varphi\rangle$ when $\mu$ is a probability measure. When $\mu=\mathbb{P}$ on $(\Omega,\mathscr{F},(\mathscr{F}_t)_t,\mathbb{P})$, the expectation is simply denoted by $\mathbb{E}$. \\
$H(\nu|\mu)$ & The relative entropy (or Kullback-Leibler divergence) between two measures $\mu,\nu$, see Defintion \cref{def:entropyfisher_summary}.\\ 
$\langle \mu,\varphi\rangle$ & The integral of a measurable function $\varphi$ with respect to a measure $\mu$. \\
$\mathrm{Law}(X)$ & The law of a random variable $X$ as an element of $\pb(E)$ where $X$ takes its value in the space $E$.\\
$(\Omega,\mathscr{F},(\mathscr{F}_t)_t,\mathbb{P})$ & A filtered probability space. Unless otherwise stated, all the random variables are defined on this set. The expectation is denoted by $\E$. \\
$\sigma(X^1,X^2,\ldots)$ & The $\sigma$-algebra generated by the random variables $X^1, X^2,\ldots$. \\
$T_\#\mu$ & The pushforward of the measure $\mu$ on a set $E$ by the measurable map $T:E\to F$. This is a measure on the set $F$ defined by $T_\#\mu(\mathscr{A}) = \mu(T^{-1}(\mathscr{A}))$ for any measurable set $\mathscr{A}$ of $F$. \\
$\|\cdot\|_{\mathrm{TV}}$ & The Total Variation (TV) norm for measures. \\
$W_p$ & The Wasserstein-$p$ distance between probability measures (see Definition \cref{def:wasserstein}).\\
$X\sim\mu$ & It means that the law of the random variable $X$ is $\mu$. \\ 
$(\mathsf{X}_t)_t$ or $(\mathsf{Z}_t)_t$ & The canonical process on the path space $D(I,E)$ defined by $\mathsf{X}_t(\omega)=\omega(t)$. \\
$(\mathbf{X}^N_t)^{}_t$ or $(\mathbf{Z}^N_t)^{}_t$ & The canonical process on the product space $D(I,E)^N$ with components $\mathbf{X}^N_t=(\mathsf{X}^1_t,\ldots,\mathsf{X}^N_t)$. 
\end{tabularx}

\subsubsection*{Systems of particles and operators}

\noindent\begin{tabularx}{\linewidth}{@{}lX}
$E$ & The state space of the particles, assumed to be at least a Polish space. \\
$f^N_t$ & The $N$-particle distribution in $\pb(E^N)$ at time $t\geq0$.\\
$f^{k,N}_t$ & The $k$-th marginal of $f^N_t$.\\
$f^N_I$ & The $N$-particle distribution on the path space in $\pb(D(I,E^N))$ or $\pb(C(I,E^N))$ for a time interval $I=[0,T]$. We identify $D(I,E^N)\simeq D(I,E)^N$. \\
$f_t$ & The limit law in $\pb(E)$ at time $t\geq0$.\\
$f_I$ & The limit law on the path space in $\pb(D(I,E))$ or $\pb(C(I,E))$. \\
$F^N_t$ & The law of the empirical process in $\pb(\pb(E))$ at time $t\geq0$.\\
$F^{\mu,N}_I$ & The weak pathwise law of the empirical process in $\pb(D(I,\pb(E)))$ on the time interval $I=[0,T]$. \\
$F^N_I$ & The strong pathwise law of the empirical process in $\pb(\pb(D(I,E)))$ on the time interval $I=[0,T]$. \\
$\mathcal{L}_N$ & The $N$-particle generator acting on (a subset of) $C_b(E^N)$. \\
$\mathcal{L}^N$ & The $N$-particle generator acting on $\pb(E^N)$ defined as the formal adjoint of $\mathcal{L}_N$. \\
$L\diamond_i\varphi_N$ & The action of an operator $L$ on (a subset of) $C_b ( E )$ against the $i$-th variable of a function $\varphi_N$ in $C_b ( E^N )$, defined as the function in (a subset of) $C_b ( E^N )$
$L \diamond_i \varphi_N : (x^1,\ldots,x^N) \mapsto L [ x \mapsto \varphi_N ( x^1 , \ldots , x^{i-1} , x , x^{i+1} , \ldots , x^N ) ] ( x^i )$. The definition readily extends to the case of an operator $L^{(2)}$ acting on $C_b(E^2)$ and two indexes $i<j$ in which case we write $L^{(2)}\diamond_{ij}\varphi_N$. \\
$(\mathcal{X}^N_t)^{}_t$ & The $N$-particle process, with components
$\mathcal{X}^N_t=(X^{1,N}_t,\ldots,X^{N,N}_t)\in E^N$. Often we write $X^{i,N}_t\equiv X^i_t$ and $(\mathcal{X}^N_t)^{}_t\equiv\mathcal{X}^N_{[0,T]}$. \\
$(\mathcal{Z}^N_t)^{}_t$ & An alternative notation for the $N$-particle process with $\mathcal{Z}^N_t=(Z^{1,N}_t,\ldots,Z^{N,N}_t)$. Often used for Boltzmann particle systems or kinetic systems. 
\end{tabularx}

\section{Summary of the first part}\label[II]{sec:summary}

\subsection{Particle systems, chaos and propagation of chaos}

The object of the present review is the study of large systems of interacting particles. Throughout this work, a particle system is defined as a Markov process $(\mathcal{X}^N_t)^{}_{t\in I}$ with values in $E^N$ where $E$ is a Polish space, $N$ is the number of particles and $I=[0,T]$, $T\in(0,+\infty]$ is a time interval. Throughout this review, we use the notation $\mathcal{X}^N_t =~(X^{1,N}_t,\ldots,X^{N,N}_t)$ for the particle system and we write $X^i_t\equiv X^{i,N}_t$ without the $N$ superscript for the $i$-th particle when no confusion is possible. 

From the theory of Markov processes (see Appendix~\cref{appendix:timeinhomoegeneousmarkov}), the probability distribution of the particle system at time $t$ denoted by $f^N_t\in \pb(E^N)$ satisfies the (weak) \emph{Liouville equation} 
\begin{equation}\label[IIeq]{eq:liouville_summary}\forall \varphi_N\in \Dom(\mathcal{L}_N), \quad \frac{\dd}{\dd t} \langle f^N_t, \varphi_N\rangle = \langle f^N_t, \mathcal{L}_N\varphi_N\rangle,\end{equation}
where $\mathcal{L}_N$ is the infinitesimal generator of the particle system acting on a (dense) subset of test functions $\Dom(\mathcal{L}_N)\subset C_b(E^N)$. In stochastic analysis, the (richer) pathwise law $f^N_{[0,T]} \in \pb(D([0,T],E^N))$ is sometimes preferred and is characterised as the solution of a martingale problem. It means that $f^N_{[0,T]}$ is the unique probability distribution on the Skorokhod space of c\`adl\`ag functions such that for all test function $\varphi_N\in \Dom(\mathcal{L}_N)$, the process defined by: 
\[M_t^{\varphi_N} := \varphi_N(\mathbf{X}^N_t) - \varphi_N(\mathbf{X}^N_0) - \int_0^t \mathcal{L}_N\varphi_N(\mathbf{X}^N_s)\dd s,\]
is a $f^N_{[0,T]}$-martingale. In this definition, the process $(\mathbf{X}^N_t)^{}_{t\geq0}$ denotes the canonical process on $D([0,T],E^N)$ defined for any $\omega\in D([0,T],E^N)$ and any $t\geq0$ by $\mathbf{X}^N_t(\omega)=\omega(t)$. 

The particle system is assumed to be \emph{exchangeable} in the sense that $f^N_t$ (resp. its pathwise version $f^N_{[0,T]}$) is a symmetric probability distribution on $E^N$ (resp. on $D([0,T],E)^N\simeq D([0,T],E^N)$). 

This review is devoted to the notions of chaos and propagation of chaos introduced by Kac \cite{kac_foundations_1956} and defined below. 

\begin{definition}[Kac's chaos]
Let $f\in\pb ( E )$. A sequence $( f^N )_{N \geq 1}$ of symmetric probability measures on $E^N$ is said to be \emph{$f$-chaotic} when for any $k\in\N$ and any function $\varphi_k \in C_b ( E^k )$, 
\begin{equation*}
\lim_{N \to+ \infty} \langle f^N , \varphi_k \otimes 1^{\otimes N - k} \rangle = \langle f^{\otimes k} , \varphi_k \rangle.
\end{equation*}
It means that for all $k\in\N$, the $k$-th marginal satisfies $f^{k,N}\to f^{\otimes k}$ for the weak topology.
\end{definition}

From now on in this review, the initial distribution $f^N_0\in\pb(E^N)$ of the particle system is always assumed to be $f_0$-chaotic for a given $f_0\in \pb(E)$. The goal is to prove that this initial chaoticity assumption is propagated at later times as in the following definition. 

\begin{definition}[Pointwise and pathwise propagation of chaos]\label[II]{def:poc_summary}
Let $f^N_0 \in\pb(E^N)$ be the initial $f_0$-chaotic distribution of $\mathcal{X}^N_0$ at time $t=0$. 
\begin{itemize}
\item \emph{Pointwise propagation of chaos} holds towards a flow of measures $(f_t)_t\in C(I,\pb(E))$ when the law $f^N_t\in\pb(E^N)$ of $\mathcal{X}^N_t$ is $f_t$-chaotic for every time $t\in I$.
\item \emph{Pathwise propagation of chaos} holds towards a distribution $f_I\in\pb(D(I,E))$ on the path space when the law $f^N_I\in\pb\big(D(I,E)^N\big)$ of the process $\mathcal{X}^N_I$ (seen as a random element in $D(I,E)^N$) is $f_I$-chaotic.
\end{itemize}
\end{definition} 

The propagation of chaos property (pointwise or pathwise) describes the limit behavior of the particle system when the number of particles grows to infinity. It implies that any subsystem (of fixed size) of the $N$-particle system asymptotically behaves as a system of i.i.d processes with common law $f_t$ (note that the particles are always identically distributed by the exchangeability assumption). This translates the physical idea that for large systems, the correlations between two (or more) given particles which are due to the interactions become negligible. By looking at the whole system, only an averaged behavior can be observed instead of the detailed correlated trajectories of each particle. This notion of average behavior can be understood through the following characterization of the notion of chaos. The proof of this fundamental lemma can be found in the first part of the present review article, see Lemma \cref{lemma:caractchaos}, or in the classical course by Sznitman \cite[Proposition~2.2]{sznitman_topics_1991}. 

\begin{lemma} \label[II]{lemma:caractchaos_summary}
Each of the following assertions is equivalent to Kac's chaos.
\begin{enumerate}[(i)]
\item There exists $k \geq 2$ such that  $f^{k,N}$ converges weakly towards $f^{\otimes k}$.
\item The random empirical measure \[\mu_{\mathcal{X}^N} := \frac{1}{N}\sum_{i=1}^N \delta_{X^i},\] 
converges in law towards the deterministic measure $f$, where for any $N\in\N$, $\mathcal{X}^N = (X^1,\ldots,X^N)\sim f^N$. 
\end{enumerate}
\end{lemma}

The central question is therefore the description of the limit law $f_t$ which will be defined as the solution of a nonlinear PDE or of a nonlinear martingale problem. For all the models presented in this review, the goal is to compute the limit $N\to+\infty$ of any marginal $f^{k,N}_t$ of the $N$-particle distribution at any time $t>0$ or the limit of the random empirical measure $\mu_{\mathcal{X}^N_t}$. The models which are considered belong to one of the three families of models described in the next Section \cref{sec:models_summary}. 

\subsection{Models}\label[II]{sec:models_summary}

The first two families of models are called \emph{mean-field models}, they are defined by a generator of the form 
\begin{equation}\label[IIeq]{eq:Nparticlemeanfieldgenerator_summary}\mathcal{L}_N\varphi_N(\mathbf{x}^N) = \sum_{i=1}^N L_{\mu_{\mathbf{x}^N}}\diamond_i\varphi_N(\mathbf{x}^N),\end{equation}
where given a probability measure $\mu\in\pb(E)$, $L_\mu$ is the generator of a Markov process on $E$ which will be either a diffusion (Section \cref{sec:mckeanvlasov_summary}) or jump-process (Section~\cref{sec:meanfieldjump_summary}). Throughout this review, the notation $L\diamond_i\varphi_N$ denotes the function:
\[L \diamond_i \varphi_N : (x^1,\ldots,x^N) \in E^N \mapsto L [ x \mapsto \varphi_N ( x^1 , \ldots , x^{i-1} , x , x^{i+1} , \ldots , x^N ) ] ( x^i ) \in\R.\]
The third family of models are the Boltzmann models (Section \cref{sec:boltzmann_summary}). 

\subsubsection{McKean-Vlasov diffusion}\label[II]{sec:mckeanvlasov_summary} When the generator $L_\mu$ in \cref{eq:Nparticlemeanfieldgenerator_summary} is the generator of a diffusion process, the particle system is the solution the following system of SDEs
\begin{equation}\label[IIeq]{eq:mckeanvlasov_summary}
\forall i\in\{1,\ldots,N\}, \quad \dd X^{i,N}_t = b \big( X^{i,N}_t,\mu_{\mathcal{X}^N_t} \big) \dd t + \sigma \big( X^{i,N}_t,\mu_{\mathcal{X}_t} \big) \dd B^{i}_t,
\end{equation}
for $i\in\{1,\ldots,N\}$ where $(B^i_t)^{}_t$ are $N$ independent Brownian motions and the drift function $b$ and diffusion matrix are of the form 
\begin{equation*}
b : \R^d\times\pb(\R^d)\to \R^d,\quad \sigma : \R^d\times\pb(\R^d)\to\mathcal{M}_d(\R).\end{equation*}

\begin{remark}
Note that there are actually $dN$ independent one-dimensional Brownian motions. This remark may be helpful in cases where the Brownian motions in the different directions are different. In particular, for \emph{kinetic} particles defined by their positions and velocities, the noise is often added on the velocity variable only (this case is nevertheless covered by \cref{eq:mckeanvlasov_summary} with a block-diagonal matrix $\sigma$ with a vanishing block on the position variable). 
\end{remark}

The mean-field limit $N\to+\infty$ is given by the nonlinear Fokker-Planck equation 
\begin{equation}\label[IIeq]{eq:mckeanvlasov-pde_summary}
\partial_t f_t( x ) = -\nabla_x\cdot\{b( x ,f_t)f_t\}+\frac{1}{2}\sum_{i,j=1}^d \partial_{x_i}\partial_{x_j}\{a_{ij}(x,f_t)f_t \},
\end{equation}
where $a(x,\mu):=\sigma(x,\mu)\sigma(x,\mu)^\mathrm{T}$. This is the law of the \emph{nonlinear McKean-Vlasov process} $(\overline{X}_t)_t$ which solves the following nonlinear SDE: 
\begin{equation}\label[IIeq]{eq:mckeanvlasov-limit_summary}
\dd \overline{X}_t = b {\left( \overline{X}_t,f_t \right)}\dd t + \sigma {\left( \overline{X}_t,f_t \right)} \dd B_t.
\end{equation}
where $B_t$ is a Brownian motion and $f_t = \mathrm{Law}(\overline{X}_t)$. The well-posedness of \cref{eq:mckeanvlasov-limit_summary} is proved under Lipschitz assumptions on $b$ and $\sigma$ in Proposition \cref{prop:wellposednessmckeanlipschitz}. 

In many applications, the particles are rather defined by their positions and velocities, respectively denoted for the $i$-th particle at time $t$ by $(X^i_t,V^i_t)\in\R^d\times\R^d$. For instance, when $\sigma\equiv0$, a particle system ruled by the Newton equations 
\[\frac{\dd X^i_t}{\dd t} = V^i_t,\quad \frac{\dd V^i_t}{\dd t} = \sum_{j=1}^N F(X^i_t-X^j_t),\]
where $F$ is a force, can be written in the form \cref{eq:mckeanvlasov_summary} with the function $b:\R^d\times\R^d\times\mathcal{P}(\R^d\times\R^d)\to\R^d\times\R^d$ given by 
\[b((x,v),\mu) = \left(v,\int_{\R^d\times\R^d} F(x-x')\mu(\dd x',\dd v')\right).\]
In this case, the limit Equation \cref{eq:mckeanvlasov-pde_summary} is the renowned Vlasov equation which is historically one of the first and most important models in plasma physics and celestial mechanics. In the following, we will nevertheless most often consider stochastic models although some of the results still apply in this deterministic case (in particular the important Theorem \cref{thm:mckean}). For a detailed account of the Vlasov equation in this context, we refer to the review article \cite{jabin_review_2014}.  

Stochastic McKean-Vlasov systems have a wide range of applications. Some examples in physics and biology are described in the first part of this review, see Example \cref{example:diffusionprocess}. In this second part, we will treat important historical applications in physics: questions related to the granular media equation are discussed in Section~\cref{sec:gradientsystems} and particle systems with Coulomb-type interactions and other singular kernels, in particular in fluid dynamics, are described in Sections \cref{sec:jabin} and \cref{sec:vorticity}. More recently, McKean-Vlasov systems have also been used to model biological phenomena, in particular self-organized swarming phenomena (Sections \cref{sec:kuramoto} and~\cref{sec:flocking}). Very recently, these models have also gain attention in data sciences for the design and study of Particle Swarm Intelligence algorithms (Section \cref{sec:datasciences}).

\subsubsection{Mean-field jump process}\label[II]{sec:meanfieldjump_summary}

The $N$-particle process is defined by a generator of the form \cref{eq:Nparticlemeanfieldgenerator_summary} where given $\mu\in\pb(E)$, $L_\mu$ is the generator of a jump process of the form
\[L_\mu\varphi(x) = \lambda(x,\mu)\int_{E} \{\varphi(y)-\varphi(x)\}P_\mu(x,\dd y).\]
It describes a system of $N$ jump processes, driven by $N$ independent Poisson processes with jump rate 
\[\lambda:E\times\pb(E)\to\R_+,\,\, (x,\mu)\mapsto \lambda(x,\mu).\]
The law of jumps is prescribed by the jump measure:
\[P:E\times\pb(E)\to\pb(E),\,\, (x,\mu)\mapsto P_\mu(x,\dd y).\]

In classical kinetic theory, mean-field jump processes can be used to give a stochastic interpretation to the famous BGK equation \cite{bhatnagar_model_1954} (see Example \cref{example:bgk}). They have also recently become a basic tool for neuron models in biology (Example \cref{example:parametricjump} and Section \cref{sec:neurons}). 

\subsubsection{Boltzmann models}\label[II]{sec:boltzmann_summary}

The $N$-particle process is defined on an abstract Polish space $E$ by a generator of the form:
\begin{equation}\label[IIeq]{eq:boltzmanngenerator_summary}
\mathcal{L}_N\varphi_N = \sum_{i=1}^N L^{(1)}\diamond_i\varphi_N+\frac{1}{N}\sum_{i< j}L^{(2)}\diamond_{ij}\varphi_N,
\end{equation}
where $\varphi_N\equiv\varphi_N(z^1,\ldots,z^N)$ is a test function on the product space $E^N$. The operator $L^{(2)}$ acts on two-variable test functions and stands for binary interactions between particles. The operator $L^{(1)}$ acts on one-variable test functions and describes the individual flow of each particle (and possibly the boundary conditions). More explicitly, let us recall the notations, for $(z^1,\ldots,z^N)\in E^N$ and $i<j$, 
\begin{equation*}
    L^{(1)}\diamond_i \varphi_N(z^1,\ldots,z^n)  =  L^{(1)} \big[u \mapsto \varphi_N(z^1,\ldots,z^{i-1},u,z^{i+1},\ldots, z^N)\big](z^i)
\end{equation*}
and 
\begin{multline*}
    L^{(2)}\diamond_{ij} \varphi_N(z^1,\ldots,z^n) \\ =L^{(2)}\big[(u,v)\mapsto \varphi_N(z^1,\ldots,z^{i-1},u,z^{i+1},\ldots,z^{j-1},v,z^{j+1},\ldots, z^N)\big](z^i,z^j).
\end{multline*}
These models are called \emph{Boltzmann models} in reference to the famous Boltzmann equation of rarefied gas dynamics which is a fundamental equation for mathematicians, physicists and philosophers. It will be explained at the end of this section (see Equation \cref{eq:Boltzmannphysics_summary}) how it can be obtained as the limit of a general particle system of the form \cref{eq:boltzmanngenerator_summary}. The specificity of Boltzmann models is that the particles interact only at random times by pair and not individually with an average of all the other particles as in mean-field models. In full generality, the state space $E$ is an abstract space. In classical kinetic theory, $E=\R^d\times\R^d$ is the phase space of positions and velocities and two particles interact when they are close enough: they are said to \emph{collide} and by analogy, we will keep this terminology to refer to an interaction between two particles even in an abstract space. In addition to these pairwise interactions, each particle is also subject to an individual flow prescribed by the operator $L^{(1)}$. Typical examples in kinetic theory include
\begin{itemize}
    \item (Free transport) $L^{(1)}\varphi(x,v) = v\cdot\nabla_x\varphi$,
    \item (Space diffusion) $L^{(1)}\varphi(x,v) = \Delta_x\varphi$.
    \item (Velocity diffusion) $L^{(1)}\varphi(x,v) = \Delta_v\varphi$.
\end{itemize} 
When two particles collide, the effect of the collision is prescribed by the operator $L^{(2)}$. In kinetic theory, this operator acts on the velocity variable only but in full generality, in an abstract space $E$, it will be assumed to satisfy the following assumptions. 
\begin{assumption}\label[II]{assum:L2_summary} The operator $L^{(2)}$ satisfies the following properties. 
\begin{enumerate}[(1)]
    \item The domain of the operator $L^{(2)}$ is a subset of $C_b(E^2)$.
    \item There exist a continuous map called the \emph{post-collisional distribution}
    \[\Gamma^{(2)}:(z_1,z_2)\in E\times E \mapsto \Gamma^{(2)}(z_1,z_2,\dd z_1',\dd z_2')\in\pb(E\times E),\]
    and a symmetric function called the \emph{collision rate}
    \[\lambda : (z_1,z_2)\in E\times E\mapsto \lambda(z_1,z_2)\in \R_+,\]
    such that for all $\varphi_2\in C_b(E^2)$ and all $z_1,z_2\in E$,
    \begin{equation}\label[IIeq]{eq:generalL2_summary}L^{(2)}\varphi_2(z_1,z_2) = \lambda(z_1,z_2)\iint_{E\times E} \{\varphi_2(z_1',z_2')-\varphi_2(z_1,z_2)\}\Gamma^{(2)}(z_1,z_2,\dd z_1',\dd z_2').\end{equation}
    \item For all $z_1,z_2\in E$, the post-collisional distribution is symmetric in the sense that
    \begin{equation}\label[IIeq]{eq:symgamma2_summary}\Gamma^{(2)}(z_1,z_2,\dd z_1',\dd z_2') = \Gamma^{(2)}(z_2,z_1,\dd z_2',\dd z_1').\end{equation}
    It  ensures that the law $f^N_t$ defined by the backward Kolmogorov equation remains symmetric for all time provided that $f^N_0$ is symmetric.
    \item The function $\lambda$ is measurable on $\{(z_1,z_2)\in E^2,\,\,z_1\ne z_2\}$ and for all $z\in E$, $\lambda(z,z)=0$. 
\end{enumerate}
\end{assumption}

The assumption that $\lambda$ is a (measurable) function prevents from considering the true classical Boltzmann inhomogeneous case in kinetic theory~$\lambda( z_1 , z_2)=\delta_{x_1 = x_2}\Phi(|v_1-v_2|)$ for some nonnegative function $\Phi$ (that is, two particles collide when they are exactly at the same position), which is beyond the scope of this review (see however Section \cref{sec:lanford}). The collision rate is often assumed to be uniformly bounded 
\begin{equation}\label[IIeq]{eq:uniformboundlambda_summary}\sup_{z_1,z_2\in E}\lambda(z_1,z_2)\leq \Lambda <\infty.\end{equation}
This \emph{cutoff} assumption is unfortunately not physically relevant for many models where an infinite number of collisions may happen in finite time. However, it may serve as a first approximation which can be simulated on a computer as explained in Proposition \cref{prop:acceptreject} (see also Algorithm \cref{algo:exact}). 

The operator \cref{eq:boltzmanngenerator_summary} describes a particle system where each pair of particles interact at a rate given by the function $\lambda$ by updating the states of both particles according to the measure $\Gamma^{(2)}$. When propagation of chaos hold, the limit law $f_t$ is the solution of the \emph{general Boltzmann equation}:
\[\frac{\dd}{\dd t}\langle f_t,\varphi\rangle = \langle f_t, L^{(1)}\varphi\rangle+\langle f_t^{\otimes 2},L^{(2)}(\varphi\otimes 1)\rangle.\]
Using Assumption \cref{assum:L2_summary}, this equation can be rewritten 
\begin{multline}\label[IIeq]{eq:Boltzmannequationgeneral_summary}
\frac{\dd}{\dd t}\langle f_t,\varphi\rangle = \langle f_t, L^{(1)}\varphi\rangle\\
+\int_{E^3}\lambda(z_1,z_2)\big\{\varphi(z_1')-\varphi(z_1)\big\}\Gamma^{(2)}(z_1,z_2,\dd z_1',E)f_t(\dd z_1)f_t(\dd z_2),
\end{multline}
or in a more symmetric form, using \cref{eq:symgamma2_summary}:
\begin{multline}\label[IIeq]{eq:symmetricBoltzmannequationgeneral_summary}
\frac{\dd}{\dd t}\langle f_t,\varphi\rangle = \langle f_t, L^{(1)}\varphi\rangle\\
+\frac{1}{2}\int_{E^4}\lambda(z_1,z_2)\big\{\varphi(z_1')+\varphi(z_2')-\varphi(z_1)-\varphi(z_2)\big\}\Gamma^{(2)}(z_1,z_2,\dd z_1',\dd z_2')f_t(\dd z_1)f_t(\dd z_2).
\end{multline}
In many applications, the post-collisional distribution is explicitly given as the image measure of a known parameter space $(\Theta,\nu)$ endowed with a probability measure $\nu$ (or a positive measure with infinite mass). In this review, this particular class of models will be called \emph{parametric Boltzmann models}.  

\begin{definition}[Parametric and semi-parametric Boltzmann model]\label[II]{def:Boltzmannparammodel_summary} Let be given two measurable functions
\[\psi_1:E\times E\times\Theta\to E,\quad\psi_2:E\times E\times\Theta\to E,\]
which satisfy the symmetry assumption 
\begin{equation*}\forall (z_1,z_2)\in E^2,\quad (\psi_1,\psi_2)(z_1,z_2,\cdot)_{\#}\nu = (\psi_2,\psi_1)(z_2,z_1,\cdot)_{\#}\nu.\end{equation*}
Let the function $\psi$ be defined by
\[\psi:E\times E\times\Theta\to E^2,\,(z_1,z_2,\theta)\mapsto \big(\psi_1(z_1,z_2,\theta),\psi_2(z_1,z_2,\theta)\big).\]
A \emph{parametric Boltzmann model} with parameters $(\Theta,\psi)$ is a Boltzmann model of the form \cref{eq:boltzmanngenerator_summary} with Assumption \cref{assum:L2_summary} and a post-collisional distribution of the form: 
\[\forall (z_1,z_2)\in E^2,\quad \Gamma^{(2)}(z_1,z_2,\dd z_1',\dd z_2') = \psi(z_1,z_2,\cdot)_{\#} \nu.\]
The post-collisional distribution of a \emph{semi-parametric Boltzmann model} is of the form
\begin{equation}\label[IIeq]{eq:semiparametric_summary}\forall (z_1,z_2)\in E^2,\quad \Gamma^{(2)}(z_1,z_2,\dd z_1',\dd z_2') = \psi(z_1,z_2,\cdot)_{\#} \big(q(z_1,z_2,\theta)\nu(\dd\theta)\big),\end{equation}
where $q:E\times E\times \Theta\to\R_+$ is a fixed nonnegative function with $\int_{\Theta} q(z_1,z_2,\theta)\nu(\dd\theta) = 1$ for every $(z_1,z_2)$ in $E^2$. We will often assume that there exists $M>0$ and $q_0({\theta})$ a probability density function with respect to ${\nu}$ such that
\begin{equation}\label[IIeq]{eq:semiparambound_summary}\forall z_1,z_2\in E,\,\forall {\theta}\in{\Theta},\quad q(z_1,z_2,{\theta})\leq M q_0({\theta}).\end{equation}
\end{definition}
In the literature, the following variant of the generator \cref{eq:boltzmanngenerator_summary} is sometimes considered (see Example \cref{example:symmetrization} for more details): for $\mathbf{z}^N=(z^1,\ldots,z^N)\in E^N$ and $\varphi_N\in C_b(E^N)$, 
\begin{multline}\label[IIeq]{eq:Boltzmannnonsym_summary}\mathcal{L}_N\varphi_N(\mathbf{z}^N) = \sum_{i=1}^N L^{(1)}\diamond_i \varphi_N(\mathbf{z}^N) \\+ \frac{1}{2N}\sum_{i\ne j} \tilde{\lambda}(z^i,z^j)\int_{\tilde{\Theta}}\big\{\varphi_N\big(\mathbf{z}^N\big(i,j,\tilde{\theta}\big)\big)-\varphi_N\big(\mathbf{z}^N\big)\big\}\tilde{\nu}(\dd\tilde{\theta}),\end{multline}
where $\tilde{\lambda}:E\times E\to\R_+$, $\tilde{\Theta}$ is a parameter set endowed with a probability measure~$\tilde{\nu}$ and $\mathbf{z}^N(i,j,\tilde{\theta})$ is the $N$ dimensional vector whose $k$ component is equal to
\[{z}^k(i,j,\theta) = \left\{\begin{array}{rcl} 
z^k & \text{if} & k\ne i,j\\
\tilde{\psi}_1(z^i,z^j,\tilde{\theta}) & \text{if} & k=i\\
\tilde{\psi}_2(z^i,z^j,\tilde{\theta}) & \text{if} & k=j
\end{array}
\right.,\]
for two given functions $\tilde{\psi}_1, \tilde{\psi}_2 :E\times E\times \tilde{\Theta}\to E$. In this case, the general Boltzmann equation \cref{eq:Boltzmannequationgeneral_summary} can be re-written: 
\begin{multline}\label[IIeq]{eq:boltzmannpsitilde_summary}
\frac{\dd}{\dd t}\langle f_t,\varphi\rangle = \langle f_t, L^{(1)}\varphi\rangle+\frac{1}{2}\int_{\tilde{\Theta}\times E^2}\tilde{\lambda}(z_1,z_2)\Big\{\varphi\big(\tilde{\psi}_1(z_1,z_2,\tilde{\theta})\big)+\varphi\big(\tilde{\psi}_2(z_1,z_2,\tilde{\theta})\big)\\-\varphi(z_1)-\varphi(z_2)\Big\}\tilde{\nu}(\dd\tilde{\theta})f_t(\dd z_1)f_t(\dd z_2).    
\end{multline}
The generator \cref{eq:Boltzmannnonsym_summary} slightly differs from \cref{eq:boltzmanngenerator_summary}, because the pair $(i,j)$ is distinguished from the pair $(j,i)$. Consequently, the double sum in \cref{eq:Boltzmannnonsym_summary} runs over all indices $i,j=1,\ldots,N$ while in the sum \cref{eq:boltzmanngenerator_summary}, it runs over the indices $i<j$. The two formulations are nevertheless equivalent as shown in Example \cref{example:symmetrization}, in the first part of this review. 

The study of Boltzmann models has historically been motivated by the study of the Boltzmann equation of rarefied gas dynamics which reads (in strong form): 
\begin{multline}\label[IIeq]{eq:Boltzmannphysics_summary}
    \partial_t f_t(x,v) + v\cdot\nabla_x f_t\\ = \int_{\R^d}\int_{\mathbb{S}^{d-1}} B(v-v_*,\sigma)\Big(f_t(x,v_*')f_t(x,v')-f_t(x,v_*)f_t(x,v)\Big)\dd v_*\dd\sigma,
\end{multline}
where
\begin{equation}\label[IIeq]{eq:postcollisionsigma_summary}
v'  =  \displaystyle{\frac{v+v_*}{2}+\frac{|v-v_*|}{2}\sigma},\quad 
v_*' = \displaystyle{\frac{v+v_*}{2}-\frac{|v-v_*|}{2}\sigma},\end{equation}
The function $B:\R^d\times\mathbb{S}^{d-1}\to\R_+$ called the \emph{cross-section} is of the form
\begin{equation}\label[IIeq]{eq:collisionkernelB_summary}B(u,\sigma) = \Phi(|u|)\Sigma(\theta),\end{equation}
with $\cos\theta=\frac{u}{|u|}\cdot\sigma$, $\theta\in[0,\pi]$. Some famous cross-sections are listed below.
\begin{itemize}
    \item (Hard spheres) 
    \begin{equation}\label[IIeq]{eq:hardsphereB_summary}\Phi(|u|)=|u|,\quad \Sigma(\theta)=1.\end{equation}
    \item (Maxwell molecules) 
    \begin{equation}\label[IIeq]{eq:truemaxwellmolecules_summary}\Phi(|u|)=1,\quad \int_0^\pi\Sigma(\theta)\dd \theta = +\infty.\end{equation}
    \item (Maxwell molecules with Grad's cutoff)
    \begin{equation}\label[IIeq]{eq:maxwellmoleculescutoff_summary}\Phi(|u|)=1, \quad \int_0^\pi \Sigma(\theta)\dd \theta <+\infty.\end{equation}
\end{itemize}

In a spatially homogeneous setting, the case of bounded $\Phi$ and integrable $\Sigma$ (including Maxwell molecules with Grad's cutoff) is a parametric Boltzmann model with 
\[\psi_1(v,v_*,\theta) = v',\quad\psi_2(v,v_*,\theta) = v'_*,\]
and 
\begin{align*}\lambda(v,v_*) &= \Phi(|v-v_*|)\int_0^\pi\Sigma(\theta)\dd\theta,\\ \Gamma(v,v_*,\dd z',\dd v_*,\dd v_*') &= \psi(v,v_*,\cdot)_{\#} \left(\frac{\Sigma}{\int_0^\pi\Sigma(\theta)\dd\theta}\right).\end{align*}
Mathematically, it is often much simpler to consider a bounded $\Phi$. However, physically, only Maxwell molecules satisfy this condition and they are therefore particularly studied because of this mathematical simplicity. The case of the unbounded models \cref{eq:hardsphereB_summary} and \cref{eq:truemaxwellmolecules_summary} is more delicate, see Example~\cref{example:noncutoffmodels} and Section \cref{sec:Boltzmannclassicalmodels}. The derivation of the Boltzmann equation of rarefied gas dynamics \cref{eq:Boltzmannphysics_summary} in various cases will be discussed in Section \cref{sec:boltzmannreview}. 

In addition to these important examples, further recent applications of Boltzmann models can be found in particular in socio-economical models of wealth and opinion dynamics such as the ones described in Section \cref{sec:socioeconomicmodels}.  

\subsection{Proving propagation of chaos}\label[II]{sec:proving_summary}

Some of the classical techniques to prove propagation of chaos are gathered in section which summarizes the content of Section \cref{sec:proving}. 

\subsubsection{Coupling methods}\label[II]{sec:coupling_summary} When a SDE description of the particle system is available, the coupling method initiated by McKean \cite{mckean_propagation_1969} and Sznitman \cite{sznitman_topics_1991} consists in comparing the trajectories of the particle system with the trajectories of a system of $N$ i.i.d processes with common law $f_t$. 

\begin{definition}[Chaos by coupling the trajectories]\label[II]{def:chaosbycouplingtrajectories_summary}
Let be given a final time $T\in (0,\infty]$, a distance $d_E$ on $E$ and $p\in\N$. Propagation of chaos holds by coupling the trajectories when for all $N\in\N$ there exist
\begin{itemize}
    \item a system of particles $(\mathcal{X}^N_t)^{}_t$ with law $f^N_t\in\pb(E^N)$ at time $t\leq T$, 
    \item a system of \emph{independent} processes $\big(\overline{\mathcal{X}}{}^N_t\big){}^{}_t$ with law $f^{\otimes N}_t\in\pb(E^N)$ at time $t\leq T$, 
    \item a number $\varepsilon(N,T)>0$ such that $\varepsilon(N,T)\underset{N\to+\infty}{\longrightarrow}0$, 
\end{itemize}
such that (pathwise case)
\begin{equation}\label[IIeq]{eq:chaoscouplingpathwise_summary}\frac{1}{N}\sum_{i=1}^N \E{\left[\sup_{t\leq T}d_E\big(X^i_t,\overline{X}{}^i_t\big)^p\right]}\leq \varepsilon(N,T),\end{equation}
or (pointwise case)
\begin{equation}\label[IIeq]{eq:chaoscouplingpointwise_summary}\frac{1}{N}\sum_{i=1}^N \sup_{t\leq T} \E{\left[d_E\big(X^i_t,\overline{X}{}^i_t\big)^p\right]}\leq \varepsilon(N,T).\end{equation}
\end{definition}

Note that \cref{eq:chaoscouplingpathwise_summary} implies \cref{eq:chaoscouplingpointwise_summary}. The bound \cref{eq:chaoscouplingpointwise_summary} implies: 
\[\sup_{t\leq T} W_p\big(f^N_{t},f^{\otimes N}_{t}\big)\leq \varepsilon(N,T)\underset{N\to+\infty}{\longrightarrow}0,\]
where $W_p$ denotes the Wasserstein-$p$ distance (see Definitions \cref{def:wasserstein} and \cref{def:spaceswasserstein}) on a $\pb(E^N)$ defined for $\mu,\nu\in E^N$ by: 
\[W_p(\mu,\nu) := \inf_{\pi\in\Pi(\mu,\nu)} \left(\frac{1}{N}\sum_{i=1}^N\int_{E^N\times E^N}  |x^i-y^j|^p\pi(\dd \mathbf{x},\dd \mathbf{y})\right)^{1/p}, \]
and $\Pi(\mu,\nu)$ is the set of all probability measures on $E^N\times E^N$ with marginals $\mu$ and $\nu$. It implies the propagation of chaos in the sense of Definition \cref{def:poc_summary} since the topology induced by the Wasserstein distance is stronger than the topology of the weak convergence of probability measures (see Section \cref{sec:def}). 

Coupling techniques are widely used and many examples will be presented below. The original argument of McKean and Sznitman is presented in Section \cref{sec:mckeancoupling}. It is based on the \emph{synchronous coupling} between the particle system \cref{eq:mckeanvlasov_summary} and the system of $N$ independent SDEs: 
\[\dd \overline{X}^i_t = b(\overline{X}^i_t,f_t)\dd t + \sigma(\overline{X}^i_t,f_t)\dd B^i_t,\]
where $(B^i_t)_t$ is the \emph{same} Brownian motion as in \cref{eq:mckeanvlasov_summary}. Other coupling techniques are presented in Section \cref{sec:othercouplingsreview}. 

For Boltzmann models, we postpone the discussion to Section \cref{sec:couplingBoltzmann}. 

\subsubsection{Compactness methods}\label[II]{sec:provingcompactness_summary} Thanks to Lemma \cref{lemma:caractchaos_summary}, the propagation of chaos property is equivalent to the convergence in law of the sequence of empirical measures. A natural strategy to prove such convergence is to prove on the one hand that it is possible to extract a converging subsequence and on the other hand to prove the uniqueness of the accumulation point. Note that these properties respectively show the existence and the uniqueness of the limit problem, which can be given, depending on the point of view, by a nonlinear PDE or a nonlinear martingale problem. The uniqueness property strongly depends on the limit nonlinear problem and it is an independent problem not necessarily related to the underlying particle system. In order to extract a converging subsequence, it is important to note that the sequence of empirical measures is a sequence of measure-valued random variables and in this context, it is natural to try to apply one of the classical or less classical stochastic tightness criteria recalled in Appendix \cref{appendix:tightness}. There is however an important subtlety to keep in mind: there are actually three strictly nonequivalent points of view on the empirical measure and depending on the one chosen, it provides three different nonequivalent results. These point of view are explained in great detailed in Section \cref{sec:poc} and we briefly recall them now. 
\begin{itemize}
\item The strongest point of view, called \emph{(strong) pathwise}, considers the empirical measure as the empirical measure associated to a sequence of $N$ random processes defined in the Skorokhod space, that is, with the previous notations, the sequence $(\mu_{\mathcal{X}^N_{[0,T]}})_N$. For each $N$, the empirical measure is thus a random element $\mu_{\mathcal{X}^N_{[0,T]}}\in\pb(D([0,T],E))$ and the goal is to prove the convergence of the laws in the space $\pb(\pb(D([0,T],E)))$. 

\item The second, weaker, point of view, called \emph{functional law of large numbers}, sees the empirical measure as a measure-valued process, that is, for each $N$, a random process $t\in[0,T]\mapsto \mu_{\mathcal{X}^N_t}\in\pb(E)$, i.e. a random variable in the space $D([0,T],\pb(E))$. The goal is thus to prove the convergence of the sequence of pathwise laws in the space $\pb(D([0,T],\pb(E)))$. 

\item Finally, the weakest point of view, called \emph{pointwise} point of view studies the flow of time marginals of the law of the empirical measure process, that is the mapping $t\in[0,T]\mapsto \mathrm{Law}(\mu_{\mathcal{X}^N_t})\in\pb(\pb(E))$. This defines a deterministic sequence in the functional space $C([0,T],\pb(\pb(E)))$. 

\end{itemize}
The first proofs of the propagation of chaos using compactness methods for spatially homogeneous version of the Boltzmann model \cref{eq:Boltzmannphysics_summary} are due to Tanaka \cite{tanaka_probabilistic_1983} and Sznitman \cite{sznitman_equations_1984}. For the McKean-Vlasov diffusion and the mean-field jump model, a detailed analysis can be found in \cite{graham_stochastic_1997, meleard_asymptotic_1996}. A more recent approach which exploits the gradient-flow structure of the McKean-Vlasov diffusion is due to~\cite{carrillo_-convexity_2020}. These results will be discussed in Section \cref{sec:mckeancompactnessreview} and Section \cref{sec:martingaleboltzmannreview}. 

\subsubsection{Generator related methods}\label[II]{sec:abstractmischlermouhot_summary} When seen as measure-valued processes, the sequence of empirical measures is a sequence of Markov processes in the space $\mathcal{P}(E)$. Since a Markov process is defined by its generator, the convergence (in law) of a sequence of processes can be recast into the convergence of the sequence of their generators. Based on this idea, the seminal article of Gr\"unbaum \cite{grunbaum_propagation_1971} is based on the asymptotic analysis of the generator of the empirical measure process when $N\to+\infty$. However, since $\pb(E)$ is only a metric space with no Banach structure, the rigorous definition of the infinitesimal generator of a measure-valued processes and the notion of convergence are in this case extremely delicate. A completed and rigorous version of Gr\"unbaum's original argument is due to \cite{mischler_kacs_2013, mischler_new_2015} and is discussed in detailed in Section \cref{sec:abstractmischlermouhot}. The main result is an abstract theorem (Theorem~\cref{thm:abstractMischler}) valid for a wide range of mean-field and Boltzmann models. This strategy has been applied to the spatially homogeneous version of the Boltzmann equation~\cref{eq:Boltzmannphysics_summary} in \cite{mischler_kacs_2013} and leads to uniform in time propagation of chaos results. The main results are gathered in Section~\cref{sec:kacprogram}. 

\subsubsection{Entropy bounds}

Most of the methods already presented require at some point some regularity assumptions on the interaction, typically a Lipschitz continuity property for the functions $b$ and $\sigma$ in \cref{eq:mckeanvlasov_summary} or of $\psi$ in \cref{eq:semiparametric_summary}. However, such assumption cannot be verified in many important cases, for instance Coulomb-type or Biot and Savart interactions. To deal with such systems, a new class of methods has recently been developed, based on the notion of entropy. In the present context, the study of entropy and entropy bounds originates from the large deviation analysis of particle systems, as reviewed in Section \cref{sec:largedeviations}. Following these techniques, recent results have been obtained for singular systems and systems with low regularity in physics and biology, see Sections \cref{sec:jabin} and \cref{sec:vorticity}. On a more probabilistic side, these techniques are also strongly linked to the Girsanov transform and also lead to propagation of chaos results for very general and abstract systems, see Section  \cref{sec:chaosviagirsanov}.

\begin{definition}[Entropy, Fisher information]\label[II]{def:entropyfisher_summary}
Let $\mathscr{E}$ be a Polish space. Given two probability measures $\mu,\nu \in \pb ( \mathscr{E})$ (or more generally two measures), the relative entropy is defined by
\[ H ( \nu | \mu ) := \int_{\mathscr{E}} \frac{\dd  \nu}{\dd  \mu} \log {\left( \frac{\dd \nu}{\dd  \mu} \right)} \dd  \mu , \]
where $\dd\nu/\dd\mu$ is the Radon-Nikodym derivative. When the two measures are mutually singular, by convention, the relative entropy is set to~$+\infty$ (the same holds for the Fisher information below). If moreover $E$ is endowed with a smooth manifold structure, the Fisher information can be defined as
\[ I ( \nu | \mu ) := \int_{\mathscr{E}} {\left| \nabla \log {\left( \frac{\dd  \nu}{\dd \mu} \right)} \right|}^2 \dd  \nu, \]
with the same conventions.
\end{definition}

The following lemma links entropy bound and Kac's chaos in Total Variation norm. It is a direct consequence of the Pinsker inequality and the Csiszar inequality~\cite{csiszar_sanov_1984}.

\begin{lemma}\label[II]{lemma:entropyPinskerCsiszar}
Let $\mathscr{E}$ be a Polish space and let $f^N\in \pb(\mathscr{E}^N)$ and $f\in \pb(\mathscr{E})$. For every nonnegative integer $k \leq N$, it holds that
\begin{equation*}\frac{1}{2}\big\| f^{k,N} - f^{\otimes k} \big\|_{\mathrm{TV}}^2 \leq H {\left( f^{k,N} \big| f^{\otimes k} \right)} \leq \frac{k}{N}  H {\left( f^N | f^{\otimes N} \right)}, \end{equation*}
where $\|\cdot\|_{\mathrm{TV}}$ is the Total Variation norm (which induces a topology stronger than the topology of the weak convergence of probability measures, see Section \cref{sec:def}). 
\end{lemma}

For the McKean-Vlasov diffusion \cref{eq:mckeanvlasov_summary}, the following lemma gives a way to bound the relative entropy between the $N$-particle distribution and its mean-field limit. The first pathwise inequality is a consequence of the Girsanov theorem (see Appendix \cref{appendix:girsanov} and Lemma \cref{lemma:entropyboundgirsanov_summary}). The second one can be formally obtained by direct computations (see Lemma \cref{lemma:computeH}). 

\begin{lemma}[Pathwise and pointwise entropy bounds]\label[II]{lemma:entropyboundgirsanov_summary}
Let $T>0$ and $I=[0,T]$. For $N\in\N$, let $f^N_I\in\pb(C([0,T],(\R^{d})^N))$ be the law of the McKean-Vlasov diffusion $(\mathcal{X}^N_t)^{}_t$ defined by \cref{eq:mckeanvlasov_summary} with $b:\R^d\times\pb(\R^d)\to\R^d$ and $\sigma = I_d$,
and let $f^N_t\in\pb((\R^d)^N)$ its time marginal at time $t\in[0,T]$. Let $f_I\in \pb(D([0,T],\R^d))$ be the pathwise law of the limit nonlinear McKean-Vlasov diffusion \cref{eq:mckeanvlasov-limit_summary} and let $f_t\in\pb(\R^d)$ be its time marginal at time $t\in[0,T]$ (it is the solution of \cref{eq:mckeanvlasov-pde_summary}).
\begin{itemize}
\item For any $k\leq N$ it holds that
\begin{equation}\label[IIeq]{eq:entropyboundgirsanov_summary}H{\big(f_I^{k,N}|f_I^{\otimes k}\big)}\leq\frac{k}{2}\E{\left[\int_0^T \big|b\big({X}^1_t,\mu_{\mathcal{X}^N_t}\big)-b({X}^1_t,f_t)\big|^2\dd t\right]}.\end{equation}
\item For every $\alpha > 0$ it holds that
\begin{equation}\label[IIeq]{eq:computeH_summary}\frac{\dd}{\dd t} H \big( f^N_t | f_t^{\otimes N} \big) \leq \frac{\alpha - 1}{2} I \big( f^N_t | f_t^{\otimes N} \big) + \frac{N}{2 \alpha} \E {\left[ \big| b \big( X^1_t , \mu_{\mathcal{X}^N_t} \big) - b ( X^1_t ,f_t ) \big|^2 \right]}. \end{equation}
\end{itemize}
\end{lemma}

\subsubsection{Interaction graphs}\label[II]{sec:interactiongraph_summary}

In an abstract Boltzmann model given by the generator~\cref{eq:boltzmanngenerator_summary} in Section~\cref{sec:boltzmann_summary}, the binary interactions can be represented by graph structures. Given a trajectorial realisation of the particle system, the \emph{interaction graph} of a particle (or a group of particles) is built backward in time and retain the genealogical interactions which determine the particle at the current time (i.e. the history of the collisions). Before building graphs from particle realisations, the minimal structure of such a possible graph is detailed in the following definition.

\begin{definition}[Interaction graph] \label[II]{def:intearctiongraph}
Consider an index $i\in\{1,\ldots,N\}$ (it will stand later for the index of a particle). An interaction graph for $i$ at time $t>0$ is the data of
\begin{enumerate}
    \item a $k$-tuple $\mathcal{T}_k=(t_1,\ldots, t_k)$ of \emph{interaction times} $t>t_1>t_2>\ldots>t_k>0$,
    \item a $k$-tuple $\mathcal{R}_k=(r_1,\ldots,r_k)$ of pairs of indexes, where for $\ell\in\{1,\ldots,k\}$, the pair denoted by $r_\ell=(i_\ell,j_\ell)$ is such that $j_\ell\in\{i_0,i_1,\ldots,i_{\ell-1}\}$ with the convention $i_0=i$ and $i_\ell\in\{1,\ldots,N\}$.
\end{enumerate}
Such an interaction graph is denoted by $\mathcal{G}_i(\mathcal{T}_k,\mathcal{R}_k)$.
\end{definition}

Given a trajectorial realisation of a Boltzmann particle system, the interaction graph of the particle $i$ retains the minimal information needed to compute the state of particle~$i$ at time $t>0$. It is constructed from Definition \cref{def:intearctiongraph} as follows.
\begin{itemize}
    \item The set $(i_1,\ldots,i_k)$ is the set of indexes of the particles which interacted directly or indirectly with particle $i$ during the time interval $(0,t)$ (an indirect interaction means that the particle has interacted with another particle which interacted directly or indirectly with particle $i$) -- note that the $i_\ell$'s may not be all distinct.
    \item The times $(t_1,\ldots,t_k)$ are the times at which an interaction occurred.
    \item For $\ell\in\{1,\ldots,k\}$, the indexes $(i_\ell,j_\ell)$ are the indexes of the two particles which interacted together at time $t_\ell$. 
\end{itemize}

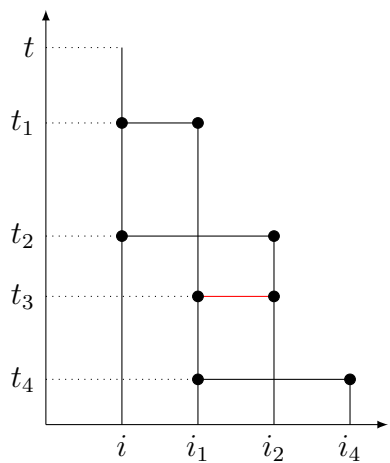
\begin{figure}
    \centering
    \begin{tikzpicture}[>=latex]
    \draw[->] (0,0)--(0,5.5);
    \draw[->] (0,0)--(4.5,0);
    
    \coordinate (t) at (0,5);
    \coordinate (t1) at (0,4);
    \coordinate (t2) at (0,2.5);
    \coordinate (t3) at (0,1.7);
    \coordinate (t4) at (0,0.6);
    
    \node[anchor=east] at (t) {$t$};
    \node[anchor=east] at (t1) {$t_1$};
    \node[anchor=east] at (t2) {$t_2$};
    \node[anchor=east] at (t3) {$t_3$};
    \node[anchor=east] at (t4) {$t_4$};
    
    \coordinate (i) at (1,5);
    \coordinate (j1) at (1,4);
    \coordinate (i1) at (2,4);
    \coordinate (j2) at (1,2.5);
    \coordinate (i2) at (3,2.5);
    \coordinate (j3) at (2,1.7);
    \coordinate (i3) at (3,1.7);
    \coordinate (j4) at (2,0.6);
    \coordinate (i4) at (4,0.6);
    
    \draw[-,color=red] (j3) -- (i3) ;
    
    \node[anchor=north] at (1,0) {$i$};
    \node[anchor=north] at (2,0) {$i_1$};
    \node[anchor=north] at (3,0) {$i_2$};
    \node[anchor=north] at (4,0) {$i_4$};

    \filldraw[black] (j1) circle (2pt);
    \filldraw[black] (i1) circle (2pt);
    \filldraw[black] (j2) circle (2pt);
    \filldraw[black] (i2) circle (2pt);
    \filldraw[black] (j3) circle (2pt);
    \filldraw[black] (i3) circle (2pt);
    \filldraw[black] (j4) circle (2pt);
    \filldraw[black] (i4) circle (2pt);
    
    \draw[-] (1,0) -- (i) ;
    \draw[-] (2,0) -- (i1) ;
    \draw[-] (3,0) -- (i2) ;
    \draw[-] (4,0) -- (i4) ;
    
    \draw[-] (j1) -- (i1) ; 
    \draw[-] (j2) -- (i2) ; 
    \draw[-] (j4) -- (i4) ;
    
    \draw[-,dotted] (t) -- (i)  ;
    \draw[-,dotted] (t1) -- (j1)  ;
    \draw[-,dotted] (t2) -- (j2)  ;
    \draw[-,dotted] (t3) -- (j3)  ;
    \draw[-,dotted] (t4) -- (j4)  ;

    \end{tikzpicture}
    \caption{An interaction graph. The vertical axis represents time. Each particle is represented by a vertical line parallel to the time axis. The index of a given particle is written on the horizontal axis. The construction is done backward in time starting from time $t$ where only particle $i$ is present. At each time $t_\ell$, if $i_\ell$ does not already belong to the graph, it is added on the right (with a vertical line which starts at $t_\ell$). The couple $r_\ell=(i_\ell,j_\ell)$ of interacting particles at time $t_\ell$ is depicted by an horizontal line joining two big black dots on the vertical line representing the particles $i_\ell$ and $j_\ell$. for instance, on the depicted graph, $r_2=(i_2,i)$. Note that at time $t_3$, $r_3=(i_1,i_2)$ (or indifferently $r_3=(i_2,i_1)$) where $i_1$ and $i_2$ were already in the system. Index $i_3$ is skipped and at time $t_4$, the route is $r_4=(i_4,i_1)$. The \emph{recollision} occurring at time $t_3$ is depicted in red.}
    \label[II]{fig:interactiongraph_summary}
\end{figure}

Following the terminology of \cite{graham_stochastic_1997}, a \emph{route} of size $q$ between $i$ and $j$ is the union of $q$ elements $r_{\ell_k}=(i_{\ell_k},j_{\ell_k})$, $k=1,\ldots,q$ such that $i_{\ell_1}=i$, $i_{\ell_{k+1}}=j_{\ell_k}$ and $j_{\ell_q}=j$. A route of size 1 (\emph{i.e} a single element $r_\ell$) is simply called a route. A route which involves two indexes which were already in the graph before the interaction time (backward in time) is called a \emph{recollision}. This construction is more easily understood with the graphical representation of an interaction graph shown on Figure \cref{fig:interactiongraph_summary}. 

The definition of interaction graphs can be extended straightforwardly starting from a group of particles instead of only one particle. This representation does not take into account the physical trajectories of the particles, it only retains the history of the interactions among a group of particles. Note that the graph is not a tree in general since the $i_\ell$'s are not necessarily distinct. It is a tree when no recollision occurs. 

The following definition extends the construction of Definition \cref{def:intearctiongraph} to the case of random parameters.

\begin{definition}[Random interaction graph] \label[II]{def:randomintearctiongraph}
Let $\Lambda>0$, $N\in\N$, $i\in\{1,\ldots,N\}$ and $t>0$. Let $(T^{m,\ell})_{1\leq k<\ell\leq N}$ be $N(N-1)/2$ independent Poisson processes with rate $\Lambda/N$. For each Poisson process $T^{m,\ell}$ we denote by $(T^{m,\ell}_n)^{}_n$ its associated increasing sequence of jump times. The sets of times $\mathcal{T}_k=(t_1,\ldots,t_k)$ and routes $\mathcal{R}_k=(r_1,\ldots,r_k)$ are defined recursively as follows. Initially, $t_0=t$ and $i_0=i$ and for $k\geq0$,
\begin{equation}\label[IIeq]{eq:recursivejumpingtimes_summary}t_{k+1} = \max_{\ell,p,n}\big\{T^{i_\ell,p}_n\,|\, T^{i_\ell,p}_n<t_k,\,\ell\leq k\big\}.\end{equation}
Then, given $(\ell,p,n)$ such that $t_{k+1}=T^{i_\ell,p}_n$, $i_{k+1}=p$ and $j_{k+1}=i_\ell$ so that $r_{k+1}=~(i_{k+1},j_{k+1})$. The procedure is stopped once the set on the right-hand side of \cref{eq:recursivejumpingtimes_summary} is empty (it happens almost surely after a finite number of iterations). The resulting interaction graph $\mathcal{G}_i(\mathcal{T}_k,\mathcal{R}_k)$ is called the random interaction graph with rate $\Lambda$ rooted on $i$ at time $t$. The definition is extended similarly starting from a finite number of indexes $(i_0,i_1,\ldots,i_k)$ instead of just $i$.  
\end{definition}

As explained before, a realisation of a Boltzmann particle system immediately gives an interaction graph for each particle. More importantly, given an interaction graph, it is possible to construct a forward realisation of a Boltzmann particle. More precisely, when the interaction graph is sampled as a random interaction graph following Definition \cref{def:randomintearctiongraph}, then the following straightforward lemma constructs a forward realisation of a stochastic process whose pathwise law is equal to $f^{1,N}_{[0,t]}$, the first marginal of the law $f^N_{[0,t]}$ of a Boltzmann particle system given by the generator~\cref{eq:boltzmanngenerator_summary} on the time interval $[0,t]$. 


\begin{lemma} \label[II]{lem:randomintearctiongraph}
Let us consider the Boltzmann setting given by Assumption \cref{assum:L2_summary} together with the uniform bound \cref{eq:uniformboundlambda_summary} on $\lambda$. Given a realisation of a random interaction graph sampled beforehand as in Definition \cref{def:randomintearctiongraph}, apply the following procedure:
\begin{enumerate}
    \item At time $t=0$, let the particles $Z^{i_\ell}_0$ be distributed according to the initial law.
    \item Between two collision times, the particles evolve according to $L^{(1)}$. 
    \item At a collision time $t_\ell$, with probability $\lambda(Z^{i_\ell}_{t_\ell^-},Z^{j_\ell}_{t_\ell^-})/\Lambda$, the new states of particles $i_\ell$ and $j_\ell$ are sampled according to \[{\left(Z^{i_\ell}_{t_\ell^+},Z^{j_\ell}_{t_\ell^+}\right)}\sim\Gamma^{(2)}{\left(Z^{i_\ell}_{t_\ell^-},Z^{j_\ell}_{t_\ell^-},\dd z_1, \dd z_2\right)}.\]
\end{enumerate}
Then the process $(Z^{i}_s)_{s\leq t}$ is distributed according to the one-particle marginal $f^{1,N}_{[0,t]}$ of the law $f^N_{[0,t]}$ of a Boltzmann particle system given by the generator \cref{eq:boltzmanngenerator_summary} on the time interval $[0,t]$.
\end{lemma}

This result will be useful later. Interaction graphs and random interaction graphs are used in \cite{lanford_time_1975} and in \cite{graham_stochastic_1997} to prove propagation of chaos by a direct control of the trajectories of the particles. This will be reviewed respectively in Sections \cref{sec:lanford} and \cref{sec:pathwisekactheorem}.

\section{McKean-Vlasov diffusion models}\label[II]{sec:mckeanreview}

Since the seminal work of McKean \cite{mckean_propagation_1969}, later extended by Sznitman \cite{sznitman_topics_1991}, a very popular method of proving propagation of chaos for mean-field systems is the synchronous coupling method (Section \cref{sec:synchronouscouplingreview}). Over the last years, some alternative coupling methods have been proposed to handle either weaker regularity or to get uniform in time estimates under mild physically relevant assumptions (Section \cref{sec:othercouplingsreview}). Alternatively to these SDE techniques, the empirical process can be studied using stochastic compactness methods \cite{sznitman_nonlinear_1984,graham_stochastic_1997}, leading to (non quantitative) results valid for mixed jump-diffusion models (Section~\cref{sec:mckeancompactnessreview}). Recent works focus on large deviation techniques, in particular the derivation of entropy bounds from Girsanov transform arguments \cite{jabin_quantitative_2018,lacker_strong_2018}, this allows interactions with a very weak regularity (Section \cref{sec:jabin}) or with a very general form (Section \cref{sec:mckeangeneralinteractions}). 
 
\subsection{Synchronous coupling}\label[II]{sec:synchronouscouplingreview}

In this section, we give several examples of the very fruitful idea of synchronous coupling presented in Section \cref{sec:coupling_summary}. The first instance of synchronous coupling that we are aware of is due to McKean himself although the most popular form of the argument is due to Sznitman. This will be discussed in Section \cref{sec:mckeancoupling}. This original argument is valid under strong Lipschitz and boundedness assumptions but it can be extended to more singular cases, as explained in Section \cref{sec:mckeantowardssingular}. Finally, in Section \cref{sec:gradientsystems}, the strategy is successfully applied to gradient systems and leads to uniform in time and convergence to equilibrium results. 

\subsubsection{McKean's theorem and beyond for Lipschitz interactions} \label[II]{sec:mckeancoupling}

The following theorem due to McKean is the most important result of this section. For a function $K : E^2 \rightarrow \R$, we recall the notation $K \star \mu(x) := \int K(x, y) \mu( \dd y)$.

\begin{theorem}[McKean]\label[II]{thm:mckean} Let the drift and diffusion coefficients in \cref{eq:mckeanvlasov_summary} be defined by 
\begin{equation}\label[IIeq]{eq:mckeanthmassum}\forall x\in \R^d,\forall \mu\in\pb(\R^d),\quad b(x,\mu) := \tilde{b}\big(x,K_1\star\mu(x)\big),\quad \sigma(x,\mu) = \tilde{\sigma}\big(x,K_2\star\mu(x)\big),\end{equation}
where $K_1:\R^d\times\R^d\to\R^m$, $K_2:\R^d\times\R^d\to\R^n$, $\tilde{b}:\R^d\times\R^m\to\R^d$ and $\tilde{\sigma}:~\R^d\times~\R^n\to\mathcal{M}_d(\R)$ are globally Lipschitz and $K_1,K_2$ are bounded. Then pathwise chaos by coupling in the sense of Definition \cref{def:chaosbycouplingtrajectories_summary} holds for any $T>0$, $p=2$, with the synchronous coupling
\begin{equation}\label[IIeq]{eq:mckeanparticlessynch}X^{i,N}_t = X^i_0 + \int_0^t\tilde{b}{\left(X^{i,N}_s,K_1\star\mu_{\mathcal{X}_s^N}{\left(X^{i,N}_s\right)}\right)}\dd s+\int_0^t \tilde{\sigma}{\left(X^{i,N}_s,K_2\star\mu_{\mathcal{X}_s^N}{\left(X^{i,N}_s\right)}\right)}\dd B^i_s,\end{equation}
and
\begin{equation}\label[IIeq]{eq:nonlinearmckeanparticlessynch}\overline{X}{}^{i,N}_t = X^i_0 + \int_0^t\tilde{b}{\left(\overline{X}{}^{i,N}_s,K_1\star f_s{\left(\overline{X}{}^{i,N}_s\right)}\right)}\dd s+\int_0^t \tilde{\sigma}{\left(\overline{X}{}^{i,N}_s,K_2\star f_s{\left(\overline{X}{}^{i,N}_s\right)}\right)}\dd B^i_s.\end{equation}
It means that the trajectories satisfy: 
\[\frac{1}{N}\sum_{i=1}^N \E{\left[\sup_{t\leq T} \big|X^i_t-\overline{X}{}^i_t\big|^2\right]}\leq \varepsilon(N,T),\]
where the convergence rate is given by 
\begin{equation}\label[IIeq]{eq:mckeanthmeps}\varepsilon(N,T) = \frac{c_1(b,\sigma,T)}{N}\e^{c_2(b,\sigma,T)T},\end{equation}
for some absolute constants $C,\tilde{C},C_\mathrm{BDG}>0$ not depending on $N,T$,
\begin{equation}\label[IIeq]{eq:mckeanthmc1}c_1(b,\sigma,T) := CT {\left(T\|K_1\|_\infty^2\|\tilde{b}\|_\mathrm{Lip}^2+C_{\mathrm{BDG}}\|K_2\|_\infty^2\|\tilde{\sigma}\|_\mathrm{Lip}^2\right)},\end{equation}
and 
\begin{equation}\label[IIeq]{eq:mckeanthmc2}c_2(b,\sigma,T) := \tilde{C}{\left(T{\left(1+\|K_1\|_\mathrm{Lip}^2\right)\|\tilde{b}\|_\mathrm{Lip}^2\,+C_{\mathrm{BDG}}\left(1+\|K_2\|_\mathrm{Lip}^2\right)}\|\tilde{\sigma}\|_\mathrm{Lip}^2\right)}.\end{equation}
\end{theorem}

We present two proofs of this result. The first one is the original proof due to McKean \cite{mckean_propagation_1969}. The second one is due to Sznitman \cite{sznitman_topics_1991}. Sznitman's proof is a slightly shorter and more general version of McKean's proof. We chose to include McKean's original argument for three reasons. First it gives an interesting and somehow unusual probabilistic point of view on the interplay between exchangeability and independence (see Section \cref{sec:infinitesystems}). This is an underlying idea for all the models presented in this review which is made very explicit in McKean's proof. Secondly, although the computations in both proofs are very much comparable, McKean's proof is philosophically an existence result while Sznitman's proof is based on the well-posedness result stated in Proposition \cref{prop:wellposednessmckeanlipschitz}. Finally, it seems that McKean's proof has been somehow forgotten in the community or is sometimes confused with Sznitman's proof which in turn has become incredibly popular. McKean's argument was first published in \cite{mckean_propagation_1967} and then re-published in \cite{mckean_propagation_1969}. Both references are not easy to find nowadays and it is probably the source of the confusion between the two proofs. 

\begin{proof}[Proof (McKean)] The originality of this proof is that the nonlinear process is not introduced initially. It appears as the limit of a Cauchy sequence of coupled systems of particles with increasing size. Let  $(B^i_t)^{}_t$, $i\geq1$ be an infinite collection of independent Brownian motions and for $N\in \N$ we recall the notation
\[\mathcal{X}^{N}_t = \big(X^{1,N}_t,\ldots,X^{N,N}_t\big)\in (\R^d)^N,\]
where $(X^{i,N}_t)^{}_t$ solves \cref{eq:mckeanparticlessynch}. The idea is to prove that the sequence (in $N$) of processes $(X^{1,N}_t)^{}_{t}$ is a Cauchy sequence in $L^2\big(\Omega,C([0,T],\R^d)\big)$ and then to identify the limit as the solution of \cref{eq:nonlinearmckeanparticlessynch}.  The proof is split into several steps. 
\medskip

\noindent\textit{\textbf{Step 1.} Cauchy estimate}
\medskip 

Let $M>N$ and let us consider the coupled particle systems $\mathcal{X}^N$ and $\mathcal{X}^M$ where the $N$ first particles in $\mathcal{X}^M$ have the same initial condition as $X^{1,N},\ldots, X^{N,N}$ and are driven by the same Brownian motions $B^1,\ldots B^N$. Using \cref{eq:mckeanparticlessynch} and the Burkholder-Davis-Gundy inequality it holds that for a constant $C_{\mathrm{BDG}}>0$,

\begin{multline}\label[IIeq]{eq:mckeanthmfirstinequality}
    \E{\left[\sup_{t\leq T} \big|X^{1,M}_t-X^{1,N}_t\big|^2\right]} \leq 2T\int_0^T \E{\left|b{\left(X^{1,M}_t,\mu_{\mathcal{X}^M_t}\right)}-b{\left(X^{1,N}_t,\mu_{\mathcal{X}^N_t}\right)}\right|}^2\dd t \\ 
    \qquad+2C_{\textrm{BDG}}\int_0^T \E{\left|\sigma{\left(X^{1,M}_t,\mu_{\mathcal{X}^M_t}\right)}-\sigma{\left(X^{1,N}_t,\mu_{\mathcal{X}^N_t}\right)}\right|}^2\dd t.
\end{multline}

For the first term on the right-hand side of \cref{eq:mckeanthmfirstinequality}, we write: 
\begin{multline}\label[IIeq]{eq:mckeanthmsecondinequality}
    \E{\left|b{\left(X^{1,M}_t,\mu_{\mathcal{X}^M_t}\right)}-b{\left(X^{1,N}_t,\mu_{\mathcal{X}^N_t}\right)}\right|}^2 \leq 2\E{\left|b{\left(X^{1,M}_t,\mu_{\mathcal{X}^M_t}\right)}-b{\left(X^{1,M}_t,\mu_{\mathcal{X}^{N,M}_t}\right)}\right|}^2\\
    \qquad+2\E{\left|b{\left(X^{1,M}_t,\mu_{\mathcal{X}^{N,M}_t}\right)}-b{\left(X^{1,N}_t,\mu_{\mathcal{X}^N_t}\right)}\right|}^2,
\end{multline}
where $\mathcal{X}^{N,M}_t=\big(X^{1,M}_t,\ldots,X^{N,M}_t\big)\in(\R^d)^N$. Each of the two terms on the right-hand side of \cref{eq:mckeanthmsecondinequality} is controlled using \cref{eq:mckeanthmassum}, the  Lipschitz assumptions and the fact that the $X^{j,M}$ are identically distributed. For the first term, expanding the square gives: 
\begin{align*}
    &\E{\left|b{\left(X^{1,M}_t,\mu_{\mathcal{X}^M_t}\right)}-b{\left(X^{1,M}_t,\mu_{\mathcal{X}^{N,M}_t}\right)}\right|}^2  \\
    & \leq \|\tilde{b}\|_\mathrm{Lip}^2\,\E\Big|\frac{1}{M}\sum_{j=1}^M K_1{\left(X^{1,M}_t,X^{j,M}_t\right)}-\frac{1}{N}\sum_{j=1}^N K_1{\left(X^{1,M}_t,X^{j,M}_t\right)}\Big|^2\\
    &\leq \|\tilde{b}\|_\mathrm{Lip}^2\, {\left(\frac{1}{M}+\frac{1}{N}-2\frac{N}{MN}\right)}\E{\left|K_1{\left(X^{1,M}_t,X^{2,M}_t\right)}\right|}^2\\
    &\quad+\|\tilde{b}\|_\mathrm{Lip}^2\,{\left(\frac{M-1}{M}+\frac{N-1}{N}-2\frac{M(N-1)}{MN}\right)}\times\\
    &\phantom{\quad+\|\tilde{b}\|_\mathrm{Lip}^2\Big(\frac{M-1}{M}+\frac{N-1}{N}-}\times\E{\left[K_1{\left(X^{1,M}_t,X^{2,M}_t\right)}\cdot K_1{\left(X^{1,M}_t,X^{3,M}_t\right)}\right]}\\
    &\leq 2\,{\left(\frac{1}{N}-\frac{1}{M}\right)}\|K_1\|_\infty^2\|\tilde{b}\|_\mathrm{Lip}^2.
\end{align*}
For the second term, the Lipschitz assumptions leads to:
\begin{align*}
    &\E{\left|b{\left(X^{1,M}_t,\mu_{\mathcal{X}^{N,M}_t}\right)}-b{\left(X^{1,N}_t,\mu_{\mathcal{X}^N_t}\right)}\right|}^2\\
    &\qquad \leq 2\|\tilde{b}\|_\mathrm{Lip}^2\, \E\Big[\big|X^{1,N}_t-X^{1,M}_t\big|^2\\
    &\phantom{\qquad \leq 2\|\tilde{b}\|_\mathrm{Lip}^2\, \E\Big[}+\Big|\frac{1}{N}\sum_{j=1}^N K_1{\left(X^{1,M}_t,X^{j,M}_t\right)}-\frac{1}{N}\sum_{j=1}^N K_1{\left(X^{1,N}_t,X^{j,N}_t\right)}\Big|^2\Big]\\
    &\qquad \leq 2\left(1+2 \|K_1\|_\mathrm{Lip}^2\right)\|\tilde{b}\|_\mathrm{Lip}^2\, \E\big|X^{1,N}_t-X^{1,M}_t\big|^2.
\end{align*}
The same estimates hold for the diffusion term on the right-hand side of \cref{eq:mckeanthmfirstinequality} with~$\sigma$ instead of $b$ and $K_2$ instead of $K_1$. Gathering everything thus leads to: 
\begin{multline*}
\E{\left[\sup_{t\leq T} \big|X^{1,M}_t-X^{1,N}_t\big|^2\right]} \\\leq {\left(\frac{1}{N}-\frac{1}{M}\right)}c_1(b,\sigma,T)+c_2(b,\sigma,T)\int_0^T \E\big|X^{1,N}_t-X^{1,M}_t\big|^2\dd t
\end{multline*}
where $c_1$ and $c_2$ are defined by \cref{eq:mckeanthmc1} and \cref{eq:mckeanthmc2}. Using (a generalisation of) Gronwall lemma, it follows that: 
\begin{equation}\label[IIeq]{eq:mckeanthmgronwall}\E{\left[\sup_{t\leq T} \big|X^{1,M}_t-X^{1,N}_t\big|^2\right]}\leq {\left(\frac{1}{N}-\frac{1}{M}\right)}c_1(b,\sigma,T)\e^{c_2(b,\sigma,T)T}.\end{equation}

\noindent\textit{\textbf{Step 2.} Cauchy limit and exchangeability}
\medskip 

The previous estimate implies that the sequence $(X^{1,N})_N$ is a Cauchy sequence in $L^2(\Omega,C([0,T],\R^d))$. Since this space is complete, this sequence has a limit denoted by $\overline{X}{}^1 \equiv (\overline{X}{}^1_t)^{}_t$. Applying the same reasoning for any $k\in \N$, there exists an infinite collection of processes $\overline{X}{}^k$, defined for each $k\geq1$ as the limit of $(X^{k,N})_N$. These processes are identically distributed and their common law depends only on $(X^i_0)^{}_{i\geq1}$ and $(B^i)^{}_{i\geq1}$ which are independent random variables. Moreover, knowing $(X^1_0,B^1)$ and for any measurable set $\mathscr{B}$, any event of the type $\{\overline{X}{}^1\in \mathscr{B}\}$  belongs to the $\sigma$-algebra of exchangeable events generated by the random variables $(X^i_0)^{}_{i\geq2}$ and $(B^i)_{i\geq2}$. Since these random variables are i.i.d, Hewitt-Savage 0-1 law (Theorem~\cref{thm:hewittsavage01}) states that this $\sigma$-algebra is actually trivial. It follows that $\overline{X}{}^1$ is a functional of~$X^1_0$ and $B^1$ only. The same reasoning applies for each $\overline{X}{}^k$ and hence the processes~$\overline{X}{}^k$ are also independent.
\medskip 

\noindent\textit{\textbf{Step 3.} Identification of the limit}
\medskip 

At this point, propagation of chaos is already proved and it only remains to identify the law of the $\overline{X}{}^k_t$ as the law of the solution of \cref{eq:nonlinearmckeanparticlessynch}. To do so, McKean defines for $i\in\{1,\ldots,N\}$ the processes 
\[\widetilde{X}^{i,N}_t = X^i_0 + \int_0^t b{\left(\overline{X}{}^i_s,\mu_{\overline{\mathcal{X}}{}^N_s}\right)}\dd s + \int_0^t \sigma{\left(\overline{X}{}^i_s,\mu_{\overline{\mathcal{X}}{}^N_s}\right)}\dd B^i_s,\]
where $\overline{\mathcal{X}}{}^N_t=(\overline{X}{}^1_t,\ldots,\overline{X}{}^N_t)$. 
From the independence of the processes and by the strong law of large numbers, the right hand side converges almost surely as $N\to+\infty$ towards the right hand side of \cref{eq:nonlinearmckeanparticlessynch} with $f_s$ being the law of $\overline{X}{}^i_s$ (which is the same for all $i$). Moreover, direct Lipschitz estimates lead to 
\[\E{\left[\sup_{t\leq T} \big|\widetilde{X}^i_t-X^i_t\big|^2\right]} \leq \frac{C}{N},\]
where $C$ is a constant which depends only on $T$, $\|K_1\|_\mathrm{Lip}$, $\|K_2\|_\mathrm{Lip}$. By uniqueness of the limit, it follows that $\overline{X}{}^i_t$ satisfies \cref{eq:nonlinearmckeanparticlessynch}. Moreover, the bound \cref{eq:mckeanthmeps} is obtained by taking the limit $M\to+\infty$ in \cref{eq:mckeanthmgronwall}.
\end{proof}

The following proof is due to Sznitman \cite{sznitman_topics_1991} in the case where $\sigma$ is constant and with $p=1$ in Definition \cref{def:chaosbycouplingtrajectories_summary}. The following (direct) adaptation to the model of Theorem \cref{thm:mckean} can be found in \cite[Proposition 2.3]{jourdain_propagation_1998}. 

\begin{proof}[Proof (Sznitman)] With a more direct approach, the strategy is to introduce both the particle system and its (known) limit given respectively by \cref{eq:mckeanparticlessynch} and \cref{eq:nonlinearmckeanparticlessynch} and to estimate directly the discrepancy between the two processes. Using the Burkholder-Davis-Gundy inequality, it holds that for a constant $C_\mathrm{BDG}>0$, 
\begin{multline}\label[IIeq]{eq:mckeanthmsznitfirstineq}
    \E{\left[\sup_{t\leq T} \big|\overline{X}{}^i_t-X^i_t\big|^2\right]}
    \leq 2T\int_0^T \E{\left|b{\left(\overline{X}{}^i_t,f_t\right)}-b{\left(X^i_t,\mu_{\mathcal{X}^N_t}\right)}\right|}^2\dd t \\
    + 2C_\mathrm{BDG}\int_0^T \E{\left|\sigma{\left(\overline{X}{}^i_t,f_t\right)}-\sigma{\left(X^i_t,\mu_{\mathcal{X}^N_t}\right)}\right|}^2\dd t. 
\end{multline}

The drift term on the right-hand side of \cref{eq:mckeanthmsznitfirstineq} is split into two terms as follows:
\begin{multline}\label[IIeq]{eq:mckeanthmsznitsecondineq}
    \E{\left|b{\left(\overline{X}{}^i_t,f_t\right)}-b{\left(X^i_t,\mu_{\mathcal{X}^N_t}\right)}\right|}^2 \leq 2\E{\left|b{\left(\overline{X}{}^i_t,f_t\right)}-b{\left(\overline{X}{}^i_t,\mu_{\overline{\mathcal{X}}{}^N_t}\right)}\right|}^2 \\
    +2\E{\left|b{\left(\overline{X}{}^i_t,\mu_{\overline{\mathcal{X}}{}^N_t}\right)}-b{\left(X^i_t,\mu_{\mathcal{X}^N_t}\right)}\right|}^2.
\end{multline}
For the first term on the right-hand side of \cref{eq:mckeanthmsznitsecondineq}, the assumption \cref{eq:mckeanthmassum} and the Lipschitz assumptions give: 
\begin{multline*}\E{\left|b{\left(\overline{X}{}^i_t,f_t\right)}-b{\left(\overline{X}{}^i_t,\mu_{\overline{\mathcal{X}}{}^N_t}\right)}\right|}^2 \leq \|\tilde{b}\|_\mathrm{Lip}^2\,\E\Big|K_1\star f_t(\overline{X}{}^i_t)-\frac{1}{N}\sum_{j=1}^N K_1(\overline{X}{}^i_t,\overline{X}{}^j_t)\Big|^2\\
= \frac{\|\tilde{b}\|_\mathrm{Lip}^2}{N^2}\,\E\Big|\sum_{j=1}^N\Big\{K_1\star f_t(\overline{X}{}^i_t)- K_1(\overline{X}{}^i_t,\overline{X}{}^j_t)\Big\}\Big|^2 .
\end{multline*}
Expanding the square, it leads to: 
\begin{align*}
    &\E{\left|b{\left(\overline{X}{}^i_t,f_t\right)}-b{\left(\overline{X}{}^i_t,\mu_{\overline{\mathcal{X}}{}^N_t}\right)}\right|}^2 \\
    &\leq \frac{\|\tilde{b}\|^2_\mathrm{Lip}}{N^2} \sum_{k,\ell=1}^N \E{\left[{\left(K_1\star f_t(\overline{X}{}^i_t)- K_1(\overline{X}{}^i_t,\overline{X}{}^k_t)\right)}\cdot {\left(K_1\star f_t(\overline{X}{}^i_t)- K_1(\overline{X}{}^i_t,\overline{X}{}^\ell_t)\right)}\right]}\\
    &\leq \frac{4\|\tilde{b}\|^2_\mathrm{Lip}\|K_1\|^2_\infty}{N}\\
    &\qquad+\frac{\|\tilde{b}\|^2_\mathrm{Lip}}{N^2}\,\sum_{k\ne \ell}\E{\left[{\left(K_1\star f_t(\overline{X}{}^i_t)- K_1(\overline{X}{}^i_t,\overline{X}{}^k_t)\right)}\cdot {\left(K_1\star f_t(\overline{X}{}^i_t)- K_1(\overline{X}{}^i_t,\overline{X}{}^\ell_t)\right)}\right]}.
\end{align*}
When $k\ne\ell$, using the fact that $\overline{X}^k_t$ and $\overline{X}^\ell_t$ are independent, $f_t$-distributed and independent of $\overline{X}^i_t$, it holds that:
\begin{align*}
    &\E{\left[{\left(K_1\star f_t(\overline{X}{}^i_t)- K_1(\overline{X}{}^i_t,\overline{X}{}^k_t)\right)}\cdot {\left(K_1\star f_t(\overline{X}{}^i_t)- K_1(\overline{X}{}^i_t,\overline{X}{}^\ell_t)\right)}\right]} \\
    &= \E\Big[ \E{\left[{\left(K_1\star f_t(\overline{X}{}^i_t)- K_1(\overline{X}{}^i_t,\overline{X}{}^k_t)\right)}\cdot {\left(K_1\star f_t(\overline{X}{}^i_t)- K_1(\overline{X}{}^i_t,\overline{X}{}^\ell_t)\right)}\Big| \overline{X}{}^i_t\right]}\Big]\\
    &= \E{\left[ \E{\left[{\left(K_1\star f_t(\overline{X}{}^i_t)- K_1(\overline{X}{}^i_t,\overline{X}{}^k_t)\right)}\Big| \overline{X}{}^i_t\right]}\E{\left[{\left(K_1\star f_t(\overline{X}{}^i_t)- K_1(\overline{X}{}^i_t,\overline{X}{}^\ell_t)\right)}\Big| \overline{X}{}^i_t\right]}\right]}\\
    &=0,
\end{align*}
To obtain the last inequality observe that since $k\ne \ell$ at least one of them is not equal to $i$, let us assume that $\ell\ne i$. Then since $\mathrm{Law}(\overline{X}{}^\ell_t)=f_t$, it holds that
\[\E{\left[{\left(K_1\star f_t(\overline{X}{}^i_t)- K_1(\overline{X}{}^i_t,\overline{X}{}^\ell_t)\right)}\Big| \overline{X}{}^i_t\right]}=0.\]
In conclusion,
\begin{equation}\label[IIeq]{eq:mckeanthmszniterror}\E{\left|b{\left(\overline{X}{}^i_t,f_t\right)}-b{\left(\overline{X}{}^i_t,\mu_{\overline{\mathcal{X}}{}^N_t}\right)}\right|}^2\leq\frac{4\|\tilde{b}\|^2_\mathrm{Lip}\|K_1\|^2_\infty}{N}.\end{equation}
For the second-term on the right-hand side of \cref{eq:mckeanthmsznitsecondineq}, the Lipschitz assumptions give:
\begin{equation}\label[IIeq]{eq:mckeanthmsznitlip}\E{\left|b{\left(\overline{X}{}^i_t,\mu_{\overline{\mathcal{X}}{}^N_t}\right)}-b{\left(X^i_t,\mu_{\mathcal{X}^N_t}\right)}\right|}^2 \leq C\|\tilde{b}\|^2_\mathrm{Lip}{\left(1+\|K_1\|^2_\mathrm{Lip}\right)}\E \big|\overline{X}{}^i_t-X^i_t\big|^2.\end{equation}
The same estimates hold when $b$ and $K_1$ are replaced by $\sigma$ and $K_2$. Gathering everything leads to: 
\begin{multline*}
\E{\left[\sup_{t\leq T} \big|\overline{X}{}^i_t-X^i_t\big|^2\right]} \leq \frac{1}{N}c_1(b,\sigma,T)+c_2(b,\sigma,T)\int_0^T \E\big|\overline{X}{}^i_t-X^i_t\big|^2\dd t \\ 
\leq \frac{1}{N}c_1(b,\sigma,T)+c_2(b,\sigma,T)\int_0^T \E{\left[\sup_{s\leq t}\big|\overline{X}{}^i_s-X^i_s\big|^2\right]}\dd t.
\end{multline*}
The conclusion follows by Gronwall lemma. 
\end{proof}

\begin{remark} \begin{enumerate}
    \item The same synchronous coupling result holds (at least) with $p=1$ (see \cite[Corollary 3.3]{andreis_mckeanvlasov_2018}) and $p=4$ (see \cite[Proposition 2.3]{jourdain_propagation_1998}) in Definition \cref{def:chaosbycouplingtrajectories_summary}. 
    \item Pointwise chaos \cref{eq:chaoscouplingpointwise_summary} in Definition \cref{def:chaosbycouplingtrajectories_summary} is a consequence of pathwise chaos \cref{eq:chaoscouplingpathwise_summary} but it can also be proved directly with the same line of argument but where the Burkholder-Davis-Gundy inequality is replaced by the It\={o} isometry. 
    \item The starting inequality (Equation \cref{eq:mckeanthmfirstinequality} in McKean's proof and Equation \cref{eq:mckeanthmsznitfirstineq} in Sznitman's proof) can be replaced by an equality using It\={o}'s lemma. This may bring a small improvement in the constants $c_1$ and $c_2$. For instance, in the common case where $\sigma$ is a constant, we can write (in Sznitman's framework), 
    \[ \big|\overline{X}{}^i_t-X^i_t\big|^2 = 2\int_0^t {\left\langle b\big(\overline{X}{}^i_s,f_s\big)-b{\left(X^i_s,\mu_{\mathcal{X}^N_s}\right)},\overline{X}{}^i_s-X^i_s\right\rangle}\dd s.\]
    And we would obtain for some constants $C,\tilde{C}>0$ (see for instance the introduction of \cite{salem_gradient_2020}):
    \[\E{\left[\sup_{t\leq T}\big|\overline{X}{}^i_t-X^i_t\big|^2\right]}\leq C\frac{\|\tilde{b}\|_\mathrm{Lip}^2\|K_1\|_\infty^2}{N}\e^{\tilde{C}\|\tilde{b}\|_\mathrm{Lip}(1+\|K_1\|_\mathrm{Lip})T},\]
    and therefore propagation of chaos holds over a time interval $T\sim \log N$. Several example will be given in the following (see in particular Theorem \cref{thm:mckeanmoment} and Theorem \cref{thm:uniformpocuniformconvex}). 
\end{enumerate}
\end{remark}

When $\sigma=I_d$, the following corollary shows that the pathwise particle system is strongly chaotic in TV norm. This result has been proved in \cite[Theorem 5.5]{malrieu_logarithmic_2001}.

\begin{corollary}[Pathwise TV chaos]\label[II]{coro:mckeantvchaos} Under the same assumptions as in McKean's theorem but with $\sigma=I_d$, for all $k<N$ it holds that
\[\big\|f^{k,N}_{[0,T]}-f^{\otimes k}_{[0,T]}\big\|_\mathrm{TV} \leq C(T)\sqrt{\frac{k}{N}}.\]
\end{corollary}

\begin{proof} By the Pinsker inequality \cref{eq:pinsker} and the inequality \cref{eq:csiszarmarginal}, it holds that 
\[\big\|f^{k,N}_{[0,T]}-f^{\otimes k}_{[0,T]}\big\|^2_\mathrm{TV}\leq 2\frac{k}{N}H(f^N_{\left[0, T\right]}|f^{\otimes N}_{\left[0, T\right]}).\]
Using \cref{eq:entropyboundgirsanov_summary}, the right-hand side is bounded by:
\[\big\|f^{k,N}_{[0,T]}-f^{\otimes k}_{[0,T]}\big\|^2_\mathrm{TV}\leq 2k\E{\left[\int_0^T\big|b\big(X^1_t,\mu_{\mathcal{X}^N_t}\big)-b(X^1_t,f_t)\big|^2\right]}.\]
By McKean's theorem, the expectation on the right-hand side is bounded by $\frac{C(T)}{N}$ and the conclusion follows. 
\end{proof}

McKean's theorem can be directly generalised to more general, yet Lipschitz, settings as we shall see in Section \cref{sec:mckeangeneralinteractionsextendingmckean}.

\subsubsection{Towards more singular interactions}\label[II]{sec:mckeantowardssingular}

The hypotheses of McKean's theorem (bounded and globally Lipschitz interactions) are most often too strong in practice. Even though there is no real hope for better results at this level of generality, many directions have been explored to weaken the hypotheses in specific cases. 
\begin{enumerate}
\item \textbf{(Moment control).} A commonly admitted idea is that propagation of chaos should also hold for only locally Lipschitz interaction functions with polynomial growth provided that moment estimates can be proved (both at the particle level and for the limiting nonlinear system). 
\item \textbf{(Moderate interaction and cut-off).} If one is mainly interested in the derivation of a singular nonlinear system, another idea is to smoothen the interaction at the particle level, for instance by adding a cutoff parameter or by convolution with a sequence of mollifiers. Such procedures typically depend on a smoothing parameter $\varepsilon$ that will go to zero. For a fixed $\varepsilon>0$ McKean's theorem gives a (quantitative) error estimate between the particle system and a smoothened nonlinear system. Then the idea is to take a smoothing parameter $\varepsilon\equiv \varepsilon_N$ which depends on $N$ such that $\varepsilon_N\to0$ as $N\to+\infty$. Taking advantage of the quantitative bound given by McKeans's theorem, the goal is to choose an appropriate $\varepsilon_N$ (usually a very slowly converging sequence) to pass to the limit directly from the smooth particle system to the singular nonlinear system.    
\end{enumerate}

In the present section, we give some examples of these ideas which naturally extend Sznitman's proof of McKean's theorem using synchronous coupling. Note that all the proofs crucially depend at some point of a well-posedness result for the nonlinear system. In practise, for singular interactions, proving propagation of chaos therefore largely depends on the considered model. Several examples for classical PDEs in kinetic theory can be found in Section \cref{sec:classicalpdekinetic}.

\subsubsection*{Moment control.} In \cite{bolley_stochastic_2011} the authors introduce some sufficient conditions on the interaction kernels $K_1$ and $K_2$ to extend the result of McKean's theorem to non globally Lipschitz bounded settings. This comes at the price of a strong assumption on the boundedness of the moments. Other examples using similar ideas will be detailed in Section \cref{sec:gradientsystems}. We first give a simple version of \cite[Theorem 1.1]{bolley_stochastic_2011} for the McKean-Vlasov model \cref{eq:mckeanvlasov_summary}. 

\begin{theorem}[\cite{bolley_stochastic_2011}]\label[II]{thm:mckeanmoment} Let us consider the McKean-Vlasov model \cref{eq:mckeanvlasov_summary}. Let $b,\sigma$ be as in McKean's Theorem \cref{thm:mckean} with $\tilde{b},\tilde{\sigma}$ globally Lipschitz and assume that there exists $\gamma>0$, $p > 0$ such that for $i=1,2$, $K_i$ satisfy for all $x,y,x',y'\in\R^d$,
\begin{equation}\label[IIeq]{eq:hyplipgrowth}\big|K_i(x,y)-K_i(x',y')\big|\leq \gamma \Big(|x-x'|+|y-y'|\Big)\Big(1+|x|^p+|y|^p+|x'|^p+|y'|^p\Big).\end{equation}
Assume there exist $\kappa>0$ and $p'\geq p$ such that for any $T>0$, Equations \cref{eq:mckeanparticlessynch} and~\cref{eq:nonlinearmckeanparticlessynch} admit solutions which verify
\begin{equation}\label[IIeq]{eq:hypexpmoment}\sup_{N}\,\sup_{t\leq T} \E\big[\e^{\kappa |X^i_t|^{p'}}\big]<+\infty,\quad \sup_{t\leq T} \E\big[\e^{\kappa |\overline{X}{}^i_t|^{p'}}\big]<+\infty.\end{equation}
and for $i=1,2$, 
\begin{equation}\label[IIeq]{eq:hypL2K}\sup_{t\leq T}\int_{\R^d\times\R^d} |K_i(x,y)|^2f_t(\dd x)f_t(\dd y)<+\infty.\end{equation}
Then for all $T>0$, there exists $C(T)>0$ such that for all $0 < t\leq T$,  
\[\E|X^i_t-\overline{X}{}^i_t|^2\leq \frac{C(T)}{N^{\e^{-Ct}}}.\]
Moreover, if the moment bound \cref{eq:hypexpmoment} holds for some $p'>p$ then for all $0<\varepsilon<1$, there exists $C(T)>0$ such that for all $t\leq T$, 
\[\E|X^i_t-\overline{X}{}^i_t|^2\leq \frac{C(T)}{N^{1-\varepsilon}}.\]
\end{theorem}

Sufficient conditions which ensure the well-posedness of \cref{eq:hypexpmoment} and
\cref{eq:hypL2K} are given by $p \leq 2$, $\tilde{b}$ and $\tilde{\sigma}$ bounded, $K_i(x,y) = \tilde{K}_i(x-y)$ with ${| \tilde{K}_i (x) |} \leq C {\left( 1 + |x| \right)}$ for $i=1,2$. This is a particular case of the more detailed result \cite[Theorem 1.2]{bolley_stochastic_2011}, see also \cite[Lemma 3.5]{bolley_stochastic_2011}).

\begin{proof}
The proof is similar to the proof of McKean's theorem using Sznitman's synchronous coupling but starting from It\=o's formula: 
\begin{align*}\frac{\dd}{\dd t}\E|\overline{X}{}^i_t-X^i_t|^2&=2\E\big\langle \overline{X}{}^i_t-X^i_t,b(\overline{X}{}^i_t,f_t)-b\big(X^i_t,\mu_{\mathcal{X}^N_t}\big)\big\rangle\\
&\quad+ 2\E\big\|\sigma(\overline{X}{}^i_t,f_t)-\sigma\big(X^i_t,\mu_{\mathcal{X}^N_t}\big)\big\|^2.\end{align*}
Then 
\begin{align*}
    \E\big\langle \overline{X}{}^i_t-X^i_t,b(\overline{X}{}^i_t,f_t)-b\big(X^i_t,\mu_{\mathcal{X}^N_t}\big)\big\rangle &= \E\big\langle \overline{X}{}^i_t-X^i_t,b(\overline{X}{}^i_t,f_t)-b\big(\overline{X}{}^i_t,\mu_{\overline{\mathcal{X}}{}^N_t}\big)\big\rangle \\ 
    &\quad+\E\big\langle \overline{X}{}^i_t-X^i_t,b\big(\overline{X}{}^i_t,\mu_{\overline{\mathcal{X}}{}^N_t}\big)-b\big(X^i_t,\mu_{\mathcal{X}^N_t}\big)\big\rangle
\end{align*}
Using Cauchy-Schwarz inequality with the same classical argument as before but replacing the boundedness of $K_1$ by \cref{eq:hypL2K} gives:
\[\E\big\langle \overline{X}{}^i_t-X^i_t,b(\overline{X}{}^i_t,f_t)-b\big(\overline{X}{}^i_t,\mu_{\overline{\mathcal{X}}{}^N_t}\big)\big\rangle \leq {\left(\E|\overline{X}{}^i_t-X^i_t|^2\right)}^{1/2}\frac{C}{\sqrt{N}}.\]
Then,
\begin{align*}
    &\E\big\langle \overline{X}{}^i_t-X^i_t,b\big(\overline{X}{}^i_t,\mu_{\overline{\mathcal{X}}{}^N_t}\big)-b\big(X^i_t,\mu_{\mathcal{X}^N_t}\big)\big\rangle \\
    &=\frac{1}{N}\sum_{i=1}^N \E\big\langle \overline{X}{}^i_t-X^i_t,b\big(\overline{X}{}^i_t,\mu_{\overline{\mathcal{X}}{}^N_t}\big)-b\big(X^i_t,\mu_{\mathcal{X}^N_t}\big)\big\rangle\\
    &\leq \frac{1}{N}\sum_{i=1}^N\E|\overline{X}{}^i_t-X^i_t|\big|b\big(\overline{X}{}^i_t,\mu_{\overline{\mathcal{X}}{}^N_t}\big)-b\big(X^i_t,\mu_{\mathcal{X}^N_t}\big)\big|\\
    &\leq \frac{\|\tilde{b}\|_\mathrm{Lip}}{N^2} \sum_{i,j=1}^N \E|\overline{X}{}^i_t-X^i_t||K_1(\overline{X}{}^i_t,\overline{X}{}^j_t)-K_1(X^i_t,X^j_t)| + \frac{\|\tilde{b}\|_\mathrm{Lip}}{N}\sum_{i=1}^N \E|\overline{X}^i_t-X^i_t|^2\\
    &\leq \frac{C\gamma\|\tilde{b}\|_\mathrm{Lip}}{N^2}\sum_{i,j=1}^N  \E\Big[\Big(|\overline{X}{}^i_t-X^i_t|^2+|\overline{X}{}^i_t-X^i_t||\overline{X}{}^j_t-X^j_t|\Big)\\
    &\phantom{\leq \frac{\gamma\|\tilde{b}\|_\mathrm{Lip}}{N^2}\sum_{i,j=1}^N  \E\Big[\Big(|\overline{X}{}^i_t-X^i_t|^2+|\overline{X}{}^i_t-X^i_t|}\times\Big(1+|\overline{X}{}^i_t|^p+|\overline{X}{}^j_t|^p+|X^i_t|^p+|X^j_t|^p\Big)\Big]\\
    &\quad =: \frac{C\gamma\|\tilde{b}\|_\mathrm{Lip}}{N^2}\sum_{i,j=1}^N \E[I_{ij}]
\end{align*}
For a given $(i,j)$ and $R>0$, the authors of \cite{bolley_stochastic_2011} define the event 
\[\mathcal{R} = \Big\{|\overline{X}{}^i_t|\leq R,|\overline{X}{}^j_t|\leq R,|{X}^i_t|\leq R,|{X}^j_t|\leq R\Big\}.\]
Then they distinguish the two cases inside the expectation:
\begin{align*}
\E[I_{ij}] &= \E[\1_{\mathcal{R}}I_{ij}] +  \E[\1_{\mathcal{R}^c}I_{ij}] \\ 
&\leq C(1+4R^p)\E|\overline{X}{}^i_t-X^i_t|^2 + \E[\1_{\mathcal{R}^c}I_{ij}]\\
&\leq C(1+4R^p)\E|\overline{X}{}^i_t-X^i_t|^2 \\
&\qquad + \left(\E[\1_{\mathcal{R}^c}]\right)^{1/2}\left(\E\Big[\big(1+|\overline{X}{}^i_t|^p+|\overline{X}{}^j_t|^p+|X^i_t|^p+|X^j_t|^p\big)^2\Big]\right)^{1/2}
\end{align*}
The probability of $\mathcal{R}^c$ is controlled by the Markov inequality, 
\begin{align*}
    \E[\1_{\mathcal{R}^c}] &\leq \E[\1_{|\overline{X}{}^i_t|>R}]+\E[\1_{|\overline{X}{}^j_t|>R}]+\E[\1_{|{X}^i_t|>R}]+\E[\1_{|{X}^j_t|>R}]\\
    &\leq C\e^{-\kappa R^{p'}}
\end{align*}
Setting $r=\kappa^{p/p'} R^p/2^{p/p'}$, it follows that
\[\E[I_{ij}]\leq C(1+r)\E|\overline{X}{}^i_t-X^i_t|^2+C\e^{-r^{p'/p}}.\]
A similar reasoning applies for the term with $\sigma$ and therefore, the function 
\[y(t) := \E|\overline{X}{}^i_t-X^i_t|^2,\]
satisfies, for every $r>0$:
\[y'(t) \leq C(1+r)y(t)+\e^{-r^{p'/p}}+\frac{C}{\sqrt{N}}\sqrt{y(t)} \leq C(1+r)y(t)+C\e^{-r^{p'/p}}+\frac{C}{N}.\]
If $p=p'$, choosing $r = \log ( 1 + 1/y(t) )$ leads to the nice differential inequality
\[y'(t) \leq Cy(t)+Cy(t) \log ( 1 + 1/y(t) ) + \frac{C}{N}, \]
and a complicated Gronwall-like argument (see e.g. \cite[Lemma 5.2.1]{chemin_fluides_1995} or \cite[Theorem 27]{dragomir_gronwall_2003}) terminates the proof. Otherwise when $p'>p$, choose $r = (\log N)^{p/p'}$, and since $y(0) = 0$, a direct integration by the classical Gronwall lemma gives
\[ y(t) \leq 2 \,\e^{CT + C T (\log N)^{p/p'} - \log N } ,\]
which concludes.
\end{proof}

The authors of \cite{bolley_stochastic_2011} write a detailed proof in the kinetic case with 
\[b(x,v,\mu) = -F(x,v) - H\star\mu(x,v),\quad \sigma(x,v,\mu) = \sqrt{2}I_d,\]
where $F,H:\R^d\times\R^d\to\R^d$ satisfy a slightly weaker assumption, namely: 
\[-\langle v-w,F(x,v)-F(x,w)\rangle \leq \gamma_1|v-w|^2,\]
and
\[|F(x,v)-F(y,v)|\leq\gamma_2\min(1,|x-y|)(1+|v|^p),\]
and similarly for $H$. They also prove \cite[Theorem 1.2]{bolley_stochastic_2011} which gives sufficient conditions on $F$ and $G$ for the well-posedness of both the particle and the nonlinear systems and such that the hypotheses of Theorem \cref{thm:mckeanmoment} are satisfied. Theorem \cref{thm:mckeanmoment} corresponds to a combination of the variant (V3), of the case given in Section 1.2.2 and of the case given in Section 1.2.3 of \cite[Theorem 1.1]{bolley_stochastic_2011}.

\subsubsection*{Moderate interaction.} In \cite{oelschlager_law_1985}, Oelschl\"ager introduced the concept of \emph{moderately interacting particles}. He studied systems of the form \cref{eq:mckeanthmassum} with a constant diffusion matrix $\sigma\equiv \sqrt{2}I_d$ and with a symmetric interaction kernel $K_1$ which \emph{depends on} $N$ as follows: \begin{equation}\label[IIeq]{eq:moderateinteractionK}\forall x,y\in\R^d,\quad K_1(x,y)\equiv K_1^N(y-x) := \frac{1}{\varepsilon_N^d}K_0\left(\frac{y-x}{\varepsilon_N}\right),\end{equation}
where $K_0:\R^d\to\R$ is a fixed symmetric radial kernel and $(\varepsilon_N)_N$ is a sequence such that $\varepsilon_N\to 0$ as $N\to+\infty$. The strength of the interaction between two particles is thus of the order $\sim \varepsilon_N^{-d}N^{-1}$. Oelschl\"ager considered the case $\varepsilon_N=N^{-\beta/d}$ with $\beta\in(0,1)$. The two extreme cases $\beta=0$ and $\beta=1$ correspond respectively to a \emph{weak interaction} of order $\sim1/N$ (actually what is usually called the mean-field scaling) and a \emph{strong interaction} of order $\sim1$ (it would be hopeless to take the limit $N\to+\infty$ in this case without further assumptions, see Section \cref{sec:Boltzmannclassicalmodels}). More generally, the term \emph{moderate interaction} refers to any situation in which $\varepsilon_N\to0$ and  $\varepsilon^{-d}_NN^{-1}=o(1)$. In this case
\[K_1^N(x,\cdot) \underset{N\to+\infty}{\longrightarrow} \delta_x,\]
in the distributional sense, which allows to recover singular \emph{purely local} interactions. 

When the diffusion matrix $\sigma\equiv\sqrt{2}I_d$ is constant, the main result of \cite[Theorem~1]{oelschlager_law_1985} is a functional law of large numbers which states the convergence of the empirical measure valued process ${\big(\mu_{\mathcal{X}^N_t}\big)}_t$ towards the deterministic singular limit $f_t$ solution of 
\[\partial_t f_t(x) = -\nabla_x\cdot\big\{\tilde{b}(x,f_t(x))f_t(x)\big\}+\Delta_x f_t.\]
We call this interaction purely local because the drift term no longer depends on the convolution $K_1\star f_t(x)$ but only on the local quantity $f_t(x)$. The strategy is roughly the same as the one explained in Section \cref{sec:martingalecompactness}. The first step is a relative compactness result in $\pb(C([0,T],\pb(\R^d)))$, the second step is the identification of the limit process which is shown to be almost surely the solution of a deterministic equation. The last step and in this case, the most difficult one, is the uniqueness of the solution of this deterministic equation. In the case of a gradient system, well-posedness results in some H\"older spaces are available in the PDE literature~\cite{ladyzenskaja_linear_1968}. 

Later, Oelschl\"ager studied the fluctuations around the limit \cite{oelschlager_fluctuation_1987} and applied these results to a multi-species reaction-diffusion system \cite{oelschlager_derivation_1989}. A pathwise extension of Oelschl\"ager's results can be found in \cite{meleard_propagation_1987}. 

The martingale approach of \cite{oelschlager_law_1985} is very restricted to the case when the diffusion matrix is equal to the identity. In the general case \cref{eq:mckeanthmassum}, the problem is revisited in \cite{jourdain_propagation_1998}. The approach is based on a careful control of the convergence rate in McKean's theorem and \emph{ad hoc} well-posedness results for the limiting purely local equation~\cref{eq:nonlinearpdemoderateiniteraction}. First note that the $L^\infty$ and Lipschitz norms of $K_1^N$ are controlled by 
\[\|K_1^N\|_\infty = \frac{C_0}{\varepsilon_N^d},\quad \|K_1^N\|_\mathrm{Lip} = \frac{C_1}{\varepsilon_N^{d+1}},\]
for some constants $C_0,C_1>0$ depending on $K_0$. We also assume that $K_2$ is of the form \cref{eq:moderateinteractionK} (possibly with another $K_0$). Thus, McKean's theorem gives for all $N$ an estimate of the form 
\begin{equation}\label[IIeq]{eq:mckeanmoderateinteraction}\E{\left[\sup_{t\leq T}|X^{i,N}_t-\overline{X}{}^{i,N}_t|^2\right]} \leq \tilde{c}_1\frac{\varepsilon_N^{-2d}}{N}\exp{\left(\tilde{c}_2\varepsilon_N^{-2(d+1)}\right)},\end{equation}
for some constants $\tilde{c}_1,\tilde{c}_2>0$ depending only on $T,K_0,\tilde{b}$ and $\tilde{\sigma}$ and where $\overline{X}{}^{i,N}_t$ satisfies 
\[\dd\overline{X}{}^{i,N}_t = \tilde{b}{\left(\overline{X}{}^{i,N}_t,K_1^N\star f_{t}^{(N)}{\left(\overline{X}{}^{i,N}_t\right)}\right)}\dd t + \tilde{\sigma}{\left(\overline{X}{}^{i,N}_t,K_2^N\star f_t^{(N)}{\left(\overline{X}{}^{i,N}_t\right)}\right)}\dd B^i_t,\]
with $f^{(N)}_t\in \pb(\R^d)$ is the law of $\overline{X}{}^{i,N}_t$. It satisfies 
\begin{multline}\label[IIeq]{eq:nonlinearpdeNmoderateiniteraction}
\partial_t f^{(N)}_t(x) = -\nabla_x\cdot{\left\{\tilde{b}{\left(x,K_1^N\star f^{(N)}_t(x)\right)}f^{(N)}_t(x)\right\}}\\
+\frac{1}{2}\sum_{i,j=1}^d \partial_{x_i}\partial_{x_j}{\left\{a_{ij}{\left(x,K_2^N\star f^{(N)}_t(x)\right)}f^{(N)}_t(x)\right\}},
\end{multline}
where $a\equiv(a_{ij}) = \tilde{\sigma}\tilde{\sigma}^\mathrm{T}$. In order to take $N\to+\infty$ in \cref{eq:mckeanmoderateinteraction}, Jourdain and M\'el\'eard~\cite{jourdain_propagation_1998} assume that $\varepsilon_N\to0$ slowly enough so that the right-hand side of~\cref{eq:mckeanmoderateinteraction} still converges to zero. A sufficient condition is 
\begin{equation}\label[IIeq]{eq:epsilonNsufficientmoderateinteraction}\varepsilon_N^{-2(d+1)} \leq \delta\log N,\end{equation}
for a small $\delta>0$. In the bound \cref{eq:mckeanmoderateinteraction}, the nonlinear process $(\overline{X}{}^{i,N}_t)^{}_t$ still depends on $N$ (through $K_1^N$ and $K_2^N$) so it is not possible to simply take the limit $N\to+\infty$. Moreover, the goal is to prove propagation of chaos towards the solution $f_t$ of the purely local PDE:
\begin{equation}\label[IIeq]{eq:nonlinearpdemoderateiniteraction}
\partial_t f_t(x) = -\nabla_x\cdot{\left\{\tilde{b}\left(x,f_t(x)\right)f_t(x)\right\}}
+\frac{1}{2}\sum_{i,j=1}^d \partial_{x_i}\partial_{x_j}{\left\{a_{ij}{\left(x,f_t(x)\right)}f_t(x)\right\}}.
\end{equation}
Well-posedness results for the PDEs \cref{eq:nonlinearpdeNmoderateiniteraction} and \cref{eq:nonlinearpdemoderateiniteraction} can be found in \cite[Section 1]{jourdain_propagation_1998}. The approach of \cite{jourdain_propagation_1998} is based on the work of \cite{ladyzenskaja_linear_1968} on parabolic PDEs. The main assumptions are the regularity of the drift and diffusion coefficients (respectively at least $C^2$ and $C^3$) and of the initial condition (at least $C^2$ with an H\"older continuous second order derivative), together with the following non-negativity assumption on~$a$: 
\[\forall x\in\R^d,\forall z\in\R^d,\forall p\in\R,\quad \langle x, (a'(z,p)p+a(z,p))x\rangle \geq 0,\]
where for $z\in\R^d$, $a'(z,p)$ denotes the derivative of $p\in\R \mapsto a(z,p)\in\mathcal{M}_d(\R)$. Then, \cite[Proposition 2.5]{jourdain_propagation_1998} shows that \cref{eq:nonlinearpdemoderateiniteraction} is well-posed, that the associated nonlinear SDE is well-posed and that the solution $\overline{X}{}^i_t$ of
\[\dd\overline{X}{}^{i}_t = \tilde{b}\Big(\overline{X}{}^{i}_t, f_{t}\big(\overline{X}{}^{i}_t\big)\Big)\dd t + \sigma\Big(\overline{X}{}^{i}_t,f_t\big(\overline{X}{}^{i}_t\big)\Big)\dd B^i_t,\]
satisfies: 
\begin{equation}\label[IIeq]{eq:nonlinearprocessmoderateinteraction}\E{\left[\sup_{t\leq T}|\overline{X}{}^{i,N}_t-\overline{X}{}^i_t|^4\right]}\leq C\varepsilon_{N}^{4\beta},\end{equation}
for some $\beta>0$. The proof of this proposition is based on PDE arguments. In particular, since the law of $\overline{X}{}^{i,N}_t$ solves \cref{eq:nonlinearpdeNmoderateiniteraction}, using Ascoli's theorem (or other compactness criteria) it is possible to extract a convergent subsequence $f^{(N)}_t\to f_t$ where $f_t$ solves \cref{eq:nonlinearpdemoderateiniteraction} with an explicit convergence rate. Combining \cref{eq:mckeanmoderateinteraction} and \cref{eq:nonlinearprocessmoderateinteraction} leads to 
\[\E{\left[\sup_{t\leq T}|X^{i,N}_t-\overline{X}_t|^2\right]} \leq {\left(C\varepsilon_N^{2\beta}+\tilde{c}_1\frac{\varepsilon_N^{-2d}}{N}\exp{\left(\tilde{c}_2\varepsilon_N^{-2(d+1)}\right)}\right)},\]
and the conclusion follows as soon as $\varepsilon_N$ satisfies \cref{eq:epsilonNsufficientmoderateinteraction}. 

Recent applications of these results can be found in \cite{chen_rigorous_2021} and \cite{diez_propagation_2020}. The reference~\cite{chen_rigorous_2021} presents a generalisation of \cite{jourdain_propagation_1998} to a multi-species system with non globally Lipschitz interactions. The article contains very detailed well-posedness results for the different systems involved. In \cite{diez_propagation_2020}, the diffusion process is replaced by a Piecewise Deterministic process on a (compact) manifold.  

\subsubsection*{Singular interactions with cutoff}

Similarly to the moderate interaction case where the goal was to approximate purely local interactions, it is possible to introduce a cutoff parameter which depends on $N$ in order to approximate interactions which do not satisfy the regularity hypotheses of McKean's theorem. The most important cases in the literature are the singular Coulomb-type interactions. These models are given either by the first order model 
\begin{equation}\label[IIeq]{eq:kellersegelparticles}\dd X^i_t = F\star \mu_{\mathcal{X}^N_t}\dd t + \sigma\dd B^i_t,\end{equation}
or by the second order kinetic system 
\begin{equation}\label[IIeq]{eq:kineticcoulombparticles}\dd X^i_t = V^i_t \dd t, \quad \dd V^i_t = \frac{1}{N}\sum_{j=1}^N F(X^i_t-X^j_t) \dd t + \sigma\dd B^i_t,  \end{equation}
where in both cases $\sigma>0$ and $F$ is a Coulomb-type force:  
\[F(x) = \xi_d\frac{x}{|x|^d},\]
with a constant $\xi_d\in \R$. This force is singular at the origin and derives from a potential $F=-\nabla_x \Phi$ with 
\begin{align*}
\Phi(x) &= \frac{\xi_d}{2}\log|x|,\quad d=2\\
\Phi(x) &= \frac{\xi_d}{(d-2)|x|^{d-2}},\quad d\geq3.
\end{align*}
In the attractive case $\xi_d>0$, the kinetic systems corresponds to the classical gravitational Newtonian dynamics. More recently, the first order system has been used to model the chemotaxis interactions between swarms of bacteria. Using the classical notations in this context and with an appropriate constant $\xi_d>0$, the formal limit $N\to+\infty$ leads to the propagation of chaos towards a limit distribution $\rho$ which satisfies the system  
\begin{subequations}\label[IIeq]{eq:parabolicelliptickellersegel}
\begin{align}
\partial_t \rho &= - \nabla \cdot(\rho\nabla c) + \frac{\sigma^2}{2}\Delta \rho,\\
-\Delta c &= \rho. \label[IIeq]{eq:parabolicelliptickellersegel_c}
\end{align}
\end{subequations}
This system is called the parabolic-elliptic Keller-Segel system. The density $\rho$ represents the (spatial) density of bacteria and $c$ the concentration of a chemical substance which is secreted by the bacteria and whose gradient drives the motion of the other bacteria. 

The repulsive case $\xi_d<0$ is classically used in plasma physics. Importantly, the kinetic case with $\sigma=0$  has also motivated the first propagation of chaos results (in a more regular setting)  by Braun and Heppp \cite{braun_vlasov_1977} and Dobrushin \cite{dobrushin_vlasov_1979}, as the formal limit equation is the renowned Vlasov equation
\begin{align*}
    &\partial_t f_t + v\cdot \nabla_x f_t + (F\star \rho_t)\cdot \nabla_v f_t = 0\\
    &\rho_t(x) = \int_{\R^d} f_t(x,v)\dd v.
\end{align*}
For a detailed review of classical and recent propagation results in this deterministic context (but possibly with random initial conditions), we refer the interested reader to the review articles \cite{jabin_review_2014, golse_dynamics_2016}, to the articles \cite{hauray_n-particles_2007,hauray_particles_2015,boers_mean_2016,lazarovici_mean_2017} or, with different techniques, to the recent articles \cite{jabin_mean_2016,serfaty_systems_2019,serfaty_mean_2020}.

In all these cases, it is not possible to directly apply McKean's theorem due to the singularity at the origin. Moreover, it has been shown by \cite[Proposition 4]{fournier_stochastic_2017} that, for the particle system associated to the Keller-Segel system, the singularity is indeed visited with nonzero probability: for any $N\geq1$, any $\xi_d>0$ and any $t_0>0$, then any solution of \cref{eq:kellersegelparticles} (if it exists) satisfies 
\[\mathbb{P}\big(\exists s\in[0,t_0],\,\,\exists i\ne j,\,\,X^i_s=X^j_s\big)>0.\]
Consequently, this raises well-posedness issues already at the particle level. A natural strategy is therefore to remove the singularity at the origin by regularizing the force using a cutoff parameter $\varepsilon>0$. For instance, in \cite{carrillo_propagation_2019}, the authors replace the force $F$ by the regularized version
\begin{equation}\label[IIeq]{eq:cutoffcarrillo}F_\varepsilon(x) = \xi_d\frac{x}{\max(|x|,\varepsilon)^d}.\end{equation}
In \cite{liu_propagation_2019,fetecau_propagation_2019}, the authors define a regularized potential $\Phi_\varepsilon = J_\varepsilon\star\Phi$ where $J_\varepsilon(x) = \varepsilon^{-d}J(\varepsilon^{-1}x)$ and $J$ is a smooth mollifier. Once a regularized particle system is defined, it is also possible to define the associated synchronously coupled system of nonlinear SDEs (which depend on the cutoff parameter). Since coupling methods typically give quantitative results, it becomes possible to take a cutoff parameter $\varepsilon\equiv\varepsilon_N$ which depends on $N$ and vanishes as $N\to+\infty$. Similarly to the moderate interaction case, it is necessary to obtain beforehand precise well-posedness results and the sharpness of the estimates will determine the size of the cutoff. In \cite{carrillo_propagation_2019}, the authors extend previous results by \cite{lazarovici_mean_2017} in the deterministic case (but starting from a random initial condition) and prove the propagation of chaos for the regularized system \cref{eq:kineticcoulombparticles} with the cutoff \cref{eq:cutoffcarrillo} and a cutoff size $\varepsilon\equiv \varepsilon_N = N^{-\delta}$, $\delta<1/d$. This result has been improved in \cite{huang_mean-field_2020}. In \cite{liu_propagation_2019,fetecau_propagation_2019}, the authors use a smooth mollifier for the regularized system \cref{eq:kellersegelparticles} with a cutoff size $\varepsilon\equiv \varepsilon_N = (\log N)^{-1/d}$.

Finally, it should be noted that there has been many recent advances on this issue using techniques that are not always based on cutoff approximations. These results will be further discussed in Section \cref{sec:jabin}, Section \cref{sec:chaosviagirsanov} and Section \cref{sec:vorticity}.

\subsubsection{Gradient systems and uniform in time estimates}\label[II]{sec:gradientsystems}

In this section, the case of McKean-Vlasov gradient systems is investigated, that is systems of the form \cref{eq:mckeanvlasov_summary} with
\begin{equation}\label[IIeq]{eq:gradientsystem}b(x,\mu) = -\nabla V(x)-\nabla W\star\mu(x),\quad \sigma(x,\mu)\equiv\sigma I_d,\quad \sigma>0,\end{equation}
where $V,W:\R^d\to\R$ are two twice continuously differentiable potentials (usually symmetric, but it will be precised each time), respectively called the confinement potential and the interaction potential. The law of the corresponding nonlinear McKean-Vlasov process satisfies the famous granular media equation:
\begin{equation}\label[IIeq]{eq:granularmedia}\partial f_t = \frac{\sigma^2}{2}\Delta f_t + \nabla\cdot(f_t\nabla(V+W\star f_t)).\end{equation}
For the modelling details, we refer the reader to
\cite{benedetto_kinetic_1997,benedetto_non-maxwellian_1998} who first derived this equation. The granular media equation has been studied analytically in \cite{carrillo_kinetic_2003,carrillo_contractions_2006} and later in \cite{bolley_uniform_2013}. The fundamental question, which also motivates this section, is the long-time asymptotic of the solution, in particular the existence of stationary solutions and the convergence to equilibrium. The probabilistic counterpart of the granular media equation is the nonlinear McKean-Vlasov process \cref{eq:mckeanvlasov-limit_summary} with $b,\sigma$ given by \cref{eq:gradientsystem}. The long-time behaviour of this process is not simpler than the direct study of \cref{eq:granularmedia} but this probabilistic approach strongly suggests to consider the (linear) McKean particle system \cref{eq:mckeanvlasov_summary} as a starting point, the idea being to replace the nonlinearity in dimension $d$ by a linear system of particles in high dimension $dN$. Since the behaviour of linear diffusion systems is well-established, this may be simpler provided that it is possible to prove convergence results with rates independent of the dimension. In a series of works reviewed in this section, it has been shown that quantitative convergence to equilibrium for the nonlinear system may follow from the study of the particle system. The crucial result is the uniform in time propagation of chaos. In this section we review some results in this sense under various convexity assumptions on the potentials. Note that uniform in time propagation of chaos is strongly linked to the uniqueness of a stationary measure for the nonlinear process. Uniform in time propagation of chaos may not hold as soon as the nonlinear system admits more than one stationary measure (in the cases studied below, this is a consequence of the fact that the particle system admits a unique equilibrium). In general uniform in time propagation of chaos and the existence of a unique stationary measure for the nonlinear process hold simultaneously. We start by stating the main theorem of this section which is due to Malrieu \cite{malrieu_logarithmic_2001}. 
 
\begin{theorem}[Uniform in time propagation of chaos \cite{malrieu_logarithmic_2001}]\label[II]{thm:uniformpocuniformconvex}
Let $\mathcal{X}^N_t$ be the particle system \cref{eq:mckeanparticlessynch} and $\overline{\mathcal{X}}^N_t$ be the synchronously coupled nonlinear system~\cref{eq:nonlinearmckeanparticlessynch}. Let $b,\sigma$ be given by~\cref{eq:gradientsystem}, and let us assume that there exist $\beta > 0$ and $p \geq 1$ such that $V,W$ satisfy the following properties.
\begin{itemize}
    \item $V$ is $\beta$-uniformly convex: 
    \[\forall x,y\in\R^d,\quad \langle x-y,\nabla V(x)-\nabla V(y)\rangle\geq \beta |x-y|^2.\]
    \item $W$ is symmetric and convex:
        \[\forall x,y\in\R^d,\quad \langle x-y,\nabla W(x)-\nabla W(y)\rangle\geq 0.\]
    \item $\nabla W$ is locally Lipschitz and has polynomial growth of order $p$.
\end{itemize}
Let the initial law $f_0\in\pb_{2p}(\R^d)$ have bounded moments of order $2p$. Then there exists a constant $C>0$ depending only on $\beta$ and $p$ such that 
\begin{equation}\label[IIeq]{eq:uniformpocuniformconvex}\sup_{t\geq0} \frac{1}{N}\sum_{i=1}^N \E|X^i_t-\overline{X}{}^i_t|^2\leq \frac{C}{N}.\end{equation}
\end{theorem}

All the well-posedness results for both the particle system and the nonlinear process are proved in \cite[Section 2]{cattiaux_probabilistic_2008}. The proof of Theorem \cref{thm:uniformpocuniformconvex} is given below. This is an extension of Sznitman's proof of McKean's theorem by synchronous coupling to the case of unbounded interactions. In a one-dimensional setting, a similar result is proved in \cite[Theorem 3.1]{benachour_nonlinear_1998-1}. It has been adapted to the current setting in \cite[Theorem 3.3]{malrieu_logarithmic_2001}. The uniform convexity of $V$ allows a uniform in time control of the trajectories. To deal with unbounded interactions, the following lemma will be needed to control the moments of the nonlinear system uniformly in time (see also \cite[Proposition 3.10]{benachour_nonlinear_1998-1} and \cite[Corollary 2.3, Proposition 2.7]{cattiaux_probabilistic_2008}). 

\begin{lemma}[Moment bound]\label[II]{lemma:uniformpocmoment} Let $(\overline{X}_t)_t$ be the nonlinear McKean-Vlasov process~\cref{eq:mckeanvlasov-limit_summary} with $b$ and $\sigma$ given by \cref{eq:gradientsystem}. Let $V$ be $\beta$-uniformly convex for a constant $\beta>0$ and let $W$ be symmetric and convex. Then for every $p \geq 1$ such that $f_0\in\pb_{2p}(\R^d)$, it holds that
\[\sup_{t\geq 0}\E|\overline{X}_t|^{2p} < +\infty.\]
\end{lemma}

\begin{proof}
It\={o}'s formula gives: 
\begin{align*}
    |\overline{X}_t|^{2p} &= |\overline{X}_0|^{2p} + 2p\int_0^t \big\langle |\overline{X}_s|^{2(p-1)}\overline{X}_s,-\nabla V(\overline{X}_s)-\nabla W\star f_s(\overline{X}_s)\big\rangle \dd s\\
    &\quad+\sigma^2 dp\int_0^t |\overline{X}_s|^{2(p-1)}\dd s + 2p(p-1)\sigma^2\int_0^t |\overline{X}_s|^{2(p-1)}\dd s \\ 
    &\quad + \sigma\int_0^t\big\langle 2p|\overline{X}_s|^{2(p-1)}\overline{X}_s,\dd B_s\big\rangle.
\end{align*}
Taking the expectation and then the time-derivative
\begin{align*}
\frac{\dd}{\dd t} \E|\overline{X}_t|^{2p}&=
-2p \E{\left[ |\overline{X}_t|^{2(p-1)} \big\langle \overline{X}_t - 0,\nabla V(\overline{X}_t)-\nabla V(0)\big\rangle\right]}\\
&\quad-2p\E{\left[ |\overline{X}_t|^{2(p-1)} \big\langle \overline{X}_t,\nabla V(0)\big\rangle\right]}\\
&\quad+p\sigma^2(d+2(p-1)) \E|\overline{X}_t|^{2(p-1)}\\
&\quad-2p \E\big\langle |\overline{X}_t|^{2(p-1)}\overline{X}_t,\nabla W\star f_t(\overline{X}_t)\big\rangle \\ 
&\leq-2p \beta \E |\overline{X}_t|^{2p}+2p|\nabla V(0)|\E|\overline{X}_t|^{2p-1}\\
&\quad+p\sigma^2(d+2(p-1)) \E|\overline{X}_t|^{2(p-1)}\\
&\quad-2p \E\big\langle |\overline{X}_t|^{2(p-1)}\overline{X}_t,\nabla W\star f_t(\overline{X}_t)\big\rangle.
\end{align*}
where the inequality follows from the uniform convexity of $V$.
Let now $(\overline{Y}_t)_t$ be an independent copy of $(\overline{X}_t)_t$, so that
\[ \E\big\langle |\overline{X}_t|^{2(p-1)}\overline{X}_t,\nabla W\star f_t(\overline{X}_t)\big\rangle = \E\big\langle |\overline{X}_t|^{2(p-1)}\overline{X}_t,\nabla W(\overline{X}_t-\overline{Y}_t)\big\rangle. \]
Since $W$ is symmetric, $\nabla W$ is odd, leading to
\[ \E\big\langle |\overline{X}_t|^{2(p-1)}\overline{X}_t,\nabla W\star f_t(\overline{X}_t)\big\rangle = - \E\big\langle |\overline{X}_t|^{2(p-1)}\overline{Y}_t,\nabla W(\overline{X}_t-\overline{Y}_t)\big\rangle, \]
using the fact $\overline{X}_t$ and $\overline{Y}_t$ are independent and have the same law. Summing these two expressions and using the convexity of $W$ gives
\begin{align*}
    \E\big\langle |\overline{X}_s|^{2(p-1)}\overline{X}_s,\nabla W\star f_s(\overline{X}_s)\big\rangle &= \frac{1}{2}\E\big\langle |\overline{X}_s|^{2(p-1)}(\overline{X}_s-\overline{Y}_s),\nabla W(\overline{X}_s-\overline{Y}_s)\big\rangle\\
    &\geq 0.
\end{align*}
Denote the moment of order $2p$ by $\mu_{2p}(t) := \E|\overline{X}_t|^{2p}$. Then it holds that
\[ \frac{\dd}{\dd t} \mu_{2p}(t) \leq -\lambda(p) \mu_{2p}(t) +c_1(p) \mu_{2(p-1)}(t) + c_2(p)\mu_{2p-1}(t),\]
where $\lambda(p),c_1(p),c_2(p)>0$ only depend on $p$ and $\beta$. Since $\mu_{2p}(0)=0$, the conclusion follows by integrating this differential Gronwall-like inequality, after noticing that for all $\varepsilon>0$ and for all exponent $q< 2p$, there exists a constant $K>0$ such that for all $x\in \R^d$,
\[|x|^{q}\leq K+\varepsilon |x|^{2p}.\]
\end{proof}

\begin{proof}[Proof (of Theorem \cref{thm:uniformpocuniformconvex})] The proof proceeds similarly as in Sznitman's approach but the convexity assumptions are used to get a better uniform in time control of the trajectories. The starting point is It\={o}'s formula: 
\begin{align*}
    |X^i_t-\overline{X}{}^i_t|^2 &= -2\int_0^t \big\langle X^i_s-\overline{X}{}^i_s,\nabla V(X^i_s)-\nabla V(\overline{X}{}^i_s)\big\rangle\dd s \\
    &\quad -2 \int_0^t\big\langle X^i_t-\overline{X}{}^i_s, \nabla W\star\mu_{\mathcal{X}^N_s}(X^i_s)-\nabla W\star f_t(\overline{X}{}^i_s)\big\rangle\dd s
\end{align*}
Taking the expectation, and then differentiating
\begin{align*}    
    \frac{\dd}{\dd t} \E |X^i_t-\overline{X}{}^i_t|^2 
    &\leq -2\beta \E |X^i_t-\overline{X}{}^i_t|^2 \\
    &\quad -2 \E \big\langle X^i_t-\overline{X}{}^i_t, \nabla W\star \mu_{\mathcal{X}^N_t}(X^i_t)-\nabla W\star \mu_{\overline{\mathcal{X}}{}^N_t}(\overline{X}{}^i_t)\big\rangle \\
    &\quad - 2 \E \big\langle X^i_t-\overline{X}{}^i_t, \nabla W\star \mu_{\overline{\mathcal{X}}{}^N_t}(\overline{X}{}^i_t)-\nabla W\star f_t(\overline{X}{}^i_t)\big\rangle,
\end{align*}
where the uniform convexity assumption on $V$ is used and the introduction of the term $\nabla W\star \mu_{\overline{\mathcal{X}}{}^N_s}$ in the second term on the right-hand side is forced as in the proof of McKean's theorem. For the second term on the right-hand side of the last inequality, summing over $i$ leads to 
\begin{align*}
    &\sum_{i=1}^N \E \big\langle X^i_t-\overline{X}{}^i_t, \nabla W\star \mu_{\mathcal{X}^N_t}(X^i_t)-\nabla W\star \mu_{\overline{\mathcal{X}}{}^N_t}(\overline{X}{}^i_t)\big\rangle \\
    &= \frac{1}{N}\sum_{i,j=1}^N \E \big\langle X^i_t-\overline{X}{}^i_t, \nabla W(X^i_t-X^j_t)-\nabla W(\overline{X}{}^i_t-\overline{X}{}^j_t)\big\rangle \\ 
    &=\frac{1}{N}\sum_{i\leq j} \E \big\langle (X^i_t-X^j_t)-(\overline{X}{}^i_t-\overline{X}{}^j_t),\nabla W(X^i_t-X^j_t)- \nabla W(\overline{X}{}^i_t-\overline{X}{}^j_t)\big\rangle \\
    &\geq 0,
\end{align*}
where the convexity and symmetry of $W$ are used. 

Then after summing over $i=1,\ldots,N$ and dividing by $N$, the Cauchy-Schwarz inequality for the last term gives: 
\begin{align*}
    \frac{\dd}{\dd t} \frac{1}{N}\sum_{i=1}^N \E|X^i_t-\overline{X}{}^i_t|^2 &\leq -\frac{2\beta}{N}\sum_{i=1}^N \E|X^i_t-\overline{X}{}^i_t|^2 \\
    &\qquad+\frac{2}{N}\sum_{i=1}^N \left(\E|X^i_t-\overline{X}{}^i_t|^2\right)^{1/2} r^i_t \dd t,
\end{align*}
where 
\[r^i_t = {\left(\E \big|\nabla W\star \mu_{\overline{\mathcal{X}}{}^N_t}(\overline{X}{}^i_t)-\nabla W\star f_t(\overline{X}{}^i_t)\big|^2\right)}^{1/2}.\]
As in Sznitman's proof, it holds that 
\begin{align*} |r^i_t|^2 = \E \big|\nabla W\star \mu_{\overline{\mathcal{X}}{}^N_t}(\overline{X}{}^i_t)-\nabla W\star f_t(\overline{X}{}^i_t)\big|^2 = \frac{C}{N^2} \sum_{j=1}^N \E \big| \nabla W(\overline{X}{}^i_t-\overline{X}{}^j_t)\big|^2,
\end{align*}
since the processes $\overline{X}{}^i$ are independent. Using the polynomial growth of $\nabla W$ and Lemma \cref{lemma:uniformpocmoment}, it follows that there exists a constant $C_p$ depending on $p$ only such that 
\[|r^i_t|^2 \leq \frac{C_p}{N}.\]
Finally, by exchangeability, it holds that
\begin{align*}
    \frac{\dd}{\dd t} \frac{1}{N}\sum_{i=1}^N \E|X^i_t-\overline{X}{}^i_t|^2 &\leq -\frac{2\beta}{N}\sum_{i=1}^N \E|X^i_t-\overline{X}{}^i_t|^2 \\&\quad+\frac{2C_p}{\sqrt{N}}\left(\frac{1}{N}\sum_{i=1}^N \E|X^i_t-\overline{X}{}^i_t|^2\right)^{1/2}.
\end{align*}
Thus, setting
\[y(t):=\left(\frac{1}{N}\sum_{i=1}^N \E|X^i_t-\overline{X}{}^i_t|^2\right)^{1/2},\]
it holds that
\[y'(t) \leq -\beta y(t)+\frac{C_p}{\sqrt{N}},\]
and since $y(0)=0$, the conclusion follows by integrating this Gronwall-like differential inequality. 
\end{proof}

As a corollary, we state the main application of this theorem which is the exponentially fast convergence to equilibrium of the nonlinear process. Once again, this result is proved in \cite{malrieu_logarithmic_2001}. 

\begin{corollary}\label[II]{coro:uniformpocuniformconvex} Let $f_0$, $V$ and $W$ satisfy the same assumptions as in Theorem \cref{thm:uniformpocuniformconvex}, with $\sigma=\sqrt{2}$ for simplicity. Let us also assume that $\int f_0\log f_0 <\infty$.
Then the following properties hold true 
\begin{enumerate}
    \item \textbf{(Entropic chaos).} There exists $C>0$ such that for every $N \geq 1$
    \[\sup_{t\geq 0} H(f^N_t | f_t^{\otimes N}) \leq C.\]
    \item \textbf{(Convergence to equilibrium for the nonlinear process).} There exists a unique $\mu_\infty\in\pb(\R^d)$ and a constant $C>0$ such that for all $t>0$,
    \[\|f_t-\mu_\infty\|_{\mathrm{TV}}\leq C\e^{-\beta t/2}.\]
Note that in both statements, the constant $C>0$ depends on $f_0$.
\end{enumerate}
\end{corollary}

\begin{proof}[Proof (sketch)]
\begin{enumerate}
    \item The first property is proved in \cite[Proposition 3.13]{malrieu_logarithmic_2001} and follows from a log-Sobolev inequality satisfied by $f_t$. More generally, it is possible to use the general bound given by Lemma \cref{lemma:computeH} with $\alpha = \frac{1}{2}$ in \cref{eq:computeH_summary}. Thanks to the Bakry-Emery criterion (Proposition \cref{prop:bakryemery}), it can be shown that there exists $\lambda>0$ such that $f_t$ (and thus $f_t^{\otimes N}$) satisfies LSI($\lambda$) for every $N \geq 1$, i.e.
    \[ -\frac{1}{4} I \left( f^N_t | f_t^{\otimes N} \right) \leq - \frac{\lambda}{2} H \left( f^N_t | f_t^{\otimes N} \right), \]
    see \cite[Proposition 3.12]{malrieu_logarithmic_2001}. Then the quantity,
    \[\frac{1}{N} \sum_{j=1}^N \nabla W \left( X^i_t - X^j_t \right) - \nabla W \star f_t \left( X^i_t \right) \] 
    is controlled by $\sum_{i=1}^N\E \left| X^i_t - \overline{X}{}^i_t \right|^2$ and the square moments of $\overline{X}{}^i_t$ which are both bounded uniformly in time by Theorem \cref{thm:uniformpocuniformconvex}. Reporting in \cref{eq:computeH_summary}, this eventually gives a constant $C > 0$ such that
    \[ \frac{\dd }{\dd t} H \left( f^N_t | f_t^{\otimes N} \right) \leq - \frac{\lambda}{2} H \left( f^N_t | f_t^{\otimes N} \right) + C \]
    and the conclusion follows by integrating this Gronwall-like differential inequality.

    \item This is the content of \cite[Theorem 3.18]{malrieu_logarithmic_2001}. The existence and uniqueness of $f_\infty$ is proved for instance in \cite[Theorem 2.2]{benedetto_non-maxwellian_1998}. To get a quantitative convergence bound, the idea is to introduce the particle system as a pivot: 
    \begin{equation}\label[IIeq]{eq:ergodicgradient3}\|f_t-f_\infty\|_{\mathrm{TV}} \leq\|f_t-f^{1,N}_t\|_{\mathrm{TV}}+\|f^{1,N}_t-\mu^{1,N}_\infty\|_\mathrm{TV}+\|\mu^{1,N}_\infty-f_\infty\|_\mathrm{TV},\end{equation}
    where $\mu^{1,N}_\infty$ is the first marginal of the probability measure $\mu^N_\infty$ with density
    \begin{equation}\label[IIeq]{eq:gibbsmeasureparticlesgradient}
    \mu^N_\infty(\dd\mathbf{x}) \propto \exp{\left(-\sum_{i=1}^N V(x^i)-\frac{1}{2N}\sum_{i,j=1}^N W(x^i-x^j)\right)}\dd\mathbf{x}.\end{equation}
    Note that $\mu^N_\infty$ is the unique invariant measure of the $N$-particle process. The first and third terms on the right-hand side of \cref{eq:ergodicgradient3} are bounded by $K/\sqrt{N}$ using the first property thanks to the Pinsker inequality. The second term on the right-hand side of \cref{eq:ergodicgradient3} is bounded by $K\sqrt{N}\e^{-\beta t}$ using a classical application of the Bakry-Emery criterion (Proposition \cref{prop:bakryemery}). Thus, 
    \[\|f_t-f_\infty\|_\mathrm{TV}\leq \frac{K}{\sqrt{N}}+K\sqrt{N}\e^{-\beta t},\]
    and the conclusion follows by taking $N$ of the order of $\e^{\beta t/2}$. 
\end{enumerate}
\end{proof}

\begin{remark}[Invariant measures and phase transitions]\label[II]{rem:phasetransitionslogsobuniformpoc}
    The fact that the $N$-particle system admits the unique invariant measure \cref{eq:gibbsmeasureparticlesgradient} for any choice of potentials $V,W$ is a very important and noticeable property. On the contrary, the limit equation~\cref{eq:granularmedia} may have more than one stationary solution and the unicity in Corollary \cref{coro:uniformpocuniformconvex} is ensured by the strong convexity assumptions on the potentials. When the limit equation has more than one stationary solution, the system is said to undergo a phase transition. In a recent work \cite{delgadino_phase_2021}, the relation between phase transitions, uniform in time propagation of chaos and log-Sobolev inequalities is explored for McKean-Vlasov gradient systems~\cref{eq:gradientsystem}. It is shown that in the absence of phase transition then uniform in time propagation of chaos is equivalent to the non degeneracy as $N\to+\infty$ of the constant in the log-Sobolev inequality satisfied by the Gibbs measure~\cref{eq:gibbsmeasureparticlesgradient} of the $N$-particle system. This work is based on the gradient-flow framework which will be discussed in Section \cref{sec:gradientflows}. 
\end{remark}

It is also possible to go beyond Theorem \cref{thm:uniformpocuniformconvex} and prove concentration inequalities by using log-Sobolev inequalities for the $N$-particle law with constants independent of $N$. These questions will be discussed in Section \cref{sec:concentrationineqgradient}.

The uniform convexity assumption is generally understood as too strong to cover cases of physical interest. Some extensions of Theorem \cref{thm:uniformpocuniformconvex} with weaker convexity assumptions are discussed below. 

\begin{enumerate}[(a)]

\item \textbf{No confinement.} The key assumption is the uniform convexity of $V$ (the confinement potential) which allows a uniform in time control of the trajectories. In \cite{carrillo_kinetic_2003}, the authors studied analytically the granular media equation which corresponds to the law of the nonlinear system when $V=0$. However, at the particle level, it has been shown in \cite{benachour_nonlinear_1998-1} and \cite[Section 4]{malrieu_logarithmic_2001}, \cite[Section 2]{malrieu_convergence_2003} that propagation of chaos does not hold uniformly in time. This is unfortunate as it annihilates any hope of studying the long-time behaviour of the nonlinear system with a probabilistic point of view as in the case when $V$ is uniformly convex. Nevertheless, Malrieu \cite{malrieu_convergence_2003} showed that uniform in time propagation of chaos does hold for the system defined by 
\begin{equation}\label[IIeq]{eq:projectedsystem}Y^i_t = X^i_t -\frac{1}{N}\sum_{j=1}^N X^j_t,\end{equation}
which is the projection of the particle system on the set 
\[\mathcal{M} := \left\{\mathbf{x}\in \R^N,\,\sum_{j=1}^N x^j = 0\right\}.\]
The proof proceeds similarly as before but requires the potential $W$ to be uniformly convex (and not only convex as in Theorem \cref{thm:uniformpocuniformconvex}). It also requires a uniform in time control of the moments of the nonlinear system, proved in dimension one in \cite[Proposition 3.10]{benachour_nonlinear_1998-1} and more generally in \cite[Lemma 5.2]{malrieu_convergence_2003}. Details can be found in \cite[Theorem 5.1]{malrieu_convergence_2003} as well as a probabilistic proof of the convergence to equilibrium for the granular media equation \cite[Theorem 6.2]{malrieu_convergence_2003}. 

\item \textbf{Non uniformly convex potentials.} In Theorem \cref{thm:uniformpocuniformconvex} and in the case $V=~0$ in \cite{malrieu_convergence_2003}, at least one of the potentials has to be uniformly convex. This condition is relaxed in \cite{cattiaux_probabilistic_2008} and replaced by the $C(A,\alpha)$-condition already introduced in \cite{carrillo_kinetic_2003}: there exists $A,\alpha>0$ such that for all $0<\varepsilon<1$, 
\[\forall x,y\in\R^d,\quad \langle x-y,\nabla W(x)-\nabla W(y)\rangle\geq A\varepsilon^{\alpha}(|x-y|^2-\varepsilon^2).\]
This condition is weaker than uniform convexity and includes important cases such as $W(x)=|x|^{2+\alpha}$. Uniform in time propagation of chaos holds either for the particle system $X^i_t$ when $V$ satisfies the $C(A,\alpha)$ condition or for the projected system $Y^i_t$ \cref{eq:projectedsystem} when $V=0$ and $W$ satisfies the $C(A,\alpha)$ condition. In both cases, the convergence rate obtained in \cite[Theorem 3.1]{cattiaux_probabilistic_2008} is $N^{-1/(\alpha+1)}$ instead of $N^{-1}$ in Theorem \cref{thm:uniformpocuniformconvex}.

\item \textbf{Convexity outside a ball of confinement and large diffusion.} As already explained, uniform in time propagation is strongly linked to the existence of a unique stationary solution to the nonlinear equation \cref{eq:granularmedia}. It has been proved in \cite{herrmann_non-uniqueness_2010,tugaut_convergence_2013,tugaut_phase_2014} that such uniqueness does not hold in general without a convexity assumption. However, uniqueness may hold even in non convex settings provided that the diffusion $\sigma$ is large enough and with the assumption of convexity outside a ball of confinement. This includes important cases such as double-well potentials. Convergence to equilibrium for the nonlinear system is studied in particular in \cite{tugaut_convergence_2013,tugaut_phase_2014,bolley_uniform_2013}. Extending these results at the particle level has been the subject of many recent works. To prove uniform in time propagation of chaos, new coupling approaches, which go beyond the traditional synchronous coupling, have been developed. They will be discussed in more details in the following sections. Let us mention in particular the reflection coupling method \cite{durmus_elementary_2020} (Section \cref{sec:reflectioncoupling}) and the optimal coupling approach of \cite{salem_gradient_2020,del_moral_uniform_2019} (Section~\cref{sec:optimalcoupling}). 
\end{enumerate}

We end this section by reviewing some cases which go beyond the gradient setting. 

\subsubsection*{More general diffusion matrices.} Taking a general diffusion matrix $\sigma\equiv\sigma(x,\mu)$ would add two terms in It\={o}'s formula in the proof of Theorem \cref{thm:uniformpocuniformconvex}:
\[\int_0^t \|\sigma(X^i_s,\mu_{\mathcal{X}^N_s})-\sigma(\overline{X}{}^i_s,f_s)\|^2\dd s\]
and
\[2\int_0^t\big\langle X^i_s-\overline{X}{}^i_s,(\sigma(X^i_s,\mu_{\mathcal{X}^N_s})-\sigma(\overline{X}{}^i_s,f_s))\dd B^i_s\big\rangle. \] 
The same proof of uniform in time propagation of chaos would work for globally Lipschitz $\sigma$ with a Lipschitz constant $L>0$ which is sufficiently small with respect to $\beta$.

\subsubsection*{Non-gradient systems.} The proof does not really depend on the form of the drift. To get uniform in time propagation of chaos, more general interactions can be considered provided that they satisfy the same convexity and growth assumptions satisfied by $\nabla V$ and $\nabla W$. The gradient system setting seems more natural to study convergence to equilibrium properties as already discussed. However, similar results than the ones presented in this section but in a very general, yet restrictive, framework can be found for instance in \cite{veretennikov_ergodic_2006}. See also \cite{mishura_existence_2020,wang_distribution_2018} for additional weak and strong well-posedness results on the corresponding nonlinear process. 

\subsubsection*{Kinetic systems.} These ideas can be extended to the case of a kinetic system $\mathcal{Z}^N_t = \big((X^1_t,V^1_t),\ldots(X^N_t,V^N_t)\big)\in (\R^d\times\R^d)^N$ defined by the $N$ coupled SDE: 
\begin{equation}\label[IIeq]{eq:particlevfp}
\left\{
\begin{array}{rcl}
    \dd X^i_t &=& V^i_t \dd t\\
    \dd V^i_t &=& -F(V^i_t)\dd t-G(X^i_t)\dd t-H\star\mu_{\mathcal{X}^N_t}(X^i_t)\dd t + \sigma\dd B^i_t
\end{array}
\right.,
\end{equation}
where $F,G,H:\R^d\to\R^d$ are respectively called the friction force, the exterior confinement force and the interaction force and $\mu_{\mathcal{X}^N_t}$ denotes the $x$-marginal of $\mu_{\mathcal{Z}^N_t}$, so that 
\[H\star\mu_{\mathcal{X}^N_t}(X^i_t) = \frac{1}{N}\sum_{j=1}^N H(X^i_t-X^j_t).\]
The corresponding nonlinear McKean-Vlasov process is obtained by replacing the empirical measure of the particle system by the law $f_t(x,v)\dd x\dd v$ of the nonlinear process which is the solution of the famous Vlasov-Fokker-Planck equation: 
\begin{equation}\label[IIeq]{eq:vfp}
    \partial_t f_t + v\cdot\nabla_x f_t - H\star\rho[f_t](x)\cdot\nabla_v f_t = \frac{\sigma^2}{2}\Delta_v f_t + \nabla_v\cdot\big((F(v)+G(x))f_t\big),
\end{equation}
where $\rho[f_t](x) := \int_{\R^d}f_t(x,v)\dd v$. 

\begin{theorem}[\cite{bolley_trend_2010}]\label[II]{thm:kineticuniformpoc} In \cref{eq:particlevfp}, assume that the forces satisfy the following properties.
\begin{itemize}
    \item There exists $\alpha,\alpha'>0$ such that for all $v,w\in\R^d$, 
    \[|F(v)-F(w)|\leq\alpha |v-w|,\quad \langle v-w,F(v)-F(w)\rangle\geq \alpha'|v-w|^2.\]
    \item There exists $\beta,\delta>0$ such that for all $x,y\in\R^d$, 
    \[G(x) = \beta x + \tilde{G}(x),\quad |\tilde{G}(x)-\tilde{G}(y)|\leq \delta |x-y|.\]
    \item There exists $\gamma>0$ such that for all $x,y\in\R^d$, 
    \[|H(x)-H(y)|\leq \gamma |x-y|.\]
    Then there exists $\varepsilon_0>0$ such that if $0\leq \gamma,\delta<\varepsilon_0$, then there exists a constant $C>0$ such that
    \[\sup_{t\geq0} \frac{1}{N}\sum_{i=1}^N \E{\left[|X^i_t-\overline{X}{}^i_t|^2+|V^i_t-\overline{V}^i_t|^2\right]}\leq\frac{C}{N}.\]
\end{itemize}
\end{theorem}

The proof of Theorem \cref{thm:kineticuniformpoc} again follows from a classical synchronous coupling. However, the standard approach would not give uniform in time estimates (it would only be a special instance of McKean's theorem in a Lipschitz setting which do not take advantage of the form of the interactions). The idea of \cite[Theorem 1.2]{bolley_trend_2010} is to introduce a new metric on the state space $E=\R^d\times\R^d$ which is equivalent to the usual Euclidean metric but for which some dissipativity can be recovered. Namely, the authors show that there exist $a,b,c>0$, such that the following expression defines a positive definite quadratic form on $\R^d\times\R^d$
\[Q(x,v) = a|x|^2+b\langle x,v\rangle + c|v|^2,\]
and which satisfies 
\[\frac{\dd}{\dd t}\E[Q(X^i_t - \overline{X}{}^i_t,V^i_t,\overline{V}^i_t)] \leq - \E[|X^i_t-\overline{X}{}^i_t|^2+|V^i_t-\overline{V}^i_t|^2] + \frac{C}{N},\]
from which the result follows. 

\begin{remark} This approach is strongly inspired by the hypocoercivity methods \cite{villani_hypocoercivity_2009}. In fact, in the same article \cite[Theorem 1.1]{bolley_trend_2010} the authors also show the exponential convergence to equilibrium of the nonlinear process, using a synchronous coupling method (between two nonlinear processes) and a perturbed Euclidean metric. This extends a classical result of Villani \cite[Theorem 56]{villani_hypocoercivity_2009} to a non-compact setting but for a weaker distance (the Wasserstein distance). Note that unlike \cite{malrieu_logarithmic_2001}, convergence to equilibrium for the Vlasov-Fokker-Planck equation follows only from its nonlinear stochastic interpretation but does not use its particle approximation. The same method could also be applied to the granular media equation \cite[Remark 2.2]{bolley_trend_2010}.
\end{remark}

Although very general, a drawback of Theorem \cref{thm:kineticuniformpoc} is that it only works for close to linear confinement force and small interactions. When the forces derive from potentials, similarly to Theorem \cref{thm:uniformpocuniformconvex}, it becomes possible to prove stronger results by using the explicit expression of the equilibria of the particle system (which is not known in general). The following theorem due to Monmarch\'e \cite[Theorem 3]{monmarche_long-time_2017} considers the uniformly convex case. 

\begin{theorem}[\cite{monmarche_long-time_2017}]\label[II]{thm:kineticgradientuniformpoc} In \cref{eq:particlevfp}, assume that the following properties hold. 
\begin{itemize}
    \item There exists $\gamma>0$ such that for all $v\in\R^d$,
    \[F(v) = -\gamma v.\]
    \item There exists a smooth potential $V:\R^d\to\R$ with bounded derivatives of order larger than 2 and such that for all $x\in\R^d$,
    \[G(x) = \nabla V(x).\]
    Moreover, $V$ is uniformly convex in the sense that there exists $c_1>0$ such that $\nabla^2 V \geq c_1$.
    \item There exists a smooth symmetric potential $W:\R^d\to\R$ with bounded derivatives of order larger than 2 and such that for all $x\in\R^d$,
    \[H(x) = \nabla W(x).\]
    Moreover, there exists a constant $c_2<\frac{1}{2}c_1$ such that $\nabla^2 W \geq -c_2$.
\end{itemize}
Let $f_0\in\pb_2(\R^d)$ admit a smooth density in $L\log L$. Then there exists $\alpha>0$ and $C>0$ such that 
    \[\sup_{t\geq 0}W_2\big(f^{1,N}_t,f_t\big)\leq \frac{C}{N^\alpha},\]
    and the same estimate also holds in total variation norm. 
\end{theorem}

Within this setting, the $N$-particle process admits a unique stationary distribution given by its density: 
\[\mu^N_\infty(\dd \mathbf{x},\dd \mathbf{v}) \propto \exp{\left(-\frac{2\gamma}{\sigma^2}{\left(\sum_{i=1}^N V(x^i)+\frac{1}{2N}\sum_{i,j=1}^N W(x^i-x^j)+\frac{1}{2}\sum_{i=1}^N |v^i|^2\right)}\right)}\dd\mathbf{x}\dd\mathbf{v}.\]
one of the main results of \cite[Theorem 1]{monmarche_long-time_2017} is the exponential decay of the relative entropy for the $N$-particle process with a rate which does not depend on $N$, namely there exist $C,\chi>0$ such that
\begin{equation}\label[IIeq]{eq:monmarcheentropy}H\big(f^N_t|\mu^N_\infty\big) \leq C\e^{-\chi t}H\big(f_0^N|\mu^N_\infty\big).\end{equation}
Combined with McKean's theorem, as in \cite{malrieu_logarithmic_2001}, it is then possible to prove the exponential convergence towards equilibrium for the nonlinear process \cite[Lemma~8, Proposition 13]{monmarche_long-time_2017}. Namely, there exist $\mu_\infty(\dd x, \dd v)\in \pb_2(\R^d\times\R^d)$ and $C>0$ such that 
\begin{equation}\label[IIeq]{eq:monmarchecvequilibrium}W^2_2(f_t,\mu_\infty)\leq C \e^{-\chi t}.\end{equation}
Combining this long-time estimate with the short-term bound given by McKean's theorem, it is possible to improve the propagation of chaos result to get a uniform in time convergence. For $t\leq \varepsilon \log N$, McKean's theorem already gives two constants $C,b>0$ such that
\[W_2\big(f^{1,N}_t,f_t\big)\leq \frac{C}{N^{1/2-b\varepsilon}}.\]
Then for $t\geq \varepsilon \log N$, using the normalised distance on $E^N$ (see Definition \cref{def:spaceswasserstein}), 
\begin{align*}
    W_2 \big(f^{1,N}_t,f_t\big) & \leq W_2\big(f^N_t,f_t^{\otimes N}\big)\\
    &\leq W_2\big(f^N_t,\mu^N_\infty\big)+W_2\big(\mu^N_\infty,\mu^{\otimes N}_\infty\big)+W_2\big(\mu^{\otimes N}_\infty,f_t^{\otimes N}\big)\\
    &\leq C{\left(\frac{1}{N^{\varepsilon\chi}}+\frac{1}{N^{1/2}}\right)},
\end{align*}
The first and third terms on the right-hand side of the second line are bounded by $CN^{-\varepsilon\chi}$ using \cref{eq:monmarcheentropy} and \cref{eq:monmarchecvequilibrium}. The second term is bounded by \cite[Lemma~8]{monmarche_long-time_2017}. Theorem \cref{thm:kineticgradientuniformpoc} follows by taking $\varepsilon=(\chi+2b)^{-1}$. Note that
unlike the previous theorems in this section, the final uniform in time estimate is not at the level of the trajectories. In the work of Malrieu, convexity is used to prove uniform in time propagation of chaos and to prove that the $N$-particle law satisfies a log-Sobolev inequality. Since the previous argument relies on the classical McKean's theorem, convexity is only used to obtain the bound \cref{eq:monmarcheentropy}. In a recent work \cite{guillin_uniform_2020}, Guillin and Monmarch\'e have used the results of \cite{guillin_uniform_2019} to remove the convexity assumptions, allowing a broader class of potentials, notably potentials which are convex outside a ball of confinement. Finally, even more recently, Guillin, Le Bris and Monmarch\'e \cite{guillin_convergence_2021} have used a completely different technique, the reflection coupling method discussed below (Section \cref{sec:reflectioncoupling}), to further weaken the assumptions on the potentials. 

\subsection{Other coupling techniques}\label[II]{sec:othercouplingsreview}

In this section, we review some of the main results obtained by the other types of couplings presented in Section \cref{sec:coupling}.  

\subsubsection{Reflection coupling for uniform in time chaos}\label[II]{sec:reflectioncoupling}

Let us consider a gradient system of the form \cref{eq:gradientsystem}. Following the work of \cite{eberle_reflection_2016,eberle_quantitative_2019} on reflection coupling (see Section \cref{sec:reflectioncouplingproving} for a general presentation), we first state the following technical lemma which is the cornerstone of \cite{eberle_reflection_2016,durmus_elementary_2020}.

\begin{lemma}[\cite{eberle_reflection_2016,durmus_elementary_2020}]\label[II]{lemma:reflfectionexistsf}
Assume that $V$ is such that there exists a continuous function $\kappa:[0,+\infty)\to\R$ satisfying $\liminf_{r\to+\infty}\kappa(r)>0$ and
\begin{equation}\label[IIeq]{eq:reflectionVkappa}\forall x,y\in\R^d,\quad\big\langle x-y,\nabla V(x)-\nabla V(y)\big\rangle\geq \frac{\sigma^2}{2}\kappa(|x-y|)|x-y|^2.\end{equation}
Then there exists an increasing $C^2$ concave function $f:[0,+\infty)\to[0,+\infty)$ with $f(0)=0$ and $f'\leq 1$ such that 
\begin{equation}\label[IIeq]{eq:reflectiondf}d_f:(x,y)\in\R^d\times\R^d\longmapsto f(|x-y|)\end{equation}
defines a distance on $\R^d$ and which satisfies for all $r\geq0$,
\begin{equation}\label[IIeq]{eq:refelctionfeq}f''(r)-\frac{1}{4}r\kappa(r)f'(r)\leq -\frac{c_0}{2\sigma^2}f(r),\end{equation}
for a constant $c_0>0$. 
\end{lemma}

The proof of Lemma \cref{lemma:reflfectionexistsf} can be found in \cite[Section 2.1]{durmus_elementary_2020} which follows closely the framework introduced in \cite[Section 2.1]{eberle_reflection_2016}. The function $f$ and the constant $c_0$ have an explicit but somehow not particularly enlightening expression as a function of $\kappa$. Their construction is nevertheless motivated and detailed in \cite[Section 4]{eberle_reflection_2016} (see also Section \cref{sec:reflectioncouplingproving}). 

The condition \cref{eq:reflectionVkappa} on $V$ implies the existence of two constants $m_V>0$ and $M_V\geq0$ such that for all $x,y\in\R^d$,
\[\big\langle x-y,\nabla V (X)-\nabla V(y)\big\rangle \geq m_V|x-y|^2-M_V.\]
This implies uniform convexity outside a ball and thus allows non globally convex settings, the prototypical example being the double well potential
\[V(x) = |x|^4-a |x|^2.\]
We now state the main theorem of this section, it is due to \cite{durmus_elementary_2020}. 

\begin{theorem}[\cite{durmus_elementary_2020}] Let us consider the McKean-Vlasov particle system \cref{eq:mckeanvlasov_summary} with $b,\sigma$ given by \cref{eq:gradientsystem}. Let $V$ be such that there exist a function $f$ and a constant $c_0$ given by Lemma \cref{lemma:reflfectionexistsf}. Assume that the interaction potential $W$ is symmetric, that $\nabla W$ is Lipschitz for the distance induced by $f$ with a Lipschitz constant $\eta>0$ and that there exists $M_W\geq0$ such that $\nabla^2 W\geq -M_W$. Let the $N$ particles be initially i.i.d with law $\tilde{f}_0\in\pb(\R^d)$. Then for all $t\geq0$, it holds that: 
\[W_{1,d_f}\big(f^N_t,f_t^{\otimes N}\big) \leq \e^{-2(c_0-\eta)t}W_{1,d_f}(\tilde{f}_0,f_0)+\frac{C(c_0,\eta)}{\sqrt{N}},\]
where we recall that $W_{1,d_f}$ denotes both the Wasserstein-1 distance on $\pb(\R^d)$ for the distance $d_f$ defined by \cref{eq:reflectiondf} and the Wasserstein-1 distance on $\pb((\R^d)^N)$ for the normalised distance induced by $d_f$.
\end{theorem}

\begin{proof}[Proof (sketch)] The strategy is to use a componentwise reflection coupling in $(\R^d)^N$ between a particle system $\mathcal{X}^N_t$ and a system $\overline{\mathcal{X}}{}^N_t$ of independent nonlinear McKean-vlasov processes. Since the reflection coupling badly behaves on the diagonal, \cite{eberle_reflection_2016,durmus_elementary_2020} introduced the following interpolation between reflection and synchronous coupling:
\begin{align*}
    \dd \overline{X}{}^i_t &= -\nabla V(\overline{X}{}^i_t)\dd t - \nabla W\star f_t(\overline{X}{}^i_t)\dd t + \sigma \Big\{\phi^\delta(E^i_t)\dd B^i_t + \big(1-\phi^\delta(E^i_t)\big)\dd \tilde{B}^i_t\Big\}\\
    \dd {X}^i_t &= -\nabla V({X}^i_t)\dd t - \nabla W\star \mu_{\mathcal{X}^N_t}({X}^i_t)\dd t \\
    &\qquad+ \sigma \Big\{\phi^\delta(E^i_t)\big(I_d-2e^i_t(e^i_t)^\mathrm{T}\big)\dd B^i_t + \big(1-\phi^\delta(E^i_t)\big)\dd \tilde{B}^i_t\Big\},
\end{align*}
where $(B^i_t)^{}_t$ and $(\tilde{B}^i_t)^{}_t$ are $2N$ independent Brownian motions, where
\[E^i_t := \overline{X}{}^i_t - X^i_t,\quad e^i_t := E^i_t/|E^i_t|,\]
and where $\phi^\delta : \R^d\to\R$ is a Lipschitz function such that 
\[\phi^\delta(x) = \left\{ 
\begin{array}{rcl} 
1 & \text{if} & |x|\geq \delta\\
0 & \text{if} & |x|\leq \delta/2
\end{array}
\right.,\]
for a parameter $\delta>0$ (ultimately $\delta\to0$). It is also assumed that $\mathcal{X}^N_0 \sim \tilde{f}_0^{\otimes N}$ and $\overline{\mathcal{X}}{}^N_0\sim f_0^{\otimes N}$ are optimally coupled for the distance $W_{1,d_f}$. 

Using \cite[Lemma 7]{durmus_elementary_2020}, It\=o's formula gives:
\begin{align*}
    \dd f(|E^i_t|) = \Big(f'(|E^i_t|)C^i_t+2\sigma^2 f''(|E^i_t|)\phi^\delta(E^i_t)^2\Big)\dd t+ f'(|E^i_t|)A^i_t\dd t + \dd M^i_t,
\end{align*}
where 
\[C^i_t = -\langle \nabla V(\overline{X}{}^i_t)-\nabla V(X^i_t), e^i_t\rangle,\]
and $(A^i_t)^{}_t$ is an adapted stochastic process such that 
\[A^i_t \leq \left|\nabla W\star f_t(\overline{X}{}^i_t)-\nabla W \star \mu_{\mathcal{X}^N_t}(X^i_t)\right|,\]
and $M^i_t$ is a martingale. As usual the drift term is split into two parts, the ``non-interacting'' part which involves $C^i_t$ and the ``interacting'' part which involves $A^i_t$. Similarly to the proof of Theorem \cref{thm:uniformpocuniformconvex}, Durmus et al. take advantage of the non-interacting part to get a uniform in time control and they treat the interacting part as a perturbation. The main difference is that, thanks to \cref{eq:refelctionfeq}, the reflection coupling gives a better control, namely it holds that: 
\begin{align*}
    f'(|E^i_t|)C^i_t+2\sigma^2 f''(|E^i_t|)\phi^\delta(E^i_t)^2 \leq -2c f(|E^i_t|) + \omega(\delta) + 2c_0f(\delta),
\end{align*}
where $\omega(r) = \sup_{s\in[0,r]} s\kappa(s)^{-}$. The interacting part is controlled as usual by forcing the introduction of a nonlinear term: 
\[A^i_t \leq \left|\nabla W \star \mu_{\overline{\mathcal{X}}{}^N_t}(\overline{X}{}^i_t)-\nabla W \star \mu_{\mathcal{X}^N_t}(X^i_t)\right| + \left|\nabla W \star \mu_{\overline{\mathcal{X}}{}^N_t}(\overline{X}{}^i_t)-\nabla W\star f_t(\overline{X}{}^i_t)\right|.\]
Using the hypotheses on $W$, the fact that $f$ is increasing and defines a distance, it holds that  
\[\left|\nabla W \star \mu_{\overline{\mathcal{X}}{}^N_t}(\overline{X}{}^i_t)-\nabla W \star \mu_{\mathcal{X}^N_t}(X^i_t)\right| \leq \eta f(|E^i_t|) +  \frac{\eta}{N}\sum_{j=1}^N f(|E^j_t|).\]
Then, as in the proof of McKean's theorem (Theorem \cref{thm:mckean}), since it holds that
\[\E\left[\nabla W(\overline{X}^i_t-\overline{X}^j_t)\,\Big|\,\overline{X}^i_t\right] = \nabla W\star f_t(\overline{X}^i_t),\]
and since $\nabla W(0)=0$ and $\nabla W$ is Lipschitz with a constant $\eta$, by the Cauchy-Schwarz inequality and the independence of the processes $\overline{X}^i_t$, it holds that 
\[\E\left|\nabla W \star \mu_{\overline{\mathcal{X}}{}^N_t}(\overline{X}{}^i_t)-\nabla W\star f_t(\overline{X}{}^i_t)\right| \leq \frac{C\eta}{\sqrt{N}},\]
where the constant $C>0$ only depends on a uniform in time moment control on $f_t$ (similar to Lemma \cref{lemma:uniformpocmoment} and proved in \cite[Lemma 8]{durmus_elementary_2020}). Therefore, 
\[\frac{1}{N}\sum_{i=1}^N\E A^i_t \leq \frac{2\eta}{N}\sum_{i=1}^N \E f(|E^i_t|) + \frac{C\eta}{\sqrt{N}}.\]
Using that $f'\leq 1$, this yields
\[\frac{\dd }{\dd t} \frac{1}{N} \sum_{i=1}^N \E f(|E^i_t|) \leq -2\frac{c_0-\eta}{N}\sum_{i=1}^N \E f(|E^i_t|) + \omega(\delta)+2c_0f(\delta)+\frac{C\eta}{\sqrt{N}}.\]
The conclusion follows by integrating this differential Gronwall-like inequality and by letting $\delta\to0$, since it holds that $\lim_{\delta\to0^+}\omega(\delta)=0$ and $f(0)=0$. 
\end{proof}

\begin{remark}[Sticky boundary and sticky coupling]
At the informal level, one could take $\delta=0$ so that the synchronous coupling only acts when $E^i_t =0$, i.e. when the two processes coincide. Since the drifts do not coincide, one could expect that the processes immediately become different and that this only happens for a Lebesgue-nul set of times. However this is not true in general. This is closely linked with the (quite strange) fact the solution to the SDE on $[0,+\infty)$ with sticky boundary behavior at $0$ 
\[ \dd r_t = a\dd t + \mathbbm{1}_{r_t>0} \dd B_t , \] 
with $a > 0$, indeed spends some (Lebesgue-positive) time at $0$, hence the term \emph{sticky}. This fact has been shown by Watanabe \cite{watanabe_stochastic_1971} and recently by Eberle and Zimmer \cite{eberle_sticky_2019}. We also mention that this latter work introduce the notion of \emph{sticky coupling} between two diffusion processes with different drifts. Using this coupling, the distance between the two processes is controlled by the solution of a SDE on $[0,+\infty)$ which is \emph{sticky} at 0. Based on this idea, the question of the long-time behaviour of nonlinear McKean-Vlasov processes and of the uniform in time propagation of chaos has been recently revisited in \cite{durmus_sticky_2022}.  
\end{remark}

The fact that the result holds for the distance $W_{1,d_f}$ may seem unsatisfactory compared to Theorem \cref{thm:uniformpocuniformconvex} and its extensions which hold in $W_2$ distance or Theorem~\cref{thm:BGVconcentrationineq} which holds in $W_1$ distance. This directly comes from the somehow ad hoc estimate \cref{eq:refelctionfeq}. The result has been recently improved in \cite{liu_long-time_2021}, where instead of the function $f$, the authors consider the function $h$ solution of the following Poisson equation
\[4h''(r)+r\kappa(r)h'(r) = -r,\quad r>0.\]
Using the same reflection coupling strategy, the authors obtain a similar uniform in time propagation of chaos result in $W_1$ distance, in both the pathwise and pointwise settings, see \cite[Theorem 2.9]{liu_long-time_2021}. This article also contains many results in $W_1$ distance regarding concentration inequalities and explicit exponential rates of convergence towards equilibrium (independent of $N$) for the particle system. The choice of the function $h$ avoids some technicalities in the definition of the function $f$ (see \cite{eberle_quantitative_2019} and Lemma \cref{lemma:reflfectionexistsf}). The setting is quite general so we do not give all the details here. We only mention that it applies to cases where $V$ is convex outside a ball but has potentially many wells. An assumption on $\nabla^2 W$ is also made in order to prevent phase transitions (which would forbid uniform in time estimates), see \cite[Section 2.2]{liu_long-time_2021}. 

Finally, we also mention the recent extension \cite{guillin_convergence_2021} of the reflection coupling method to kinetic systems, that is for systems of the form \cref{eq:particlevfp} in the setting of Theorem \cref{thm:kineticgradientuniformpoc} but with relaxed convexity assumptions. This work generalizes \cite{monmarche_long-time_2017,guillin_uniform_2020}, in particular, it does not require the knowledge of the invariant measure of the particle system and allows a broader class of potentials. 

\subsubsection{Chaos via Glivenko-Cantelli}\label[II]{sec:holding}

In a recent article \cite{holding_propagation_2016},
Holding proves pathwise chaos on finite time intervals using a coupling on vector-fields instead of particles for a system of the form \cref{eq:mckeanvlasov_summary} with Assumption \cref{eq:mckeanthmassum} and constant diffusion matrix $\sigma = I_d$. The interaction kernel $K_1$ is not assumed to be Lipschitz-continuous, it can be only H\"older-continuous (with and exponent larger than $ 2/3$ for kinetic systems).  There is no strong assumption on the initial data. The argument is based on a new Glivenko-Cantelli theorem for vector fields which moves the need for regularity properties from the SDE system onto the limit equation.

Holding introduces the random measure $\widetilde{f}^{b_N}_t\in\pb(E)$, solution of the equation 
\begin{equation}\label[IIeq]{eq:holdingpde}\partial_t\widetilde{f}^{b_N}_t + \nabla_x\cdot{\left(b\big(x,\mu_{\mathcal{X}^N_t}\big)\widetilde{f}^{b_N}_t\right)} = \frac{1}{2}\Delta_x\widetilde{f}^{b_N}_t,\quad \widetilde{f}^{b_N}_0 = f_0.\end{equation}
From a SDE point of view, the coupling introduced by \cite{holding_propagation_2016} is given by the following exchangeable particle system, defined conditionally on the random vector-field $b_N:=b ( \cdot , \mu_{\mathcal{X}^N_s})$ 
\[ \widetilde{X}^{i,N}_t =  \widetilde{X}^{i,N}_0 + \int_0^t b {\left( \widetilde{X}^{i,N}_s , \mu_{\mathcal{X}^N_s} \right)} \dd s + \widetilde{B}^{i}_t, \]
where $\widetilde{B}^i_t$ are independent Brownian motions, which are also independent of $\mathcal{X}^N$. The starting point of the proof is then similar to the one of McKean's theorem but with the splitting step: 
\begin{equation}\label[IIeq]{eq:splittingholding} \E {\left[ W_2^2 {\left( \mu_{\mathcal{X}^N_t} , f_t  \right)} \right]} \leq 2 \E {\left[ W_2^2 {\left( \mu_{\mathcal{X}^N_t} , \widetilde{f}^{b_N}_t \right)} \right]} + 2 \E {\left[ W_2^2 {\left( \widetilde{f}^{b_N}_t , f_t \right)} \right]}. \end{equation}
The main difference with \cref{eq:mckeanthmsznitsecondineq} is that $\widetilde{f}^{b_N}_t$ replaces $\mu_{\overline{\mathcal{X}}{}^N_t}$. 

Since the diffusion matrix is constant in \cref{eq:mckeanvlasov_summary}, the interaction between the particles is only prescribed by the vector field $b_N := b(\cdot,\mu_{\mathcal{X}^N_s})$. One of the main ideas in \cite{holding_propagation_2016} is to see this vector field as a random variable (with values in a suitable H\"older space). Heuristically, if a random vector field $\tilde{b}$ with the same law as $\tilde{b}_N$ is given independently of everything, then the particle system \eqref{eq:mckeanvlasov_summary} where the term $b(\cdot,\mu_{\mathcal{X}^N_s})$ is replaced by $\tilde{b}$ simply becomes a system of independent processes since the interaction term has been removed. In other words, conditionally on the random vector-field $b_N:=b ( \cdot , \mu_{\mathcal{X}^N_s})$, the particles $X^{i,N}_t$ are $\widetilde{f}^{b_N}_t$-distributed random variables. To see how this observation can help to bound the first term on the right-hand side of~\cref{eq:splittingholding}, given a fixed (smooth) vector field $b:\R^d\to\R^d$ (random or not), we denote by $\mu_{\mathcal{X}^N_t|b}$ the empirical measure of the $N$-particle system \cref{eq:mckeanthmassum} where the drift is replaced by this fixed $b$. Note that $\mu_{\mathcal{X}^N_t|b_N}=\mu_{\mathcal{X}^N_t}$. Similarly we denote by $\widetilde{f}^b_t$ the solution of~\cref{eq:holdingpde} with $b_N$ replaced by $b$. Then the first term on the right-hand side of \cref{eq:splittingholding} reads:
\begin{equation}\label[IIeq]{eq:conditionalvectorfield}\E {\left[ W_2^2 {\left( \mu_{\mathcal{X}^N_t} , \widetilde{f}^{b_N}_t \right)} \right]}=\E {\left[ \E{\left[W_2^2 {\left( \mu_{\mathcal{X}^N_t|b_N} , \widetilde{f}^{b_N}_t \right)} \big| b_N\right]} \right]}.\end{equation}
Since $b_N$ cannot be easily controlled, in order to bound the inside conditional expectation, the author of \cite{holding_propagation_2016} simply takes the supremum of the Wasserstein distance over all vector fields in a suitable space. More precisely, one of the main result \cite[Corollary 2.3]{holding_propagation_2016} is a generalized Glivenko-Cantelli theorem for SDE which gives the explicit bound: 
\[ \mathbb{E} {\left[ \sup_{b \in \mathcal{B}} \sup_{0 \leq t\leq T}  W_2^2 {\left( \mu_{\mathcal{X}^N_t|b} , \widetilde{f}^{b}_t \right)}  \right]} \leq \varepsilon(N,T), \]
where $\mathcal{B}$ is a subset of H\"older regular vector fields and $\varepsilon(N,T)\to 0$ is an explicit polynomial rate of convergence. Taking the supremum in \cref{eq:conditionalvectorfield} over all the vector fields in $\mathcal{B}$, this controls the first term on the right-hand side of \cref{eq:splittingholding}. To conclude, the control of the second term on the right-hand side of \cref{eq:splittingholding} shall result from stability estimates on the solution of the limit equation with respect to its vector-field parameter $b$.

With a traditional synchronous coupling, the Lipschitz regularity of $b$ is used to control the particle system and a crude $L^\infty$ estimate is used to control the error term which depends on the limit equation. With this approach, the need for regularity on $b$ is weakened by the generalized Glivenko-Cantelli theorem and the control of the error term can take advantage of the regularity properties of the limit equation. This idea is successfully applied in \cite{holding_propagation_2016} to various first-order and kinetic systems, at the pathwise level. 

\subsubsection{Optimal coupling and WJ inequality}\label[II]{sec:optimalcoupling} 

This section is devoted to the analytical coupling approach of \cite{salem_gradient_2020} described in the introductory Section \cref{sec:optimalcouplingintro}. In this section, this approach is mainly applied to gradient-systems but, as explained in \cite[Section 2.2]{salem_gradient_2020}, it also allows to recover, at the level of the laws, many of the results obtained by synchronous or reflection coupling for more general McKean-Vlasov systems. The strategy originated from the earlier works \cite{bolley_convergence_2012,bolley_uniform_2013} where the author prove the convergence to equilibrium of the solution of respectively the linear Fokker-Planck equation and the nonlinear granular media equation. The strategy is adapted and carried out at the particle level in \cite{salem_gradient_2020} and in \cite{del_moral_uniform_2019} to prove at the same time the convergence to equilibrium and the propagation of chaos in a non globally convex setting. 
\medskip

In this section, we recall (see Definition \cref{def:spaceswasserstein}) that $\widetilde{W}_2$ denotes the non-normalised Wasserstein distance on $E^N$ defined by $\widetilde{W}_2^2(f^N,g^N) = NW_2^2(f^N,g^N)$ for $f^N,g^N\in\pb(E^N)$. 
\medskip

The starting point of the argument is the following observation. For McKean-Vlasov systems \cref{eq:mckeanvlasov_summary} with linear interaction functions of the form
\[b(x,\mu) \equiv \int_{\R^d}b(x,y)\mu(\dd y), \quad \sigma(x,\mu)\equiv \sigma I_d,\quad \sigma>0,\]
the laws $f^N_t$ and $f^{\otimes N}_t$ are both absolutely continuous solutions of continuity equations in $\R^{dN}$. It is therefore possible to compute the dissipation rate in $\widetilde{W}_2$ distance between them using a result which originates from the theory of gradient-flows \cite[Theorem 23.9]{villani_optimal_2009}, namely it holds that: 
\begin{align}\label[IIeq]{eq:derivW2}
&\frac{\dd}{\dd t}\frac{1}{2}\widetilde{W}_2^2\big(f^N_t, f^{\otimes N}_t\big) \nonumber\\
&= \int_{\R^{dN}} \Big\langle \mathbf{b}^N(\mathbf{x}^N) - \frac{\sigma^2}{2}\nabla \log f^N_t(\mathbf{x}^N),\nabla \psi^{N\star}_t(\mathbf{x}^N)-\mathbf{x}^N\Big\rangle f^N_t\big(\dd \mathbf{x}^N\big)\nonumber \\
&\quad- \int_{\R^{dN}} \Big\langle \Big(b\star f_t - \frac{\sigma^2}{2}\nabla \log f_t\Big)^{\otimes N}(\mathbf{x}^N),\nabla \psi^{N}_t(\mathbf{x}^N)-\mathbf{x}^N\Big\rangle f^{\otimes N}_t\big(\dd \mathbf{x}^N\big),
\end{align}
where
\[\mathbf{b}^N : \mathbf{x}^N\in \R^{dN}\mapsto \left(\frac{1}{N}\sum_{i=1}^N b(x^1,x^i),\ldots,\frac{1}{N}\sum_{i=1}^N b(x^N,x^i)\right)\in\R^{dN},\]
and $\psi^N_t$ is the maximizing Kantorovich potential between $f^{\otimes N}_t$ and $f^N_t$ given by Brenier's theorem \cite[Theorem 9.4]{villani_optimal_2009} and defined by 
\begin{equation}\label[IIeq]{eq:optimalW2Brenier}
    \widetilde{W}_2^2\big(f^N_t,f^{\otimes N}_t\big) = \int_{\R^{dN}} | \nabla \psi^N_t(\mathbf{x}^N)-\mathbf{x}^N|^2f_t^{\otimes N}\big(\dd \mathbf{x}^N\big),
\end{equation}
that is, the coupling $(\nabla \psi^N_t)_\# f^{\otimes N}_t = f^N_t$ is optimal for the $\widetilde{W}_2$ distance. The relation \cref{eq:derivW2} can be obtained by a formal derivation of \cref{eq:optimalW2Brenier}, the rigorous proof is the content of \cite[Theorem 23.9]{villani_optimal_2009}. The cornerstone of \cite{salem_gradient_2020} is the following proposition, which gives an explicit bound for the right-hand side of \cref{eq:derivW2}. For now on we fix $\sigma(x,\mu)=\sqrt{2}I_d$ for simplicity. 

\begin{proposition}\label[II]{prop:salemineq} Given a symmetric probability measure $g^N\in \pb_2((\R^d)^N)$ and $\mu\in~\pb(\R^d)$, Salem introduces the quantity: 
\begin{multline}\label[IIeq]{eq:salemJ}
    \mathcal{J}\big(g^N|\mathbf{b}^N,\nu^{\otimes N}\big):= \int_{\R^{dN}}(\Delta \psi^N(\mathbf{x}^N)+\Delta\psi^{N\star}(\nabla\psi^N)-2dN)\nu^{\otimes N}\big(\dd \mathbf{x}^N\big)\\
    +\frac{1}{N}\sum_{i,j=1}^N\int_{\R^{2dN}} \langle b(\nabla_i\psi^N(\mathbf{x}^N),\nabla_j\psi^N(\mathbf{x}^N))-b(x^i,x^j),\nabla_i\psi^N(\mathbf{x}^N)-x^i\rangle \nu^{\otimes N}\big(\dd\mathbf{x}^N\big), 
\end{multline}
where $\psi^N$ is the maximizing Kantorovich potential such that $(\nabla\psi^N)_\#\nu^{\otimes N}= g^N$. Assume that the vector fields $(\mathbf{b}^N - \nabla\log f^N_s)^{}_{s\geq0}$ and $(b\star f_s-\nabla\log f_s)_{s\geq0}$ are locally Lipschitz and satisfy for any $t\geq0$,
\begin{equation}\label[IIeq]{eq:hypsalemineq}\int_0^t \int_{\R^{dN}} |\mathbf{b}^N - \nabla\log f^N_s|^2 \dd f^N_s \dd s + \int_0^t\int_{\R^{d}} |b\star f_s-\nabla\log f_s|^2\dd f_s\dd s <+\infty.\end{equation}
Then for all $\eta>0$ and all $0\leq r<t$, it holds that 
\begin{multline}\label[IIeq]{eq:saleminequality}\widetilde{W}^2_2(f_t^{\otimes N},f^N_t) \leq \widetilde{W}^2_2(f_r^{\otimes N},f^N_r) -2\int_r^t \mathcal{J}(f^N_s|\mathbf{b}^N,f_s^{\otimes N})\dd s\\ + \eta\int_r^t \widetilde{W}^2_2(f_s^{\otimes N},f^N_s)\dd s + \eta^{-1}\int_r^t \mathcal{F}_N(b,f_s)\dd s,
\end{multline}
where
\begin{equation}\label[IIeq]{eq:salemF}
\mathcal{F}_N(b,f_s) = \frac{1}{N^2}\sum_{i=1}^N\sum_{j\ne i}^N \int_{\R^{dN}}|b(x^i,x^j)-b\star f_s(x^i)|^2 f_s^{\otimes N}\big(\dd \mathbf{x}^N\big).
\end{equation}
\end{proposition}
The proof is detailed in \cite[Proposition 1]{salem_gradient_2020}. Under mild local Lipschitz assumptions on $b$, the functional $\mathcal{F}_N$ can be easily bounded uniformly in $N$. The whole point is therefore to find a good control of the functional $\mathcal{J}$. Two main ideas are given.
\begin{itemize}
    \item First, it is possible to prove (see \cite[Lemma 3.2]{bolley_convergence_2012}):
    \begin{equation}\label[IIeq]{eq:partJgeq0}\int_{\R^{dN}}(\Delta \psi^N_s(x)+\Delta\psi^{N*}_s(\nabla\psi^N_s)-2dN)\mu_s^{\otimes N}(\dd x)\geq0.\end{equation}
    From this crude estimate, one can just neglect the corresponding term in \cref{eq:saleminequality} and retrieve all the results based on synchronous coupling (in particular McKean's theorem and Theorem \cref{thm:uniformpocuniformconvex}). 
    \item More generally, in order apply the Gronwall lemma in \cref{eq:saleminequality}, it is desirable to bound $\mathcal{J}$ from below by a $W_2$ distance. This lead \cite{bolley_convergence_2012} and later \cite{salem_gradient_2020,del_moral_uniform_2019} to introduce the WJ inequality. In this context, a probability measure $\nu\in\pb_2(\R^d)$ is said to satisfy a symmetric WJ($\kappa$) inequality for a constant $\kappa>0$ when for all symmetric probability measure $g^N\in\pb_2((\R^d)^N)$, it holds that
    \begin{equation}\label[IIeq]{eq:WJineq}\kappa \widetilde{W}_2^2\big(g^N,\nu^{\otimes N}\big)\leq \mathcal{J}\big(g^N|\mathbf{b}^N,\nu^{\otimes N}\big).\end{equation}
    For a gradient system which possesses a unique stationary measure $\mu_\infty$, \cite[Proposition 3]{salem_gradient_2020} shows that $\mu_\infty$ satisfies a WJ($\kappa$) inequality. 
\end{itemize}

The main results \cite[Theorem 2.2, Corollary 1]{salem_gradient_2020} are summarised in the following theorem. 

\begin{theorem}[\cite{salem_gradient_2020}]\label[II]{thm:salem} Assume that $\sigma(x,\mu)=\sqrt{2}I_d$ and $b(x,\mu)\equiv b\star \mu(x)$ where
\[b(x,y) = -\nabla V(x) -\varepsilon \nabla W(x-y),\]
with $V(x)=|x|^4-a|x|^2$ and $W(x)=-|x|^2$, where $a,\varepsilon>0$. Let $f^N_0 \in \pb_6(\R^{dN})\cap L\log L(\R^{dN})$. Then there exist $a_0>0$ and $\varepsilon_0>0$ such that if $a<a_0$ and $\varepsilon<\varepsilon_0$, then the nonlinear McKean-Vlasov equation has a unique stationary solution $\mu_\infty\in\pb_2(\R^d)$ and there exist two constants $C,\alpha>0$ such that (for the normalised Wasserstein distance):
\[\forall t \geq 0, \quad W_2^2\big(f^N_t,\mu_\infty^{\otimes N}\big) \leq W_2^2\big(f^N_0,\mu_\infty^{\otimes N}\big)\e^{-\alpha t}+\frac{C}{N}.\]
Moreover if $f_0\in \pb_6(\R^d)\cap L\log L(\R^d)$ and $f^N_0 = f_0^{\otimes N} $ then there exists $\beta\in(0,1)$ such that
\[\sup_{t \geq 0} W_2^2\big(f_t^{\otimes N},f^N_t\big)\leq CN^{-\beta}.\]
\end{theorem}

\begin{proof}[Proof (summary)]
Let us summarise the main steps of the proof. 
\begin{enumerate}
    \item As usual, some a priori bounds are needed. In \cite[Lemma 4.1]{salem_gradient_2020}, the potentials are shown to satisfy an explicit property of convexity at infinity as well as an explicit $L^\infty$ bound near the origin. Then it is possible to prove classical moment estimates which ensure that if the initial conditions have sufficiently many moments, then the moments of any order of both $f^N_t$ and $f_t$ are uniformly bounded in time. 
    \item The fundamental property is stated in \cite[Proposition 3]{salem_gradient_2020}. First, by \cite[Proposition 4.4 (iii)]{bolley_uniform_2013}, given potentials which satisfy \cite[Lemma 4.1]{salem_gradient_2020}, there exists a stationary solution $\mu_\infty$ of the nonlinear McKean-Vlasov equation. Such a measure is a minimizer of the free energy of the system. Then for $a,\varepsilon$ sufficiently small, such a measure $\mu_\infty$ is shown to satisfy a symmetric WJ($\kappa$) inequality \cref{eq:WJineq} for some $\kappa>0$. This implies the uniqueness of $\mu_\infty$. 
    \item In order to apply Proposition \cref{prop:salemineq}, it is necessary to check the assumption \cref{eq:hypsalemineq}, which again follows from the preliminary bounds derived in \cite[Lemma 4.1, Lemma 4.2]{salem_gradient_2020}. Then since $\mathcal{F}_N$ can be bounded uniformly in $N$ by a constant $C(a,\varepsilon)$, the inequality \cref{eq:saleminequality} applied with the stationary measure $\mu_\infty$ gives for any $\eta>0$ and any $0\leq r<t$:
    \begin{equation}\label[IIeq]{eq:salemgronwall}\widetilde{W}_2^2\big(\mu^{\otimes N}_\infty,f_t^N\big)\leq \widetilde{W}_2^2\big(\mu_\infty^{\otimes N},f^N_r\big)-(\kappa-\eta)\int_r^t \widetilde{W}_2^2\big(\mu_\infty^{\otimes N},f^N_s\big)\dd s + \frac{(t-r)C(a,\varepsilon)}{\eta}.\end{equation}
    Since this holds for any $r<t$, this implies that $\widetilde{W}_2^2\big(\mu^{\otimes N}_\infty,f_t^N\big)$ is differentiable in $t$ and satisfies the differential inequality
    \begin{equation}\label[IIeq]{eq:salemgronwalldifferential}\frac{\dd}{\dd t}\widetilde{W}_2^2\big(\mu^{\otimes N}_\infty,f_t^N\big)\leq -(\kappa-\eta)\widetilde{W}_2^2\big(\mu_\infty^{\otimes N},f^N_t\big) + \frac{C(a,\varepsilon)}{\eta}.\end{equation}
    The Gronwall lemma finally gives the first point of Theorem \cref{thm:salem}.
    \item The above point gives an optimal convergence rate towards the stationary measure $\mu_\infty$. To control the distance to $f^{\otimes N}_t$ at any time $t>0$, the classical strategy is to use on the one hand the exponential convergence of $f_t$ towards $\mu_\infty$ to control the long-time behaviour and on the other hand the non uniform in time McKean's theorem to control the short time behaviour. First using \cite[Theorem 23.9]{villani_optimal_2009} and the WJ($\kappa$) inequality satisfied by $\mu_\infty$, it holds that
    \[W_2^2(f_t,\mu_\infty) \leq W_2^2(f_0,\mu_\infty)\e^{-\kappa t}.\]
    Since the inequality is preserved by tensorization, the triangle inequality yields
    \begin{equation}\label[IIeq]{eq:salemlongtime}W_2^2\big(f^N_t,f^{\otimes N}_t\big)\leq W_2^2\big(f^N_0,\mu_\infty^{\otimes N}\big)\e^{-C(a,\varepsilon)t} + \frac{C}{N}+W_2^2(f_0,\mu_\infty)\e^{-\kappa t}.\end{equation}
    Moreover, for $f_0^N=f_0^{\otimes N}$, using \cref{eq:saleminequality} and \cref{eq:partJgeq0} (or equivalently, McKean's theorem), it holds that for all $t\geq0$, 
    \begin{equation}\label[IIeq]{eq:salemshorttime}W_2^2\big(f^N_t,f_t^{\otimes N}\big)\leq \frac{C\e^{C(a,\varepsilon,\eta)t}}{N}.\end{equation}
    Choosing $T_N = \delta \log N$, the result follows by combining \cref{eq:salemlongtime} for $t\geq T_N$ and~\cref{eq:salemshorttime} for $t\leq T_N$. 
\end{enumerate}
\end{proof}

\begin{remark}
We have slightly changed the original formulation of inequalities \cref{eq:saleminequality} and \cref{eq:salemgronwall} which are stated with $r=0$ only in \cite{salem_gradient_2020} but they actually hold with any $r$. However, as pointed to us by an anonymous reviewer, the Gronwall lemma with negative parameters does not hold true in integral form and the differential inequality \cref{eq:salemgronwalldifferential} is needed to conclude. Moreover, in the inequality on top of page 5750 in \cite{salem_gradient_2020}, the measure $\mu_t$ should read $\mu_\infty$.  
\end{remark}

A similar result is obtained in \cite[Theorem C, Corollary D]{del_moral_uniform_2019} also by means of a WJ inequality but in the equivalent case where $\sigma$ is taken large enough. The authors consider a broader class of potentials, though the main assumption remains convexity of $V$ outside a ball of confinement. In fact, it seems that the result of \cite{salem_gradient_2020} holds for potentials which satisfy \cite[Lemma 4.1]{salem_gradient_2020}, which is very similar to \cite[Assumptions (A-1)-(A-10)]{del_moral_uniform_2019}. The main difference with \cite{salem_gradient_2020} is that the authors do not derive the general inequality \cref{eq:saleminequality} but prove an improved version of McKean's theorem (using a synchronous coupling) in the case of non independent initial conditions, see \cite[Proposition B]{del_moral_uniform_2019}. Both approaches are motivated by \cite{bolley_convergence_2012,bolley_uniform_2013}. More precisely, they are based on \cite[Proposition 3.4]{bolley_convergence_2012} which gives a criterion for an invariant measure $\mu_\infty$ to satisfy a WJ inequality. This leads to the equivalent results \cite[Proposition 3]{salem_gradient_2020} and \cite[Proposition 2.3]{del_moral_uniform_2019}.

\subsection{Compactness methods for mixed systems and gradient flows}\label[II]{sec:mckeancompactnessreview}

The content of this section develops the compactness arguments briefly introduced in Section \cref{sec:provingcompactness_summary} (see also Section \cref{sec:provingcompactness}) in two cases. Section \cref{sec:martingalecompactness} focuses on the functional law of large numbers and the strong pathwise empirical propagation of chaos via martingale arguments (see Definition \cref{def:empiricalchaos}). Section \cref{sec:gradientflows} uses the gradient-flow formulation to prove a pointwise empirical propagation of chaos result for gradient systems.  

\subsubsection{Pathwise chaos via martingale arguments}\label[II]{sec:martingalecompactness}

We first state the assumptions on the generator of the particles process. In all this section we consider a dense separable subspace of the set of test functions $\mathcal{F} \subset C_b ( E )$ such that $\| . \|_{\infty} \leq C_{\mathcal{F}} \| . \|_{\mathcal{F}}$. We assume that $\mathcal{F}$ is contained in the domain of $L_\mu$ for all $\mu\in\pb(E)$ and $\mathcal{F}^{\otimes N}\subset \Dom(\mathcal{L}_N)$. 

\begin{assumption}[Mean-field generator and initial well-posedness] \label[II]{assum:initialdefined}
The generator of the process $( \mathcal{X}^N_t )^{}_{t \geq 0}$ is of the mean-field type \cref{eq:Nparticlemeanfieldgenerator_summary} and the associated martingale problem is wellposed. Moreover the initial law $f^N_0\in\pb(E^N)$ satisfies the moment bound: 
\[\sup_{N}\E\big|X^{i,N}_0\big|^2 <+\infty.\]
\end{assumption}

Since $L_{\mu}$ can involve any differential operator with no homogeneous term and any integral jump operator, this generator covers the case of the McKean-Vlasov diffusion and of the mean-field jump processes. It can also be a mixed jump-diffusion generator. 

\begin{assumption}[Bounds on the limit generator] \label[II]{assum:genebound}
There exists a constant $C_L > 0$ such that
\[\forall x\in E,\,\, \forall \varphi \in \mathcal{F}, \quad \sup_N \sup_{\mu \in \hat{\pb}_N ( E )} {\left\{\big| L_{\mu} \varphi(x )
\big|^2 + \Gamma_{L_{\mu}} ( \varphi , \varphi )(x)\right\}} \leq C_{L} \big( 1 + |x|^2 \big), \]
where the carr\'e du champ operator is defined for regular test functions by
\[ \Gamma_{L_{\mu}} ( \varphi , \psi ) := L_{\mu} {\left[ \varphi \psi \right]} - \varphi L_{\mu} \psi - \psi L_{\mu} \varphi. \]
\end{assumption}

The main consequence of Assumption \cref{assum:genebound} is to ensure the uniform control of the second moment on any interval $[0,T]$ (weaker assumptions could thus be sufficient in specific cases): 
\begin{equation}\label[IIeq]{eq:momentboundmartingale}
    \E{\left[\sup_{t\leq T} \big|X^{1,N}_t\big|^2\right]}\leq C_T{\left(1+\E\big|X^{1,N}_0\big|^2\right)},\quad \E{\left[\sup_{t\leq T} \big|M^{1,N}_t\big|^2\right]}\leq C_T{\left(1+\E\big|X^{1,N}_0\big|^2\right)},
\end{equation}
where $X^{1,N}_t= X^{1,N}_0 + M^{1,N}_t + A^{1,N}_t$ is the semimartingale decomposition of $X^{1,N}_t$ (see Appendix \cref{appendix:Dsemimartingales}). This is proved in \cite[Lemma 3.2.2]{joffe_weak_1986}. It relies on the use of Gronwall lemma in It\=o's formula: the bound on the generator controls the integral term and the bound on the carr\'e du champ operator controls the martingale part (see also Proposition \cref{thm:joffecriterion}). For the jump and diffusion processes, Assumption~\cref{assum:genebound} holds under the usual global Lipschitz assumptions which also ensure the well-posedness of both the particle process and the nonlinear system. We also recall that for the mean-field jump process 
\[\Gamma_{L_\mu}(\varphi,\varphi)(x) = \int_E [\varphi(y)-\varphi(x)]^2 P_\mu(x,\dd y),\]
and for the McKean-Vlasov diffusion, 
\[\Gamma_{L_\mu}(\varphi,\varphi)(x) = 2 \big( \nabla\varphi(x) \big)^{\mathrm{T}} a(x,\mu) \nabla\varphi(x). \]

The main difference between the functional law of large numbers (Theorem \cref{thm:martingaleweakpathwiseempiricalchaos}) and the strong pathwise empirical propagation of chaos result (Theorem \cref{thm:martingalestrongpathwiseempiricalchaos}) will be the assumption on the limit law. 

\subsubsection*{Functional law of large numbers.}

We first prove a functional law of large numbers, that is the convergence of the sequence of 
\begin{equation}\label[IIeq]{eq:FmuNlaw}F^{\mu,N}_{[0,T]} = \mathrm{Law}{\left(\big(\mu_{\mathcal{X}^N_t}\big)_{0\leq t\leq T}\right)}\in \pb\big(D([0,T],\pb(E))\big),\end{equation}
which means that the empirical process is seen as a random c\`adl\`ag measure-valued process $t \mapsto \mu_{\mathcal{X}^N_t}$, where $(\mathcal{X}^N_t)^{}_t$ is the $N$-particle process given by Assumption \cref{assum:initialdefined}. Two additional assumptions are needed.

\begin{assumption}[Limit continuity] \label[II]{assum:continu} 
The generator $L$ satisfies:
\begin{itemize}
    \item For every $\varphi$ in $\mathcal{F}$, $(x,\mu) \mapsto L_{\mu} \varphi (x)$ is a bounded continuous function. 
    \item For every $\mu$ in $\pb(E)$, $\varphi \mapsto L_{\mu} \varphi$ is a $C_b(E)$-valued continuous mapping.
\end{itemize}
In the first point the topology on $\pb(E)$ is the one induced by the weak convergence of probability measures (cf. Definition \cref{def:convergenceproba}). More precisely, a sequence of probability measures $(\mu_N)_N$ is said to converge towards $\mu$ when $\langle \mu_N,\varphi\rangle\to\langle \mu,\varphi\rangle$ for any test function $\varphi\in C_b(E)$. For the second point, the topology is the topology of the uniform convergence. 
\end{assumption}

This assumption is satisfied in particular for generators which are differential or integral operators with continuous integrable coefficients. This assumption is necessary to take the limit within an equation, instead of using direct c\`adl\`ag characterizations.
The last assumption concerns the limit law. 

\begin{assumption}[Limit uniqueness]  \label[II]{assum:wellposedweakpathwise}
For every $T > 0$ and any $f_0\in\pb(E)$, the limit nonlinear weak PDE
\begin{equation}\label[IIeq]{eq:weakpdemeanfield_summary}
\forall \varphi\in\mathcal{F},\quad \frac{\dd}{\dd t}\langle f_t,\varphi\rangle = \langle f_t, L_{f_t}\varphi\rangle,     
\end{equation}
has at most one unique solution in $C ( [ 0,T ] , \pb ( E ) )$.
\end{assumption}

Note that existence is not required as it will be included in the following propagation of chaos result. However the uniqueness assumption at the limit is very strong, because uniqueness is closely linked with convergence of approximate sequences (such as particle systems). For physically relevant systems, uniqueness is often the most difficult issue.
 
\begin{theorem}[Functional law of large numbers]\label[II]{thm:martingaleweakpathwiseempiricalchaos}
Let $(f^N_0)_N$ be an initial $f_0$-chaotic sequence and let $(\mathcal{X}^N_{t})^{}_t$ be the $E^N$-valued $N$-particle process given by Assumption \cref{assum:initialdefined} with initial distribution $f^N_0$. Assume that Assumptions \cref{assum:initialdefined}, \cref{assum:genebound}, \cref{assum:continu} and~\cref{assum:wellposedweakpathwise} hold true. Then the nonlinear weak PDE \cref{eq:weakpdemeanfield_summary} is well-posed and its solution $(f_t)_t\in C([0,T],\pb(E))$ satisfies:
\[F^{\mu,N}_{[0,T]} \underset{N \to + \infty}{\longrightarrow} \delta_{\left( f_t \right)_{0 \leq t \leq T}} \in \pb(D([0,T],\pb(E)).\]
where $F^{\mu,N}_{[0,T]}$ is the law of the measure-valued empirical process defined by \cref{eq:FmuNlaw}. 
\end{theorem}

To prove this theorem, we will follow a method which can be found in \cite{merino-aceituno_isotropic_2016} and which we adapt to the more abstract present framework. 

\begin{proof} The proof is split into several steps: using (the general) It\=o's formula, we start with some preliminary computations in the linear case which will be used to prove a tightness result on the weak pathwise law $F^{\mu,N}_{[0,T]}$. Then we identify the limit points by controlling the stochastic remainder. 

\medskip
\noindent\textit{\textbf{Step 1.} Some preliminary computations for linear test functions.}
\medskip 

Let us consider a one-particle test functions $\varphi\in\mathcal{F}$ and let us define the average $N$-particle test function 
\[\bar{\varphi}_N : \mathbf{x}^N \mapsto \langle \mu_{\mathbf{x}^N} , \varphi \rangle.\]
By Assumption \cref{assum:initialdefined}, it holds that
\begin{equation}\label[IIeq]{eq:averagedgenerator} \mathcal{L}_N \bar{\varphi}_N \big( \mathbf{x}^N \big) = \sum_{i=1}^N \frac{1}{N} L_{\mu_{\mathbf{x}^N}} \varphi ( x^i ) = \left\langle \mu_{\mathbf{x}^N}, L_{\mu_{\mathbf{x}^N}} \varphi \right\rangle, \end{equation}
so that It\=o's formula gives
\begin{equation}\label[IIeq]{eq:averagedequation}
\big\langle \mu_{\mathcal{X}^N_t} , \varphi \big\rangle =  \big\langle \mu_{\mathcal{X}^N_0} , \varphi \big\rangle + \int_0^t  \big\langle \mu_{\mathcal{X}^N_s} , L_{\mu_{\mathcal{X}^N_s} } \varphi \big\rangle \dd  s + M^{N,\varphi}_t,
\end{equation}
where $M^{N,\varphi}_t$ is a martingale. Using Assumption \cref{assum:initialdefined} again, the carr\'e du champ operator reads:
\[ \Gamma_{\mathcal{L}_N} ( \bar{\varphi}_N , \bar{\varphi}_N ) \big( \mathbf{x}^N \big) = \sum_{i=1}^N {\left[ L_{\mu_{\mathbf{x}^N} } \diamond_i [ {\bar{\varphi}_N}^2 ] \big( \mathbf{x}^N \big) - 2 \big\langle \mu_{\mathbf{x}^N} , L_{\mu_{\mathbf{x}^N} } \varphi \big\rangle L_{\mu_{\mathbf{x}^N} } \diamond_i \bar{\varphi}_N \big( \mathbf{x}^N \big) \right]}. \]
Since $L_{\mu_{\mathbf{x}^N}}$ is linear and vanishes on constant functions, one obtains for any index $i\in\{1,\ldots,N\}$, 
\[N^2 L_{\mu_{\mathbf{x}^N} } \diamond_i [ {\bar{\varphi}_N}^2 ] \big( \mathbf{x}^N \big) = L_{\mu_{\mathbf{x}^N} } [ \varphi^2 ] (x^i) + 2 {\left( \sum_{ j \neq i  } \varphi ( x^j ) \right)} L_{\mu_{\mathbf{x}^N} } \varphi ( x^i ), \]
and
\[2 N^2 \big\langle \mu_{\mathbf{x}^N} , L_{\mu_{\mathbf{x}^N} } \varphi \big\rangle L_{\mu_{\mathbf{x}^N} } \diamond_i \bar{\varphi}_N \big( \mathbf{x}^N \big) = 2 \varphi ( x^i ) L_{\mu_{\mathbf{x}^N} } \varphi ( x^i) + 2 {\left( \sum_{  j \neq i  } \varphi ( x^j) \right)} L_{\mu_{\mathbf{x}^N} } \varphi ( x^i). \]
We conclude that: 
\begin{equation}\label[IIeq]{eq:averagedcarreduchamp} \Gamma_{\mathcal{L}_N} ( \bar{\varphi}_N , \bar{\varphi}_N ) \big( \mathbf{x}^N \big) = \frac{1}{N} {\left\langle \mu_{\mathbf{x}^N} , \Gamma_{L_{\mu_{\mathbf{x}^N} }} \left( \varphi , \varphi \right) \right\rangle}. \end{equation}
The right-hand side goes to $0$ as $N \to + \infty$ thanks to \cref{eq:momentboundmartingale} and Assumptions \cref{assum:initialdefined} and~\cref{assum:genebound}.

\medskip
\noindent\textit{\textbf{Step 2.} Tightness of the sequence $\big( F^{\mu,N}_{[0,T]} \big)^{}_{N \geq 1}$.}
\medskip 

We follow the method of \cite{merino-aceituno_isotropic_2016}. The tightness is proved using Jakuboswki's criterion (Theorem \cref{thm:jakubowski}).
\begin{enumerate}[(i)] 
\item We first prove that for any $\varepsilon>0$, there exists a compact set $K_\varepsilon\subset\pb(E)$ such that 
\[\forall t \in [0,T], \quad \mathbb{P}\big(\mu_{\mathcal{X}^N_t}\in K_\varepsilon\big)>1-\varepsilon.\]
Since for every $M > 0$ and $x_0\in E$, the set
\[ {\left\{ \nu \in \pb ( E ),  \,\, \int_E d^2_E ( x , x_0 ) \nu ( \dd x ) \leq M \right\}} \]
is compact for the weak topology on $\pb \left( E \right)$, it is enough to prove that the uniform $L^2$ moment bound on $f^{1,N}_0$ is propagated on $[ 0 , T ]$ uniformly in $N$. Thanks to Assumption \cref{assum:genebound}, this is the content of \cref{eq:momentboundmartingale}. 

\item The set of linear functions on $\pb(E)$ $\Phi:\nu\mapsto\langle \nu,\varphi\rangle$, $\varphi\in \mathcal{F}$ separates points and is closed under addition. We therefore fix $\varphi\in \mathcal{F}$ and we prove the tightness of the laws in $\pb(D([0,T],\R))$ of the real-valued process $\big(\langle \mu_{\mathcal{X}^N_t},\varphi\rangle\big)_t$. To do that, we use Aldous criterion (Theorem \cref{thm:aldous}) and we use the decomposition \cref{eq:averagedequation}. Since the process is bounded, the first condition is automatically satisfied. Then, let us fix two $\mathscr{F}_t$-adapted stopping times $\tau_1\leq \tau_2\leq\tau_1+\theta$ for a fixed $\theta>0$. On the one hand, by Doob's optional sampling theorem, we have: \[\E\big|M^{N,\varphi}_{\tau_2}-M^{N,\varphi}_{\tau_1}\big|^2 = \E{\left[\big|M^{N,\varphi}_{\tau_2}\big|^2-\big|M^{N,\varphi}_{\tau_1}\big|^2\right]} = \E{\left[\int_{\tau_1}^{\tau_2} \dd \big\langle M^{N,\varphi}\rangle_t\right]}.\]
Using Lemma \cref{lemma:semiMartVariation} and \cref{eq:averagedcarreduchamp}, we deduce that:
\[\E\big|M^{N,\varphi}_{\tau_2}-M^{N,\varphi}_{\tau_1}\big|^2 \leq \E{\left[\int_{\tau_1}^{\tau_2} \Gamma_{\mathcal{L}_N} ( \bar{\varphi}_N , \bar{\varphi}_N ) \big( \mathcal{X}^N_t \big) \dd t\right]}\leq C_\varphi\frac{\theta}{N}.\]
On the other hand, using \cref{eq:averagedgenerator}, Assumption \cref{assum:genebound} and \cref{eq:momentboundmartingale}, one gets (up to changing the constant)
\[ \E {\left[ {\left( \int_{\tau_1}^{\tau_2} \mathcal{L}_N \bar{\varphi}_N {\left( \mathcal{X}^N_t \right)} \dd t \right)}^2 \right]} \leq C_{\varphi} \theta^2. \]
Formula \cref{eq:averagedequation} therefore leads to
\[ \E {\left[ \big\langle \mu_{\mathcal{X}^N_{\tau_2}} - \mu_{\mathcal{X}^N_{\tau_1}} , \varphi \big\rangle^2 \right]} \leq C_{\varphi} {\left[ \theta^2 + \frac{\theta}{N} \right]}.\]
We conclude using the Markov inequality that the conditions of Aldous criterion are fulfilled.
\end{enumerate}

\medskip

\noindent\textit{\textbf{Step 3.} Skorokhod representation for limit points and well-posedness.}

\medskip

\noindent For any $T >0$, the sequence $( F^{\mu,N}_{[0,T]} )_{N \geq 1}$ is thus relatively compact for the weak topology on $\pb ( D ( [ 0 , T ] , \pb ( E ) ) )$. Let $\pi$ be a limit point. Skorokhod  representation theorem provides then a probability space $\Omega$ on which a realisation of $\mu_{\mathcal{X}^N_t}$ converges almost surely (up to an extraction which we do not relabel) towards a $\pi$-distributed $D ( [ 0 , T ] , \pb ( E ) )$-valued random variable $( \bar{f}_t )_{0 \leq t \leq T}$, such that a.s. $\bar{f}_0 = f_0$ thanks to the initial chaos assumption. We want to prove that $\overline{f}_t$ is almost surely a solution of \cref{eq:weakpdemeanfield_summary}. Using Assumption \cref{assum:wellposedweakpathwise}, we will deduce that this PDE is well-posed and that $\pi$ is the Dirac mass at this solution. Using the BDG inequality, it holds that:
\[ \mathbb{E} {\left[ \sup_{0 \leq t \leq T} {\left( M^{N,\varphi}_t \right)}^2 \right]} \leq 4 \mathbb{E} \big[ [ M^{N,\varphi} ]_T \big] = 4 \mathbb{E} \big[ \langle M^{N,\varphi} \rangle_T \big],\] 
where we have used that $( [ M^{N,\varphi} ]_t - \langle M^{N,\varphi} \rangle_t )_{0 \leq t \leq T}$ is a martingale. Using lemma~\cref{lemma:semiMartVariation} and Step $1$ we conclude that
\[ \mathbb{E} {\left[ \sup_{0 \leq t \leq T} {\left( M^{N,\varphi}_t \right)}^2 \right]} \leq \frac{4}{N} \E{\left[\int_0^T {\Big\langle \mu_{\mathcal{X}^N_t} , \Gamma_{L_{\mu_{\mathcal{X}^N_t} }} \left( \varphi , \varphi \right) \Big\rangle} \dd t\right]}\underset{N\to+\infty}{\longrightarrow} 0, \]
where we have used Assumption \cref{assum:genebound}. Up to extracting once more, we can assume that the above $L^2$ convergence is almost sure:
\begin{equation} \label[IIeq]{eq:AsMartVanish}
\sup_{0 \leq t \leq T} M^{N,\varphi}_t  \underset{N\to+\infty}{\longrightarrow}0,\quad\text{a.s.}
\end{equation}
By the first part of Assumption \cref{assum:continu} (continuity with respect to $\mu$), we can take the limit in \cref{eq:averagedequation} and we obtain by dominated convergence that for all $\varphi\in\mathcal{F}$,
\[ \forall t \in \left[ 0,T \right], \quad \langle \bar{f}_t , \varphi \rangle =  \langle \bar{f}_0 , \varphi \rangle + \int_0^t  \langle \bar{f}_s , L_{\bar{f}_s } \varphi \rangle \dd s\quad\text{a.s.}  \]
To recover the limit equation, one needs to invert the ``$\forall \varphi \in \mathcal{F}$'' term and the ``almost surely'' mention. To do that, let us consider a dense countable subset $( \varphi^n )_n$ of $\mathcal{F}$ (it exists because $E$ is a Polish space). The previous steps tells that for each $\varphi^n$, the set of issues in $\Omega$ such that the equality $\frac{\dd}{\dd t} \langle \bar{f}_t , \varphi^n \rangle = \langle \bar{f}_t , L_{\bar{f}_t} \varphi^n \rangle$ does not hold for some $0 \leq t \leq T$ is negligible. By countable union, the set of issues such that this equality does not hold for every $0 \leq t \leq T$ for any of the $\varphi^n$ is still negligible. We then use the continuity with respect to $\varphi$ from Assumption \cref{assum:continu} to conclude by density that $( \bar{f}_t )_{0 \leq t \leq T}$ almost surely solves 
\[ \forall \varphi \in \mathcal{F},\,\, \forall t \in \left[ 0,T \right], \quad \frac{\dd}{\dd t} \langle \bar{f}_t , \varphi \rangle = \langle \bar{f}_t , L_{\bar{f}_t} \varphi \rangle. \]

Theorem \cref{thm:continuitysara} now proves that $t \mapsto \bar{f}_t$ is almost surely continuous: indeed the vanishing of jumps directly stems from the decomposition \cref{eq:averagedequation} together with Equation~\cref{eq:AsMartVanish}, as required by Theorem \cref{thm:continuitysara} (note this condition is reminiscent from Aldous criterion in Step $2$). This shows that any $\pi$-distributed random function is almost surely a solution of \cref{eq:weakpdemeanfield_summary}. Since this solution is unique by Assumption \cref{assum:wellposedweakpathwise}, this shows the well-posedness of \cref{eq:weakpdemeanfield_summary} and proves that  $\pi = \delta_{\left( f_t \right)_{0 \leq t \leq T}}$ where $f_t$ is the unique solution of \cref{eq:weakpdemeanfield_summary}. 
\end{proof}

\begin{example} In addition to the historical works \cite{oelschlager_martingale_1984,gartner_mckean-vlasov_1988} already mentioned, this method has been recently applied in \cite{merino-aceituno_isotropic_2016} in a coagulation-fragmentation model leading to the 4-wave kinetic equation and in \cite{diez_propagation_2020} for a mean-field PDMP on a manifold leading to a BGK equation. This approach also works to prove moderate interaction results \cite{oelschlager_law_1985}. This proof remains true for Boltzmann molecules, in which case the first step (which corresponds to Lemma \cref{lemma:quadraticmeanfield}) has to be replaced by Lemma~\cref{lemma:quadraticboltzmann}.
\end{example}

\begin{remark}[The need for quadratic estimates] This proof may seem surprising because only one-particle test functions on $E$ are considered even though it leads to  a convergence result on random measure-valued process. The quadratic estimates actually lie in the computation of the quadratic variation of the martingale $M^{N,\varphi}_t$ in Step 2 and in the control of the carr\'e du champ operator \cref{eq:averagedcarreduchamp}. This last computation is a special case of the more general result in Lemma \cref{lemma:quadraticmeanfield} about the behaviour of the generator for polynomial test functions of order two. Namely, taking a test function $\varphi_2=\varphi^1\otimes\varphi^2\in\mathcal{F}^{\otimes 2}$ and denoting by 
\[\forall \nu \in\pb(E),\quad R_{\varphi^1\otimes\varphi^2}(\nu) = \langle \nu^{\otimes 2}, \varphi_2\rangle,\]
the associated polynomial function on $\pb(E)$, it holds that: 
\begin{multline*}\mathcal{L}_N {\left[ R_{\varphi^1 \otimes \varphi^2} \circ \boldsymbol{\mu}_N \right]}\big(\mathbf{x}^N\big) = R_{L_{ \mu_{\mathbf{x}^N} } \varphi^1 \otimes \varphi^2} \big( \mu_{\mathbf{x}^N}\big) + R_{\varphi^1 \otimes L_{ \mu_{\mathbf{x}^N} } \varphi^2} \big( \mu_{\mathbf{x}^N}\big) \\+ \frac{1}{N} {\left\langle \mu_{\mathbf{x}^N} , \Gamma_{L_{\mu_{\mathbf{x}^N}}} ( \varphi^1 , \varphi^2 ) \right\rangle}, \end{multline*}
and the carr\'e du champ estimate \cref{eq:averagedcarreduchamp} stems from that since
\begin{multline} \label[IIeq]{eq:mean-fieldcarre}
\Gamma_{\mathcal{L}_N} {\left( \bar{\varphi}^1_N , \bar{\varphi}^2_N \right)}\big(\mathbf{x}^N\big) = \mathcal{L}_N {\left[ R_{\varphi^1 \otimes \varphi^2} \circ \boldsymbol{\mu}_N \right]}\big(\mathbf{x}^N\big)\\ - R_{L_{ \mu_{\mathbf{x}^N} } \varphi^1 \otimes \varphi^2} \big( \mu_{\mathbf{x}^N}\big) - R_{\varphi^1 \otimes L_{ \mu_{\mathbf{x}^N} } \varphi^2} \big( \mu_{\mathbf{x}^N}\big),
\end{multline}
thanks to the mean-field property $\mathcal{L}_N \bar{\varphi}_N = \left\langle \mu_{\mathbf{x}^N} , L_{ \mu_{\mathbf{x}^N} } \varphi \right\rangle$.
Note that purely one particle-related methods are not possible, because the weak convergence of $f^{k,N}_t$ characterizing Kac's chaos has to hold at least with $k \geq 2$ (Lemma \cref{lemma:caractchaos_summary}). 
\end{remark}

\subsubsection*{Strong pathwise empirical chaos.} For the strong pathwise result, the goal is to prove the convergence of the sequence of 
\begin{equation}\label[IIeq]{eq:strongpathwiseFNlaw}F^{N}_{[0,T]} = \mathrm{Law}{\left(\mu_{\mathcal{X}^N_{[0,T]}}\right)}\in \pb\big(\pb(D([0,T],E))\big),\end{equation}
which means that the empirical process is seen as a random empirical measure on the path space $D([0,T],\pb(E))$, to which belongs each component of the $N$-particle process $\mathcal{X}^N_{[0,T]}$ given by Assumption \cref{assum:initialdefined}.

The proof of the following theorem can be found in \cite{graham_stochastic_1997,meleard_asymptotic_1996} and relies on the classical and powerful framework described in \cite{joffe_weak_1986}. This technique has also been used by Sznitman \cite{sznitman_equations_1984} for Boltzmann models (see Section \cref{sec:martingaleboltzmannreview}). The starting point is a strong uniqueness result for the limit martingale problem. 

\begin{assumption}[Uniqueness for the limit martingale problem] \label[II]{assum:uniquemart}
Given an initial value $f_0\in\pb(E)$, there exists at most one probability distribution on the Skorokhod space $f_{[0,T]}\in\pb(D([0,T],E))$ such for all all $\varphi\in\mathcal{F}$,
\begin{equation*} M^\varphi_t := \varphi(\mathsf{X}_t)-\varphi(\mathsf{X}_0) - \int_0^t L_{f_s}\varphi(\mathsf{X}_s)\dd s,\end{equation*}
is a $f_{[0,T]}$-martingale, where $(\mathsf{X}_t)_t$ is the canonical process and for $t\geq0$, $f_t := (\mathsf{X}_t)_\# f_I$.
\end{assumption}

Note once more that this remains a strong assumption and that uniqueness for the limit system is often the hardest property to prove for physical systems. However, existence is not needed as it is included in the following theorem. 

\begin{theorem}[Strong pathwise empirical chaos]\label[II]{thm:martingalestrongpathwiseempiricalchaos}
Let $(f^N_0)_N$ an initial $f_0$-chaotic sequence and let $\mathcal{X}^N_{[0,T]}\in D([0,T],E^N)$ be the $N$-particle process given by Assumption \cref{assum:initialdefined} with initial distribution $f^N_0$. Assume that Assumptions \cref{assum:initialdefined}, \cref{assum:genebound}, \cref{assum:continu} and \cref{assum:uniquemart} hold true. Then the nonlinear mean-field martingale problem stated in Assumption \cref{assum:uniquemart} is well-posed and its solution $f_{[0,T]}\in\pb(D([0,T],E))$ satisfies:
\[F^N_{[0,T]}\underset{N\to+\infty}{\longrightarrow} \delta_{f_{[0,T]}}\in \pb(\pb(D([0,T],E))),\]
where $F^N_{[0,T]}$ is the pathwise empirical law defined by \cref{eq:strongpathwiseFNlaw}.
\end{theorem}

\begin{proof} The first step is to show the tightness of the sequence $\big(F^N_{[0,T]}\big)_N$ in the space $\pb \big( \pb ( D ( [ 0 , T ] , E ) ) \big)$. 
\medskip

\noindent\textit{\textbf{Step 1.} Tightness.}
\medskip

Thanks to the exchangeability and Lemma \cref{lemma:empiricaltightness}, it is sufficient to prove the tightness of the sequence
\[F^{1,N}_{[0,T]}= \mathrm{Law} {\left( X^{1,N}_{[ 0 ,T ]} \right)}\in \pb(D([0,T],E)). \]
The process $( X^{1,N}_t )_{0 \leq t \leq T}$ can be characterized as a $D$-semimartingale (see Definition \cref{def:semiMartD}) thanks to Assumption \cref{assum:initialdefined} by taking $\varphi_N=\varphi\otimes 1^{\otimes(N-1)}$ as a test function, given a one-particle test function $\varphi\in \mathcal{F}$. It implies that
\[  {M}^{\varphi,1,N}_t  := \varphi \big( X^{1,N}_t \big) - \varphi \big( X^{1,N}_0 \big) - \int_0^t L_{\mu_{\mathcal{X}^N_s}} \varphi \big( X^{1,N}_s \big) \text{ds}, \]
is a martingale. The Joffe-Metivier criterion \cref{thm:joffecriterion} can then be applied: Assumption \cref{assum:genebound} implies the tightness of $\big(F^{1,N}_{[0,T]}\big)_{N \geq 1}$. Moreover, using Lemma \cref{lemma:semiMartVariation} and Assumption \cref{assum:initialdefined} the predictable quadratic variation is given by 
\begin{align*} \big\langle  {M}^{\varphi,1,N}  \big\rangle_t &= \int_0^t \Gamma_{\mathcal{L}_N} \left( \varphi \otimes 1^{\otimes(N-1)} , \varphi  \otimes 1^{\otimes(N-1)} \right) \left( \mathcal{X}^N_s \right) \dd s \\&= \int_0^t \Gamma_{L_{\mu_{\mathcal{X}^N_s}}} \left( \varphi , \varphi \right) \left( X^{1,N}_s \right) \dd s.\end{align*}
Similarly, for $k\leq N$, taking $\varphi_N=1^{\otimes(k-1)}\otimes\varphi\otimes 1^{\otimes(N-k)}$, the following process is a martingale: 
\[{M}^{\varphi,k,N}_t  := \varphi \big( X^{k,N}_t \big) - \varphi \big( X^{k,N}_0 \big) - \int_0^t L_{\mu_{\mathcal{X}^N_s}} \varphi \big( X^{k,N}_s \big) \text{ds}. \]
The predictable cross variation can be computed the same way taking $\varphi_N = \varphi \otimes \psi \otimes 1^{\otimes(N-2)}$, 
\begin{equation}\label[IIeq]{eq:crossbracketsstrongpathwise} \big\langle {M}^{\varphi,1,N},{M}^{\psi,2,N} \big\rangle_t = \int_0^t \Gamma_{\mathcal{L}_N} ( \varphi \otimes 1^{\otimes(N-1)} , 1\otimes \psi \otimes 1^{\otimes(N-2)} ) \big( \mathcal{X}^N_s \big) \dd s = 0. \end{equation}
It will be useful for Step $3$.
\medskip

\noindent\textit{\textbf{Step 2.} Skorokhod representation for limit points and continuity points.}
\medskip

Let $\pi\in\pb\big(\pb(D([0,T],E))\big)$ be a limit point of 
$\big(F^N_{[0,T]}\big)_{N \geq 1}$. Using Skorokhod representation theorem, it is possible to consider a probability space and a $\pi$-distributed random variable $f_{[0,T]}\in\pb(D([0,T],E))$ such that (for the weak topology):
\[\mu_{\mathcal{X}^N_{[0,T]}} \underset{N\to+\infty}{\longrightarrow} f_{[0,T]}\quad \text{a.s.}\]
Consider now $n \geq 1$ with some positive real numbers $s_1 \leq \ldots \leq s_n \leq s < t$ and some functions $\varphi,\varphi^1, \ldots , \varphi^n \in \mathcal{F}$ and let us consider the function:
\begin{multline*} 
F_{s_1,\ldots , s_n , s , t} : \nu \in \pb ( D ( [ 0 , T ] , E ) ) \\
\mapsto {\left\langle \nu ,  {\left( \varphi ( \mathsf{X}_t ) - \varphi ( \mathsf{X}_s ) - \int_s^t L_{\nu_r} \varphi ( \mathsf{X}_r ) \dd r \right)} \varphi^1 ( \mathsf{X}_{s_1} ) \ldots \varphi^n ( \mathsf{X}_{s_n} ) \right\rangle}\in\R,
\end{multline*}
where $\nu_r=(\mathsf{X}_r)_\#\nu\in\pb ( E )$ denotes the $r$-marginal of $\nu$. Thanks to Assumption~\cref{assum:continu}, the map $F_{s_1,\ldots , s_n , s , t}$ would be continuous if the coordinates maps $\mathsf{X} \mapsto \mathsf{X}_t$ were continuous. However these maps are not continuous in general for the Skorokhod topology. For $u$ in $\R_+$, consider the event
\[ A_u := \Big\{ { Q \in \pb(D([0,T],E)) \, : \, Q\big( \{ X \in D([0,T],E) \, : \, |\Delta X_u| > 0 \} \big) > 0 } \Big\}. \]
The map $F_{s_1,\ldots , s_n , s , t}$ will thus be $\pi$-a.s. continuous when $s_1,\ldots , s_n , s , t$ are taken out of the set
\[ J := \{ u \in \R_+ ,\, \, \pi (A_u) > 0 \}, \]
Adapting a proof from \cite{graham_stochastic_1997}, let us show that $J$ is at most countable. The key idea is that given $k\geq1$, \emph{a c\`adl\`ag function $X$ on a compact time-interval admits a finite numbers of jumps with amplitudes bigger than $1/k$}. Let us denote by $\mathcal{J}(X,1/k,[0,k])$ the number of jumps of $X$ with amplitude $|\Delta X_t | > 1/k$ for $t\in[0,k]$. Define then for $m\geq1$, 
\begin{multline*}
A_{u}^{k,m} := \Big\{Q \in \pb(D([0,T],E)) \, : \\ \, Q\big( \{ X \in D([0,T],E) \, : \, |\Delta X_u| > 1/k \text{ and }\mathcal{J}(X,1/k,[0,k]) \leq mk \} \big) > 1/k  \Big\}.
\end{multline*}
Moreover, the following properties hold. 
\begin{itemize}
    \item The sequence ${\left( \bigcup_{m \geq 1} A_{u}^{k,m} \right)}_{k \geq 1}$ is non-decreasing (for the set inclusion) in $k$.
    \item For a fixed $k\geq1$, the sequence $( A_{u}^{k,m} )_{m \geq 1}$ is non-decreasing in $m$.
    \item The set $A_u$ can be decomposed as
    \[ A_u = \bigcup_{k \geq 1} \bigcup_{m \geq 1} A_{u}^{k,m}. \]
\end{itemize}
The monotonic convergence of probability measures thus gives
\[ \pi(A_u) = \lim_{k \to + \infty} \pi{\left( {\bigcup_{m \geq 1} A_{u}^{k,m}} \right)} = \lim_{k \to + \infty} \lim_{m \to + \infty} \pi{\left( {A_{u}^{k,m}} \right)}. \]
Introducing
\[ J^{k,m} := \big\{ u \in [0,k],\, \, \pi{\left( {A_{u}^{k,m}} \right)} > 1/k \big\}, \]
the same trick leads to
\[ J = \bigcup_{k \geq 1} \bigcup_{m \geq 1} J^{k,m}. \]
Let us now prove that $J^{k,m}$ is finite. If it were not, there would exist a sequence $(u_n)_{n \geq 1}$ of pairwise distinct numbers in $[0,k]$ such that
\[ \forall n \geq 1, \quad \pi{\left(  A_{u_n}^{k,m} \right)} > 1/k. \]
We apply now (the consequence of) Lemma \cref{lem:Intersection} to $(A_{u_n}^{k,m})_{n \geq 1}$ in the probability space $\Omega =\pb(D([0,T],E))$ for $P = \pi$: for every $n \geq 1$, there exists an intersection involving $n$ of the $A_{u_i}$ which has positive $\pi$-measure, and this leads to the existence of integers $i^n_1<\ldots<i^n_n$ such that 
\begin{multline*} \pi\Big(\Big\{ Q \in \pb(D([0,T],E)) \, : \forall 1 \leq j \leq n, \\  Q{\left( \big\{ X \in D([0,T],E): \, |\Delta X_{u_{i^n_j}}| > \frac{1}{k} \text{ and  } \mathcal{J}(X,1/k,[0,k]) \leq mk \big\} \right)} > \frac{1}{k} \Big\}\Big) > 0. \end{multline*}
The same reasoning can be applied within the probability, considering a probability measure $Q \in \pb(D([0,T],E))$ such that for all $j\in\{1,\ldots,n\}$, 
\[Q{\left( \big\{ X \in D([0,T],E) \, : \, |\Delta X_{u_{i^n_j}}| > 1/k \text{ and  } \mathcal{J}(X,1/k,[0,k]) \leq mk \big\} \right)} > 1/k, \]
and applying Lemma \cref{lem:Intersection} with $P = Q$ and $\Omega = D([0,T],E)$ to the events
\[ {\left(\big\{ X \in D([0,T],E) \, : \, |\Delta X_{u_{i^n_j}}| > 1/k \text{ and  } \mathcal{J}(X,1/k,[0,k]) \leq mk \big\}\right)}_{1 \leq j \leq n}. \]
Since $n$ can be taken arbitrarily large, this allows to consider an arbitrary large intersection of these events which has $Q$-positive measure. This is contradictory since the number of jumps with amplitude bigger than $1/k$ allowed on $[0,k]$ is at most $mk$.

This proves the finiteness of $J^{k,m}$ for any $k,m \geq 1$, so $J$ is at most countable by countable union. This implies the $\pi$-almost sure continuity of $F_{s_1,\ldots , s_n , s , t}$ for $s_1,\ldots , s_n , s , t$ outside an at most countable set $D_{\pi}$. Outside of this set
\[ F_{s_1,\ldots , s_n , s , t} {\left( \mu_{\mathcal{X}^N_{\left[ 0 ,T \right]}} \right)} \underset{N\to+\infty}\longrightarrow  F_{s_1,\ldots , s_n , s , t} {\left( f_{\left[ 0 ,T \right]} \right)}. \]
Note this argument is more general and can be adapted to Boltzmann models as in \cite{graham_stochastic_1997}, see also Section \cref{sec:martingaleboltzmannreview}. 
\medskip

\noindent\textit{\textbf{Step 3.} Identifying the limit points using the martingale problem.}
\medskip

To recover $\langle \pi , F_{s_1,\ldots , s_n , s , t} \rangle$, it is now sufficient to take the expectation. Using the Cauchy-Schwarz inequality and then Fatou's lemma, it holds that for $s_1,\ldots , s_n , s , t$ outside of $D_{\pi}$
\begin{align*}
\langle \pi , | F_{s_1,\ldots , s_n , s , t} | \rangle^2 &\leq \langle \pi , F_{s_1,\ldots , s_n , s , t}^2 \rangle \\
&\leq \lim_{N \to \infty} \E_{f_{[ 0 ,T ]}^N} {\left[ F_{s_1,\ldots , s_n , s , t}^2 {\left( \mu_{\mathcal{X}^N_{[ 0 ,T ]}} \right)} \right]} \\
&= \E_{f_{[ 0 ,T ]}^N} {\left[ {\left\{ \frac{1}{N} \sum_{i=1}^N  \big(  {M}^{\varphi,i,N}_t -  {M}^{\varphi,i,N}_s \big) \varphi^1 ( \mathsf{X}^{1,N}_{s_1} ) \ldots \varphi^n ( \mathsf{X}^{1,N}_{s_n} ) \right\}}^2 \right]} \\
&= \frac{1}{N} \E_{f_{[ 0 ,T ]}^N} {\left[ \Big\{ \big(  {M}^{\varphi,1,N}_t -  {M}^{\varphi,1,N}_s \big) \varphi^1 ( \mathsf{X}_{s_1} ) \ldots \varphi^n ( \mathsf{X}_{s_n} ) \Big\}^2 \right]} \\
& \quad + \frac{N-1}{N} \E_{f_{[ 0 ,T ]}^N}\Big[\big( {M}^{\varphi,1,N}_t - {M}^{\varphi,1,N}_s \big) \big({M}^{\varphi,2,N}_t - {M}^{\varphi,2,N}_s\big)\times  \\
&\phantom{\quad + \frac{N-1}{N} \E_{f_{[ 0 ,T ]}^N}} \times\varphi^1 ( \mathsf{X}^{1,N}_{s_1} ) \ldots \varphi^n ( \mathsf{X}^{1,N}_{s_n} ) \varphi^1 ( \mathsf{X}^{2,N}_{s_1} ) \ldots \varphi^n ( \mathsf{X}^{2,N}_{s_n} )\Big]. 
\end{align*}
Assumption \cref{assum:genebound} ensures that $\mathcal{M}^{\varphi,1,N}_t$ is bounded in $L^2$ by the carr\'e du champ vector, so that the first term on the right-hand side vanishes as $N \to +\infty$. For the second one, we write:
\[ \E {\left[ {M}^{\varphi,1,N}_t {M}^{\varphi,2,N}_s \big| \sigma {\big( ( \mathsf{X}_r )_{0 \leq r \leq s} \big)} \right]} =  {M}^{\varphi,1,N}_s {M}^{\varphi,2,N}_s. \]
Then, since the cross-brackets \cref{eq:crossbracketsstrongpathwise} are equal to zero, taking the expectation leads to: 
\[\E {\left[ {M}^{\varphi,1,N}_t {M}^{\varphi,2,N}_s \right]}=0.\]
So the second term is actually equal to zero. 

This proves $F_{s_1,\ldots , s_n , s , t}$ is $0$ $\pi$-almost surely: this holds for every $0 \leq s_1 \leq \ldots \leq s_n \leq s < t$ outside the countable set $D_{\pi}$, and every $\varphi,\varphi^1,\ldots,\varphi^n$ in $\mathcal{F}$. By density of $\mathcal{F}$ in $C_b(E)$ and $\R_+ \setminus D_{\psi}$ in $\R_+ $, this allows to replace $\varphi^1 ( \mathsf{X}_{s_1} ) \ldots \varphi^n ( \mathsf{X}_{s_n} )$ by any $\sigma {\left( ( \mathsf{X}_r )_{0 \leq r \leq s} \right)}$-measurable function to obtain
\begin{multline*} \E {\left[ \varphi ( \mathsf{X}_t ) - \varphi ( \mathsf{X}_0 ) - \int_0^t L_{f_r} \varphi ( \mathsf{X}_r ) \dd r \Big| \sigma {\big( ( \mathsf{X}_r )_{0 \leq r \leq s} \big)} \right]} \\= \varphi ( \mathsf{X}_s ) - \varphi ( \mathsf{X}_0 ) - \int_0^s L_{f_r} \varphi ( \mathsf{X}_r ) \dd r, \end{multline*}
for $\pi$-almost every pathwise law $f_{[ 0 , T ]}$, every $\varphi$ in $\mathcal{F}$ and every $s,t$ outside the countable set $D_{\pi}$. A limit $f_{[ 0 , T ]}$-distributed process being  c\`adl\`ag,  this is sufficient to prove that $\pi$-almost every pathwise law $f_{[ 0 , T ]}$ solves the martingale problem of Assumption \cref{assum:uniquemart}. Consequently, this proves existence for this problem and since uniqueness holds, the problem is well-posed and $\pi$ has to be a Dirac measure $\delta_{f_{[ 0 , T ]}}$, which concludes the proof. 
\end{proof}

\begin{example} In \cite{chiang_mckean-vlasov_1994}, the argument is reversed: Theorem \cref{thm:martingalestrongpathwiseempiricalchaos} states only an existence result which is then used to prove the strong uniqueness result using a synchronous coupling argument. Propagation of chaos follows. This allows to treat the case of McKean-Vlasov diffusions with more general interaction functions. 
\end{example}

\subsubsection{Gradient systems as gradient flows}\label[II]{sec:gradientflows}

In this section we consider McKean-Vlasov gradient systems \cref{eq:mckeanvlasov_summary} with: 
\[b(x,\mu) = -\nabla V(x) -\nabla W\star\mu(x),\quad \sigma=\sqrt{2}I_d.\]

The following theorem states that the McKean-Vlasov gradient systems can be characterised as gradient flows at the three levels of description: the nonlinear solution of the limit equation, the $N$-particle distribution and the $\pb(\pb(E))$-valued curve inherited from the nonlinear semigroup generated by the limit PDE. The notion of gradient-flow is recalled in Section \cref{sec:gradientflowsintro}. 

\begin{theorem}[McKean-Vlasov as gradient flows]\label[II]{thm:mvgradientflows} Let $f_0\in\pb_4(\R^d)$ and $f^N_0\in\pb^\mathrm{sym}_4(\R^{dN})$ admit a density. Let $V,W$ be respectively a confinement potential and an interaction potential which are both bounded below, $\lambda$-convex for some $\lambda\in\R$. Assume also that $W$ is symmetric and satisfies the doubling condition
\[\exists C>0,\,\,\forall x,y\in\R^d,\quad W(x+y)\leq C(1+W(x)+W(y)).\]
\begin{enumerate}
    \item  In $\pb_2(\R^d)$, let the energy $\mathcal{F}$ be defined by: 
    \[\mathcal{F}(\rho) :=  \int_{\R^d}\rho(x)\log\rho(x)\dd x + \int_{\R^d} V(x)\rho(x)\dd x + \frac{1}{2}\int_{\R^d} W(x-y)\rho(x)\rho(y)\dd y,\]
    whenever $\rho$ has a density with respect to the Lebesgue measure and $\mathcal{F}(\rho)=+\infty$ otherwise. Then there exists a unique $2\lambda$-gradient flow $f_t$ for $\mathcal{F}$ such that $\lim_{t\downarrow0}f_t = f_0$ in $\pb_2(\R^d)$. Moreover $f_t$ is a weak distributional solution of the nonlinear McKean-Vlasov-Fokker-Planck equation \cref{eq:mckeanvlasov-pde_summary}.
    \item In $\pb_2(\R^{dN})$, let the energy $\mathcal{F}^N$ be defined by: 
    \begin{multline*}\mathcal{F}^N(\rho^N) := \frac{1}{N}\int_{\R^{dN}}\rho^N(\mathbf{x}^N)\log\rho^N(\mathbf{x}^N)\dd\mathbf{x}^N + \frac{1}{N}\sum_{i=1}^N\int_{\R^{dN}}V(x^i)\rho^N(\dd\mathbf{x}^N)\\+\frac{1}{2N^2}\sum_{i\ne j}\int_{\R^{dN}}W(x^i-x^j)\rho^N(\dd\mathbf{x}^N).\end{multline*}
    Then there exists a unique $3\lambda$-gradient flow $f^N_t$ for $\mathcal{F}^N$ such that $\lim_{t\downarrow0}f^N_t = f^N_0$ in $\pb_2(\R^d)$. Moreover $f^N_t$ is a weak distributional solution of the $N$-particle Liouville equation \cref{eq:liouville_summary}. 
    \item In $\pb_2(\pb_2(\R^d))$, let the energy $\mathcal{F}^\infty$ be defined by
    \[\mathcal{F}^\infty(\pi) := \int_{\pb_2(\R^d)}\mathcal{F}(\rho)\pi(\dd\rho).\]
    Then there exists a unique $3\lambda$-gradient flow $\pi_t$ for $\mathcal{F}^\infty$ such that $\lim_{t\downarrow0}\pi_t = \pi_0:=\lim_{N\to+\infty}F^N_0\in\pb_2(\pb_2(\R^d))$. Moreover $\pi_t$ is explicitely given by 
    \[\pi_t = (\overline{S}_t)_\#\pi_0,\]
    where $\overline{S}_t : \pb_2(\R^d)\to\pb_2(\R^d)$ is the nonlinear semi-group generated by the McKean-Vlasov-Fokker-Planck equation \cref{eq:mckeanvlasov-pde_summary} in the sense that the solution of this PDE is given by $f_t = \overline{S}_t(f_0)$ (see also Section \cref{sec:limitsemigroup} and Section \cref{appendix:timeinhomoegeneousmarkov}). 
\end{enumerate}
\end{theorem}

The first two points are classical, see \cite[Theorem 6.31]{pajot_lecture_2014} or \cite[Chapter 11]{ambrosio_gradient_2008}. The third point is proved in \cite[Lemma 19]{carrillo_-convexity_2020}. Within this setting, propagation of chaos is proved in \cite[Theorem 2]{carrillo_-convexity_2020}. 

\begin{theorem}[\cite{carrillo_-convexity_2020}] Under the same assumptions as in Theorem \cref{thm:mvgradientflows} and with the same notations, for all $T>0$ it holds that
\[\lim_{N\to+\infty} \sup_{t\in[0,T]} W_2(f^N_t,\pi^N_t) = 0,\]
where $\pi^N_t$ is the $N$-th moment measure of $\pi_t$ defined by:
\[\pi^N_t = \int_{\pb_2(\R^{d})} \rho^{\otimes N} \pi_t(\dd\rho) = \int_{\pb_2(\R^d)} \overline{S}_t(\rho)^{\otimes N}\pi_0(\dd\rho).\]
In particular if $f^N_0$ is $f_0$-chaotic, then
\[\lim_{N\to+\infty} \sup_{t\in[0,T]} W_2(f^N_t,f_t^{\otimes N}) = 0.\]
\end{theorem}

The key result is \cite[Lemma 13]{carrillo_-convexity_2020}. It is based on Ascoli's theorem in the space $C\big([0,T],\pb_2(\pb_2(\R^d))\big)$ and states that there exists $\pi=(\pi_t)_t\in C([0,T],\pb_2(\pb_2(\R^d)))$ such that the law of the empirical process $F^N_t := \mathrm{Law}(\mu_{\mathcal{X}^N_t})$ satisfies
\[\sup_{t\in[0,T]}\mathcal{W}_2^2\big(F^N_t,\pi_t)\underset{N\to+\infty}{\longrightarrow}0,\]
up to extracting a subsequence and where $\mathcal{W}_2\equiv W_{2,W_2}$ is the Wasserstein-2 distance on $\pb_2(\pb_2(\R^d))$ for the $W_2$ distance on $\pb_2(\R^d)$ (see Definition \cref{def:spaceswasserstein}). Using the fact that the push-forward by the empirical measure map is an isometry for the Wasserstein distance, it is possible to prove that this convergence is equivalent to the convergence of the $N$-particle distribution: 
\[\sup_{[0,T]} W_2^2\big(f^N_t,\pi_t^{N}\big) \underset{N\to+\infty}{\longrightarrow}0.\]
See for instance \cite[Lemma 10]{carrillo_-convexity_2020} or \cite[Theorem 5.3]{hauray_kacs_2014}. Once a converging subsequence is extracted, the limit $\pi$ is identified as the unique gradient flow with energy $\mathcal{F}^\infty$ by passing to the limit in the Evolution Variational Inequality \cref{eq:evi} which characterises the gradient-flow with energy $\mathcal{F}^N$ using a $\Gamma$-convergence result \cite[Lemma 16]{carrillo_-convexity_2020}. 

\subsection{Entropy bounds with very weak regularity}\label[II]{sec:jabin}

In this section, the problem is to weaken the regularity assumptions of Theorem \cref{thm:mckean} for the McKean-Vlasov diffusion with coefficients:
\begin{equation}\label[IIeq]{eq:bsectionjabingeneral}b(x,\mu) = \tilde{b}(x,K\star\mu(x)),\quad \sigma = I_d,\end{equation}
where $\tilde{b}:\R^n\to\R^d$ is still assumed to be Lipschitz but the interaction kernel $K:\R^d\times\R^d\to\R^n$ has a very weak regularity. Among the methods that are introduced in Section \cref{sec:proving_summary}, the entropy-based methods are particularly adapted to handle weak regularity. From a probabilistic point of view, the relative entropy functional (Definition \cref{def:entropyfisher_summary}) naturally arises as the rate function of a large deviation principle and entropy bounds are classically obtained as an application of Girsanov theorem which does not require any particular regularity assumptions (see for instance Lemma \cref{lemma:entropyboundgirsanov_summary}). The content of this section will be based on the entropy methods introduced in \cite{jabin_mean_2016,jabin_quantitative_2018} and which can be seen as an analytical counterpart of these observations. The main object of study will therefore be the Liouville equation (rather than the system of SDEs) for which it is possible to define a notion of entropy solution which is well adapted to the context (see Definition \cref{def:entropysolution} below). For the limit solution $f$ of \cref{eq:mckeanvlasov-pde_summary}, things are easier because it is possible to propagate the regularity of $f_0$ and it is therefore possible to assume that $f$ can be taken very regular. The starting point is the evolution equation \cref{eq:computeH_summary} satisfied by $H(f^N_t|f^{\otimes N}_t)$. 

\begin{remark}
We would like to emphasize the importance of the Girsanov theorem as the underlying idea although this section contain purely analytical arguments. A probabilistic pathwise version of the results presented in this section which are directly based on Girsanov theorem can be found in \cite{lacker_strong_2018} and \cite{jabir_rate_2019}. These works focus more on the ability to take an abstract general interaction function $b$ rather than on regularity questions. We will discuss these aspects in Section \cref{sec:chaosviagirsanov}. 
\end{remark}

\subsubsection{An introductory example in the \texorpdfstring{$L^\infty$}{bounded} case}\label[II]{sec:jabinintro} 

As an introductory example to the work of \cite{jabin_quantitative_2018}, let us first start with the case where $K\in L^\infty(\R^d)$ (that is, compared to McKean's theorem, the Lipschitz and continuity assumptions on $K$ are removed). The following computations are essentially formal but the ideas will be used in a rigorous framework in the next paragraph. In particular, we assume that $f_t$ is regular enough so that $\log f_t$ can be taken as a test function in the weak Liouville equation \cref{eq:liouville_summary}. Using the entropy dissipation relation on $f^N_t$ which defines the notion of entropy solution, the computations in the proof Lemma \cref{lemma:computeH} which leads to~\cref{eq:computeH_summary} can be fully justified. With $\alpha=1$ in the conclusion, we recall that we obtained (in integrated form): 
\begin{equation}\label[IIeq]{eq:computeHintrojabin}H\big(f^N_t|f^{\otimes N}_t\big) \leq H\big(f^N_0|f^{\otimes N}_0\big)+N\int_0^t \E_{f^N_s}{\left[\big| b\big(X^1_s,\mu_{\mathcal{X}^N_s}\big)-b(X^1_s,f_s)\big|^2\right]} \dd s.\end{equation}
The goal is to find a uniform bound (in $N$) for the expectation on the right-hand side in terms of $H(f^N_s|f^{\otimes N}_s)$. Gronwall lemma will then gives a bound on the entropy and propagation of chaos will follow by Lemma \cref{lemma:entropyPinskerCsiszar}. Note that this quantity is not very far from \cref{eq:mckeanthmszniterror} in the proof of McKean's theorem. The main difference is that the expectation on the right-hand side is an expectation with respect to $f^N_t$ instead of an expectation with respect to $f^{\otimes N}_t$. Of course the latter is more amenable as it allows to use the very simple but efficient argument of Sznitman based on the law of large number and which uses only the boundedness of $K$. The next idea is thus a change of measure argument which is the content of \cite[Lemma 1]{jabin_quantitative_2018}: for all $\eta>0$ and all $\varphi_N\in L^\infty(E^N)$,
\begin{equation}\label[IIeq]{eq:changeofmeasure} \int_{E^N} \varphi_N \big( \mathbf{x}^N \big) f^N_t \big( \dd\mathbf{x}^N \big) \leq \frac{1}{\eta N} {\left( H \big( f^N_t | f_t^{\otimes N} \big) + \log \int_{E^N} \e^{\eta N \varphi_N \big( \mathbf{x}^N \big)} f^{\otimes N}_t  \big(\dd \mathbf{x}^N \big) \right)}.\end{equation}
This identity is a straightforward rewriting of $H(f^N_t|u) \geq 0$ (which is always true) for the probability density $u:=\e^{\eta N \varphi_N} f^{\otimes N}_t / \int_{E^N} \e^{\eta N \varphi_N } f^{\otimes N}_t$. Using this relation gives: 
\begin{align*}H\big(f^N_t|f^{\otimes N}_t\big) &\leq H\big(f^N_0|f^{\otimes N}_0\big)+\frac{1}{\eta}\int_0^t H\big(f^N_s|f^{\otimes N}_s\big)\dd s
\\
&\quad+\frac{1}{\eta}\int_0^t\log\int_{E^N}\exp{\left(\eta N\big| b\big(x^1,\mu_{\mathbf{\mathbf{x}}^N}\big)-b(x^1,f_s)\big|^2\right)}f^{\otimes N}_s\big(\dd\mathbf{x}^N\big)\dd s.\end{align*}
Expanding the square and using \cref{eq:bsectionjabingeneral} as in the proof of McKean's theorem leads to: 
\[H\big(f^N_t|f^{\otimes N}_t\big) \leq H\big(f^N_0|f^{\otimes N}_0\big)+\frac{1}{\eta}\int_0^t H\big(f^N_s|f^{\otimes N}_s\big)\dd s+\frac{1}{\eta}\int_0^t \log Z_N \dd s,\]
with 
\[Z_N := \int_{E^N} \exp{\left(\frac{\eta\|\tilde{b}\|^2_\mathrm{Lip}}{N}\sum_{i,j=1}^N\psi(x^1,x^i)\psi(x^1,x^j)\right)}f_s^{\otimes N}\big(\dd \mathbf{x}^N\big),\]
where
\begin{equation}\label[IIeq]{eq:cancellationspsi}\psi(x,y) = K(x,y)-K\star f_s(x).\end{equation}
The goal is to prove that $Z_N$ is bounded; the conclusion will then follow by Gronwall lemma. Note that there is still the cancellation 
\[\forall x\in E,\quad \int_{E} \psi(x,y)f_s(\dd y) = 0,\]
but it is not possible to use it directly as in the proof of McKean's theorem because now, this quantity appears inside the exponential. Note however that $Z_N$ can be seen as the partition function of a Gibbs measure with a potential which, up to the first variable which plays a special role, is very much reminiscent of a polynomial potential of order two in Theorem \cref{thm:ldpgibbspolynomial}. The second assertion in Theorem \cref{thm:chaosfromLDP} precisely implies that $Z_N$ is bounded. However, in this context, there is a way to bound $Z_N$ more directly (for $\eta$ small enough): this is the content of \cite[Theorem~3]{jabin_quantitative_2018}. The proof is based on the series expansion: 
\begin{multline*}\exp{\left(\frac{\eta\|\tilde{b}\|^2_\mathrm{Lip}}{N}\sum_{i,j=1}^N\psi(x^1,x^i)\psi(x^1,x^j)\right)} \\= \sum_{k=0}^{+\infty}\frac{1}{k!}\frac{\eta^k\|\tilde{b}\|^{2k}_\mathrm{Lip}}{N^k}{\left(\sum_{i,j=1}^N\psi(x^1,x^i)\psi(x^1,x^j)\right)}^{k}.\end{multline*}
Then, by expanding the power term Jabin and Wang recover polynomial terms in $\psi$ and by separating the terms with $k<N$ from the ones with $k>N$, they use combinatorial arguments to identify the right cancellations (using \cref{eq:cancellationspsi}) which lead to the conclusion. In conclusion, there exists a constant $C>0$ such that 
\[H\big(f^N_t|f^{\otimes N}_t\big) = H\big(f^N_0|f^{\otimes N}_0\big)+\frac{1}{\eta}\int_0^t H\big(f^N_s|f^{\otimes N}_s\big)\dd s+Ct,\]
and the result follows.

\subsubsection{With \texorpdfstring{$W^{-1,\infty}$}{less than bounded} kernels} 

In \cite{jabin_quantitative_2018}, the above arguments are presented in a completely rigorous framework in the fully linear case 
\[b(x,\mu) = F(x) + K\star\mu(x),\quad \sigma = I_d,\]
where the state space is the $d$-dimensional torus $E=\mathbb{T}^d$. The force term $F$ is implicitly regular (to ensure that $f$ can be taken regular) but the interaction kernel $K:\mathbb{T}^d\to\mathbb{T}^d$ is less than bounded, it is assumed to be an element of the following functional space. Although it didn't produce a quantitative estimate, a similar idea had been used in \cite{fournier_propagation_2014}.

\begin{definition} A vector field $K$ such that $\int_{\mathbb{T}^d} K = 0$ is said to belongs to $\dot{W}^{-1,\infty}(\mathbb{T}^d)$ when there exists a matrix field $V$ in $L^\infty(\mathbb{T}^d)$ such that $K=\nabla\cdot V$. The definition extends similarly to scalar functions. 
\end{definition}

The regularity on $K$ is extremely weak. It includes the case $K\in L^\infty$ which is the original framework of \cite{jabin_mean_2016} but it is also possible to consider singular kernels and in particular the Biot-Savart kernel in dimension 2: 
\[K(x) = \alpha\frac{x^\perp}{|x|^2}+K_0(x),\]
where $x^\perp$ is the rotation of $x\in\mathbb{R}^2$ by $\pi$ and $K_0$ is a correction which makes $K$ periodic. Other examples of relevant kernels include collision-like kernels where two particles interact when they are exactly at a given distance. We refer the interested reader to \cite[Section 1.3]{jabin_quantitative_2018} and to the end of this section for further examples. It is not easily possible to construct SDE solutions of the particle system with this weak regularity, Jabin and Wang thus introduce the following notion of entropy solution for the solution of the Liouville equation.

\begin{definition}[Entropy solution]\label[II]{def:entropysolution} A probability density $f^N_t\in L^1(\mathbb{T}^{dN})$ for a time $t\in[0,T]$ is an entropy solution to the Liouville equation \cref{eq:liouville_summary} when it solves \cref{eq:liouville_summary} in the sense of distributions and for almost every $t\leq T$, 
\begin{align}\label[IIeq]{eq:entropysolutionineq}
&\int_{\mathbb{T}^{dN}} f^N_t\log f^N_t + \frac{1}{2}\int_0^t\int_{\mathbb{T}^{dN}} \frac{|\nabla f^N_t|^2}{f^N_t}\,\,\dd s \nonumber \\
&\leq \int_{\mathbb{T}^{dN}} f^N_0\log f^N_0 -\frac{1}{N}\sum_{i,j=1}^N \int_0^t\int_{\mathbb{T}^{dN}} {\left(\nabla\cdot F(x^i)+\nabla\cdot K(x^i-x^j)\right)}f^N_t\big(\mathbf{x}^N\big)\dd \mathbf{x}^N\dd s.
\end{align}
\end{definition}

It is much easier to prove that there exists an entropy solution, this typically comes from a regularisation argument with a smoothened kernel \cite[Proposition~1]{jabin_quantitative_2018}. The entropy dissipation inequality \cref{eq:entropysolutionineq} classically comes from a formal derivation of the entropy $\int f^N_t\log f^N_t$ and here it is taken as a definition. For the limit equation~\cref{eq:mckeanvlasov-pde_summary}, one can ask for a stronger regularity as in the main theorem \cite[Theorem 1]{jabin_quantitative_2018} stated below. 

\begin{theorem}[Pointwise McKean-Vlasov, $\dot{W}^{-1,\infty}$ kernel \cite{jabin_quantitative_2018}]\label[II]{thm:jabin} Assume that $\nabla\cdot F\in L^\infty(\mathbb{T}^d)$ and that $K\in \dot{W}^{-1,\infty}(\mathbb{T}^d)$ with $\nabla\cdot K \in \dot{W}^{-1,\infty}(\mathbb{T}^d)$. Let $f^N_t$ be an entropy solution of the Liouville equation in the sense of Definition \cref{def:entropysolution}. Assume the limit law satisfies $f\in L^\infty([0,T],W^{2,p}(\mathbb{T}^d))$ for any $p<\infty$ and $\inf f>0$. Then the following entropy bound holds:
\begin{equation}\label[IIeq]{eq:boundjabin}H\big(f^N_t|f^{\otimes N}_t\big) \leq \e^{Ct}{\left(H\big(f^N_0|f^{\otimes N}_0\big)+1\right)},\end{equation}
where $C>0$ depends on $d$, the derivative bounds on $K$, $F$, $f$ and the initial condition. 
\end{theorem}

We sketch the main arguments of the proof in the case $F=0$ for simplicity. The starting point is as before the computations which lead to Lemma \cref{lemma:entropyboundgirsanov_summary} (for the details, see the proof of Lemma \cref{lemma:computeH}) but with a much finer analysis based on the divergence form of the kernel $K=\nabla\cdot V$. 

\begin{proof}[Proof (main ideas)] The computations of Lemma \cref{lemma:computeH} become fully rigorous with the notion of entropy solution and the regularity assumptions on $f$ \cite[Lemma 2]{jabin_quantitative_2018}. Carrying on the computations up to the last step \cref{eq:laststepcomputeH} and using \cref{eq:entropysolutionineq}, Jabin and Wang obtained the following inequality (in integrated form): 
\begin{align}
H\big(f^N_t|f^{\otimes N}_t\big) &\leq H\big(f^N_0|f^{\otimes N}_0\big)\nonumber\\
&-\sum_{i=1}^N \int_0^t\int_{\mathbb{T}^{dN}} \big(K\star\mu_{\mathbf{x}^N}(x^i)-K\star f_s(x^i)\big)\nabla \log f_s^{\otimes N}\big(\mathbf{x}^N\big)f^N_s\big(\dd\mathbf{x}^N\big)\dd s\nonumber\\
&-\sum_{i=1}^N \int_0^t\int_{\mathbb{T}^{dN}} \big(\nabla\cdot K\star\mu_{\mathbf{x}^N}-\nabla\cdot K\star f_s\big) f^N_s\big(\dd\mathbf{x}^N\big)\dd s\nonumber\\
&- \frac{1}{2}\int_0^t I\big(f^N_s|f^{\otimes N}_s\big)\dd s. \label[IIeq]{eq:computeHjabin}
\end{align}
In Lemma \cref{lemma:computeH}, the terms involving $K$ are handled by using the Young inequality. Here, owing to the assumptions on $K$, Jabin and Wang use the decomposition 
\[K = \overline{K}+\widetilde{K},\]
where $\overline{K}=\nabla\cdot V\in \dot{W}^{-1,\infty}(\mathbb{T}^d)$ with $\nabla\cdot \overline{K}=0$, $V\in L^\infty(\mathbb{T}^d)$ and $\widetilde{K}\in L^\infty$. The term involving $\widetilde{K}$ is slightly more technical because of the divergence term but it can be handled following the same ideas than the ones used for $\overline{K}$ and leads to the same conclusion (see \cite[Lemma 4]{jabin_quantitative_2018}). We skip the computations and focus on $\overline{K}$ (this is \cite[Lemma 3]{jabin_quantitative_2018}). By integration by parts, it holds that: 
\begin{align*}
    &\sum_{i=1}^N \int_{\mathbb{T}^{dN}} \big(\overline{K}\star\mu_{\mathbf{x}^N}(x^i)-\overline{K}\star f_s(x^i)\big)\nabla \log f_s^{\otimes N}\big(\mathbf{x}^N\big)f^N_s\big(\dd\mathbf{x}^N\big)\\
    &\quad= \sum_{i=1}^N\int_{\mathbb{T}^d} \big(V\star\mu_{\mathbf{x}^N}(x^i)-V\star f_t(x^i)\big):\nabla_{x^i}f^{\otimes N}_s \nabla_{x^i}{\left(\frac{f^N_s}{f^{\otimes N}_s}\right)}^{\mathrm{T}}\dd\mathbf{x}^N\\
    &\quad\quad+\sum_{i=1}^N\int_{\mathbb{T}^{dN}}\big(V\star\mu_{\mathbf{x}^N}(x^i)-V\star f_t(x^i)\big):\frac{\nabla^2_{x^i}f^{\otimes N}_s}{f^{\otimes N}_s}f^N_s\big(\dd\mathbf{x}^N\big)\\
    &\quad=:A(s)+B(s).
\end{align*}
The two terms $A$ and $B$ are of different nature. For the first one, it is possible to use the similar trick as the one at the end of the proof of Lemma \cref{lemma:computeH}. Using Cauchy-Schwarz inequality and Young inequality, for any $\gamma>0$, 
\begin{align*}
A(s)&\leq \frac{\gamma}{2}I\big(f^N_s|f^{\otimes N}_s\big)+\frac{C}{2\gamma} \sum_{i=1}^N \int_{\mathbb{T}^{dN}}\big|V\star\mu_{\mathbf{x}^N}(x^i)-V\star f_s(x^i)\big|^2f^N_t\big(\dd\mathbf{x}^N\big)\\
&=\frac{\gamma}{2}I\big(f^N_s|f^{\otimes N}_s\big)\\
&\qquad+\frac{C}{2\gamma}\sum_{i=1}^N\sum_{\alpha,\beta=1}^d \int_{\mathbb{T}^{dN}}\big|V_{\alpha,\beta}\star\mu_{\mathbf{x}^N}(x^i)-V_{\alpha,\beta}\star f_t(x^i)\big|^2f^N_s\big(\dd\mathbf{x}^N\big)\\
&=:\frac{\gamma}{2}I\big(f^N_s|f^{\otimes N}_s\big)+\frac{C}{2\gamma}\sum_{i=1}^N\sum_{\alpha,\beta=1}^d A^i_{\alpha,\beta}(s)
\end{align*}
where the constant $C>0$ comes from the bounds on $f$ and $V_{\alpha,\beta}$ are the coordinates of $V$. Choosing the appropriate $\gamma$ will cancel the Fisher information term in \cref{eq:computeHjabin}. It remains to bound the terms $A^i_{\alpha,\beta}(s)$ and $B(s)$ (the term which involves $\widetilde{K}$ in~\cref{eq:computeHjabin} would give analogous terms). As in the conclusion of Lemma \cref{lemma:computeH} and \cref{eq:computeHintrojabin}, since they are observables of the particle system and it is possible to use the change of measure identity \cref{eq:changeofmeasure}. For each $A^i_{\alpha,\beta}(s)$, since $V\in L^\infty(\mathbb{T}^d)$ it will give exactly the same kind of terms as at the beginning of this section. They can be bounded uniformly in $N$ using \cite[Theorem 3]{jabin_quantitative_2018}. For $B(s)$, the change of measure identity~\cref{eq:changeofmeasure} yields: 
\[B(s) \leq \frac{1}{\eta}H\big(f^N_s|f^{\otimes N}_s\big)+\frac{1}{\eta} \log Z_N,\]
where $Z_N$ is of the form
\[ Z_N =  \int_{\mathbb{T}^{dN}} \exp {\left[ \eta N G {\left( \mu_{\mathbf{x}^N} \right)} \right]} f_s^{\otimes N} \big( \dd \mathbf{x}^N \big),\]
with $G:\mu\mapsto\langle\mu\otimes\mu,\phi_2\rangle$ is polynomial function of order two. Namely: 
\[G(\mu) = \langle \mu\otimes\mu, \phi_2\rangle,\]
where
\[\phi_2(x,z) =(V(x-z)-V\star f_s(x)):\frac{\nabla^2 f_s(x)}{f_s(x)}.\]
If $V$ were continuous, Theorem \cref{thm:chaosfromLDP} would say that $\lim_{N\to+\infty}\e^{Nm_0}Z_N$, exists and is finite for a computable $m_0\geq0$, which is more than what is needed here. However, in this case $V$ is only bounded. The authors thus introduce a ``modified law of large numbers'' \cite[Theorem 4]{jabin_quantitative_2018} which implies that $Z_N$ is bounded by a universal constant. The proof of \cite[Theorem 4]{jabin_quantitative_2018} follows similar but much more difficult combinatorial arguments as the ones in the proof of \cite[Theorem 3]{jabin_quantitative_2018}. It is based on a fine use of the two cancellations: 
\[\int_{\mathbb{T}^d} \phi_2(x,z)f_s(\dd z) = 0, \quad \int_{\mathbb{T}^d} \phi_2(x,z)f_s(\dd x) = 0,\]
for all $x,z\in\mathbb{T}^d$. It also needs $L^p$ bounds on $\phi_2$ which depend on the regularity of~$f_s$. The final bound \cref{eq:boundjabin} then follows from Gronwall lemma as before. 
\end{proof}

We conclude this section with some additional remarks and extensions of Theorem \cref{thm:jabin}.

\begin{enumerate}
    \item It is interesting to see how the tricky combinatorial results \cite[Theorem 3, Theorem 5]{jabin_quantitative_2018} can lead to the desired law of large numbers: an insightful use of exchangeability allows to remove extra continuity assumptions.
    \item A reminiscent pattern is the compromise between regularity whether on the initial equation through coefficients, or on the limit process by strong well-posedness result. Here very weak regularity is taken for the particle process, but strong regularity on the limit measure is required. This is in a sense, the opposite of what is done in Section \cref{sec:mckeantowardssingular}. 
    \item The setting of Theorem \cref{thm:jabin} is in fact more general as it also allows a diffusion coefficient $\sigma_N$ which depends on $N$ (we took $\sigma=1$). The behaviour is different depending on whether $\sigma_N\geq\sigma_0>0$ (non-degenerate case) or $\sigma_N\to 0$ (vanishing diffusion case). The first case would add an additional term which depends on $|\sigma_N-\sigma|$ in the final bound \cref{eq:boundjabin}. The vanishing diffusion case is handled by \cite[Theorem 2]{jabin_quantitative_2018} under slightly stronger regularity assumptions on $K$. 
    \item The kinetic case (in $\mathbb{R}^d$) is the original one investigated in \cite{jabin_mean_2016}. The modified law of large number \cite[Theorem 2]{jabin_mean_2016} analogous to \cite[Theorem 4]{jabin_quantitative_2018} is slightly simpler because of the symplectic structure of the system.  
    \item Recent extensions concern gradient systems with an interaction kernel of the form $K=-\nabla W$. The analysis in \cite{bresch_mean-field_2019,serfaty_systems_2019,duerinckx_mean-field_2016} is based on a new modulated free energy which includes in its definition the Gibbs equilibrium measures of the particle and nonlinear systems. 
    \item Guillin, Le Bris and Monmarch\'e \cite{guillin_uniform_2021} have recently shown that Theorem~\cref{thm:jabin} can be made uniform in time for a broad class of singular kernels which includes the Biot-Savart kernel on the torus and the 2D vortex model (see Section \cref{sec:vorticity}). To obtain this result, the authors use a control of the relative entropy by the Fisher information which appears in \cref{eq:computeHjabin}.
    \item The entropy methods described here could be referred as \emph{global} entropy methods because a bound on the global entropy $H(f^N_t|f^{\otimes N}_t)$ (i.e. with $N$ particles)  gives a \emph{local} bound for the marginals of lower order $k$ (by Lemma \cref{lemma:entropyPinskerCsiszar}). In a recent work \cite{lacker_hierarchies_2021}, Lacker has developed \emph{local entropy methods} to prove directly a bound on the $k$-particle relative entropy
    \[H(f^{k,N}_t|f^{\otimes k}_t) = \mathcal{O}((k/N)^2),\]
    for any $k\leq N$ (which also implies the bound $\|f^{k,N}_t-f^{\otimes k}_t\|_{\mathrm{TV}} = \mathcal{O}{\left(\frac{k}{N}\right)}$ by the Pinsker inequality). The approach is based on a kind of BBGK hierarchy for the family of $k$-particle relative entropies. 
\end{enumerate}

\subsection{Concentration inequalities for gradient systems}\label[II]{sec:concentrationineqgradient}

In this section, we make a step forward after propagation of chaos and briefly state two large deviation results for gradient systems. The first one is a weaker result which follows from Theorem \cref{thm:uniformpocuniformconvex} and the Bakry-Emery criterion. The second result is stronger but requires a significant amount of work which will not be detailed here. 

The Bakry-Emery criterion (Proposition \cref{prop:bakryemery}) is applied to McKean-Vlasov gradient systems in Malrieu \cite{malrieu_logarithmic_2001} to obtain concentration inequalities at the particle level. For each observable $\varphi$, it provides a quantitative estimate in both $N$ and $t$ of the deviation between the $N$-particle system and its McKean-Vlasov limit. When the latter converges as $t\to+\infty$ towards its unique invariant measure $\mu_\infty$ (see Corollary \cref{coro:uniformpocuniformconvex}), this also provides confidence interval for the convergence of the $N$-particle system towards $\mu_\infty$. The following theorem summarises the results of \cite{malrieu_logarithmic_2001}.

\begin{theorem}[Concentration inequalities for gradient systems] \label[II]{malrieuCI}
Let us consider the model \cref{eq:mckeanvlasov_summary} with the coefficients \cref{eq:gradientsystem} and a potential $V$ which is $\beta$-uniformly convex. Let $f_0$ satisfy a log-Sobolev inequality with constant $\lambda_0$ (see Section \cref{concentrationI}) and assume that the $N$-particles are initially i.i.d. with common law $f_0$. The following properties hold under the same assumptions as in Theorem \cref{thm:uniformpocuniformconvex}.
\begin{enumerate}
\item There exists $C>0$ such that for all $\varepsilon\geq0$, $N \geq 2$, $t \geq 0$,
\begin{equation}\label[IIeq]{eq:concentrationmalrieut}\sup_{\| \varphi \|_{\mathrm{Lip}} \leq 1} \mathbb{P} {\left( {\left| \frac{1}{N} \sum_{i=1}^N \varphi {\left( X^i_t \right)} - \int_E \varphi(x)  f_t(\dd x) \right|} \geq \varepsilon  + \sqrt{\frac{C}{N}} \right)} \leq 2 \e^{-\frac{ N \lambda_t \varepsilon^2}{2}}. \end{equation}
\item There exists $C>0$ such that for all $\varepsilon\geq0$,, $N \geq 2$, $t \geq 0$,
\begin{equation}\label[IIeq]{eq:concentrationmalrieuinfty}\sup_{\| \varphi \|_{\mathrm{Lip}} \leq 1} \mathbb{P} {\left( {\left| \frac{1}{N} \sum_{i=1}^N \varphi {\left( X^i_t \right)} - \int_E \varphi(x)\mu_\infty(\dd x) \right|} \geq \varepsilon  + \sqrt{\frac{C}{N}} + C \e^{-\beta t} \right)} \leq 2 \e^{-\frac{ N \lambda_t \varepsilon^2}{2}},\end{equation}
\end{enumerate}
where $\lambda_t>0$ is bounded from below and above and will be given in the proof. 
\end{theorem}

\begin{proof}[Proof (sketch)] 
A straightforward computation (see \cite[Lemma 3.5]{malrieu_logarithmic_2001}) shows that the $N$-particle system satisfies the Bakry-Emery criterion (Proposition \cref{prop:bakryemery}) with constant $\beta$. Then, if $f_0$ satisfies $LSI( \lambda_0 )$, \cite[Corollary 3.7]{malrieu_logarithmic_2001} shows that the one-particle distribution $f^{1,N}_t$ satisfies a log-Sobolev inequality with constant $\lambda_t$ such that
\[ \frac{1}{\lambda_t} = \frac{1 - \e^{- 2 \beta t}}{\beta} + \frac{\e^{- 2 \beta t}}{\lambda_0}.\]
Thanks to Lemma \cref{concledoux}, it implies that $f^{1,N}_t$ is concentrated around its mean with an explicit error estimate. The first property therefore follows from the uniform in time bound \cref{eq:uniformpocuniformconvex}. Then, the distance between $f_t$ and $\mu_\infty$ can be quantified in Wasserstein distance: 
\begin{align*}
    W_2(f_t,\mu_\infty) &\leq W_2(f_t,f_t^{1,N})+W_2(f_t^{1,N},\mu_\infty^{1,N})+W_2(\mu_\infty^{1,N},\mu_\infty)\\
    &\leq \frac{C}{\sqrt{N}}+\sqrt{\frac{C}{N}H(f_t^N|\mu^N_\infty)}\\
    &\leq \frac{C}{\sqrt{N}}+C\e^{-\beta t},
\end{align*}
where the first and third terms on the right-hand side of the first line are bounded by $CN^{-1/2}$ by \cref{eq:uniformpocuniformconvex}, the second term is controlled by the relative entropy by the Talagrand inequality \cref{eq:talagrand} and the last line follows as in the proof Corollary \cref{coro:uniformpocuniformconvex}. Letting $N\to+\infty$ leads to
\[W_2(f_t,\mu_\infty)\leq C \e^{-\beta t}.\]
The second property \cref{eq:concentrationmalrieuinfty} thus follows by inserting this last bound in \cref{eq:concentrationmalrieut} (since the Wasserstein-2 distance controls the Wasserstein-1 distance and using Proposition~\cref{prop:dualitywasserstein}).
\end{proof}

Theorem \cref{malrieuCI} quantifies how the empirical measure $\mu_{\mathcal{X}^N_t}$ is close from its limit (in $N$ and $t$) for the distance $D_1$ given by \cref{eq:D0expression}. The distance $D_1$ is dominated by the stronger Wasserstein distance (both metrize the weak topology). Related stronger results have been shown by~\cite{bolley_quantitative_2006} using different techniques, based on the quantitative version of Sanov theorem given by Theorem \cref{thm:PQS}. Note that compared to Malrieu's results \cref{eq:concentrationmalrieut} and \cref{eq:concentrationmalrieuinfty}, the goal is to interchange the supremum and the probability (thanks to the Monge-Kantorovich duality formula Proposition \cref{prop:dualitywasserstein}). This comes at the price of stronger assumptions and with an eventually worse rate of convergence. The following theorem summarises the results of \cite[Theorem 2.9 and Theorem~2.12]{bolley_quantitative_2006}.

\begin{theorem}[Pointwise $W_1$ concentration inequalities]\label[II]{thm:BGVconcentrationineq} Let us consider the model \cref{eq:mckeanvlasov_summary} with coefficients \cref{eq:gradientsystem}, assume that there exist some constants $\beta,\gamma,\gamma'\in\R$ such that the potentials $V,W$ satisfy
\[\nabla^2 V\geq \beta I_d,\quad \gamma I_d\leq \nabla^2 W\leq \gamma' I_d\]
and
\[\forall x\in\R^d,\,\,\forall a>0,\quad |\nabla V(x)|=\mathcal{O}\big(\e^{a|x|^2}\big).\]
Assume that the initial data admits a finite square exponential moment:
\[\exists \alpha_0>0,\quad \int_{\R^d} \e^{\alpha_0 |x|^2} f_0(\dd x)<+\infty.\]
Then the following properties hold. 
\begin{enumerate}
    \item For all $T > 0$, there exists $\lambda, C > 0$ such that for all $\varepsilon>0$, there exists $N_\varepsilon$ such that for $N\geq N_\varepsilon$ : 
    \[\mathbb{P} {\left( \sup_{0 \leq t \leq T} W_1 \big( \mu_{\mathcal{X}^N_t} , f_t \big) > \varepsilon \right)} \leq C {\left( 1 + T \varepsilon^{-2} \right)} \e^{-\lambda N \varepsilon^2}.\]
    \item In the uniformly convex case $\beta>0$ and $\beta+2\gamma>0$, there exists $\lambda, C, T_0, \varepsilon_0 > 0$ such that for all $\varepsilon>0$, there exists $N_\varepsilon$ such that for $N\geq N_\varepsilon$ :
    \[\sup_{t \geq T_0 \log (\varepsilon_0 / \varepsilon)} \mathbb{P} {\left( W_1 \big( \mu_{\mathcal{X}^N_t} , \mu_\infty \big) > \varepsilon \right)} \leq C {\left( 1 + \varepsilon^{-2} \right)} \e^{-\lambda N \varepsilon^2},\]
    where $\mu_\infty$ is the unique invariant measure of the nonlinear McKean-Vlasov system. 
\end{enumerate}
\end{theorem}

A pathwise generalisation is done in \cite{bolley_quantitative_2010} in the case of a bounded time interval. 

\begin{theorem}[Pathwise $W_1$ concentration inequality]
With the same assumptions as in Theorem \cref{thm:BGVconcentrationineq}, for all $T > 0$, there exist $\lambda, C > 0$ such that for all $\varepsilon>0$, there exists $N_\varepsilon$ such that for $N\geq N_\varepsilon$ : 
    \[\mathbb{P} {\left( W_1 \big( \mu_{\mathcal{X}^N_{[0,T]}} , f_{[0,T]} \big) > \varepsilon \right)} \leq C {\left( 1 + T \varepsilon^{-2} \right)} \e^{-\lambda N \varepsilon^2},\]
where $W_1$ denotes the Wasserstein-1 distance on the path space $C([0,T],\R^d)$ (see Definition \cref{def:spaceswasserstein}). 
\end{theorem}

\subsection{General interactions}\label[II]{sec:mckeangeneralinteractions}

In this section, we discuss some results in the very general case of a McKean-Vlasov diffusion of the form
\begin{equation}\label[IIeq]{eq:mckeanvlasovgeneralbsigma} b:\R^d\times\pb(\R^d)\to \R^d,\quad \sigma:\R^d\times\pb(\R^d)\to \mathcal{M}_d(\R),\end{equation}
without assuming any particular form for these functions. 

\subsubsection{Extending McKean's theorem}\label[II]{sec:mckeangeneralinteractionsextendingmckean}

When $b$ and $\sigma$ are Lipschitz for the Wasserstein distance, then McKean's theorem and its proof can be easily extended. 

\begin{theorem}\label[II]{thm:mckeangeneral} In \cref{eq:mckeanvlasov_summary}, let the drift and diffusion coefficients \cref{eq:mckeanvlasovgeneralbsigma} satisfy the following Lipschitz bound for all $(x,y)\in E^2$ and $(\mu,\nu)\in\pb(E)^2$: 
\[\max\Big(|b(x,\mu)-b(y,\nu)|,|\sigma(x,\mu)-\sigma(y,\nu)|\Big) \leq L\Big(|x-y|+W_2(\mu,\nu)\Big).\]
Assume that $f_0\in\pb_q(E)$ for some $q>2$. Then pathwise propagation of chaos in the sense of Definition \cref{def:chaosbycouplingtrajectories_summary} holds for any $T>0$, with $p=2$ and with the synchronous coupling introduced in Theorem \cref{thm:mckean}. The convergence rate is given by
\[\varepsilon(N,T) = C(b,\sigma,T)\beta(N),\]
where $C(b,\sigma,T)>0$ is a constant depending only on $b,\sigma,q$ and $T$ and $\beta(N)$ is given by \cite[Theorem 1]{fournier_rate_2015} : 
\[\beta(N) = \left\{\begin{array}{ll}
N^{-1/2}+N^{-(q-2)/q} & \text{if}\,\,d<4\,\,\text{and}\,\,q\ne 4\\ 
N^{-1/2}\log(1+N)+N^{-(q-2)/q} & \text{if}\,\,d=4\,\,\text{and}\,\,q\ne 4\\
N^{-2/d}+N^{-(q-2)/q}&\text{if}\,\,d>4\,\,\text{and}\,\,q\ne d/(d-2)
\end{array}
\right.
\]
\end{theorem}

\begin{proof}[Proof (sketch)] We follow the same line of arguments of Sznitman's proof. The main change is that \cref{eq:mckeanthmszniterror} should be replaced by 
\[\E{\left|b{\left(\overline{X}{}^i_t,f_t\right)}-b{\left(\overline{X}{}^i_t,\mu_{\overline{\mathcal{X}}{}^N_t}\right)}\right|}^2\leq L \E W^2_2{\left(\mu_{\overline{\mathcal{X}}{}^N_t},f_t\right)} \leq C(T)\beta(N),\]
where the last inequality (with a constant $C(T)>0$) comes from \cite[Theorem 1]{fournier_rate_2015} since the $\overline{X}{}^k_t$ are independent and using a uniform moment bound on $[0,T]$. The inequality \cref{eq:mckeanthmsznitlip} still holds (with a different constant) thanks to the straightforward inequality 
\[\E W_2^2{\left(\mu_{\overline{\mathcal{X}}{}^N_t},\mu_{\mathcal{X}^N_t}\right)}\leq \frac{1}{N}\sum_{j=1}^N \E\big|\overline{X}{}^j_t-X^j_t\big|^2 = \E\big|\overline{X}{}^i_t-X^i_t\big|^2,\]
for any $i\in\{1,\ldots,N\}$ by symmetry. The rest of the proof proceeds as before. 
\end{proof}

The proof of Theorem \cref{thm:mckeangeneral} is also detailed very concisely but precisely in \cite[Section 1]{carmona_lectures_2016}.

\begin{remark}[Completeness and exchangeability]
It may also be interesting to try to adapt McKean's argument (Section \cref{sec:mckeancoupling}) to the setting of Theorem \cref{thm:mckeangeneral}. Most of the proof remains unchanged, the main difficulty (which arises just after \cref{eq:mckeanthmsecondinequality}) is the control the quantity 
\[\E W_2^2{\left(\mathcal{\mu}_{\mathcal{X}^{N,M}_t},\mu_{\mathcal{X}^M_t}\right)},\]
that is, we need to control the Wasserstein distance between two empirical measures with different numbers of samples. To do that, we can mimic the proof of the Hewitt-Savage Theorem \cref{thm:quantitativehewitt} in \cite[Theorem 5.1]{hauray_kacs_2014} and replace the Wasserstein distance by a Sobolev norm $H^{-s}$ (Definition \cref{def:sobolevnorm}). Under some moment assumptions, it defines a distance which is equivalent to the Wassertein distances \cite[Lemma 2.1]{hauray_kacs_2014}. Taking advantage of the polynomial structure property stated in Lemma \cref{lemma:sobolevpolynomial}, it is shown in Proposition \cref{mckeanG} that: 
\[\E\big\|\mu_{\mathcal{X}^{N,M}_t}-\mu_{\mathcal{X}^M_t}\big\|_{H^{-s}}^2 \leq 2\|\Phi_s\|_{\infty} \left(\frac{1}{N}-\frac{1}{M}\right).\]
\end{remark}

As a general rule, if $b,\sigma$ are globally Lipschitz for a Wasserstein metric, then it is possible to extend any result obtained by (synchronous) coupling. The price to pay is a possibly bad convergence rate, in particular with respect to the dimension $d$. Since the convergence rate typically comes from the quantitative Glivenko-Cantelli theorem \cite{fournier_rate_2015} which is sharp in general, it seems hard to obtain better results with this technique. One can also readily check that the approach of Section \cref{sec:gradientsystems} based on It\=o's formula can be applied under convexity assumptions, for instance when 
\[b(x,\mu) = -\nabla V(x) + b_0(x,\mu),\]
where $V$ is convex and $b_0:\R^d\times\pb(\R^d)\to\R^d$ is globally Lipschitz. Following these ideas, the most general and comprehensive article that we are aware of is \cite{andreis_mckeanvlasov_2018}. The authors use the synchronous coupling method to prove pathwise propagation of chaos in various Lipschitz and non Lipschitz cases for a mixed jump-diffusion model with simultaneous jumps (see Example \cref{eq:simultaneousjumps}). Because of the jump interactions, the authors work in a more amenable $L^1$ framework (the results are stated for the $W_1$ distance). Compared to the $L^2$ framework of Theorem \cref{thm:mckeangeneral} this brings some additional technicalities regarding the diffusion part but it does not modify the argument. See also \cite{graham_mckean-vlasov_1992} for an earlier work on jump-diffusion models in a $L^2$ framework but using martingale arguments similar to \cite{sznitman_equations_1984,sznitman_nonlinear_1984}. 

Finally, the globally Lipschitz framework of \cite{andreis_mckeanvlasov_2018} has recently been weakened in \cite{erny_well-posedness_2021} where the author proves the well-posedness and the propagation of chaos for general jump-diffusion McKean models with local Lipschitz coefficients but with an additional assumption about bounded exponential moments. This result is reminiscent from \cite{bolley_stochastic_2011} (see Section \cref{sec:mckeantowardssingular}).  

\subsubsection{Chaos via Girsanov theorem}\label[II]{sec:chaosviagirsanov}

When $\sigma=I_d$ (or more generally when $\sigma$ is non singular and does not depend on the measure argument), under a Lipschitz assumption on the drift, it is also possible to prove strong pathwise chaos in TV norm as in Corollary \cref{coro:mckeantvchaos}, via a Girsanov transform argument. When the drift is Lipschitz in Wasserstein distance, this follows immediately from Theorem \cref{thm:mckeangeneral}, Lemma \cref{lemma:entropyboundgirsanov_summary} and \cite{fournier_rate_2015} (this extends Corollary \cref{coro:mckeantvchaos}).

A recent strategy improves this idea without requiring the preliminary propagation of chaos result which holds only with strong Lipschitz assumptions. The following theorem is a weakened version of \cite[Theorem 2.6]{lacker_strong_2018}.

\begin{theorem}[\cite{lacker_strong_2018}]\label[II]{thm:mckeangirsanov} Fix $T>0$ and $I=[0,T]$. In \cref{eq:mckeanvlasov_summary}, assume that $\sigma=I_d$, that $b$ is bounded and that $b(x,\cdot)$ is Lipschitz for the total variation norm uniformly in $x$. Then for all $k\in\N$ it holds that
\[\lim_{N\to+\infty}H\big(f^{k,N}_I|f^{\otimes k}_I\big)=0.\]
\end{theorem}

This result relies of course on Lemma \cref{lemma:entropyboundgirsanov_summary}. The strategy of \cite{lacker_strong_2018} is then to use a crude large deviation principle to show that the right-hand side of \cref{eq:entropyboundgirsanov_summary} goes to zero as $N\to+\infty$. The key argument is the following result: there exists a constant $C>0$ which depends only on $b$ such that for all measurable open neighbourhood of $f_I$, 
\[\limsup_{N\to+\infty} \frac{1}{N} \log\mathbb{P}{\left(\mu_{\mathcal{X}^N_I}\notin U\right)}=\e^{-CT}\inf_{\nu\notin U} H(\nu|f_I).\]
This result is a kind of Sanov theorem obtained by a change of measure argument from the classical Sanov theorem applied to an i.i.d sequence of $f_I$-distributed random variables. This is \cite[Theorem 2.6 (1)]{lacker_strong_2018}. This result implies that $\Phi(\mu_{\mathcal{X}^N_I})\to \Phi(f_I)$ in probability for all bounded continuous measurable $\Phi$ on $\pb(C([0,T],\R^d))$ (see \cite[Remark 2.8]{lacker_strong_2018}). The conclusion follows by noting that the right-hand side of \cref{eq:entropyboundgirsanov_summary} is precisely an observable of this form. The detailed proof is actually written in a much more general setting than \cref{eq:mckeanvlasovgeneralbsigma}, since it is assumed that $b$ and $\sigma$ are of the form:
\[b:[0,T]\times C([0,T],\R^d)\times \pb(C([0,T],\R^d))\to\R^d,\, \sigma: [0,T]\times C([0,T],\R^d)\to \mathcal{M}_d(\R),\]
that is they depend on the time argument and on the full pathwise trajectories of the particles (instead of their local in time state). The diffusion matrix is assumed to be invertible everywhere and does not depend on the measure argument. The power of Girsanov theorem is precisely that despite this level of generality, the argument is not much modified and the proof remains relatively short. The main change is maybe the more careful look at the topology (since we work fully on the path space) and the questions of measurability which are discussed in \cite[Sections 2.1 and 2.2]{lacker_strong_2018}. Various well-posedness results for the particle and the nonlinear systems within this setting are also presented. 

An important example of such generalized framework is the probabilitstic interpretation of the parabolic-parabolic Keller-Segel model 
\begin{subequations}\label[IIeq]{eq:parabolicparabolickellersegel}
\begin{align}
    \partial_t \rho(t,x)  &= -\chi\nabla\cdot(\rho\nabla c) + \frac{1}{2}\Delta \rho \\ 
    \partial_t c(t,x) &= -\lambda c + \rho + \frac{1}{2}\Delta c,
\end{align}
\end{subequations}
where $\xi,\lambda>0$. Compared to the parabolic-elliptic Keller-Segel model \cref{eq:parabolicelliptickellersegel}, the elliptic equation \cref{eq:parabolicelliptickellersegel_c} is replaced by a parabolic equation which models the diffusivity and evaporation of the chemical substance. In \cite{talay_new_2020,tomasevic_new_2021}, the authors proposed the following probabilistic interpretation of \cref{eq:parabolicparabolickellersegel}, which takes the form of a nonlinear non-Markovian McKean-Vlasov process
\begin{subequations}\label[IIeq]{eq:paraparakellersegelmckeanvlasov}
\begin{align}
    \dd \overline{X}_t &= b(t,\overline{X}_t)\dd t + \Big\{\int_0^t (K_{t-s}\star \rho_s)(\overline{X}_t)\dd s\Big\} + \dd B_t, \\
    c(t,x) &= \e^{-\lambda t}(g_t\star c_0)(x) + \int_0^t \rho_{t-s}\star \e^{-\lambda s} g_s(x) \dd s, 
\end{align}
\end{subequations}

where $\rho_t = \mathrm{Law}(\overline{X}_t)$, $c_0$ is an initial concentration and 
\[K_t(x) = \chi\e^{-\lambda t}\nabla g_t,\quad b(t,x) = \chi\e^{-\lambda t} \nabla(c_0\star g_t),\quad g_t(x) = \frac{1}{(2\pi t)^{d/2}}\e^{-\frac{|x|^2}{2t}}.\]
In this context, the natural mean-field particle version of \cref{eq:paraparakellersegelmckeanvlasov} can be obtained formally by taking $N$ independent Brownian motions and by replacing the density $\rho_t$ by the empirical measure of the particle system. The singularity of $K_0$ and the non-Markovian setting raise many issues and the rigorous mean-field limit is proved in dimension one only in \cite{jabir_mean-field_2018} using a Girsanov transform argument. This argument has also been applied in \cite{tomasevic_propagation_2020} in a Markovian setting similar to \cref{eq:kellersegelparticles} but with a time dependent force $F\equiv F(t,x)$ which belongs to the space $L^q_\mathrm{loc}([0,+\infty),L^p(\R^d))$, with exponents $p,q\in(2,\infty)$ such that $\frac{d}{p}+\frac{2}{q}<1$.

A drawback of the previous results is that there are not quantitative (as they rely on a large deviation principle or on compactness arguments). A sharper analysis of the Girsanov transform argument is presented in \cite[Theorem 2.1]{jabir_rate_2019} and leads to the same kind of result with a quantitative optimal rate of convergence. The argument is very probabilistic and can be understood as the probabilistic counterpart of \cite{jabin_quantitative_2018} (see Section \cref{sec:jabin}). The assumptions are taken to ensure a fine control of the computations in Girsanov theorem and may not be easily interpreted within our usual setting but various detailed applications to more usual forms of McKean-Vlasov diffusion are presented, for instance the case with only bounded coefficients (as in Section \cref{sec:jabinintro}).

\subsubsection{Other techniques}\label[II]{sec:mckeangeneralinteractionsothers}

It turns out that it quickly becomes quite challenging to go beyond the nice globally Lipschitz setting. Depending on the chosen topology, even seemingly simple linear cases such as
\[b(x,\mu) = K\star\mu(x),\quad K:\R^d\times\R^d\to\R^d,\]
can become problematic: if $K$ is unbounded, even if it has a linear growth, then $b$ is not continuous any more for the weak topology. In addition to the continuity, a sufficient set of assumptions under which well-posedness and propagation of chaos can be proved are given in \cite[Section 5]{gartner_mckean-vlasov_1988}. We reproduce it below. 

\begin{assumption}[\cite{gartner_mckean-vlasov_1988}] Given $p\geq2$ and $R>0$, let us define
\[\pb_{p,R}(\R^d) := \left\{\mu\in\pb(\R^d),\,\,\int_{\R^d} |x|^p\mu(\dd x)\leq R\right\},\]
endowed with the topology induced by the weak topology on $\pb(\R^d)$ (in the sense of Definition \cref{def:convergenceproba}). Assume that $\pb_p(\R^d)$ is equipped with the ``inductive topology'' defined by: $\mathscr{A}\subset\pb_p(\R^d)$ is open if and only if $\mathscr{A}\cap\pb_{p,R}(\R^d)$ is open in $\pb_{p,R}(\R^d)$ for each $R>0$. Assume that there exists $p\geq2$ such that 
\[b:\R^d\times\pb_p(\R^d)\to\R^d,\quad \sigma:\R^d\times\pb_p(\R^d)\to\mathcal{M}_d(\R)\]
are continuous and that $\sigma(x,\mu)$ is invertible for all $(x,\mu)\in\R^d\times\pb_p(\R^d)$. 
Assume that there exists $C>0$ and $C_R>0$ for each $R>0$ such that $b,\sigma$ satisfy the following properties. 
\begin{itemize}
    \item \textbf{(Coercivity and growth).} For all $\mu\in\pb(\R^d)$ with compact support
    \[\int_{\R^d}{\left[(p-1) \|\sigma(x,\mu)\|^2+2\langle x,b(x,\mu)\rangle\right]}|x|^{p-2}\mu(\dd x)\leq C{\left(1+\int_{\R^d} |x|^p\mu(\dd x)\right)},\]
    and for all $R>0$, $\mu\in\pb_{p,R}(\R^d)$, $x\in\R^d$,
    \[\|\sigma(x,\mu)\|^2+2\langle x,b(x,\mu)\rangle \leq C_R(1+|x|^2).\]
    \item \textbf{(Monotonicity).} For all $R>0$, for all $\mu,\nu\in\pb_{p,R}(\R^d)$ and for any coupling $\Pi\in\pb(\R^d\times\R^d)$ between $\mu,\nu$, 
    \[\iint_{\R^d\times\R^d} \left[\|\sigma(x,\mu)-\sigma(y,\nu)\|^2+2\langle x-y,b(x,\mu)-b(y,\nu)\rangle\right]_+ \Pi(\dd x,\dd y)\leq C_R,\]
    and
    \begin{multline*}\iint_{\R^d\times\R^d} \left(\|\sigma(x,\mu)-\sigma(y,\nu)\|^2+2\langle x-y,b(x,\mu)-b(y,\nu)\rangle\right) \Pi(\dd x,\dd y)\\\leq C_R\iint_{\R^d\times\R^d}|x-y|^2\Pi(\dd x,\dd y).\end{multline*}
\end{itemize}
\end{assumption}

\begin{remark} Note that the inductive topology on $\pb_p(\R^d)$ is not so far from the topology induced by the $W_p$ distance. Actually, from \cite[Proposition B.3]{gartner_mckean-vlasov_1988}, a sequence $(\mu_n)_n$ in $\pb_p(\R^d)$ converges towards $\mu$ for the inductive topology if and only if 
\[\mu_n\to\mu,\quad\sup_n \int_{\R^d}|x|^p\mu_n(\dd x)<+\infty,\]
where the convergence is the weak convergence.  A slightly simpler set of assumptions expressed in the space $(\pb_p(\R^d),W_p)$ is given for instance in \cite[Section 2]{wang_distribution_2018}. See also the recent \cite{mishura_existence_2020}. Note however that the inductive topology can also be defined when the bound on the $p$-th moment is replaced by a bound on $\langle \mu,\varphi\rangle$ for a fixed nonnegative continuous test function $\varphi$ on $\R^d$, usually called a Lyapunov function. The main results of \cite{gartner_mckean-vlasov_1988} are proved within this generalised setting. Additional topological details are given in \cite[Appendix B]{gartner_mckean-vlasov_1988}.
\end{remark}

The very detailed article of G\"artner \cite{gartner_mckean-vlasov_1988} proves (weak) pathwise well-posedness and propagation of chaos using martingale arguments. This extends earlier works due to Funaki \cite{funaki_certain_1984} (for the well-posedness of the nonlinear system only) and L\'eonard \cite{leonard_loi_1986}. For further works using martingale and compactness arguments, let us also mention \cite{chiang_mckean-vlasov_1994} for a slightly weakened Lipschitz assumption and \cite{dawson_stochastic_1995} for a generalised case where the particles depend on possibly correlated Brownian motions. Note that in this last case, propagation of chaos does not always hold and the empirical measure process converges weakly towards a (non-deterministic) measure-valued process. 

While propagation of chaos has never stopped being an active field of research, it seems that, regarding the case of very general interaction functions \cref{eq:mckeanvlasovgeneralbsigma}, the work of G\"artner has long stayed one of the most, if not the only, complete and general result. Almost three decades later, this question enjoyed a sudden resurgence of interest, motivated mainly on the one hand by biological models (in particular neuron models) and on the other hand by the theory of mean-field games. In addition to the aforementioned works \cite{andreis_mckeanvlasov_2018,carmona_lectures_2016,lacker_strong_2018}, we will conclude this section with some recent directions of research which originate in the mean-field games community. Note that due to the (necessary) higher degree of technicality, we will not enter into much details. Classical references on the mean-field games theory include \cite{cardaliaguet_master_2019,cardaliaguet_notes_2010,carmona_probabilistic_2018-1,carmona_probabilistic_2018}.

\begin{itemize}
    \item In \cite{chassagneux_weak_2019}, the authors prove a very neat bound of the form 
    \[{\left|\Phi(f_t)-\E\Phi\big(\mu_{\mathcal{X}^N_t}\big)\right|} = \sum_{j=1}^{k-1}\frac{C_j}{N^j} + \mathcal{O}{\left(\frac{1}{N^k}\right)},\]
    where $\Phi:\pb_2(\R^d)\to\R$, the constants $C_j$ do not depend on $N$ and $k$ depends on the regularity of $\Phi$, $b$ and $\sigma$. In this context, regularity means differentiability in the Wasserstein space $(\pb_2(\R^d),W_2)$. As we have already seen in Section \cref{sec:limitgenerator} regarding \cite{mischler_new_2015,mischler_kacs_2013}, defining a differential calculus on the space of measures is not an easy task. The framework detailed in \cite[Section 2]{chassagneux_weak_2019} is based on the notion of ``linear functional derivatives'' and ``L-derivatives'' introduced in \cite{cardaliaguet_master_2019}. Note that the authors still assume at least a uniform bound on the diffusion matrix but also that $b$ and $\sigma$ are globally Lipschitz for the $W_2$ distance. But contrary to the results obtained using the Glivenko-Cantelli theorem \cite{fournier_rate_2015}, the constants $C_j$ do not depend on the dimension. In fact, the framework of \cite{chassagneux_weak_2019} is also applicable to the static case of $N$ $\mu$-distributed i.i.d random variables $\mathcal{X}^N$ and thus it provides explicit convergence rate of $\E\Phi(\mu_{\mathcal{X}^N})$ towards $\Phi(\mu)$ for smooth observables on $\pb_2(\R^d)$. The above result in both the static and McKean-Vlasov cases is obtained when $\Phi$ is ``$(2k+1)$-times differentiable with respect to the functional derivative''. 
    \item In \cite{chaudru_de_raynal_well-posedness_2019} (see also \cite{chaudru_de_raynal_strong_2020}) the authors revisit the question of the well-posedness of the martingale problem associated to McKean-Vlasov equations with general interactions and relate this question to the study of a class of (linear) parabolic type PDEs on the Wasserstein space (the backward Kolmogorov equation with source term and terminal condition). In the subsequent work \cite{chaudru_de_raynal_backward_2019}, the problem is investigated at the particle level which provides (quantitative) propagation of chaos results concerning the trajectories of the particles, the convergence of their distribution and the convergence of the emprirical measure process. The results hold when $b$ and $\sigma$ are bounded, H\"older continuous in space and with two bounded and H\"older continuous linear functional derivatives in the measure argument and when $\sigma$ is also uniformly elliptic. The strategy is also linked to the notion of regularization by noise and the Zvonkin transform, see \cite{zvonkin_transformation_1974,veretennikov_strong_1981}. 
\end{itemize}

\section{Boltzmann models}\label[II]{sec:boltzmannreview}

The probabilistic treatment of the Boltzmann model has been initiated by Kac in the seminal article \cite{kac_foundations_1956}. The original treatment of Kac model (Example \cref{example:kacmodels}) is based on the continuity of the generator $\mathcal{L}_N$ on the space of test functions $(C_b(E^N),\|\cdot\|_\infty)$. The arguments have been later generalised \cite{carlen_kinetic_2013-1} for a wider class of models under boundedness assumptions at the pointwise level (Section~\cref{sec:kactheorem}). A pathwise generalisation of Kac's theorem is due to \cite{graham_stochastic_1997} (Section \cref{sec:pathwisekactheorem}). Many physical models (for instance \cref{eq:hardsphereB_summary} and \cref{eq:truemaxwellmolecules_summary}) do not fit into this framework because of the strong boundedness assumption on the collision rate. To prove more general results, we will first discuss the historical stochastic martingale arguments \cite{tanaka_probabilistic_1983,sznitman_equations_1984} (Section \cref{sec:martingaleboltzmannreview}) and then three historical arguments which have recently been brought up to date and completed: first the SDE and coupling method due to Murata \cite{murata_propagation_1977} (Section \cref{sec:couplingBoltzmann}); then the pointwise study of the generator of the empirical process initiated by Gr\"unbaum \cite{grunbaum_propagation_1971} (Section \cref{sec:kacprogram}); finally, we briefly present Lanford's approach \cite{lanford_time_1975} on the deterministic hard-sphere system (Section~\cref{sec:lanford}). 

\subsection{Kac's theorem via series expansions}\label[II]{sec:kactheorem}

The following theorem, originally due to Kac, is the most important result of this section.  

\begin{theorem}[Kac]\label[II]{thm:kac} Let $(f_0^N)^{}_N$ be a sequence of symmetric probability measures on $E^N$ which is $f_0$-chaotic for a given $f_0\in\pb(E)$. Let $(\mathcal{Z}^N_t)^{}_t$ be the $N$-particle process with initial law $f^N_0$ and with generator 
\[\mathcal{L}_N\varphi_N = \frac{1}{N}\sum_{i<j}L^{(2)}\diamond_{ij}\varphi_N,\]
with $L^{(2)}$ given by Assumption \cref{assum:L2_summary} together with the uniform bound \cref{eq:uniformboundlambda_summary} on the interaction rate $\lambda$. Let $s\in\N$, $s\geq1$, and let $\varphi_s\in C_b(E^s)$ be a test function. Then for any time $t>0$ there exists $f_t\in\pb(E)$ such that
\[\E{\big[\varphi_s{\big(\mathcal{Z}^{s,N}_t\big)}\big]}\underset{N\to+\infty}{\longrightarrow}\langle f_t^{\otimes s},\varphi_s\rangle.\]
where we recall that $\mathcal{Z}^{s,N}_t$ denotes the process in $E^s$ extracted from the $s$ first components of $\mathcal{Z}^N_t$. Moreover $f_t$ is a weak measure solution of the general Boltzmann equation \cref{eq:Boltzmannequationgeneral_summary}.
\end{theorem}

We present two proofs of this theorem. Both are based on the explicit solution of the Liouville equation given by a series expansion. The first proof works at the level of observables. The second proof is slightly shorter but also requires a $L^1$ framework to work at the level of the laws (forward Kolmogorov point of view). The first proof is due to Kac \cite{kac_foundations_1956} for a one-dimensional caricature of a Maxwellian gas. The arguments are generalised in \cite{carlen_kinetic_2013}. Our presentation is also inspired by the work of McKean \cite{mckean_exponential_1967}. The second proof is the probabilistic version of Lanford's approach on the deterministic hard-sphere system (see Section \cref{sec:lanford}). The bound \cref{eq:uniformboundlambda_summary} and the fact that the interactions are delocalised considerably simplify the proof. The detail of the proof can be found in \cite{pulvirenti_kinetic_1996}. 

\begin{proof}[Proof (at the level of the observables)]
Since the operator $\mathcal{L}_N$ is bounded for the $\|\cdot\|_\infty$ norm, the exponential series $\e^{t\mathcal{L}_N}$ is convergent and it holds that:
\begin{equation}\label[IIeq]{eq:backwardseriess}\E{\big[\varphi_s{\big(\mathcal{Z}^{s,N}_t\big)}\big]} = \sum_{k=0}^{+\infty}\frac{t^k}{k!}\langle f_0^N, \mathcal{L}^k_N\varphi_s\rangle.\end{equation}
The strategy is to apply the dominated convergence theorem to pass to the limit in this series. The crucial observation is that the series converges for $t$ small enough, uniformly in $N$. Using only the continuity estimate 
\begin{equation}\label[IIeq]{eq:continuityLN}\|\mathcal{L}_N\varphi_s\|_\infty\leq C(\Lambda)N\|\varphi_s\|_\infty,\end{equation}
would give the convergence on a time interval $t<1/(NC(\Lambda))$ and it would not be possible to take the limit $N\to+\infty$. However, when $s\geq1$ is fixed, better estimates are available which are summarised in the following lemma. The basic idea is to split the general term of the series into two parts \cref{eq:kactheoremsplittingsNNsN}, one of order $1/N$ which vanishes when $s$ is fixed and a leading term of order one which converges and which will give the desired limit. 
\begin{lemma} Let us consider the linear operator $\mathbf{D}$ on $C_b(E^\infty):=\cup_{\ell\geq0}C_b(E^\ell)$ defined for $\varphi_s\in C_b(E^s)$ by: 
\[(\mathbf{D}\varphi_s)(z^1,\ldots,z^s,z^{s+1}):=\sum_{i=1}^s (L^{(2)}\diamond_{i,s+1}(\varphi_s\otimes 1))(z^1,\ldots,z^s,z^{s+1}).\]
Note that since $C_b(E^\ell)\subset C_b(E^{\ell+1})$ by the inclusion $\varphi_s\mapsto \varphi_s\otimes 1$, the space $C_b(E^\infty)$ is actually a vector space. The following properties hold. 
\begin{enumerate}[(1)]
    \item For all $k,s$ such that $k+s\leq N$, 
    \begin{equation}\label[IIeq]{eq:kactheoremsplittingsNNsN}\langle f_0^N,\mathcal{L}_N^k\varphi_s\rangle = \frac{u_{s,k}(\varphi_s)}{N}+\alpha^{(s,k)}_N\langle f_0^{s+k,N},\mathbf{D}^k\varphi_s\rangle, \end{equation}
    where $u_{s,k}(\varphi_s)$ satisfies 
    \begin{equation}\label[IIeq]{eq:bounduskphis}|u_{s,k}(\varphi_s)|\leq C(\Lambda)^k\|\varphi_s\|_\infty \frac{(s+k-2)!}{(s-1)!}\sum_{\ell=0}^{k-1}(s+\ell)^2,\end{equation}
    where $C(\Lambda)$ is the constant in \cref{eq:continuityLN} (in particular, it does not depend on $N$ nor on $k$), and
    \begin{equation}\label[IIeq]{eq:alphaskN}\alpha^{(s,k)}_N := \frac{(N-s)\ldots(N-s-k+1)}{N^k}.\end{equation}
    \item There exists $t_0>0$ which depends only on $s$ and $\Lambda$ such that the series \cref{eq:backwardseriess} converges absolutely, uniformly in $N$ and $t\in[0,t_0]$.
    \item For each $k\geq1$, it holds that 
    \begin{equation}\label[IIeq]{eq:termbytermbackward}\langle f^{N}_{0},\mathcal{L}^k_N\varphi_s\rangle\underset{N\to+\infty}{\longrightarrow} \langle f_0^{\otimes (s+k)}, \mathbf{D}^k\varphi_s\rangle.\end{equation}
\end{enumerate}
\end{lemma}

The second point is proved in \cite[Lemma 3.1]{carlen_kinetic_2013}. The only difference is that in our setting, we have to take into account the constant $\Lambda$. Their proof is based on an estimate similar to \cref{eq:bounduskphis} obtained by a combinatorial argument which does not use the splitting \cref{eq:kactheoremsplittingsNNsN}. The third point is essentially the content of \cite[Lemma 3.3]{carlen_kinetic_2013}. We give an alternative proof here based on the properties of the operator $\mathbf{D}$ which was introduced by McKean \cite{mckean_exponential_1967}. 

\begin{proof} Let us start from the following observation: for all $z^1,\ldots,z^N\in E$, 
\begin{equation}\label[IIeq]{eq:LNphis}\mathcal{L}_N\varphi_s(z^1,\ldots,z^N) = \frac{s}{N}\mathcal{L}_s\varphi_s(z^1,\ldots,z^s) + \frac{1}{N}\sum_{\ell=s+1}^N (\mathbf{D}\varphi_s)(z^1,\ldots,z^s,z^\ell).
\end{equation}
Note that $\mathcal{L}_N\varphi_s$ is a function of $N$ variables but it can be written as the sum of $s$ functions of $s$ variables and $(N-s)$ functions of $(s+1)$ variables. By symmetry we deduce that:
\begin{equation}\label[IIeq]{eq:kaclemmak1}\langle f^{N}_{0},\mathcal{L}_N\varphi_s\rangle = \frac{s}{N}\langle f^{s,N}_{0},\mathcal{L}_s\varphi_s\rangle +\frac{N-s}{N}\langle f^{s+1,N}_{0}, \mathbf{D}\varphi_s\rangle.\end{equation}
Moreover, the following continuity estimates hold for all $s\geq1$, 
\begin{equation}\label[IIeq]{eq:continuityLsD}\|\mathcal{L}_s\varphi_s\|_\infty \leq C(\Lambda) s\|\varphi_s\|_\infty,\quad \|\mathbf{D}\varphi_s\|_\infty \leq C(\Lambda) s\|\varphi_s\|_\infty,\end{equation}
where $C(\Lambda)$ depends only on $\Lambda$.
\begin{enumerate}[(1)]
\item The first point is proved by induction on $k\leq N$. The case $k=0$ is the initial chaoticity assumption and the case $k=1$ immediately follows from \cref{eq:kaclemmak1} and~\cref{eq:continuityLsD}. Let us assume the result for $k\geq1$ and let us take $s\in\N$ such that $s+k+1\leq N$. Using \cref{eq:LNphis}, by exchangeability it holds that
\[\langle f^N_0,\mathcal{L}_N^{k+1}\varphi_s\rangle = \frac{s}{N}\langle f^N_0,\mathcal{L}_N^k(\mathcal{L}_s\varphi_s)\rangle+\frac{N-s}{N}\langle f^N_0,\mathcal{L}_N^k(\mathbf{D}\varphi_s)\rangle.\]
Since $\mathcal{L}_s\varphi_s$ is a function of $s$ variables and $\mathbf{D}\varphi_s$ is a function of $(s+1)$ variables with $(s+1)+k\leq N$, the induction hypothesis for each of the two terms on the right-hand side gives:
\begin{align*}
    \langle f^N_0,\mathcal{L}_N^{k+1}\varphi_s\rangle &= \frac{s}{N}\left(\frac{u_{s,k}(\mathcal{L}_s\varphi_s)}{N}+\alpha^{(s,k)}_{N}\langle f^{s+k,N}_0,\mathbf{D}^k(\mathcal{L}_s\varphi_s)\rangle\right)\\
    &\quad+\frac{N-s}{N}\left(\frac{u_{s+1,k}(\mathbf{D}\varphi_s)}{N}+\alpha^{(s+1,k)}_N\langle f_0^{s+1+k,N},\mathbf{D}^{k+1}\varphi_s\rangle\right).
    \end{align*}
First we note that:
\[\alpha^{(s,k+1)}_N = \frac{N-s}{N}\alpha_N^{(s+1,k)},\]
as expected. Then we set: 
\begin{equation}\label[IIeq]{eq:uskplus1}u_{s,k+1}(\varphi_s) := \frac{s}{N}u_{s,k}(\mathcal{L}_s\varphi_s)+{s}\alpha^{(s,k)}_N\langle f^{s+k,N}_0,\mathbf{D}^k(\mathcal{L}_s\varphi_s)\rangle+\frac{N-s}{N}u_{s+1,k}(\mathbf{D}\varphi_s).\end{equation}
The induction hypothesis \cref{eq:bounduskphis} can be used again to bound $u_{s,k+1}(\varphi_s)$. First we note that 
\[\frac{(s+k-2)!}{(s-1)!}\sum_{\ell=0}^{k-1}(s+\ell)^2\leq \frac{(s+k-1)!}{s!}\sum_{\ell=0}^{k-1}(s+1+\ell)^2.\]
Thus using the continuity bounds \cref{eq:continuityLsD} and the induction hypothesis \cref{eq:bounduskphis}, we deduce that: 
\begin{equation}\label[IIeq]{eq:boundsnNsN}\frac{s}{N}|u_{s,k}(\mathcal{L}_s\varphi_s)|+\frac{N-s}{N}|u_{s+1,k}(\mathbf{D}\varphi_s)|\leq C(\Lambda)^{k+1}\|\varphi_s\|_\infty \frac{(s+k-1)!}{(s-1)!}\sum_{\ell=0}^{k-1}(s+1+\ell)^2.\end{equation}
Moreover, it holds that $\alpha^{(s,k)}_N\leq 1$, so using \cref{eq:continuityLsD} again leads to
\begin{equation}\label[IIeq]{eq:boundalphaknlanglerangle}{s}\alpha^{(s,k)}_N\langle f^{s+k,N}_0,\mathbf{D}^k(\mathcal{L}_s\varphi_s)\rangle\leq C(\Lambda)^{k+1}s^2\frac{(s+k-1)!}{(s-1)!}\|\varphi_s\|_\infty.\end{equation}
Reporting \cref{eq:boundsnNsN} and \cref{eq:boundalphaknlanglerangle} into \cref{eq:uskplus1} finally gives:
\[|u_{s,k+1}(\varphi_s)|\leq C(\Lambda)^{k+1}\|\varphi_s\|_\infty \frac{(s+k-1)!}{(s-1)!}\sum_{\ell=0}^{k}(s+\ell)^2,\]
which concludes the proof of the first point. 
\item Let us split the series \cref{eq:backwardseriess} into two parts, the first one for $k=0,\ldots,N-s$ and the second one for $k\geq N-s+1$. For the second part, we use the crude estimate: 
\[\|\mathcal{L}_N^k\varphi_s\|_\infty \leq C(\Lambda)N^k\|\varphi_s\|_\infty.\]
Then using Stirling's formula, the series 
\[\sum_{k=N-s+1}^{+\infty} \frac{(C(\Lambda)t)^k}{k!}N^{k}\leq \sum_{k=N-s+1}^{+\infty} \frac{(C(\Lambda)t)^k}{k!}(k+s-1)^{k},\]
is convergent for $t<\frac{1}{2C(\Lambda)\e}$. Then using the first point it holds that: 
\[\sum_{k=0}^{N-s} \frac{t^k}{k!}\langle f^N_0,\mathcal{L}_N^k\varphi_s\rangle\leq \frac{1}{N}\sum_{k=0}^{N-s}\frac{t^k}{k!} u_{s,k}(\varphi_s) + \sum_{k=0}^{N-s}\frac{t^k}{k!}\langle f^{s+k,N}_0,\mathbf{D}^k\varphi_s\rangle.\]
From \cref{eq:bounduskphis}, the following elementary estimate holds for $k\geq1$: 
\begin{align*}
\frac{t^k}{k!}|u_{s,k}(\varphi_s)|&\leq (C(\Lambda)t)^k\|\varphi_s\|_\infty\frac{(s+k-2)!}{k!(s-1)!}k(s+k-1)^2 \\
&\leq (C(\Lambda)t)^k\|\varphi_s\|_\infty \binom{s+k-2}{s-1}(s+k-1)^2\\
&\leq (C(\Lambda)t)^k\|\varphi_s\|_\infty \e^{s-1}\left(1+\frac{k-1}{s-1}\right)^{s-1}(s+k-1)^2.
\end{align*}
It follows that the series whose general term is $(t^k/k!)u_{s,k}(\varphi_s)$ is absolutely convergent uniformly in $N$ for $t$ small enough. Similarly, for the series whose general term is bounded by
\begin{align*}
    \frac{t^k}{k!}|\langle f^{s+k,N}_0,\mathbf{D}^k\varphi_s\rangle| &\leq (C(\Lambda)t)^k\|\varphi_s\|_\infty \binom{s+k-1}{s-1} \\ 
    &\leq(C(\Lambda)t)^k\|\varphi_s\|_\infty \e^{s-1}\left(1+\frac{k}{s-1}\right)^{s-1},
\end{align*}
the same conclusion holds. This concludes the proof of the second point. 
\item This follows immediately from the first point, the fact that $\alpha^{(s,k)}_N\to 1$ as $N\to~+\infty$ and the initial chaoticity assumption. 
\end{enumerate}
\end{proof}

Once the lemma is proved, it follows that for any $t<t_0$ there exists a family of probability measures $(f^{(s)}_{t})_s$ on $E^s$ such that 
\begin{equation}\label[IIeq]{eq:convergeEvarphisZN}
\E{\big[\varphi_s{\big(\mathcal{Z}^{s,N}_t\big)}\big]}\underset{N\to+\infty}{\longrightarrow}\langle f_t^{(s)}, \varphi_s\rangle:=\sum_{k=0}^{+\infty}\frac{t^k}{k!}\langle f_0^{\otimes(s+k)},\mathbf{D}^k\varphi_s\rangle.\end{equation}
It remains to prove that $f_t^{(s)}=f_t^{\otimes s}$ where $f_t=f_t^{(1)}$. The following argument is due to McKean \cite{mckean_exponential_1967} who noted that the operator $\mathbf{D}$ is a derivation in the sense that for any $s_1,s_2\in\N$
\begin{equation*}\mathbf{D}(\varphi_{s_1}\otimes\varphi_{s_2}) = \mathbf{D}\varphi_{s_1}\otimes \varphi_{s_2}+\varphi_{s_1}\otimes \mathbf{D}\varphi_{s_2}.\end{equation*}
Leibniz rule therefore implies that for any $s_1+s_2=s$ and $\varphi_{s_1}\in C_b(E^{s_1})$, $\varphi_{s_2}\in C_b(E^{s_2})$, 
    \begin{align*}
        \langle f^{(s)}_t, \varphi_{s_1}\otimes\varphi_{s_2}\rangle &= \sum_{k=0}^{+\infty} \frac{t^k}{k!}\sum_{
        \ell=0}^k \binom{k}{\ell} \langle f_0^{\otimes(s_1+s_2+k)}, \mathbf{D}^{\ell}\varphi_{s_1}\otimes \mathbf{D}^{k-\ell}\varphi_{s_2}\rangle \\ 
        &= \sum_{k=0}^{+\infty} \frac{t^k}{k!}\sum_{
        \ell=0}^k \binom{k}{\ell} \langle f_0^{\otimes(s_1+\ell)}, \mathbf{D}^\ell\varphi_{s_1}\rangle\langle f_0^{\otimes(s_2+k-\ell)}, \mathbf{D}^{k-\ell}\varphi_{s_2}\rangle \\ 
        &= \sum_{k=0}^{+\infty} \frac{t^k}{k!}\langle f_0^{\otimes(s_1+k)}, \mathbf{D}^k\varphi_{s_1}\rangle\sum_{\ell=0}^{+\infty}\frac{t^\ell}{\ell!}\langle f_0^{\otimes(s_2+\ell)}, \mathbf{D}^\ell\varphi_{s_2}\rangle \\
        &=\langle f_t^{(s_1)}, \varphi_{s_1}\rangle\langle f_t^{(s_2)},\varphi_{s_2}\rangle.
    \end{align*}
From which it follows that $f^{(s)}_t = f_t^{(s_1)}\otimes f_t^{(s_2)}$ and therefore $f_t^{(s)}=f_t^{\otimes s}$.
Then, by absolute convergence of all the series, 
it is possible to differentiate with respect to time and directly check that $f_t$ is a weak-measure solution of the Boltzmann equation: for a test function $\varphi\in C_b(E)$, 
\begin{align*}
    \frac{\dd}{\dd t}\langle f_t,\varphi\rangle &= \sum_{k=0}^{+\infty} \frac{t^k}{k!} \langle f_0^{\otimes (k+2)},\mathbf{D}^{k+1}\varphi\rangle\\
    &=\sum_{k=0}^{+\infty} \frac{t^k}{k!} \langle f_0^{\otimes (k+2)},\mathbf{D}^{k}[\mathbf{D}\varphi]\rangle\\
    &= \langle f_t^{\otimes 2}, \mathbf{D}\varphi\rangle, 
\end{align*}
where the last line follows from \cref{eq:convergeEvarphisZN} with $s=2$. Finally, since $t_0$ does not depend on the initial condition, the same reasoning applies on $[t_0,2t_0]$ and so on and therefore the result holds for any $t>0$.
\end{proof}

\begin{remark}[Convergence rate] Although we did not write it in the statement of the theorem, it can be seen from the proof (Equation \cref{eq:bounduskphis}) that for any $\varphi_s$, $\E[\varphi_s(\mathcal{Z}^{s,N}_t)]$ converges towards $\langle f_t^{\otimes s}, \varphi_s\rangle$ at rate $1/N$ (with a constant which depends on $\varphi_s$ and~$s$). This rate is optimal since it implies 
\[\E\big|\varphi_s\big(\mathcal{Z}^{s,N}_t\big)-\langle f_t^{\otimes s},\varphi_s\rangle\big|^2 = \mathcal{O}(1/N).\]
\end{remark}

Within this approach, the limit $f_t$ is defined weakly and the above proof is actually a proof of existence of a weak-measure solution of the Boltzmann equation. The dual proof follows the same arguments at the level of the laws. For simplicity, we present it in a $L^1$ framework and follow closely the arguments of \cite{pulvirenti_kinetic_1996}.

\begin{proof}[Proof (Forward point of view)] Let the initial law $f_0^N\in L^1(E^N)$ be in $L^1(E^N)$, for all $N\in\N$. We denote by $f_{t}^{s,N}$ the $s$-marginal of the law of the particle system at time $t>0$. By integrating the Liouville equation (in strong form) with respect to the variables $s+1$ to $N$, the $s$-th marginal $f^{s,N}_t$ is shown to satisfy the famous BBGKY hierarchy (see also Section \cref{sec:finitesystems}):
\begin{equation}\label[IIeq]{eq:forwardBBGKYBoltzmann}\partial_t f_t^{s,N} = \frac{s}{N}\mathcal{L}^s f_t^{s,N}+\frac{N-s}{N}\mathcal{C}_{s,s+1}f_t^{s+1,N},\end{equation}
where the operator $\mathcal{C}_{s,s+1}:\pb(E^{s+1})\to \pb(E^s)$ is defined as the dual of $\mathbf{D}$ restricted to $C_b(E^s)$, for $f^{(s+1)}\in \pb(E^{s+1})$ and $\varphi_s\in C_b(E^s)$, 
\[\langle \mathcal{C}_{s,s+1}f^{(s+1)}, \varphi_s\rangle := \langle f^{(s+1)}, \mathbf{D}\varphi_s\rangle.\]
Equation \cref{eq:forwardBBGKYBoltzmann} can be re-written using Duhamel's formula:  
\[f_t^{s,N} = \mathbf{T}_N^{(s)}(t) f_0^{(s,N)}+\frac{N-s}{N}\int_0^t \mathbf{T}_N^{(s)}(t-t_1)\mathcal{C}_{s,s+1}f_{t_1}^{s+1,N}\dd t_1,\]
where $\mathbf{T}_N^{(s)}$ is the Markov semi-group acting on $\mathcal{P}(E^s)$ generated by $\frac{s}{N}\mathcal{L}^s$. Iterating this formula gives an explicit formula for the solution of \cref{eq:forwardBBGKYBoltzmann}, namely: 
\begin{multline}\label[IIeq]{eq:BBGKYDuhamel}f_{t}^{s,N} = \sum_{k=0}^{+\infty} \alpha_N^{(s,k)} \int_0^t\int_0^{t_1}\ldots\int_0^{t_{k-1}}\mathbf{T}_N^{(s)}(t-t_1)\mathcal{C}_{s,s+1}\mathbf{T}_N^{(s+1)}(t_1-t_2)\mathcal{C}_{s+1,s+2}\ldots\\
\mathcal{C}_{s+k-1,s+k}\mathbf{T}_N^{(s+k)}(t_k)f_0^{s+k,N}\dd t_1\ldots\dd t_k,\end{multline}
where $\alpha^{(s,k)}_N$ is given by \cref{eq:alphaskN}. Just as in the previous proof, the goal is to show from this series expansion that it is possible take the limit $N\to+\infty$ in the series and that the limit defines a $f_t$-chaotic family where $f_t$ solves the Boltzmann equation. The strategy is again to show the uniform convergence of the series for small $t$ and then the term-by-term convergence. The uniform convergence of the series is straightforward in a $L^1$ framework since the operator $\mathbf{T}_N^{(s)}$ is an isometry in $L^1(E^s)$ and for all $s\geq1$: 
\[\forall f^{(s+1)}\in L^1(E^{s+1}),\quad\|\mathcal{C}_{s,s+1}f^{(s+1)}\|_{L^1(E^s)}\leq sC(\Lambda)\|f^{(s+1)}\|_{L^1(E^{s+1})}.\]
Thus the series of the $L^1$ norms are bounded by, 
\[\sum_{k=0}^{+\infty}s(s+1)\ldots(s+k-1) \frac{t^k}{k!} C(\Lambda)^k\leq \sum_{k=0}^{+\infty} {(2C(\Lambda) t)^k},\]
and uniform convergence in $L^1$ holds for $t<1/(2C_\Lambda)$. Assume that it is possible to prove the existence and uniqueness of the solution $f_t$ of the Boltzmann equation, as an element of $C([0,t_0],L^1(E))$ (typically by a fixed point method). Then a direct computation shows that starting from $f_0^{\otimes s}$ the function $f_t^{\otimes s}$ satisfies: 
\begin{equation}\label[IIeq]{eq:DuhamelBoltzmannhomogeneous}f_t^{\otimes s} = \sum_{k=0}^{+\infty}\frac{t^k}{k!} \mathcal{C}_{s,s+1}\mathcal{C}_{s+1,s+2}\ldots
\mathcal{C}_{s+k-1,s+k}f_{0}^{\otimes(s+k)}.\end{equation}
Each term of this series is exactly the limit in $L^1$ of the corresponding term in the series \cref{eq:BBGKYDuhamel} since 
\[\forall f^{(s)}\in L^1(E^s),\quad\|(\mathbf{T}_N^{(s)}-Id)f^{(s)}\|_{L^1(E^s)}\underset{N\to+\infty}{\longrightarrow}0.\]
The proof can be terminated by iterating the argument for all $t>0$ as in the previous proof. 
\end{proof}

\subsubsection*{Adding one-particle individual flows} To conclude this section, we briefly explain how to extend the result to the more general case where each particle also has an individual flow given by an operator $L^{(1)}$. As in Section \cref{sec:boltzmann_summary}, let us consider a $N$-particle system defined by the operator 
\[\mathcal{L}_N\varphi_N = \sum_{i=1}^N L^{(1)}\diamond_i \varphi_N + \frac{1}{N}\sum_{i<j}L^{(2)}\diamond_{ij}\varphi_N,\]
where $L^{(2)}$ satisfies the assumptions of Theorem \cref{thm:kac} and where we assume that the operator $L^{(1)}$ generates a continuous Markov semi-group acting on a sufficiently large subset of $C_b(E)$. For instance and as explained at the beginning of Section~\cref{sec:boltzmann_summary}, this extension is particularly important as it includes the case of kinetic particles in $E=\R^d\times\R^d$ where the particles are subject to the free-transport between the collisions. Note however that, as mentioned in Section \cref{sec:boltzmann_summary}, the boundedness assumption on the rate of collision in Theorem \cref{thm:kac} does not allow the physically important case where two particles collide only when they are exactly at the same position, since in this case, the collision rate is a Dirac delta and thus unbounded.

In the first approach, the proof is exactly the same with $\mathbf{D}$ replaced by $\mathbf{D}+\mathbf{S}$ where $\mathbf{S}$ is the linear operator on $C_b(E^\infty)$ defined by 
\[\forall \varphi_s\in C_b(E^s),\quad \mathbf{S}\varphi_s = \sum_{i=1}^s L^{(1)}\diamond_i \varphi_s.\]
The exponential formula \cref{eq:backwardseriess} does not converge when $\mathbf{S}$ is not continuous for the $L^\infty$ norm, which includes many interesting case such as free transport or diffusion. However, when $\mathbf{S}$ generates a backward semi-group $\mathbf{T}$ on $C_b(E^\infty)$ which is continuous for the $L^\infty$ norm, one can write 
\[\langle f_{t}^N,\varphi_s\rangle = \langle f_{0}^{s,N}, \mathbf{T}(t)\varphi_s\rangle+\int_0^t \frac{\dd}{\dd t_1}\langle f_{t_1}^{N}, \mathbf{T}(t-t_1)\varphi_s\rangle\dd t_1.\]
A direct computation shows that 
\[\frac{\dd}{\dd t_1}\langle f_{t_1}^{N}, \mathbf{T}(t-t_1)\varphi_s\rangle = \langle f_{t_1}^{N}, \mathcal{L}_N^B \mathbf{T}(t-t_1)\varphi_s\rangle,\]
where $\mathcal{L}_N^B = \frac{1}{N}\sum_{i<j}L^{(2)}\diamond_{ij}$. Iterating this formula, one gets the backward series expansion: 
\begin{multline*}\langle f^N_t, \varphi_s\rangle = \sum_{k=0}^{+\infty}\int_0^t\int_0^{t_1}\ldots\int_0^{t_{k-1}}\langle f_0^N, \mathbf{T}(t_k)\mathcal{L}_N^B\mathbf{T}(t_{k-1}-t_k)\ldots\\
\ldots\mathcal{L}_N^B\mathbf{T}(t-t_1)\varphi_s\rangle\dd t_1\ldots\dd t_{k}.
\end{multline*}
Tedious combinatorial arguments lead to the term-by-term convergence: 
\begin{multline*}
\langle f_0^N, \mathbf{T}(t_k)\mathcal{L}_N^B\mathbf{T}(t_{k-1}-t_k)\ldots
\ldots\mathcal{L}_N^B\mathbf{T}(t-t_1)\varphi_s\rangle\\
\underset{N\to+\infty}{\longrightarrow}\langle f_0^{\otimes (s+k)}, \mathbf{T}(t_k)\mathbf{D}\mathbf{T}(t_{k-1}-t_k)\ldots
\ldots\mathbf{D}\mathbf{T}(t-t_1)\varphi_s\rangle.
\end{multline*}
Note that when the exponential series $\mathbf{T}(t)=\e^{t\mathbf{S}}$ converges, then:
\begin{multline*}\int_0^t\int_0^{t_1}\ldots\int_0^{t_{k-1}}\langle f_0^{\otimes (s+k)}, \mathbf{T}(t_k)\mathbf{D}\mathbf{T}(t_{k-1}-t_k)\ldots
\ldots\mathbf{D}\mathbf{T}(t-t_1)\varphi_s\rangle\dd t_1,\ldots \dd t_{k} \\= \frac{t^k}{k!} \langle f_{0}^{\otimes(s+k)}, (\mathbf{D}+\mathbf{T})^k\varphi_s\rangle.\end{multline*}

With the second approach, the non-homogeneous case is thoroughly detailed in~\cite{pulvirenti_kinetic_1996}. The main difference with the proof in the homogeneous case is that Equation~\cref{eq:DuhamelBoltzmannhomogeneous} should be replaced by 
\begin{multline*}f_{t}^{\otimes s} = \sum_{k=0}^{+\infty} \int_0^t\int_0^{t_1}\ldots\int_0^{t_{k-1}}\mathbf{T}_\infty^{(s)}(t-t_1)\mathcal{C}_{s,s+1}\mathbf{T}_\infty^{(s+1)}(t_1-t_2)\mathcal{C}_{s+1,s+2}\ldots\\
\mathcal{C}_{s+k-1,s+k}\mathbf{T}_\infty^{(s+k)}(t_k)f_{0}^{\otimes(s+k)}\dd t_1\ldots\dd t_k.\end{multline*}
where $\mathbf{T}_\infty^{(s)}$ is the Markov semi-group generated by $\sum_{i=1}^s L^{(1)}\diamond_i$. The domination part is similar to the homogeneous case and the term-by-term convergence becomes
\[\|(\mathbf{T}_N^{(s)}-\mathbf{T}_\infty^{(s)})f^{(s)}\|_{L^1(E^s)}\to0.\]

\subsection{Pathwise Kac's theorem via random interaction graphs}\label[II]{sec:pathwisekactheorem}

Under the same (strong) hypotheses of Kac's theorem, a  more powerful result is due to Graham and M\'el\'eard \cite{graham_stochastic_1997,meleard_asymptotic_1996}. The proof follows a completely different strategy and relies on a trajectorial representation of the process based on the notion of interaction graphs presented in the introductory Section \cref{sec:interactiongraph_summary}. Kac's theorem states a pointwise result, the following theorem works at the pathwise level. 

\begin{theorem} Let $\mathcal{L}_N$ be of the form \cref{eq:boltzmanngenerator_summary} with Assumption \cref{assum:L2_summary} and let us assume the uniform bound \cref{eq:uniformboundlambda_summary}. Let $T>0$ be a fixed time and $I=[0,T]$. Let $f_I^N\in\pb(D([0,T],E)^N)$ be the pathwise law (with initial time marginal $f_0^{\otimes N}$) of the $N$-particle system defined by $\mathcal{L}_N$ and denote by 
$f_I^{s,N}\in \pb(D([0,T],E)^s)$ its $s$-marginal for $s\in\N$. Then the following properties hold. 
\begin{enumerate}[(i)]
    \item There is propagation of chaos in total variation norm: there exists $C>0$ such that for any $s\in\N$ :  
    \begin{equation}\label[IIeq]{eq:pathwisechaosTV}\big\|f_I^{s,N}-\big(f_I^{1,N}\big)^{\otimes s}\big\|_{\mathrm{TV}}\leq Cs(s-1)\frac{\Lambda T+\Lambda^2T^2}{N},\end{equation}
    where the $\mathrm{TV}$ norm is the $s$-dimensional total variation norm.
    \item There exist $C>0$ and a probability measure $f_I\in\pb(D([0,T],E))$ such that 
    \[\big\|f_I^{1,N}-f_I\big\|_{\mathrm{TV}} \leq \frac{C\e^{\Lambda T}}{N},\]
    moreover $f_I$ solves the nonlinear Boltzmann martingale problem with initial time marginal $f_0$ (see Definition~\cref{def:nonlinearboltzmannmartingaleproblem}). 
    \item Let $(\mathcal{Z}^N_t)^{}_t$ be a particle process with law $f^N_I$. Then for all $\Phi\in C_b(D([0,T],E))$, 
    \[\E\big|\big\langle {\mu}_{\mathcal{Z}^N_I}-f_I, \Phi \big\rangle\big|^2 = \mathcal{O}(1/N).\]
\end{enumerate}
\end{theorem}
The main result is the propagation of chaos in total variation norm with an explicit convergence rate. The other properties follow more easily so we focus on the first point. 
\begin{proof}[Proof (sketch)]
The proof is based on the observation that given an interaction graph $\mathcal{G}_i(\mathcal{T}_k,\mathcal{R}_k)$ (as defined in Definition \cref{def:intearctiongraph}
), it is possible to construct a (forward) trajectorial representation of the process $(Z^i_t)^{}_t$ of particle $i$ on $[0,T]$. To do so, the particles at time $t=0$ $(Z^i_0,Z^{i_1}_0,\ldots,Z^{i_k}_0)$ are distributed according to $f_0^{k'+1,N}$, where $k'$ is the number of distinct indices $i_1,\ldots,i_k$. At each $t_\ell\in\mathcal{T}_k$, the two corresponding particles collide according to the chosen interaction mechanism and between two collisions, the particles evolve independently according to $L^{(1)}$.

According to Lemma \cref{lem:randomintearctiongraph}: if a random interaction graph (see Definition \cref{def:randomintearctiongraph}) is first sampled with rate $\Lambda$ and rooted on $i$ at time $T$, then $(Z^i_t)^{}_t$ is distributed according to $f_I^{1,N}\in\pb(D([0,T],E))$.  

Now let be given two indexes $(i,j)$ and $\mathcal{G}^N_{ij}$ the random interaction graph with rate $\Lambda$ rooted on $(i,j)$ at time $T$. Starting from either $i$ or $j$ and following the graph backward in time, it is possible to extract two interaction subgraphs, denoted respectively by $\mathcal{G}_{ij}^{i,N}$ for the subgraph rooted on $i$ and $\mathcal{G}_{ij}^{j,N}$ for the subgraph rooted on $j$. Two cases may happen: either $\mathcal{G}_{ij}^{N}$ is a connected graph or $\mathcal{G}_{ij}^N$ has two (disjoint) connected components given by the two subgraphs $\mathcal{G}_{ij}^{i,N}$ and $\mathcal{G}_{ij}^{j,N}$. We denote by $\mathscr{A}^N_{ij}$ the event ``$\mathcal{G}_{ij}^{N}$ is a connected graph''. Conditionally on the complementary event $(\mathscr{A}^N_{ij})^c$, the processes $Z^i=(Z^i_t)^{}_t$ and $Z^j=(Z^j_t)^{}_t$ are independent since their trajectorial representations depend on two disjoints sets of independent random variables. Moreover, $\mathrm{Law}(Z^i)=\mathrm{Law}(Z^j)=f_I^{1,N}$ and $\mathrm{Law}(Z^i,Z^j)=f_I^{2,N}$. Therefore, it is sufficient to look at laws conditionally on $\mathscr{A}^N_{ij}$:
\[f^{2,N}_I-f^{1,N}_I\otimes f^{1,N}_I = \Big(\mathrm{Law}(Z^i,Z^j|\mathscr{A}^N_{ij})-\mathrm{Law}(Z^i|\mathscr{A}^N_{ij})\otimes \mathrm{Law}(Z^j|\mathscr{A}^N_{ij})\Big)\PP(\mathscr{A}^N_{ij}),\]
and it holds that
\[\|f_I^{2,N}-f_I^{1,N}\otimes f_I^{1,N}\|_\mathrm{TV}\leq2\PP(\mathcal{A}^N_{ij}).\]
The question of propagation of chaos is thus reduced to the computation of the probability of sampling a connected graph. This probability can be bounded by:
\[\PP(\mathscr{A}^N_{ij})\leq\sum_{q=1}^{+\infty} Q^N_q(T),\]
where $Q^N_q(T)=\PP(\mathscr{Q}^N_q(T))$ and $\mathscr{Q}^N_q(T)$ denotes the event ``there is a route of size $q$ joining $i$ and $j$ on $[0,T]$'' (we recall that a route is simply a path in the interaction graph, see Section \cref{sec:interactiongraph_summary}
for the precise definition), as depicted on the figure below:

\begin{figure}[H]
    \centering
    \begin{tikzpicture}[>=latex]
    \draw[->] (0,0) -- (0,3.4);
    \draw[->] (0,0) -- (8,0);
    
    \draw[-,dotted] (0,2.8) -- (7.5,2.8);
    
    \node[anchor=east] at (0,2.8) {$T$};
    \node[anchor=south] at (0.5,2.8) {$i$};
    \node[anchor=south] at (7,2.8) {$j$};
    
    \coordinate (i) at (0.5,2.8);
    \coordinate (n2) at (0.5,2.2);
    \coordinate (n23) at (1.25,2.2);
    \coordinate (n3) at (2,2.2);
    \coordinate (n4) at (2,1.4);
    \coordinate (n45) at (2.5,1.4);
    \coordinate (n5) at (3,1.4); 
    \coordinate (n6) at (3,0.9);
    \coordinate (n7) at (4,0.9); 
    \coordinate (n8) at (4,0.3);
    \coordinate (n89) at (4.7,0.3);
    \coordinate (n9) at (5.4,0.3); 
    \coordinate (n10) at (5.4,1.7);
    \coordinate (n1011) at (6.2,1.7);
    \coordinate (n11) at (7, 1.7); 
    \coordinate (j) at (7,2.8); 
    
    \node[anchor=south] at (n23) {$r_{\ell_1}$};
    \node[anchor=south] at (n45) {$r_{\ell_2}$};
    \node[anchor=south] at (n89) {$r_{\ell_{q-1}}$};
    \node[anchor=south] at (n1011) {$r_{\ell_q}$};
    
    \filldraw[black] (n2) circle (2pt); 
    \filldraw[black] (n3) circle (2pt);
    \filldraw[black] (n4) circle (2pt);
    \filldraw[black] (n5) circle (2pt);
    \filldraw[black] (n8) circle (2pt);
    \filldraw[black] (n9) circle (2pt);
    \filldraw[black] (n10) circle (2pt);
    \filldraw[black] (n11) circle (2pt);
    
    \draw[-] (i) -- (n2) ; 
    \draw[-] (n2) -- (n3) ; 
    \draw[-] (n3) -- (n4);
    \draw[-] (n4) -- (n5);
    \draw[-] (n5) -- (n6);
    \draw[-,dotted] (n6) -- (n7);
    \draw[-,dotted] (n7) -- (n8);
    \draw[-] (n8) -- (n9);
    \draw[-] (n9) -- (n10);
    \draw[-] (n10) -- (n11);
    \draw[-] (n11) -- (j); 
    
    \end{tikzpicture}
    \caption{A route of size $q$ between $i$ and $j$. The chain of interactions which links $i$ and $j$ are depicted by horizontal lines as explained in Section \cref{sec:interactiongraph_summary}.}
\end{figure}

Clearly, 
\[{Q}^N_1(T) = 1-\exp{\left(-\frac{\Lambda}{N}T\right)}\leq \frac{\Lambda T}{N},\]
since this event is equal to $\{\inf_n T^{i,j}_n<T\}$. Then for $q\geq2$, to construct a route of size $q$ it is necessary to first construct a route of size 1 from either $i$ or $j$ and then a route of size $q-1$ from the new index created to the other index $i$ or $j$. Since branching happens with a rate bounded by $\Lambda$, it holds that
\[{Q}^N_q(T) \leq \int_0^T {Q}^N_{q-1}(T-t)2\Lambda\exp\left(-2\Lambda t\right)\dd t=Q^N_{q-1}\star \e_{2\Lambda}(T),\]
where $\e_{2\Lambda}$ is the density of the exponential law with parameter $2\Lambda$. Therefore 
\[Q^N_q(T)\leq Q^N_1\star\e_{2\Lambda}^{\star(q-1)}(T),\]
and a direct computation shows that
\[\sum_{q=1}^{+\infty}{Q}^N_q(T)\leq C\frac{\Lambda T+(\Lambda T)^2}{N}.\]
The same reasoning extends for any interaction graph rooted on an arbitrary number of particles and gives the estimate \cref{eq:pathwisechaosTV}. This ends the proof of the first point of the theorem. The remaining steps are sketched below.

\begin{enumerate}
    \item With a similar reasoning, it is possible to prove that the law of any particle converges towards the law $f_I$ of the process constructed on a limit Boltzmann tree with rate $\Lambda$. To do so, the argument is based on an estimate on the probability that there is a recollision in the sampled random graph. Since as $N\to+\infty$ the number of branches is bounded (of the order $\e^{\Lambda T}$), and that the Poisson processes have rate $\Lambda/N$ it holds that 
    \[\|f^{1,N}_I-f^{}_I\|_{\mathrm{TV}}\leq C\frac{\e^{\Lambda T}}{N}.\]
    \item Since the convergence holds in total variation norm, the empirical measure process converges in probability and in law towards $f_I$. 
    \item It remains to prove that the law $f_I$ satisfies the nonlinear martingale problem. As in the McKean-Vlasov case (see Section \cref{sec:mckeancompactnessreview}), it can be proved by passing to the limit in the martingale problem satisfied by the $N$-particle system (which is possible thanks to the previous step). 
\end{enumerate}
We refer the reader to \cite{meleard_asymptotic_1996} for the details of the proof. 
\end{proof}

\subsection{Martingale methods}\label[II]{sec:martingaleboltzmannreview}

The probabilistic treatment of the spatially homogeneous version of the Boltzmann equation of rarefied gas dynamics \cref{eq:Boltzmannphysics_summary} and the question of proving propagation of chaos via martingale techniques has been initiated by \cite{tanaka_probabilistic_1983}. Such techniques lead to very powerful results as they only rely on abstract compactness criteria which apply on the path space. A drawback of the approach is that it does not provide any rate of convergence. The framework is briefly explained in the introductory Section \cref{sec:provingcompactness_summary}. The paradigmatic proof of strong pathwise empirical chaos is due to Sznitman \cite{sznitman_equations_1984}. The strategy is quite general, it does not restrict to Boltzmann-like models and can be applied to various models, in particular diffusion or jump models. A complete example in the case of McKean-Vlasov diffusion with jumps is shown in Section \cref{sec:martingalecompactness}. In this section we make some comments specific to Boltzmann models and state the final result of~\cite{sznitman_equations_1984}. Then we extend the functional law of large numbers (Theorem \cref{thm:martingaleweakpathwiseempiricalchaos}) proved in the mean-field case to general Boltzmann models. 

\subsubsection*{Strong pathwise empirical chaos.}
Sznitman \cite{sznitman_equations_1984} considers Boltzmann parametric models as given by Definition \cref{def:Boltzmannparammodel_summary} in $E=\R^d$ and such that there exists a function $\psi_0:E\times E\times \Theta\to E$ which satisfies for all $z_1,z_2\in E$ and $\theta\in\Theta$,
\[\psi_0(z_1,z_2,\theta)=\psi_1(z_1,z_2,\theta)=\psi_2(z_2,z_1,\theta),\]
that is, $\psi(z_1,z_2,\theta)=(\psi_0(z_1,z_2,\theta),\psi_0(z_2,z_1,\theta))$.

The assumptions on the interaction function $\psi_0$ are as follows.
\begin{assumption}\label[II]{assum:sznitmanpsi} There exists a continuous function $m:E\to\R_+$ such that $m\geq 1$, $\lim_{|z|\to+\infty}m(z)=+\infty$ and such that the interaction function $\psi_0$ and the interaction rate $\lambda$ satisfy:
\begin{enumerate}[(i)]
    \item for all $z_1,z_2\in E$ and all $\theta\in\Theta$,
    \[m(\psi_0(z_1,z_2,\theta))+m(\psi_0(z_2,z_1,\theta))\leq m(z_1)+m(z_2),\]
    \item there exists some real $p$ with $0\leq p\leq1$, such that for all $z_1,z_2\in E$ 
    \[\lambda(z_1,z_2)\leq m(z_1)^p+m(z_2)^p.\]
\end{enumerate}
\end{assumption}
In most cases, the function $m$ is a polynomial function of the form $m(z)=1+|z|^k$ and the above assumptions are thus mostly used to control the moments of the particle system or of the limiting equation which is often a crucial in Boltzmann models. Sznitman uses the martingale characterisation of the $N$-particle system.

\begin{assumption}\label[II]{assum:boltzmannmartingalestrongpathwise}
For any $T\in[0,+\infty]$, well-posedness holds true for the martingale problem associated to the $N$-particle system (Definition \cref{def:martingaleproblemparticles}) supplemented with the condition: for all $t>0$,
\begin{multline*}\int_{D(\R_+,E)^N}{\left[m(\mathsf{Z}^{1,N}_t)+\ldots+m(\mathsf{Z}^{N,N}_t)\right]} f^N_I\big(\dd\mathbf{Z}^N\big) \\\leq \int_{E^N} {\left[m(z^1)+\ldots+m(z^N) f^N_0\big(\dd\mathbf{z}^N\big)\right]}.\end{multline*}
\end{assumption}
 
The main result \cite[Theorem 3.3]{sznitman_equations_1984} is the following. 
\begin{theorem}\label[II]{thm:sznitmanboltzmann} Let us assume that Assumptions \cref{assum:sznitmanpsi} and \cref{assum:boltzmannmartingalestrongpathwise} hold true. Let $f_0\in\pb(E)$ and let $(f^N_0)^{}_N$ a sequence of $f_0$-chaotic probability measures on $E^N$. Assume that
\begin{enumerate}[(i)]
    \item there exists $C>0$ such that for all $N\geq1$, $\frac{m(Z^1)+\ldots+m(Z^N)}{N}\leq C$ $f^N_0$-almost surely, 
    \item $\sup_N\int_{E^N} m(z)^{1+p}  f^{1,N}_0(\dd z)<+\infty$.
\end{enumerate}
Then the laws $f^N_I\in \pb(D(\R_+,E)^N)$ are $f_I$-chaotic where $f_I\in\pb(D(\R_+,E))$ satisfies the nonlinear Boltzmann martingale problem (Definition \cref{def:nonlinearboltzmannmartingaleproblem}) supplemented with the condition 
\[\forall T>0,\quad\sup_{t\leq T} \int_{D(\R_+,E)}m(\mathsf{Z}_t)\dd f_I(\dd\mathsf{Z})<+\infty.\]
\end{theorem}

The theorem states the usual pathwise propagation of chaos result. It is obtained as a consequence of the strong pathwise empirical propagation of chaos. This setting includes the case of the hard-sphere cross-section.

\subsubsection*{Functional law of large numbers.}

Wagner \cite{wagner_functional_1996} proves a functional law of large numbers for Boltzmann parametric models of the form \cref{eq:Boltzmannnonsym_summary} with $L^{(1)}\ne 0$, adding some individual flow to particles. The proof is based on compactness arguments and a pointwise martingale characterisation of the particle system. The nonlinear process is defined by a series expansion reminiscent from Kac theorem (Theorem~\cref{thm:kac}). 

To conclude this section, we wish now to briefly discuss the extension of the method of Theorem \cref{thm:martingaleweakpathwiseempiricalchaos} to Boltzmann-type collision systems. The first part of Assumption \cref{assum:initialdefined} has to be replaced by

\begin{assumption}[Boltzmann generator]\label[II]{assum:initialboltzmannmartingale}
The generator of the process $(\mathcal{X}^N_t)^{}_t$ is a Boltzmann generator \cref{eq:boltzmanngenerator_summary} with $L^{(1)}=0$ and $L^{(2)}$ which satisfies Assumption \cref{assum:L2_summary}. Moreover the associated martingale problem (Definition \cref{def:martingaleproblemparticles}) is wellposed and the initial distribution satisfies the second assumption of Theorem \cref{thm:sznitmanboltzmann} with $m(z) = z$ and some $p>0$. 
\end{assumption}

We define the symmetrized version of $L^{(2)}$
\begin{align*}
 L^{(2)}_{\mathrm{sym}} \varphi_2 ( z_1,z_2 ) &= \frac{L^{(2)} \varphi_2 ( z_1,z_2 ) + L^{(2)} \varphi_2 ( z_2,z_1 )}{2} \\
 &= \frac{\lambda ( z_1,z_2 )}{2} \iint_{E^2} \Big\{\varphi_2 ( z'_1,z'_2 ) + \varphi_2 ( z'_2,z'_1 ) \\
& \phantom{abcdefghijklmnopq} - \varphi_2 ( z_1,z_2 ) - \varphi_2 ( z_2,z_1 )\Big\} \Gamma^{(2)} ( z_1' , z_2' \dd z_1' , \dd z_2'). 
\end{align*}
This implies $L^{(2)}_{\mathrm{sym}} [ \varphi^1 \otimes \varphi^2 ] = L^{(2)}_{\mathrm{sym}} [ \varphi^2 \otimes \varphi^1 ]$ for every  $\varphi^1 , \varphi^2 \in \mathcal{F}$. For the limit generator, given $\mu\in\pb(E)$, we define $L_{\mu}$ as
\begin{equation}\label[IIeq]{eq:Lmuboltzmann} \forall \varphi \in \mathcal{F},\,\, \forall x \in E,  \quad L_{\mu} \varphi ( x ) := {\left\langle \mu , L^{(2)}_{\mathrm{sym}} [ \varphi \otimes 1 ] ( x , \cdot )  \right\rangle} = {\left\langle \mu , L^{(2)}_{\mathrm{sym}} [ \varphi \otimes 1 ] ( \cdot , x )  \right\rangle}, \end{equation}
and equivalently $\varphi \otimes 1$ can be taken instead of $1 \otimes \varphi$ in the above definition. With this definition, the general Boltzmann equation \cref{eq:symmetricBoltzmannequationgeneral_summary} can be rewritten as in the mean-field case: 
\begin{equation}\label[IIeq]{eq:BoltzmannequationLmu}\forall \varphi\in\mathcal{F},\quad \frac{\dd}{\dd t}\langle f_t,\varphi\rangle = \langle f_t,L_{f_t}\varphi\rangle.\end{equation}
we recall the notation: 
\[\forall\nu\in\pb(E),\quad R_{\varphi^1\otimes\varphi^2}(\nu) := \langle \nu^{\otimes 2},\varphi^1\otimes\varphi^2\rangle,\]
for the polynomial function on $\pb(E)$ associated to $\varphi_2=\varphi^1\otimes\varphi^2\in\mathcal{F}^{\otimes 2}$. We will need the following quadratic estimate: 

\begin{lemma}[Quadratic estimate for Boltzmann collisions]\label[II]{lemma:quadraticestimateboltzmann}
The quadratic estimates reads 
\begin{multline*}
\mathcal{L}_N {\left[ R_{\varphi^1 \otimes \varphi^2} \circ \boldsymbol{\mu}_N \right]}\big(\mathbf{x}^N\big) = R_{L_{ \mu_{\mathbf{x}^N} } \varphi^1 \otimes \varphi^2} \big(\mu_{\mathbf{x}^N}\big) + R_{\varphi^1 \otimes L_{ \mu_{\mathbf{x}^N} } \varphi^2} \big(\mu_{\mathbf{x}^N}\big)\\+ \frac{1}{N} R_{L^{(2)}_{\mathrm{sym}} [ \varphi^1 \otimes \varphi^2 ]} \big(\mu_{\mathbf{x}^N}\big) 
+ \frac{1}{N} {\left\langle \mu_{\mathbf{x}^N} , \Gamma_{L_{\mu_{\mathbf{x}^N}}} \left( \varphi^1 , \varphi^2 \right)  \right\rangle}.
\end{multline*}
\end{lemma}

\begin{proof} See Lemma \cref{lemma:quadraticboltzmann} in the appendix.
\end{proof}

Compared to the mean-field case \cref{eq:mean-fieldcarre}, a correcting crossed-term appears for Boltzmann collisions, but this term can be handled in the same way by Assumption~\cref{assum:genebound}. One can eventually state the propagation of chaos theorem.

\begin{theorem}[Functional law of large numbers for Boltzmann models]
Let us assume that Assumptions \cref{assum:initialboltzmannmartingale}, \cref{assum:genebound}, \cref{assum:continu} and \cref{assum:wellposedweakpathwise} hold true for $L_\mu$ given by \cref{eq:Lmuboltzmann}. Then the weak Boltzmann equation \cref{eq:BoltzmannequationLmu} is wellposed and weak pathwise empirical propagation of chaos towards its solution holds for the Boltzmann model on every time interval $[0,T]$. 
\end{theorem}

\begin{proof}[Proof (sketch)] The proof is exactly the same as the one in the mean-field case (Theorem \cref{thm:martingaleweakpathwiseempiricalchaos}). The mean-field property reads this time
\[ \mathcal{L}_N \bar{\varphi}_N \big( \mathbf{x}^N \big) = {\left\langle \mu_{\mathbf{x}^N} \otimes \mu_{\mathbf{x}^N} , L^{(2)}_{\mathrm{sym}} {\left[ \varphi \otimes 1 \right]} \right\rangle} = {\left\langle \mu_{\mathbf{x}^N} , L_{\mu_{\mathbf{x}^N}} \varphi \right\rangle}, \] 
and It\=o's formula can be written the same way to complete Step $1$. The control of the carr\'e du champ is provided by Lemma \cref{lemma:quadraticestimateboltzmann} above. Step $2$ and Step $3$ are identical provided that $L_{ \mu }$ satisfies the boundedness continuity and uniqueness assumptions.
\end{proof}

\subsection{SDE and coupling} \label[II]{sec:couplingBoltzmann}

In this section, we continue the discussion started at the end of Section \cref{sec:boltzmannparametricmodels} and we prove propagation of chaos for a class of Boltzmann parametric models (Definition \cref{def:Boltzmannparammodel_summary}) using a coupling argument based on a SDE representation of the particle system. The main theorem of this section is due to Murata \cite{murata_propagation_1977} in the particular case of the 2D true Maxwellian molecules (non-cutoff). The technique of the proof has recently been revisited in \cite{cortez_quantitative_2016, cortez_quantitative_2018}. The proof in this section globally follows the same presentation as in \cite{murata_propagation_1977} although we sometimes use modernised optimal transport arguments taken from \cite{cortez_quantitative_2016}. The classical nonlinear SDE representation of the Boltzmann equation originally due to Tanaka \cite{tanaka_probabilistic_1978} for \cref{eq:Boltzmannphysics_summary} can be found in the proof. Let us first recall the setting of Definition \cref{def:Boltzmannparammodel_summary}: we take $E=\R^d$, we assume that the collision rate is constant $\lambda\equiv\Lambda$ and that the post-collisional distribution $\Gamma^{(2)}$ is of the following form: for any $\varphi_2\in C_b(E^2)$,
\begin{equation}\label[IIeq]{eq:gamma2sdecoupling}\iint_{E\times E} \varphi_2(z_1',z_2')\Gamma^{(2)}(z_1,z_2,\dd z_1',\dd z_2') = \int_{\Theta} \varphi_2(\psi_1(z_1,z_2,\theta),\psi_2(z_1,z_2,\theta))\nu(\dd\theta),\end{equation}
with $(\psi_1,\psi_2)(z_1,z_2,\cdot)_{\#}\nu = (\psi_2,\psi_1)(z_2,z_1,\cdot)_{\#}\nu$. We make the following reasonable Lipschitz and growth assumptions. 
\begin{assumption}\label[II]{assum:boltzmannsdeassum} The interaction functions $\psi_1,\psi_2$ satisfy the following properties. 
\begin{enumerate}[(i)]
    \item (Lipschitz). There exists a function $L\in L^1_\nu(\Theta)$ such that for $i=1,2$, 
    \[\forall (\theta,z_1,z_2,z_1',z_2')\in\theta\times E^4,\quad |\psi_i(z_1,z_2,\theta)-\psi_i(z_1',z_2',\theta)|\leq L(\theta)(|z_1-z_1'|+|z_2-z_2'|).\]
    \item (Linear growth). There exists a function $M\in L^1_\nu(\Theta)$ such that for $i=1,2$, 
    \[\forall (\theta,z_1,z_2)\in \Theta\times E\times E,\quad |\psi_i(z_1,z_2,\theta)|\leq M(\theta)(1+|z_1|+|z_2|).\]
\end{enumerate}
\end{assumption}

\begin{remark}
One can alternatively assume (it is maybe more classical) that: 
\[\forall (z_1,z_2,z_1',z_2')\in E^4,\quad \int_\Theta |\psi_i(z_1,z_2,\theta)-\psi_i(z_1',z_2',\theta)|\nu(\dd\theta)\leq C(|z_1-z_1'|+|z_2-z_2'|),\]
for a constant $C>0$ and similarly for the linear growth assumption.
\end{remark}

\begin{remark}
It should be noted that the Lipschitz assumption does \emph{not} hold true for 3D Maxwell molecules. However, it holds true for 2D Maxwell molecules which is the original setting of Murata's proof. This non-Lipschitz issue in the 3D case has been encompassed by Tanaka in \cite{tanaka_probabilistic_1978}.
\end{remark}

Under the assumption of linear growth, it follows easily using Gronwall lemma that the moments of all order are exponentially controlled for the nonlinear process.
\begin{lemma}\label[II]{lemma:boundmomentboltzmannsde} For all $p\geq1$, there exist $C(p)>0$ such that for all $t>0$, 
\[\int_E |z|^p f_t(\dd z) \leq \left(\int_E |z|^p f_0(\dd z)\right)\e^{C(p)t}.\]
\end{lemma}
Without loss of generality (up to redefining a process with fictitious collisions), we also assume that the interaction rate is a constant $\Lambda$ and for all $\theta\in\Theta,z\in\R^d$, $\psi_i(z,z,\theta)=z$. A system of stochastic differential equations corresponding to the particle system is given by:
\begin{equation}\label[IIeq]{eq:particleboltzmannsde}Z^i_t = Z^i_0 + \sum_{j\ne i}\int_0^t \int_\Theta\int_{\{0,1\}} a\big(Z^i_{s^-},Z^j_{s-},\theta,\sigma\big)\mathcal{N}_{ij}(\dd s,\dd\theta,\dd \sigma).\end{equation}
where
\[a(z_1,z_2,\theta,\sigma) = (1-\sigma)\psi_1(z_1,z_2,\theta)+\sigma\psi_2(z_2,z_1,\theta)-z_1.\]
For all $i,j$, $\mathcal{N}_{ij}$ is a Poisson random measure on $\R_+\times\Theta\times\{0,1\}$ with intensity $\frac{\Lambda}{N}\dd t\nu(\dd\theta)\dd\sigma$, where $\dd\sigma$ is the uniform measure on $\{0,1\}$. We also assume that for all $i,j$, the Poisson measure satisfy:
\[\mathcal{N}_{ij}=\check{\mathcal{N}}_{ji},\]
where for a Poisson measure $\mathcal{N}$ on $\R_+\times\Theta\times\{0,1\}$ with intensity $\frac{\Lambda}{N}\dd t\nu(\dd\theta)\dd\sigma$, we write 
\[\check{\mathcal{N}}(\mathscr{B}) = \mathcal{N}(\check{\mathscr{B}}),\]
where given a measurable set $\mathscr{B}\subset \R_+\times\Theta\times\{0,1\}$,
\[\check{\mathscr{B}} := \{(t,\theta,\sigma)\,|\,(t,\theta,1-\sigma)\in\mathscr{B}\}.\]
Classical results and classical references on this type of SDEs can be found in Appendix \cref{appendix:poisson}. The main result of this section is the following coupling estimate. 

\begin{theorem}\label[II]{thm:boltzmannsdecoupling}
Let $T>0$. Let $f_0\in\pb_1(E)$ and $(Z^i_0)_{1\leq i\leq N}$ be $N$ independent initial random variables with common law $f_0$. Let us assume that $\mathcal{L}_N$ is of the form \cref{eq:boltzmanngenerator_summary} with $\Gamma^{(2)}$ given by \cref{eq:gamma2sdecoupling} and $\lambda$ being a constant $\Lambda$, together with Assumption \cref{assum:boltzmannsdeassum}. Then there exist
\begin{itemize}
\item a $N$-particle system $\mathcal{Z}^N_t$ with law $f^N_t$,
\item $N$ nonlinear processes $\overline{\mathcal{Z}}{}^N_t$ which are independent and identically distributed with common law $f_t$ solution of the Boltzmann equation \cref{eq:Boltzmannequationgeneral_summary}, 
\end{itemize}
which satisfy the following property: there exists a constant $C(T)>0$ depending only on $T$ and the constants in Assumption \cref{assum:boltzmannsdeassum} such that for any constant $\eta<(2d+1)^{-1}$ and for all $1\leq i\leq N$ it holds that
\[\sup_{0\leq t \leq T}\E{\big|Z^i_t-\overline{Z}{}^i_t\big|}\leq C(T)\left(\frac{1}{N^\eta}+\frac{i}{N}\right).\]
\end{theorem}

\begin{remark}
Note that the particles defined by the processes $(\mathcal{Z}^N_t)^{}_t$ and $(\overline{\mathcal{Z}}^N_t)^{}_t$ are jointly constructed and numbered so that the coupling bound holds true for the~$i$-th particle for any fixed $i$. Although the bound depends on the particle numbering, it does not contradict the exchangeability: the random variables $(Z^1_t,\ldots,Z^N_t)$ are exchangeable and $(\overline{Z}^1_t,\ldots,\overline{Z}^N_t)$ are i.i.d. but it does not imply that the random variables $(Z^1_t-\overline{Z}^1_t,\ldots,Z^N_t-\overline{Z}^N_t)$ are exchangeable (and indeed they are not). Nevertheless, the coupling bound still provides propagation of chaos for any block of size $k$ (and even for $k\equiv k(N)\to+\infty$ with $N$ provided that $k=o(N)$). Namely, for any fixed $1\leq k\leq N$, Theorem \cref{thm:boltzmannsdecoupling} provides a coupling bound for the $k$ first particles (which are $f^{k,N}_t$-distributed) which implies that
\[\sup_{0\leq t\leq T} W_1(f_t^{k,N},f_t^{\otimes k})\leq C(T)\left(\frac{1}{N^\eta}+\frac{k}{N}\right).\]
As usual the case $k=2$ is sufficient. This coupling method requires more subtle arguments than the analog for McKean-Vlasov processes (Theorem \cref{thm:mckean}) where the starting point can be \emph{any} particle system with $N$ arbitrary independent Brownian motions. In the present case, in order to prove the desired coupling bound, the analogous Poisson random measures are constructed in the proof and define a specific particle numbering. 
\end{remark}

\begin{proof} Following Murata's work, the proof is split into several steps. The first step is devoted to the construction of the particle system. In the second step, the particle system is coupled with a system of independent nonlinear SDEs \emph{\`a la} Tanaka. The third step introduces an intermediate system of non independent processes which is used as a pivot between the particle system and the nonlinear system. In the fourth and fifth steps we use the coupling to derive explicit error estimates and we conclude the proof.
\medskip

\noindent\textit{\textbf{Step 1.} Construction of a particle system.}
\medskip

Following Murata's work, let us define $N^2$ independent Poisson random measures $\overline{\mathcal{N}}_{ij}$, indexed by $1\leq i, j\leq N$, on $\R_+\times\Theta\times\{0,1\}\times(0,\frac{1}{N}]$ with intensity $\Lambda\dd t\nu(\dd\theta)\dd\sigma\dd\alpha$. We consider the following filtration: 
\[\mathscr{F}_t=\sigma\Big(Z^i_0,\, \mathcal{N}_{ij}(\mathscr{B}),\,\,1\leq i,j\leq N,\,\mathscr{B}\,\,\text{measurable subset}\Big).\]
We define: 
\[\mathcal{N}_{ij} = \left\{
\begin{array}{rcl}
\overline{\mathcal{N}}_{ij} & \text{if} & i\leq j\\
\check{\overline{\mathcal{N}}}_{ji} & \text{if} & j>i
\end{array}
\right.,\]
so that $\mathcal{N}_{ij}=\check{\mathcal{N}}_{ji}$. We write
\[\mathcal{N}_{ij}(\dd s,\dd\theta,\dd\sigma) \equiv \int_{\alpha\in(0,1/N]}\mathcal{N}_{ij}(\dd s,\dd\theta,\dd\sigma,\dd\alpha),\]
so that $\mathcal{N}_{ij}(\dd s,\dd \theta)$ is a random Poisson measure on $\R_+\times\Theta$ with intensity $\frac{\Lambda}{N}\dd t\nu(\dd \theta)$. With this choice of Poisson measures, let $(\mathcal{Z}^N_t)^{}_t$ be the $\mathscr{F}_t$-adapted particle system given by Equation \cref{eq:particleboltzmannsde}. We can write:
\begin{equation}\label[IIeq]{eq:particlebotlzmannsdemu}Z^i_t = Z^i_0 + \int_0^t \int_\Theta\int_{\{0,1\}}\int_0^1 a\big(Z^i_{s^-},Z^\mu_{s^-},\theta,\sigma\big)\mathcal{N}_{i}(\dd s,\dd\theta,\dd\sigma,\dd \alpha),\end{equation}
where for each $\omega\in\Omega$, $t\in[0,T]$ and $\alpha\in[0,1]$, we define 
\[Z^\mu_t(\omega,\alpha) := \sum_{j=1}^N \1_{(\frac{j-1}{N},\frac{j}{N}]}(\alpha)Z^j_t(\omega),\]
and for $\mathscr{B}$ a measurable subset of $\R_+\times\Theta\times\{0,1\}\times[0,1]$, we define the Poisson random measure:
\[\mathcal{N}_i(\mathscr{B}) := \sum_{j=1}^N \mathcal{N}_{ij}(\mathscr{B}_j),\]
where
\[\mathscr{B}_j := \Big\{(t,\theta,\sigma,\alpha)\in\R_+\times\Theta\times\{0,1\}\times(0,1/N]\,\big|\,(t,\theta,\sigma,\alpha+(j-1)/N)\in\mathscr{B}\Big\}.\]
The key observation is the following: for each $\omega\in\Omega$, $Z^\mu_t(\omega,\alpha)$ is a $E$-valued process constructed on the probability space $([0,1],\dd\alpha)$ such that the $\alpha$-law of $Z^\mu_s(\omega)$ is $\hat{\mu}_{\mathcal{Z}^N_t}(\dd z)$. In the following, we call $\alpha$-random variable a random variable constructed on $([0,1],\dd\alpha)$, its law is called the $\alpha$-law and we denote by $\E_\alpha$ the expectation on this space. 
\medskip

\noindent\textit{\textbf{Step 2.} Construction of a nonlinear system and coupling.}
\medskip

First let us define the random Poisson measures on $\R_+\times\Theta\times\{0,1\}\times[0,1]$: 
\[\overline{\mathcal{N}}_i(\mathscr{B}) = \sum_{j=1}^N \overline{\mathcal{N}}_{ij}(\mathscr{B}_j).\]
They are independent. In \cite{tanaka_probabilistic_1978}, Tanaka introduced the following stochastic version of the Boltzmann equation:
\begin{equation}\label[IIeq]{eq:tanakasde}\overline{Z}{}^i_t = Z_0^i + \int_0^t\int_\Theta\int_{\{0,1\}}\int_0^1a\big(\overline{Z}{}^i_{s^-},{Y}_{s^-},\theta,\sigma\big)\overline{\mathcal{N}}_{i}(\dd s,\dd\theta,\dd\sigma,\dd \alpha),\end{equation}
where for each $t$ and $\omega$, $\overline{Y}_t(\omega,\alpha)$ is a $E$-valued $\alpha$-random variable with $\alpha$-law $\mathrm{Law}(\overline{Z}{}^i_t)$. It can be checked that the $\overline{Z}{}^i_t$ are independent and identically distributed with common law $f_t$ the solution of the Boltzmann equation. Note that as in the McKean-vlasov case, this defines a class of processes given by a SDE which depends on the own law of the process. 

Note that the above nonlinear processes are already coupled with the particle system \cref{eq:particleboltzmannsde} through the Poisson random measures and the initial condition. We go one step further by choosing an appropriate process $Y$ which couples optimally the solution of the Boltzmann equation and the emprirical measure of the particle system. We take the process $Y$ given by the following key lemma. 

\begin{lemma}[Optimal empirical coupling]\label[II]{lemma:boltzmannoptimalcoupling}
There exists a process $Y=Y_t(\omega,\alpha)$ such that
\begin{enumerate}[(i)]
    \item $(Y_t)_t$ is $\mathscr{F}_t$-predictable
    \item For each $t$ and $\omega$, the $\alpha$-law of $Y_t(\omega)$ is $f_{t^-}$. 
    \item For each $t$ and $\omega$, 
    \[\E_\alpha\big[\big|Z_t^\mu(\omega)-Y_t(\omega)\big|\big] = W_1\big(\mu_{\mathcal{Z}^N_t},f_t\big).\]
\end{enumerate}
\end{lemma}

\begin{proof} Using \cite[Corollary 5.22]{villani_optimal_2009}, we know that there exists a measurable mapping 
\[\R_+\times \Omega\to \pb(E\times E),\,\,(t,\omega)\mapsto \pi_{t,\omega},\]
such that $\pi_{t,\omega}$ is an optimal transfer plan between $\mu_{\mathcal{Z}^N_t}$ and $f_t$. Let us define for $j\in\{1,\ldots,N\}$ and $\mathscr{B}$ a measurable subset of $E$,
\[G^j_{t,\omega}(\mathscr{B}) = \frac{\pi_{t,\omega}(\mathscr{B}\times\{Z^j_t\})}{\pi_{t,\omega}(E\times\{Z^j_t\})}=:\pi_{t,\omega}\big(\mathscr{B}\times\{Z^j_t\}|E\times\{Z^j_t\}\big).\]
Using a randomization lemma there exists an $\alpha$-random variable $g^j_{t,\omega}(\alpha)$ on the probability space $\big([0,\frac{1}{N}],N\dd\alpha\big)$ such that the $\alpha$-distribution of $g^j_{t,\omega}$ is equal to  $G^j_{t,\omega}$. Then, let us define for $\alpha\in[0,1]$,
\[Y_t(\omega,\alpha) :=\sum_{j=1}^N \1_{I_j}(\alpha)g^j_{t,\omega}\left(\alpha-\frac{j-1}{N}\right),\]
where $I_j = [(j-1)/N,j/N]$. Then one can check that
\[\PP_\alpha\big(\{Y_t(\alpha)\in\mathscr{B}\}\cap\{Z^\mu_t(\alpha)=Z^j_t\}\big) = \pi_{t,\omega}(\mathscr{B}\times\{Z^j_t\}),\]
which concludes the proof. 
\end{proof}

The third property (optimal coupling) and the above proof are exactly the content of \cite[Lemma 3]{cortez_quantitative_2016}. Murata was obviously not aware of the optimal transport results that we used but he managed to prove the existence of a coupling which is optimal up to an arbitrary $\varepsilon>0$ which is enough for the rest of the argument.  

Note that with this choice of $Y$, it is not clear anymore whether the nonlinear processes \cref{eq:tanakasde} remain independent. Fortunately they are, as stated in the following lemma. 

\begin{lemma} The processes $(\overline{\mathcal{Z}}{}^N_t)^{}_t$ satisfy the following properties. 
\begin{enumerate}
    \item They are well defined $\mathscr{F}_t$-adpated processes
    \item They are identically distributed and their law is a weak measure solution of the Boltzmann equation \cref{eq:Boltzmannequationgeneral_summary}.  
    \item They are independent. 
\end{enumerate}
\end{lemma}

\begin{proof}[Proof (sketch)]
The first two properties follow from Tanaka's construction \cite{tanaka_probabilistic_1978} which are summarised in Murata's article \cite[Theorem 4.1 and Theorem 4.2]{murata_propagation_1977}. The independence is proved in \cite[Lemma 6.4]{murata_propagation_1977} (see also the proof of \cite[Lemma~6]{cortez_quantitative_2016}). The idea is to prove (using elementary martingale properties) the independence of the measures defined by
\[\overline{\mathcal{N}}_i^{\#}(\mathscr{B}) := \int_{\R_+\times \Theta\times\{0,1\}\times[0,1]} \1_\mathscr{B}(s,\theta,\sigma,Y_s(\omega,\alpha))\overline{\mathcal{N}}_i(\dd s,\dd \theta,\dd\sigma,\dd\alpha) \]
for any measurable subset $\mathscr{B}\subset \R_+\times\Theta\times\sigma\times E$. 
\end{proof}

\noindent\textit{\textbf{Step 3.} An intermediate process.}
\medskip

At this point, we have defined $N$ couples of processes $(Z^i,\overline{Z}{}^i)$ with the correct laws and the nonlinear processes are independent. We are exactly in the good position to prove the theorem. To carry out the proof let us notice that there are actually two couplings. In addition to the optimal coupling defined by Lemma \cref{lemma:boltzmannoptimalcoupling}, there is also a coupling between the jump times and between the jump random variables given by the Poisson measures $\mathcal{N}_i$ and $\overline{\mathcal{N}}_i$ which are not independent. As in Murata's proof, we separate these two sources of discrepancy by writing: 
\begin{equation}\label[IIeq]{eq:triangleZZbarZtilde}\E|Z^i_t-\overline{Z}{}^i_t|\leq \E\big|Z^i_t-\widetilde{Z}^i_t\big|+\E\big|\widetilde{Z}^i_t-\overline{Z}{}^i_t\big|,\end{equation}
where the process $\widetilde{Z}^i_t$ is defined by:
\[\widetilde{Z}^i_t = Z_0^i + \int_0^t\int_\Theta\int_{\{0,1\}}\int_0^1a\big(\widetilde{Z}^i_{s^-},{Y}_{s^-},\theta,\sigma\big)\mathcal{N}_{i}(\dd s,\dd\theta,\dd\sigma,\dd \alpha).\]
Note that the processes $\widetilde{Z}^i$ are exchangeable but not independent. In \cite{cortez_quantitative_2016}, Cortez and Fontbona consider only (an analog of) these processes and introduce later the nonlinear processes. Nevertheless, these intermediate processes propagate chaos: this follows from the following coupling bounded which holds true for the $k$ first processes $\widetilde{Z}^1,\ldots,\widetilde{Z}^k$ with any fixed $k$. This result can be found in both works, \cite[Lemma~6]{cortez_quantitative_2016} and \cite[Lemma 6.5]{murata_propagation_1977}.

\begin{lemma}\label[II]{lemma:tildeZbarZ} There exists a constant $C(T)>0$ depending only on $T$ and the constants in Assumption \cref{assum:boltzmannsdeassum} such that for all $i\leq N$,
\begin{equation}\label[IIeq]{eq:boundZbarZtilde}\sup_{0\leq t\leq T}\E\big|\widetilde{Z}^i_t-\overline{Z}{}^i_t\big|\leq C(T)\frac{i}{N}.\end{equation}
\end{lemma}

\begin{proof} Let $i\leq N$. Using the definition of the Poisson random measures $\mathcal{N}_i$ and $\overline{\mathcal{N}}_i$, we write: 
\begin{align*}
    &\big|\overline{Z}{}^i_t-\widetilde{Z}^i_t\big| \\ 
    &\leq\int_0^t\int_\Theta\int_{\{0,1\}}\int_{\frac{i-1}{N}}^{1} \big| a\big(\overline{Z}{}^i_{s^-},{Y}_{s^-},\theta,\sigma\big)-a\big(\widetilde{Z}^i_{s^-},{Y}_{s^-},\theta,\sigma\big)\big|\overline{\mathcal{N}}_i(\dd s,\dd\theta,\dd\sigma,\dd\alpha)\\
    &\quad+\sum_{j=1}^{i-1}\int_0^t\int_\Theta\int_{\{0,1\}}\int_{0}^{1/N} \big|a\big(\overline{Z}{}^i_{s^-},{Y}_{s^-}\big(\omega,\alpha+\frac{j-1}{N}\big),\theta,\sigma\big)\big|\overline{\mathcal{N}}_{ij}(\dd s,\dd\theta,\dd\sigma,\dd\alpha)\\
    &\quad+\sum_{j=1}^{i-1}\int_0^t\int_\Theta\int_{\{0,1\}}\int_{0}^{1/N} \big|a\big(\widetilde{Z}^i_{s^-},{Y}_{s^-}\big(\omega,\alpha+\frac{j-1}{N}\big),\theta,1-\sigma\big)\big|\overline{\mathcal{N}}_{ji}(\dd s,\dd\theta,\dd\sigma,\dd\alpha)
\end{align*}
Using the assumptions on the functions $\psi_1$ and $\psi_2$ and since $i\leq k$, it leads after taking the expectation to: 
\begin{multline*}
    \E\big|\overline{Z}{}^i_t-\widetilde{Z}^i_t\big|\leq c_1\int_0^t \E\big|\overline{Z}{}^i_s-\widetilde{Z}^i_s\big|\dd s + c_2\frac{i}{N} \int_0^t \E\Big[\big|\overline{Z}{}^i_s|+|\widetilde{Z}^i_s\big|\Big]\dd s\\
    +2c_3\sum_{j=1}^{i-1}\int_0^t \E\left[\int_{(j-1)/N}^{j/N} |Y_s(\omega,\alpha)|\dd\alpha\right]\dd s
\end{multline*}
With the notations of the proof of Lemma \cref{lemma:boltzmannoptimalcoupling}, one can see that 
\[\E\left[\int_{(j-1)/N}^{j/N} |Y_s(\omega,\alpha)|\dd\alpha\right] = 
\E\int_E |z|\pi_{t,\omega}(\dd z\times\{Z^j_s\}|E\times \{Z^j_t\}).\]
By exchangeability, we see that this expression is independent of $j$ and since the $\alpha$-law of $Y_s(\omega)$ is $f_s$ for any $(s,\omega)$, this expression is thus equal to
\[\E\left[\int_{(j-1)/N}^{j/N} |Y_s(\omega,\alpha)|\dd\alpha\right] = \frac{1}{N}\int_{E} |z| f_s(\dd s).\]
The conclusion thus follows from Gronwall lemma and Lemma \cref{lemma:boundmomentboltzmannsde}.
\end{proof}

\noindent\textit{\textbf{Step 4.} Coupling bound.}
\medskip

Let us now focus on the estimate of the first term on the right-hand side of \cref{eq:triangleZZbarZtilde}. We write for any $i\leq N$, 
\begin{align*}
    \big|Z^i_t-\widetilde{Z}^i_t\big| &\leq \int_0^t\int_\Theta\int_{\{0,1\}}\int_0^1 \big|a\big(Z^i_{s^-},Z^\mu_{s^-}(\omega,\alpha),\theta,\sigma\big)\\
    &\phantom{\leq \int_0^t\int_\Theta\int_{\{0,1\}}\int_0^1}\qquad-a\big(\widetilde{Z}^i_{s^-},Y_{s}(\omega,\alpha),\theta,\sigma\big)\big|
    \mathcal{N}_i(\dd s,\dd\theta,\dd\sigma,\dd\alpha) \\ 
    &\leq\int_0^t\int_\Theta\int_{\{0,1\}}\int_0^1 \Big\{(1+L(\theta))\big|Z^i_{s^-}-\widetilde{Z}^i_{s^-}\big|\\
    &\phantom{\leq\int_0^t\int_\Theta\int_{\{0,1\}}\int_0^1}\qquad+L(\theta)\big|Z^\mu_{s^-}(\omega,\alpha)-Y_{s}(\omega,\alpha)\big|\Big\}\mathcal{N}_i(\dd s,\dd\theta,\dd\sigma,\dd\alpha).
\end{align*}
Taking the expectation gives a constant $M>0$ such that 
\begin{align*}
    \E\big|Z^i_t-\tilde{Z}^i_t\big|&\leq M \int_0^t\E\left[\big|Z^i_{s^-}-\widetilde{Z}^i_{s^-}\big|+\int_0^1 \big|Z^\mu_{s^-}(\omega,\alpha)-Y_s(\omega,\alpha)\big|\dd \alpha\right]\dd s\\
    &\leq M\int_0^t \E{\left[\big|Z^i_{s^-}-\tilde{Z}^i_{s^-}\big|+W_1\big(\mu_{\mathcal{Z}^N_s},f_s\big)\right]}\dd s\\
    &\leq M\int_0^t \E{\Big[\big|Z^i_{s^-}-\widetilde{Z}^i_{s^-}\big|+\frac{1}{N}\sum_j \big|Z^j_{s^-}-\widetilde{Z}^j_{s^-}\big|+W_1\big(\mu_{\widetilde{\mathcal{Z}}^N_s},f_s\big)\Big]}\dd s
\end{align*}
where the second inequality is actually an equality and comes from the optimal coupling property (Lemma \cref{lemma:boltzmannoptimalcoupling}) and the third inequality follows from the triangle inequality and: 
\[W_1\big(\mu_{\mathcal{Z}^N_s},\mu_{\widetilde{\mathcal{Z}}^N_s}\big)\leq \frac{1}{N}\sum_{j=1}^N \big|Z^j_{s^-}-\widetilde{Z}^j_{s^-}\big|.\]
By classical arguments, we first sum this relation over $i$ and then divide by $N$ to obtain that the process $S_t := \frac{1}{N}\sum_i \E\big|Z^i_t-\widetilde{Z}^i_t\big|$ satisfies: 
\[S_t \leq M\int_0^t \E W_1\big(\mu_{\mathcal{\widetilde{Z}}^N_s},f_s\big)\dd s+2M\int_0^t S_s\dd s,\]
thus by Gronwall lemma and by exchangeability we get: 
\begin{equation}\label[IIeq]{eq:boundZZtilde}\sup_{t\leq T}\E{\big|Z^i_t-\widetilde{Z}^i_t\big|}=\sup_{t\leq T}\E{\left[S_t\right]} \leq \left(M\int_0^T \E W_1\big(\mu_{\mathcal{\widetilde{Z}}^N_s},f_s\big)\dd s\right)\e^{2M T}.\end{equation}

\noindent\textit{\textbf{Step 5.} Conclusion.}
\medskip 

It remains to estimate the quantity 
\[\E W_1\big(\mu_{\mathcal{\tilde{Z}}^N_s},f_s\big).\]
To do so, Murata proved a ``decorrelation lemma'' \cite[Lemma 6.6]{murata_propagation_1977} to directly estimate quantities of the form 
\[\E\Big[\varphi\big(\widetilde{Z}^k_t\big)\varphi\big(\widetilde{Z}^\ell_t\big)-\langle\varphi,f_t\rangle^2\Big],\]
but as noted by \cite{cortez_quantitative_2016}, we can skip these computations using a recent result on exchangeable systems (see \cite[Theorem 1.2]{hauray_kacs_2014} and Theorem \cref{thm:equivalencechaosW1}) which implies the equivalence between the different notions of chaos in Wasserstein-1 distance. Namely it holds that, 
\begin{equation}\label[IIeq]{eq:boundwassersteinztilde}\E W_1\big(\mu_{\mathcal{\widetilde{Z}}^N_s},f_s\big)\leq C\left(W_1\big(\mathrm{Law}\big(\widetilde{Z}^1_s,\widetilde{Z}^{2}_s\big),f_s^{\otimes  2}\big)+\frac{1}{N}\right)^\gamma,\end{equation}
for all $\gamma<(2d+1)^{-1}$ where the constant $C$ depends on the moment of order 1. The right-hand side is controlled by Lemma \cref{lemma:tildeZbarZ} (and the control of the moments). The conclusion thus follows by gathering \cref{eq:boundwassersteinztilde}, \cref{eq:boundZZtilde}, \cref{eq:boundZbarZtilde} and \cref{eq:triangleZZbarZtilde}. 

\end{proof}

We end this section with some additional remarks on the theorem and its proof and a few more bibliographical comments. 
\begin{enumerate}
    \item The same proof works in the case of a non constant but bounded interaction rate, with some Lipschitz conditions. In such case we  do as usual and allow  fictitious collisions. The probability of a fictitious collision can be added in the Poisson random measure. 
    \item Keeping a constant interaction rate, we have never used the fact that $\nu(\dd \theta)$ is a probability measure. The only thing  that we need is that the Lipschitz and growth functions $L(\theta)$ and $M(\theta)$ are integrable. This theoretically allows us to consider the case of non-cutoff particles when $\int_{\Theta}\nu(\dd\theta)=+\infty$. This was one of the original motivations of \cite{murata_propagation_1977} which treats the case of non-cutoff 2D Maxwell molecules. 
    \item One of the advantages of such a coupling technique is that it gives an explicit convergence rate. In our example we use crude Lipschitz and growth estimates which classically lead to a bad behaviour with time. Just as in the McKean-Vlasov case, uniform in time estimates can be obtained for specific models. An example can be found in \cite{cortez_quantitative_2016}. The authors study a ``generalised'' Kac model with linear interactions and various conservation laws (which in particular imply uniform in time control of the moments of the nonlinear law). The same method (together with an additional coupling argument) leads to quantitative uniform in time propagation of chaos for 3D Maxwell molecules (with an optimal rate) in \cite{cortez_quantitative_2018}. This latter work crucially relies on a previous work by Rousset \cite{rousset_n-uniform_2014} dealing with the uniform (in $N$) trend to equilibrium of the particle system, using coupling methods as in this subsection.
    \item Similar techniques and in particular an ``optimal coupling'' argument are also used in \cite{fournier_rate_2016, xu_uniqueness_2018} for a Nanbu system, so without binary collisions. This work illustrates the power of coupling techniques as it treats a much more difficult case than the one treated in this section. The authors managed to treat the case of hard-sphere particles (unbounded cross section) but also hard-potential particles (unbounded cross section and non integrable interaction law). For similar results in the case of binary collisions, see the recent article \cite{heydecker_kacs_2020}. 
    \item Finally, we also point out that the idea of working with an optimal coupling between the empirical measure of a particle system and its limit law also appears in an earlier work \cite{fontbona_measurability_2009} by one of the authors of \cite{cortez_quantitative_2016,cortez_quantitative_2018}. In \cite{fontbona_measurability_2009}, the authors propose a derivation of the Landau equation from a system of interacting diffusions. The stochastic interpretation of the Landau equation is given by a nonlinear SDE (in the sense of McKean) driven by a space-time white noise (instead of a classical Brownian motion in the usual McKean-Vlasov case). The associated particle system is actually better understood as a system of SDEs driven by martingale measures as described in \cite{meleard_systemes_1988}. This setting goes beyond this review and we refer the interested reader to the aforementioned articles for more details. In a sense, with modern eyes, Murata's work \cite{murata_propagation_1977} may look incomplete, essentially because it does not (could not) benefit from the recent development of optimal transport. It should be noted however that the idea of optimal coupling appeared, apparently independently, in two different contexts and several decades apart, in \cite{murata_propagation_1977} and \cite{fontbona_measurability_2009}, respectively for the derivation of the Boltzmann equation in the 70's and for the derivation of the Landau equation about 12 years ago. While coupling methods are nowadays a very important subject in the literature on particle systems, the pioneering (and maybe, in a sense, incomplete) work of Murata seems to have been largely forgotten. 
\end{enumerate}

\subsection{Some pointwise and uniform in time results in unbounded cases via the empirical process}\label[II]{sec:kacprogram}

In this section we gather some of the results obtained in \cite{mischler_new_2015} in two classical unbounded cases: the true Maxwell molecules (i.e. without cutoff) and the hard-sphere molecules, both in the spatially homogeneous setting (see Section \cref{sec:Boltzmannclassicalmodels} and Example \cref{example:spatiallyhomogeneousboltzmann}). These results are obtained via the abstract method developed in~\cite{mischler_kacs_2013,mischler_new_2015} following the seminal (incomplete) work of \cite{grunbaum_propagation_1971}. The general method is described in Section \cref{sec:abstractmischlermouhot} and Theorem \cref{thm:abstractMischler}. It reduces the problem to the careful check of five assumptions which are stated in a simple form in Assumption~\cref{assum:mischlermouhot} but which are extended and stated in a more complex form in \cite[Section 3.1]{mischler_kacs_2013} in order to treat unbounded cases and the uniform in time propagation of chaos. 

\begin{theorem}[\cite{mischler_new_2015}] 
Let $f_0\in \pb_2(\R^d)$ be compactly supported with zero momentum and finite energy:
\[\int_{\R^d} vf_0(\dd v) = 0,\quad \mathcal{E}:= \int_{\R^d} |v|^2f_0(\dd v),\]
and let $f_t$ be the solution at time $t>0$ of the spatially homogeneous version of the Boltzmann equation \cref{eq:Boltzmannphysics_summary} with initial condition $f_0$ and collision kernel $B(u,\sigma)$ given by \cref{eq:collisionkernelB_summary}.
For $N\geq 1$, the Kac sphere (or Boltzmann sphere) is defined by: 
\[\mathcal{S}^N(\mathcal{E}) := {\left\{\mathbf{v}^N\in (\R^d)^N,\,\,\frac{1}{N}\sum_{i=1}^N |v^i|^2=\mathcal{E},\,\,\sum_{i=1}^N v^i=0\right\}}.\]
Fix $T\in(0,+\infty]$. Assume that the initial $N$-particle distribution  $f^N_0$ is either tensorized $f^N_0=f_0^{\otimes N}$ or is $f_0$-chaotic and constrained on the Kac sphere $\mathcal{S}^N(\mathcal{E})$ (see~\cite[Lemma 4.4 and 4.7]{mischler_kacs_2013}).  
\begin{itemize}
    \item \textbf{(Maxwell molecules).} Let $B$ be of the form \cref{eq:truemaxwellmolecules_summary} or \cref{eq:maxwellmoleculescutoff_summary}. Then there exist a subset $\mathcal{F}\subset C_b(\R^d)$ and come constants $C(T)>0$ and $\kappa(d)>0$, which depend respectively only on $T$ and $d$, such that for any $\varphi_k\in \mathcal{F}^{\otimes k}$, $2k\leq N $, it holds that
    \[\sup_{t\leq T} \big|\big\langle f^{k,N}_t-f^{\otimes k}_t, \varphi_k\big\rangle\big|\leq \frac{C(T)k^2}{N^{\kappa(d)}}\|\varphi_k\|_{\mathcal{F}^{\otimes k}}.\]
    Moreover when $T=+\infty$ then $\kappa(d)$ is given by \cite[Theorem 1]{fournier_rate_2015} or \cite[Lemma 4.2 ]{mischler_kacs_2013}. In the cutoff case \cref{eq:maxwellmoleculescutoff_summary}, for any $T<+\infty$, the result holds with the optimal rate $\kappa(d)=\frac{1}{2}$. 
    \item \textbf{(Hard-spheres).} Let $B$ be of the form \cref{eq:hardsphereB_summary}. Then when $T<+\infty$, there exist a subset $\mathcal{F}\subset C_b(\R^d)$ and some constants $C(T)>0$ and $\alpha>0$ such that for any $\varphi_k\in \mathcal{F}^{\otimes k}$, $2k\leq N $, it holds that
    \[\sup_{t\leq T} \big|\big\langle f^{k,N}_t-f^{\otimes k}_t, \varphi_k\big\rangle\big|\leq \frac{C(T)k^2}{(1+|\log N|)^\alpha}\|\varphi_k\|_{\mathcal{F}^{\otimes k}}.\]
    Moreover if $f_0$ is instead assumed to be bounded and to have a bounded exponential moment and if $f^N_0$ is $f_0$-chaotic and constrained on the Kac sphere $\mathcal{S}^N(\mathcal{E})$, then so is the $N$-particle distribution $f^N_t$ for all $t>0$ and the previous estimate holds with $T=+\infty$. 
\end{itemize}
\end{theorem}

The results of this theorem also imply the propagation of finite and infinite dimensional Wasserstein-1 chaos as defined in Definition \cref{def:wassersteinchaos} (see \cite[Theorems~5.2 and~6.2]{mischler_kacs_2013}). The authors also prove the propagation of entropic chaos (Definition~\cref{def:entropychaotic}) for the cutoff Maxwell molecules and the hard-spheres together with the relaxation towards equilibrium with a rate independent of $N$ (see \cite[Theorem 7.1]{mischler_kacs_2013}). These results positively answer many of the conjectures raised by Kac in his seminal article \cite{kac_foundations_1956} (also known as the \emph{Kac's program in kinetic theory} \cite[Section 1.4]{mischler_kacs_2013}). In particular, it provides a \textit{``satisfactory justification of Boltzmann H-theorem''} for unbounded models (which, with a modernised terminology, corresponds to the notion of entropic chaos in the sense of Definition \cref{def:entropychaotic}).  

In the hard-sphere case, the results have recently been improved in a pathwise setting in  \cite{heydecker_pathwise_2019}. The improvement is due to a better H\"older stability result \cite[Theorem 1.6]{heydecker_pathwise_2019} for the nonlinear Boltzmann flow which improves the control of the third term on the right-hand side of \cref{eq:abstracttheoremsplitting} and leads to a polynomial convergence rate (instead of logarithmic).

\subsection{Lanford's theorem for the deterministic hard-sphere system}\label[II]{sec:lanford}

\subsubsection{The hard-sphere system}
This section is the only one which concerns a completely deterministic system called the hard-sphere system. A \emph{hard-sphere} is a spherical particle defined by its position $X^i_t$, its velocity $V^i_t$ and its diameter $\varepsilon>0$. The hard-spheres are simply subject to the free-transport but it is assumed that two hard-spheres cannot overlap. A system of $N$ hard-spheres is thus defined by the deterministic Newton equations for $i\in\{1,\ldots,N\}$, 
\begin{equation}\label[IIeq]{eq:freetransport}\frac{\dd X^i_t}{\dd t} = V^i_t,\quad \frac{\dd V^i_t}{\dd t} = 0,\end{equation}
on the domain: 
\begin{equation}\label[IIeq]{eq:domainhardspheres}\mathcal{D}_N := \big\{\mathbf{z}^N = (x^i,v^i)_{i\in\{1,\ldots,N\}}\in (\R^{d}\times\R^d)^N,\,\,\forall i\ne j,\,|x^i-x^j|\geq \varepsilon\big\}.\end{equation}
On the boundary of $\mathcal{D}_N$, that is when two particles are at a distance $\varepsilon$, the \emph{collision} of two hard-spheres is an elastic collision which preserves energy and momentum. Starting with a pair of pre-collisional velocities $(v^i,v^j)$, writing down the conservation laws leads to the following formula for the post-collisional velocities: 
\begin{equation}\label[IIeq]{eq:postcollisionnu_II}
\begin{array}{ll}
v^{i*}= v^i-\nu^{i,j}\cdot(v^i-v^j)\nu^{i,j}\vspace{0.2cm}\\
v^{j*}= v^j+\nu^{i,j}\cdot(v^i-v^j)\nu^{i,j}
\end{array},
\end{equation}
where $\nu^{i,j}:=(x^i-x^j)/|x^i-x^j|\in\mathbb{S}^{d-1}$. This formula is not the same as \cref{eq:postcollisionsigma_summary} but it can be shown that they are actually equivalent \cite[Chapter 1, Section 4.6]{villani_review_2002}, with a suitable choice of $\sigma$. Pre-collisional means that $(v^i,v^j)$ are such that $(v^i-v^j)\cdot\nu^{i,j}<0$. It can also be checked that the post-collisional velocities satisfy $(v^{i*}-v^{j*})\cdot\nu^{i,j}>0$. Note that this transformation is an involution in the sense that if $(v^i-v^j)\cdot\nu^{i,j}>0$ (that is the $v^i$ and $v^j$ are in a post-collisional configuration), then \cref{eq:postcollisionnu_II} gives the pre-collisional velocities.

For the hard-sphere system, the Liouville equation \cref{eq:liouville_summary} reduces to a simple transport equation 
\[\partial_t f^N_t + \sum_{i=1}^N v^i\cdot \nabla_{x^i}\cdot f^N_t = 0, \]
on the domain \cref{eq:domainhardspheres}. The goal is to prove the propagation of chaos when $N\to+\infty$ and $\varepsilon\equiv\varepsilon_N\to0$ with a suitable scaling. The limit distribution $f_t\equiv f_t(x,v)$ satisfies the Boltzmann equation with hard-sphere cross section, which reads in strong form: 
\begin{equation}\label[IIeq]{eq:boltzmannhardpsheres}\partial_t f_t + c\cdot\nabla_x f_t = \int_{\mathbb{S}^{d-1}\times\mathbb{R}^d} (\nu\cdot (v - v'))_+ \big(f_t(x,v^*)f_t(x,{v'}^*) - f_t(x,v)f_t(x,v')\big)\dd \nu\dd v',\end{equation}
where the $*$ notation denotes the post-collisional velocities \cref{eq:postcollisionnu_II} for the couple $(v',v)$. 

We chose to include the hard-sphere system in this review because of its historical importance. This is also at the same time one of the simplest physical model and one of the most difficult to analyse and less well understood. Rigorous analytical results are available only for short times, way too short to be physically relevant. In fact, the well-posedness of the Boltzmann equation \cref{eq:Boltzmannphysics_summary} is itself a long-standing problem of interest.

The first formal derivation of the Boltzmann equation from a system of interacting particles is due to Grad \cite{grad_principles_1958,grad_asymptotic_1963} in the scaling $N\varepsilon^{d-1}=\mathcal{O}(1)$, nowadays called the Boltzmann-Grad scaling. A few decades later, Lanford \cite{lanford_time_1975} provided the first almost complete proof of the convergence of the BBGKY hierarchy towards the Boltzmann hierarchy and thus propagation of chaos for short times for the hard-sphere system. The extension to particles interacting via short-range potentials was achieved in \cite{king_bbgky_1975}. Lanford's proof has then been improved and completed over the following years, let us cite in particular the classical references  \cite{uchiyama_derivation_1988,cercignani_mathematical_1994}. The most complete and up-to-date reference on the subject is \cite{gallagher_newton_2014} (in both the hard-sphere and short-range potentials cases). Following the seminal idea of Lanford, the very detailed proof is based on a fine analysis of the ``recollision trees'' (see also Section \cref{sec:interactiongraph_summary}). This section presents a quite general and very brief overview of Lanford's theorem and its proof. In addition to the reference article \cite{gallagher_newton_2014}, we also refer the interested reader to the reviews \cite{goncalves_derivation_2016} and \cite{golse_newton_2015}.
 
\subsubsection{The BBGKY and Boltzmann hierarchies} Before stating Lanford's theorem, we recall the notion of BBGKY hierarchy in the specific case of the hard-sphere system. As we shall see, the proof of Lanford's theorem follows roughly the same ideas as the forward point of view of Kac theorem (Theorem \cref{thm:kac}). The notion of Boltzmann hierarchy for the nonlinear limit system will also be needed. We recall the notation 
\[\mathbf{z}^s = (x^1,v^1,\ldots,x^s,v^s),\]
for a generic element of $(\R^d\times\R^d)^s$. 

\begin{definition}[mild BBGKY and Boltzmann hierarchies for hard-spheres] Let $N\in\N$ and $\varepsilon>0$. 
\begin{itemize}
    \item For each $s\in\N$, the domain of the system of $s$ hard-spheres of diameter $\varepsilon>0$ is defined by: 
    \[\mathcal{D}_s := \big\{\mathbf{z}^s\in (\R^d\times\R^d)^s,\,\,\forall i\ne j,\,|x^i-x^j|\geq \varepsilon\big\}.\]
    A set of $N$ functions $(f_t^{s,N})^{}_t\in C(\R_+,L^\infty(\mathcal{D}_s))$, $s\in\{1,\ldots,N\}$, is said to satisfy the (mild) BBGKY hierarchy when it satisfies:
    \[f_t^{s,N}(\mathbf{z}^s)=\mathbf{T}_s(t)f_0^{s,N}(\mathbf{z}^s)+\int_0^t \mathbf{T}_s(t-\tau)\mathcal{C}_{s,s+1}f_\tau^{s+1,N}(\mathbf{z}^s)\dd\tau,\]
    where $\mathbf{T}_s$ is the backward flow associated to the $s$-particle hard-sphere system (i.e. the backward flow generated by the deterministic system \cref{eq:freetransport} on the domain \cref{eq:domainhardspheres} with $s$ particles) and the collision operator $\mathcal{C}_{s,s+1}:L^\infty(\mathcal{D}_{s+1})\to L^\infty(\mathcal{D}_s)$ is defined for a test function $g^{s+1}\in L^\infty(\mathcal{D}_{s+1})$ by: 
    \begin{multline*}\mathcal{C}_{s,s+1}g^{s+1}(\mathbf{z}^s)\\:=(N-s)\varepsilon^{d-1}\sum_{i=1}^{s}\int_{\mathbb{S}^{d-1}\times\R^{d}}\nu\cdot(v^{s+1}-v^i)g^{s+1}(\mathbf{z}^s,x^i+\varepsilon\nu,v^{s+1})\dd\nu \dd v^{s+1}.\end{multline*}
    \item For each $s\in\N$, the following set is the formal limit of $\mathcal{D}_s$ when $\varepsilon\to0$:
    \[\Omega_s := \big\{\mathbf{z}^s\in (\R^d\times\R^{d})^s,\,\,\forall i\ne j,\,x^i\neq x^j\big\}.\]
    An infinite set of functions $(f^{s}_t)^{}_t\in C(\R_+,L^\infty(\mathcal{D}_s))$, indexed by $s\in\N$, is said to satisfy the (mild) Boltzmann hierarchy when it satisfies:
    \[f_t^{s}(\mathbf{z}^s)=\mathbf{S}_s(t)f_0^{s}(\mathbf{z}^s)+\int_0^t \mathbf{S}_s(t-\tau)\mathcal{C}^0_{s,s+1}f_\tau^{s+1}(\mathbf{z}^s)\dd\tau,\]
    where $\mathbf{S}_s$ is the backward free-flow associated to the $s$-particle system (which reduces to the backward flow generated by \cref{eq:freetransport} on the full space) and the collision operator $\mathcal{C}^0_{s,s+1}:L^\infty(\Omega_{s+1})\to L^\infty(\Omega_s)$ is defined for $g^{s+1}\in L^\infty(\Omega_s)$ by:
    \begin{multline*}
    \mathcal{C}^0_{s,s+1}g^{s+1}(\mathbf{z}^s)=\sum_{i=1}^{s}\int_{\mathbb{S}^{d-1}\times\R^{d}}(\nu\cdot(v^{s+1}-v^i))_+\\\times\left[g^{s+1}(\mathbf{z}^{s*},x^i,v^{(s+1)*})-g^{s+1}(\mathbf{z}^s,x^i,v^{s+1})\right]\dd\nu \dd v^{s+1},
    \end{multline*}
    where we recall that the star notation $\mathbf{z}^{s*}$ and $v^{(s+1)*}$ refers to the transformation \cref{eq:postcollisionnu_II} for the velocities (between the $i$-th and $(s+1)$-th coordinates) with angle $\nu$.  
\end{itemize}
\end{definition}
As explained in Section \cref{sec:finitesystems}, the BBGKY hierarchy can be formally derived by taking the marginals of the Liouville equation. It is slightly more technical for the hard-sphere system because of the boundary conditions, see for instance \cite[Section 4.2]{gallagher_newton_2014} and \cite{cercignani_mathematical_1994}. The Boltzmann hierarchy is the formal limit of the BBGKY hierarchy when $N\to+\infty$ in the Boltzmann-Grad limit $N\varepsilon^{d-1}\to1$. The rigorous proof of this limit is the core of Lanford's theorem. At this point, let us point out some hidden technicalities, in particular regarding the well-posedness of the two hierarchies.
\begin{itemize}
    \item A first observation is that the set of \emph{pathological} initial configurations (leading to collisions involving more than two particles or to grazing collisons) is of measure zero \cite[Proposition 4.1.1]{gallagher_newton_2014}.
    \item An unfortunate consequence of the previous observation is that all the functions that we are considering are now defined only almost everywhere. In particular it is not clear whether the collision operators $\mathcal{C}_{s,s+1}$ make sense since they involve integration over a set of zero measure (the sphere). This problem has been addressed (for the first time only) in \cite[Section 5.1]{gallagher_newton_2014}.
    \item The well-posedness of the BBGKY and Boltzmann hierarchies can be shown for short times for initial data which satisfy an energy bound given in  \cite[Theorem 6 and Theorem 7]{gallagher_newton_2014}. The main assumption is an estimate on the initial condition of the form: for almost every $x,v\in\R^d$,
    \[f_0(x,v)\leq \e^{-\mu_0-\beta_0|v|^2}.\]
\end{itemize}

\subsubsection{Lanford's theorem} The following form of Lanford's theorem is the one given in \cite[Theorem 8]{gallagher_newton_2014}.  

\begin{theorem}[Lanford] Let $(f_{0}^{s,N})_{s\leq N}$ and $(f_0^{s})_{s\geq1}$ two initial data which satisfy the well-posedness results \cite[Theorem 6 and Theorem 7]{gallagher_newton_2014} and which are \emph{admissible} in the sense that they are compatible and satisfy for all $s\in\N$ 
\begin{equation}\label[IIeq]{eq:initialcvBoltzmann}f_{0}^{(s),N}:=\int_{\R^{2d(N-s)}}\1_{\mathbf{z}^N\in\mathcal{D}_N} f_0^N(\mathbf{z}^N)\dd z^{s+1}\ldots \dd z^N\longrightarrow f_0^{s}\end{equation}
locally uniformly in $\Omega_s$ as $N\to+\infty$ in the Boltzmann-Grad limit. Let $(f_t^{s,N})$ and $(f_t^{s})$ be the solutions of the BBGKY and Boltzmann hierarchies respectively associated to the initial data $(f_{0}^{s,N})_{s\leq N}$ and $(f_0^{s})_{s\geq1}$. Then there exists a time $T>0$ such that, uniformly in $t\in [0,T]$, the following convergence in the sense of observables holds:
\[\forall s\in\N,\,\forall\varphi_s\in C_c(\R^{ds}),\quad \int_{\R^{ds}} \varphi_s(\mathbf{v}^s)\left(f^{s,N}_t(\mathbf{x}^s,\mathbf{v}^s)-f_t^{s}(\mathbf{x}^s,\mathbf{v}^s)\right)\dd \mathbf{v}^s\to0\]
locally uniformly on $\{\mathbf{x}^s\in(\R^{d})^s,\,\forall i\neq j,\,x^i\neq x^j\}$ as $N\to+\infty$ in the Boltzmann-Grad limit. 
\end{theorem}
Tensorized initial Boltzmann data $(f_0^{\otimes s})_{s\geq1}$ are admissible in the sense that there exists a BBGKY initial data which satisfy \cref{eq:initialcvBoltzmann}. This is a consequence of the Hewitt-Savage theorem (see \cite{gallagher_newton_2014}). In this case the $k$-th marginal $f^{1kN}_t$ of the BBGKY hierarchy converges towards the $k$ tensor product of the (mild) solution of the Boltzmann equation \eqref{eq:boltzmannhardpsheres}. 

We now briefly sketch the main ideas of the proof. Similarly to Kac's theorem (Theorem \cref{thm:kac}), the dominated convergence theorem is used for the iterated BBGKY and Boltzmann hierarchies (see the forward point of view of Kac's theorem). The term-by-term convergence is however way more difficult. Let us fix $s\in\mathbb{N}$. By iterating the definition of mild solution, the $s$-marginal can be written as a finite sum:
\begin{multline*}f_t^{s,N}(\mathbf{z}^s)=\sum_{k=0}^{N-s}\int_0^t \int_0^{t_1}\ldots \int_0^{t_{k-1}}\mathbf{T}_s(t-t_1)\mathcal{C}_{s,s+1}\mathbf{T}_{s+1}(t_1-t_2)\mathcal{C}_{s+1,s+2}\ldots  \\ 
\mathcal{C}_{s+k-1,s+k}\mathbf{T}_{s+k}(t_k)f_{0}^{s+k,N} \dd t_1\ldots  \dd t_k.
\end{multline*}
Compared to the initial formula, this may look more complicated but the key observation is that now, it involves only the initial condition. Of course the sum becomes infinite when $N\to+\infty$. In \cite{gallagher_newton_2014} it is therefore written directly as an infinite series, up to setting $f^{(s)}_{N,0}\equiv0$ for $s>N$. For an observable $\varphi_s$, the quantity to control is therefore: 
\[I_s(t,\mathbf{x}^s):=\sum_{k=0}^\infty I_{s,k}(t,\mathbf{x}^s),\]
where
\begin{multline*}I_{s,k}(t,\mathbf{x}^s):= \int \dd \mathbf{v}^s\varphi_s(\mathbf{v}^s)\int_0^t \int_0^{t_1}\ldots \int_0^{t_{k-1}}\mathbf{T}_s(t-t_1)\mathcal{C}_{s,s+1}\mathbf{T}_{s+1}(t_1-t_2)\mathcal{C}_{s+1,s+2}  \\ \ldots\mathcal{C}_{s+k-1,s+k}\mathbf{T}_{s+k}(t_k)f_{0}^{s+k,N} \dd t_1\ldots  \dd t_k.\end{multline*}
Similarly for the Boltzmann hierarchy, the authors of \cite{gallagher_newton_2014} define: 
\[I^0_s(t,\mathbf{x}^s):=\sum_{k=0}^\infty I^0_{s,k}(t,\mathbf{x}^s),\]
where
\begin{multline*}I^0_{s,k}(t,\mathbf{x}^s):= \int \dd \mathbf{v}^s\varphi_s(\mathbf{v}^s)\int_0^t \int_0^{t_1}\ldots \int_0^{t_{k-1}}\mathbf{S}_s(t-t_1)\mathcal{C}^0_{s,s+1}\mathbf{S}_{s+1}(t_1-t_2)\mathcal{C}^0_{s+1,s+2}  \\ \ldots\mathcal{C}^0_{s+k-1,s+k}\mathbf{S}_{s+k}(t_k)f_{0}^{s+k} \dd t_1\ldots  \dd t_k.\end{multline*}
The strategy is to use the dominated convergence theorem to prove that : 
\[\sum_{k=0}^{\infty}I_{s,k}(t,\mathbf{x}^s)\,\,\underset{N\to+\infty}{\longrightarrow}\,\, \sum_{k=0}^{\infty}I_{s,k}^0(t,\mathbf{x}^s).\]
in the Boltzmann-Grad limit and locally uniformly in $\mathbf{x}^s$.

The domination part is the easiest one (see \cite[Section 5.3]{golse_newton_2015} and \cite[Theorem~6]{gallagher_newton_2014}). The term-by-term convergence is way more technical and is based on the reformulation of the observables in terms of pseudo-trajectories. For typographical reasons, in the following definition, we change our usual convention and we write the time as an argument and not as a subscript: $Z(t)\equiv Z_t$ (it is also a usual convention for deterministic systems). In the following definition, we also use the notion of \emph{interaction tree} which is an interaction graph assumed to be without recollision, as defined in Section \cref{sec:interactiongraph_summary}.

\begin{definition}[Pseudo-trajectory]
Let $s\in\N$, $t>0$, and \[\tilde{Z}^{s,\varepsilon}(t)=\big(\tilde{X}^{1,\varepsilon}(t),\tilde{V}^{1,\varepsilon}(t),\ldots,\tilde{X}^{s,\varepsilon}(t),\tilde{V}^{s,\varepsilon}(t)\big)\in\R^{2ds}.\] Let $k\in\N$ and $\mathcal{G}_{(1,\ldots,s)}(\mathcal{T}_k,\mathcal{R}_k)$ be an interaction tree with $i_{\ell}=s+\ell$, $\ell\geq1$. Given a $k$-tuple of velocities $\mathbf{v}^k=(v^{1},\ldots,v^{k})\in\R^{ds}$ and angles ${\nu}^k=(\nu_{1},\ldots,\nu_{k})\in\mathbb{S}^{d-1}$, the BBGKY pseudo-trajectory $\tilde{Z}^\varepsilon_{s+\ell}(\tau)$ at time $\tau\geq0$, $\ell\geq0$ is defined recursively backward in time by: 
\begin{itemize}
    \item for $\tau\in(t_{\ell},t_{\ell-1}]$, $\tilde{Z}^{s+\ell-1,\varepsilon}(\tau)$ is given by the particle backward flow (with boundary conditions) starting from $\tilde{Z}^{s+\ell-1,\varepsilon}({t_{\ell-1}})$ with the convention $t_0=t$,
    \item at time $t_{\ell}^+$, a particle is adjoined to the system at position 
    $\tilde{X}^{j_{\ell},\varepsilon}({t_\ell^+})+\varepsilon\nu_\ell$ with velocity $v^\ell$, 
    \item the state of the system $\tilde{Z}^{s+\ell,\varepsilon}(t_\ell^-)$ after adjunction of the particle $s+\ell$ depends on whether the velocities $(j_\ell,s+\ell)$ at $t_\ell^+$ are pre- or post-collisional, namely we take: 
    \begin{equation*}
        \left\{
        \begin{array}{rclcl}
        \big(\tilde{V}^{j_\ell,\varepsilon}(t_\ell^-),\tilde{V}^{s+\ell,\varepsilon}(t_\ell^-)\big) &=& \big(\tilde{V}^{j_\ell,\varepsilon}(t_\ell^+),v^\ell\big) & \text{if} & \nu_\ell\cdot\big(\tilde{V}^{j_\ell,\varepsilon}(t_\ell^+)-v^\ell\big)<0 \\
        \big(\tilde{V}^{j_\ell,\varepsilon}(t_\ell^-),\tilde{V}^{s+\ell,\varepsilon}(t_\ell^-)\big) &=& \big(\tilde{V}^{j_\ell,\varepsilon}(t_\ell^+)^*,v^{\ell*}\big) & \text{if} & \nu_\ell\cdot\big(\tilde{V}^{j_\ell,\varepsilon}(t_\ell^+)-v^\ell\big)>0, \\
\end{array}
        \right.
    \end{equation*}
    where $(v^*,w^*)$ denotes the pre-collisional velocities associated to $(v,w)$ after scattering, defined by \cref{eq:postcollisionnu_II}.
\end{itemize}
When $\varepsilon=0$, $\tilde{Z}^{s+\ell,0}(\tau)$ for $\tau\in(t_\ell,t_{\ell-1}]$ is defined similarly by replacing the particle backward flow by the backward free-transport flow. The dynamical system $\tilde{Z}^{s+\ell,0}$ is called the Boltzmann pseudo-trajectory. 
\end{definition}

\begin{figure}
\centering
\begin{tikzpicture}[line cap=round,line join=round,>=latex,x=0.4501558112799655cm,y=0.4501558112799655cm,scale=1.25]
\clip(-5.472426205840941,-5.467638573091305) rectangle (16.74210429542298,8.008474716775025);
\draw [->] (16.02,-4.68) -- (-3.98,-4.68);
\draw [fill=black] (-1.9255710214520754,-4.68) circle (1.5pt);
\draw[color=black] (-1.3558586676806533,-5.0236017875921728) node {$t_1$};
\draw [dotted] (-1.9255710214520754,-4.68) -- (-1.9255710214520754,8.008474716775025);
{\draw[fill=white](-1.92,0.25721570135180016) circle (0.22507790563998275cm);
\draw [->,black] (3.7,-0.0713921493241001) -- (-1.92,0.25721570135180016);}
{\draw[fill=red](-1.92,0.25721570135180016) circle (0.045507790563998275cm);}
{\draw [pattern=dots] (-1.92,1.2831918040663135) circle (0.22507790563998275cm);
\draw [->,black] (3.7,3.74) -- (-1.92,1.2831918040663135);}
{\draw [fill=black] (3.72,-4.68) circle (1.5pt);
\draw[color=black] (4.550433053497094,-5.067566459798399) node {$t_2$};
\draw [dotted] (3.72,-4.68) -- (3.72,8.008474716775025);}
{\draw[fill=white](3.7,-0.0713921493241001) circle (0.22507790563998275cm);}
{\draw[fill=red](3.7,-0.0713921493241001) circle (0.045507790563998275cm);}
{\draw[pattern=dots](3.7,3.74) circle (0.22507790563998275cm);}
{\draw [pattern=horizontal lines,pattern color=black] (3.6790860288746425,4.74782242653284) circle (0.22507790563998275cm);}
{\draw [->,black] (9.32,-0.4) -- (3.7,-0.0713921493241001);}
{\draw [->,black] (9.32,6.0204630622466526) -- (3.679086028874643,4.74782242653284);}
{\draw [->,black] (9.32,0.6) -- (3.7,3.74);}
{
\draw [fill=black] (9.32,-4.68) circle (1.5pt);
\draw[color=black] (9.9298837628705673,-5.1307357408446475) node {$\tau$};
\draw [dotted] (9.32,-4.68) -- (9.32,8.008474716775025);}
{\draw [pattern=horizontal lines,pattern color=black] (9.32,6.0204630622466526) circle (0.22507790563998275cm);
\draw [->] (14.94,7.259523300789418) -- (9.32,6.0204630622466526);}
{\draw [fill=black] (14.94,-4.68) circle (1.5pt);
\draw[color=black] (15.5710339371667594,-5.0962284324721318) node {$t_3$};
\draw [dash pattern=on 1pt off 1pt] (14.94,-4.68) -- (14.94,8.008474716775025);
\draw [pattern=horizontal lines,pattern color=black] (14.94,7.259523300789418) circle (0.22507790563998275cm);}
{\draw[fill=white](9.32,-0.4) circle (0.22507790563998275cm);
\draw[fill=red](9.32,-0.4) circle (0.045507790563998275cm);}
{\draw[pattern=dots](9.32,0.6) circle (0.22507790563998275cm);}
{\draw (6.240445689686171,2.73442958496280824) node[anchor=north west] { \footnotesize{ \bf RECOLLISION}};}
{\draw [->] (14.94,0.9902713921493240994) -- (9.32,0.6);
\draw[pattern=dots](14.94,0.9902713921493240994) circle (0.22507790563998275cm);}
{\draw [->] (14.94,-3.54) -- (9.32,-0.4);
\draw[fill=white](14.94,-3.54) circle (0.22507790563998275cm);}
{
\draw [->,red] (14.94,-0.7286078506759005) -- (9.32,-0.4);
\draw[fill=red](14.94,-0.7286078506759005) circle (0.045507790563998275cm);}
\end{tikzpicture}
\caption{At time $t_1^+$, there is only one particle, at the same position for the BBGKY (in white) and Boltzmann (in red) pseudo-trajectories. At time $t_1^+$, a particle (dotted) is added next to the white particle in a pre-collisional way. At time $t_2^+$, a particle (dashed) is added next to the dotted particle in a post-collisional way. Due to a recollision at time $\tau\in(t_3,t_2)$, the Boltzmann and BBGKY pseudo-trajectories of the white/red particle are no longer close to each other at time $t_3$.}
\label[II]{fig:recollisionlanford}
\end{figure}
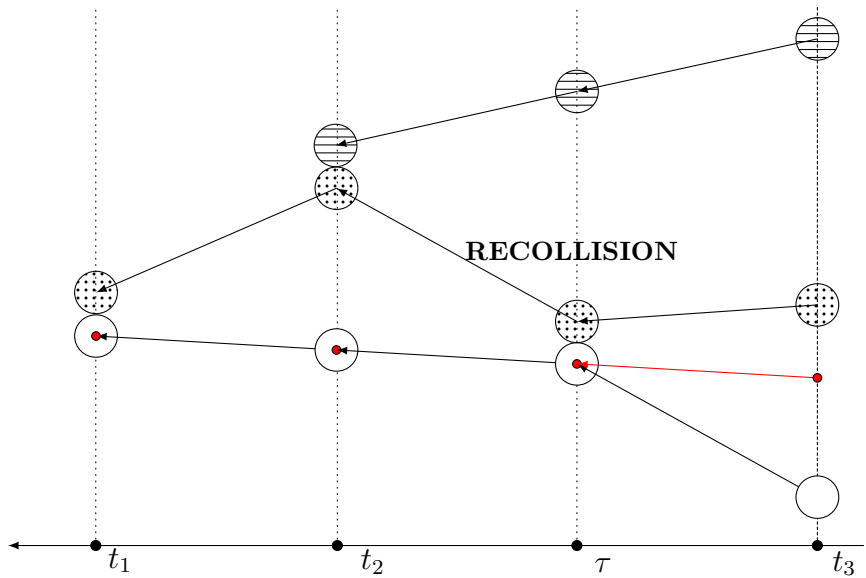

The BBGKY observable can be re-written in terms of pseudo-trajectory as:
\begin{multline}\label[IIeq]{eq:observablepseudotrajectory}I_{s,k}(t,\mathbf{x}_s) = \int \dd \mathbf{v}_s\varphi_s\int_0^t \int_0^{t_1}\ldots\\ \ldots\int_0^{t_{k-1}}\int_{\R^{dk}}\int_{(\mathbb{S}^{d-1})^k} \mathcal{A}\big(\mathcal{T}_k,\mathbf{v}^k,{\nu}^k\big) f^{s+k,N}_{0}\big(\tilde{Z}^{s+k,\varepsilon}(0)\big)\dd\mathcal{T}_k\dd\mathbf{v}^k\dd{\nu}^k,\end{multline}
where
\[\mathcal{A}(\mathcal{T}_k,\mathbf{v}^k,{\nu}^k):=\prod_{\ell=1}^k \nu_\ell\cdot\big(\tilde{V}^{j_\ell,\varepsilon}(t_\ell^+)-v^\ell\big),\quad \dd\mathcal{T}_k = \dd t_1\ldots  \dd t_k\]
and similarly for the Boltzmann observable. In order to take the Boltzmann-Grad limit $N\varepsilon^{d-1}\to1$ in \cref{eq:observablepseudotrajectory} it is necessary to prove that 
\[\tilde{Z}^{s+k,\varepsilon}\longrightarrow \tilde{Z}^{s+k,0},\]
where the pseudo-trajectories are defined taking the same initial condition at time $t$ and giving the same interaction tree and new velocities and deviation angles.
The convergence needs to be strong enough to imply the uniform convergence of the observables. The adjunction of a new particle only gives an error of size $\varepsilon$ (since it is added exactly at the position $\tilde{X}^{j_\ell,0}$ in the Boltzmann case and at a distance $\varepsilon$ in the BBGKY case) which is then transported (backward) in time and can then be controlled. The fundamental difference between the BBGKY and Boltzmann pseudo-trajectories is that BBGKY pseudo-trajectories are subject to recollisions due to the boundary conditions, that is collisions which happen between two times $t_\ell$ and which are not encoded in the collision tree $\mathcal{G}_{(1,\ldots,s)}(\mathcal{T}_k,\mathcal{R}_k)$. Such recollisons do not exist for the Boltzmann pseudo-trajectory since the particles have zero radius.
When a recollision occurs, the situation is illustrated in Figure \cref{fig:recollisionlanford}. 

The fundamental idea of Lanford is to add the particles in such a way that there is no recollision. The proof thus consists in constructing approximate observables by truncating the integration domain in \cref{eq:observablepseudotrajectory}. One of the main contributions of~\cite{gallagher_newton_2014} is an explicit control on the size of the integration domain which leads to recollisions and which is shown to converge to zero fast enough in the Boltzmann-Grad limit. This relies on several geometrical arguments detailed in \cite[Section 12]{gallagher_newton_2014}. We end this section with a few additional bibliographical comments on old and recent problems raised by Lanford's theorem. 

\begin{enumerate}
    \item The case of short-range potentials follows globally the same ideas as the hard-sphere case. On the other hand, the case of long-range interactions is mostly open. A derivation of the linear Boltzmann equation from a system of particles interacting via long-range potentials can be found in \cite{ayi_newtons_2017}. 
    \item Lanford's theorem is valid only for short times. Results on the long-time behaviour are known only for systems close to equilibrium. In the subsequent article \cite{bodineau_brownian_2016}, the authors study a system close to the equilibrium and prove that the linear Boltzmann equation can be obtained as the limit of a system of hard-spheres on a time interval growing to infinity with the number of particles. The proof is based on the same pruning procedure as in Lanford's theorem. The authors also study the motion of a \emph{tagged} particle and show that under the proper scaling, it converges towards a Brownian motion. The striking point is that this derivation starts from a purely deterministic dynamical system (although randomly initialized close to the equilibrium). See also \cite{bodineau_hard_2017} for the derivation of the Stokes-Fourier equations with a similar method. 
    \item The derivation of the Boltzmann equation has long been the source of controversies and somehow metaphysical debates around the question of time reversal and the emergence of irreversibility: the hard-sphere system obeys the Newton laws of motion which are time reversible (at the scale of the whole system) but the Boltzmann equation (which describe the evolution of a single particle) is irreversible (it is a consequence of the famous H-theorem). These quite fundamental questions were addressed in relation with Lanford's theorem already in King's thesis (see for instance the remarks at the end of \cite[Chapter 3]{king_bbgky_1975}). More recent articles also focus on a kind of large-deviation analysis and on the measurement of the \emph{size of chaos} \cite{pulvirenti_boltzmanngrad_2017,bodineau_one-sided_2018,bodineau_statistical_2020}. Note that the work \cite{bodineau_hard_2017} also partly answers the question of the emergence of irreversibility for a system close to equilibrium: despite the deterministic interactions, the tagged particle has a stochastic motion because it inherits, little by little, collision after collision, the (initial) randomness of the other particles. Its (limit) evolution is described on any time interval by the linear Boltzmann equation (because the system is initially close to equilibrium, there is no restriction on the time interval). While the deterministic particle system is reversible when taken in its whole, this is not the case from the point of view of one tagged particle.
\end{enumerate}

\section{Applications and modelling}\label[II]{sec:applications}

In Section \cref{sec:mckeanreview} and Section \cref{sec:boltzmannreview}, we have presented the prototypical application cases of the methods introduced in Section \cref{sec:proving}. In this last section, we go one step further and present a selection of mostly recent applications of these ideas to more concrete problems. Most of the examples presented are not simple direct applications of the previous results. One common and important issue (that we have already discussed) is the difficulty to handle weak regularity. Other topics which will be considered in this section include: time-discrete models which naturally arise in numerical problems, the modelling of noise and the source of stochasticity, the long-time behaviour of particle systems and their behaviour under other scaling limits. The primary objective of this section is to show through various examples how concrete modelling problems lead to these new tough theoretical questions. However, although we hope to give a panorama as faithful as possible of current research, this example-based section is by no means exhaustive and we will mainly stay at an introductory level. One important topic that will not be addressed is the theory of mean-field games. The specific question of propagation of chaos in mean-field games is discussed in great details in \cite{cardaliaguet_master_2019} and its  introduction is itself a quite complete and self-contained review on the subject. Other classical references on mean-field games include \cite{cardaliaguet_notes_2010, carmona_probabilistic_2018, carmona_probabilistic_2018-1}. 

In Section \cref{sec:classicalpdekinetic} we detail the particle interpretation of various classical PDEs in mathematical physics, with a special emphasis on the numerical consequences of these ideas. Section \cref{sec:selforganization} is devoted to a gallery of models of self-organization, mostly inspired by biological systems. Nowadays, mean-field models also have applications in data sciences, either to design more efficient algorithms or to prove their convergence; examples are given in Section \cref{sec:datasciences}. Finally, in the last Section \cref{sec:beyondpoc}, we give a glimpse on some results which go beyond the pure propagation of chaos property.

\subsection{Classical PDEs in kinetic theory: derivation and numerical methods}\label[II]{sec:classicalpdekinetic}

The derivation of classical equations in mathematical kinetic theory is the first \emph{raison d'\^etre} of the propagation of chaos theory. We have already presented some of the main examples: the Fokker-Planck equation, the BGK equation, the granular media equation or the various variants of the Boltzmann equation. In this section, we present further results for the Burgers equation (Subsection \cref{sec:burgers}), the vorticity equation (Subsection \cref{sec:vorticity}) and the Landau equation (Subsection \cref{sec:landau}). In all these cases, the propagation of chaos results derived before cannot be directly applied so we will discuss (without proof) the necessary adaptations.  

Another motivation for this section is the observation that it is usually difficult to numerically solve these kinetic equations using deterministic quadrature methods. Going back to their particle interpretation, the propagation of chaos theory naturally suggests to simulate the underlying particle system and directly use it as a basis for the approximation of the solution of the associated kinetic PDE. In the smooth case, an example is shown in Subsection \cref{sec:bossytalay}. Despite their inherent stochasticity, all the particle systems that have been studied in this review are relatively easy to simulate and modern computers can easily handle from thousands to millions of particles (it remains very far from the $\sim 10^{23}$ order of magnitude in thermodynamics, though). Note however that particle methods may suffer from a high complexity (typically quadratic in $N$), the convergence may be slow (at best $\mathcal{O}(N^{-1/2})$) and the convergence analysis may be difficult. Still, stochastic particle methods have been used with great success in particular for the Boltzmann equation, following the \emph{Direct Simulation Monte Carlo} (DSMC) methods developed by Bird in the sixties (Subsection \cref{sec:dsmc}). 

\subsubsection{Stochastic particle methods for the McKean-Vlasov model}\label[II]{sec:bossytalay}

The propagation of chaos theory for the McKean-Vlasov diffusion with smooth coefficients has been treated in Section \cref{sec:mckeancoupling}. The (quantitative) result Theorem \cref{thm:mckean} readily suggests to approximate the limit Fokker-Planck PDE by a smoothened version of the empirical measure of the particle system. In dimension one, a detailed algorithm and its convergence analysis is due to Bossy and Talay \cite{bossy_stochastic_1997}. The algorithm is based on the Fokker-Planck equation satisfied by the cumulative distribution function $V(t,x) :=~\int_{-\infty}^ x f_t(x)\dd x$, namely with the notations of Theorem \cref{thm:mckean}, 
\begin{multline}\label[IIeq]{eq:mckeancdfpde}
\partial_t V(t,x) = - {\left[\int_\R K_1(x,y)\partial_x V(t,y)\dd y\right]}\partial_x V(t,x) \\+ \frac{1}{2}\partial_x{\left\{{\left[\int_\R K_2(x,y)\partial_x V(t,y)\dd y\right]}^2\partial_x V(t,x)\right\}}.
\end{multline}
Given the a sequence of time steps $t_k = k\Delta t$ for $k\in\{0,\ldots K\}$, the particle system~\cref{eq:mckeanvlasov_summary} is approximated by a first order Euler-Maruyama scheme: 
\begin{equation}\label[IIeq]{eq:discretemckean}Y^i_{t_{k+1}} = Y^i_{t_k} + \frac{1}{N}\sum_{j=1}^N K_1(Y^i_{t_k},Y^k_{t_k})\Delta t + \frac{\sqrt{t_{k+1}-t_k}}{N}\sum_{j=1}^N K_2(Y^i_{t_k},Y^k_{t_k}){\left(G^{i}_{t_{k+1}}-G^i_{t_k}\right)},\end{equation}
where $(G^i_{t_k})_{i,k}$ are $KN$ independent standard Gaussian random variables and $Y^i_0$ are $N$ points in $\R$. Then the solution $V(t_k,\cdot)$ of \cref{eq:mckeancdfpde} is approximated at time $t_k$ by the empirical cumulative distribution function: 
\[V^N_{t_k}(x) = \frac{1}{N}\sum_{i=1}^N H{\left(x-Y^i_{t_k}\right)},\]
where $H$ is the Heaviside function $H(z) := \1_{z\geq0}$. The main results \cite[Theorem 2.1 and Theorem 2.2]{bossy_stochastic_1997} prove the convergence bounds: 
\[\max_{k=0,...,K} \E {\left\|V(t_k,\cdot)-V^N_{t_k}\right\|}_{L^1(\R)} \leq C\left(\|V(0,\cdot)-V^N_0\|_{L^1(\R)} + \frac{1}{\sqrt{N}} + \sqrt{\Delta t}\right),\]
and 
\begin{multline*}\max_{k=0,...,K} \E{\left\|f_{t_k} - \Phi_{\varepsilon}\star \mu_{\mathcal{Y}^N_{t_k}}\right\|}_{L^1(\R)} \\ \leq C{\left[\varepsilon^2 + \frac{1}{\varepsilon}\left(\|V(0,\cdot)-V^N_0\|_{L^1(\R)} + \frac{1}{\sqrt{N}} + \sqrt{\Delta t}\right)\right]},\end{multline*}
where $\Phi_\varepsilon$ is the density of the Gaussian law $\mathcal{N}(0,\varepsilon^2)$, and $C$ depends on $K$ and $\Delta t$. Similarly to Theorem \cref{thm:mckean}, the proof relies on an analogous synchronous coupling for the time discrete system~\cref{eq:discretemckean}, see \cite[Lemma 2.8]{bossy_stochastic_1997}. Still in a regular setting, we also mention that for the granular media equation, the concentration inequality of Theorem \cref{thm:BGVconcentrationineq} leads to an explicit convergence rate for the smoothened empirical measure towards the invariant measure of the nonlinear system, see \cite[Theorem 2.14]{bolley_quantitative_2006}.

For more singular kernels, general results are difficult to obtain and the ``good'' approach most often depends on the specific properties of the considered model. In the next subsections, we give a brief overview of important results for classical equations in mathematical physics. When available, we also discuss their numerical approximation via particle methods. On this last topic, a much more complete reference is the review \cite{bossy_stochastic_2005}. 

\subsubsection{The Burgers equation}\label[II]{sec:burgers}

In his seminal article \cite{mckean_propagation_1969}, McKean raised the problem of the derivation of the Burgers equation from an interacting particle system. The Burgers equation is the following one-dimensional PDE on $\R_+\times\R$:
\begin{equation}\label[IIeq]{eq:burgers}
\partial_t f_t(x) = -f_t(x)\partial_x f_t(x) + \frac{\sigma^2 }{2}\partial_{xx}^2 f_t(x). 
\end{equation}
In view of Theorem \cref{thm:mckean}, the associated particle system should be given by \cref{eq:mckeanvlasov_summary} with 
\[b(x,\mu) = \frac{1}{2}\mu(x),\quad \sigma(x,\mu)\equiv \sigma = \text{constant},\]
or equivalently, with the notations of Theorem \cref{thm:mckean}, 
\[\tilde{b}(x,y) = \frac{1}{2}y,\quad K_1(x,y) = \delta_{x,y}.\]
Clearly, $K_1$ is much too singular to satisfy the hypotheses of Theorem \cref{thm:mckean}. The main approaches to tackle the problem are the following, sorted in chronological order. 
\begin{itemize}
    \item In \cite[Theorem 3.1]{calderoni_propagation_1983}, the Dirac delta $K_1$ is approximated by a smooth function with a smoothing parameter $\varepsilon(N)$ which depends on $N$. Using McKean's quantitative approach of Theorem \cref{thm:mckean}, Calderoni and Pulvirenti show a moderate interaction result (see Section \cref{sec:mckeantowardssingular}) and prove that there exists a sequence $\varepsilon(N)$ for which the propagation of chaos result holds towards the (singular) solution of the Burgers equation. 
    \item In \cite{osada_propagation_1985}, Osada and Kotani use a more analytical approach based on the observation that the generator of the particle system with $K_1(x,y)=\delta_{x,y}$ can be written in divergence form and after a change of time, this generator can be seen as a perturbation of order $N^{-1}$ of an Ornstein-Uhlenbeck generator. The result then follows from a careful analysis of the associated $N$-particle semigroup written as a series expansion via the iterated Duhamel formula.   
    \item In \cite{sznitman_propagation_1986} (see also \cite[Chapter 2]{sznitman_topics_1991}), Sznitman replaces the deterministic drift $\delta_{X^i_t,X^j_t}\dd t$ by the symmetric local time in 0 of $X^i_t-X^j_t$. 
    \item In \cite{bossy_stochastic_1997} and in \cite{jourdain_diffusions_1997}, Bossy, Talay and Jourdain use a different particle system: they interpret the Burgers Equation \cref{eq:burgers} as the equation satisfied by the cumulative distribution function \cref{eq:mckeancdfpde} of the solution of the McKean-Vlasov equation with the kernel $K_1(x,y) = H(x-y)$ where $H$ is the Heaviside function (and a constant diffusion). Still, the kernel does not satisfies the hypothesis of McKean theorem because of the discontinuity in zero. The dedicated propagation of chaos result is proved using the strong pathwise martingale method (similar Theorem \cref{thm:martingalestrongpathwiseempiricalchaos}), see \cite[Theorem 3.2]{bossy_convergence_1996} and the generalised result \cite[Proposition 2.4]{jourdain_diffusions_1997}. In \cite[Theorem 3.1]{bossy_convergence_1996}, Bossy and Talay also prove the convergence of the particle scheme \cref{eq:discretemckean} with the same convergence rate as in the Lipschitz framework. 
    \item In \cite{lacker_strong_2018}, Lacker shows that the Burgers equation can be derived by a direct application of the generalised McKean Theorem \cref{thm:mckeangirsanov} using the Girsanov transform. 
\end{itemize}

\subsubsection{The vorticity equation and other singular kernels}\label[II]{sec:vorticity}

In dimension 2, the vorticity formulation of the Navier-Stokes equation reads
\begin{equation}\label[IIeq]{eq:vorticity}
    \partial_t w_t(x) = (K\star w_t)(x)\cdot \nabla_x w_t(x) + \nu\Delta_x w_t(x),
\end{equation}
where $\nu>0$ is called the viscosity and $K$ is the Biot and Savart kernel 
\begin{equation}\label[IIeq]{eq:biotsavart}K(x) = \frac{x^\perp}{|x|^2} = \frac{1}{|x|^2}(-x_2,x_1).\end{equation}
It is important to note that the solution $w_t$ of \cref{eq:vorticity} is not assumed to be positive so its interpretation as the law of a limit particle system is not obvious. A particle system associated to \cref{eq:vorticity} has been introduced by Chorin in \cite{chorin_numerical_1973} as a simple numerical method to solve the vorticity equation. Later, computational improvements have been proposed in \cite{greengard_fast_1987} to cope with the high complexity of the algorithm, which is quadratic in the number of particles. The idea is based on a clever specific treatment of the short and long range interactions. The method is known as the \emph{fast multipole method}. 

The mathematical treatment of Chorin algorithm and more generally the problem of the derivation of Equation \cref{eq:vorticity} from a particle system was initiated in the 80's and is still an active topic. Important progress have been made very recently. The particle system is described below and the propagation of chaos result is stated (informally). 

Let be given $N$ real-valued random variables $(m_i^N)_{i\in\{1,\ldots,N\}}$ called the \emph{circulations}. The $N$ particle system is a standard linear McKean-Vlasov system where the convolution with the drift kernel $K$ is ``weighted'' by the circulations: 
\[\dd X^i_t = \frac{1}{N}\sum_{j\neq i} m_j^N K(X^j_t-X^i_t)\dd t + \sigma\dd B^i_t,\]
where $\sigma>0$ is such that $\nu=\sigma^2/2$ and the $B^i_t$ are $N$ independent Brownian motions. Note that the circulations do not depend on time. The problem is to prove the propagation of chaos for the system $(X^i_t,m^N_i)_{i\in\{1,\ldots,N\}}$ towards the law of the nonlinear random variable $(\overline{X}_t, \overline{m})$ defined by the SDE:
\[\dd \overline{X}_t = K\star w_t(\overline{X}_t)\dd t + \sigma\dd B_t,\]
where $B_t$ is a Brownian motion and $w_t$ is the measure on $\R^2$ defined by
\[\forall \mathscr{B}\in\mathcal{B}(\R^2),\quad w_t(\mathscr{B}) = \E\big[\overline{m}\1_{\overline{X}_t \in \mathscr{B}}\big] = \iint_{\R\times\mathscr{B}} m f_t(\dd m,\dd x),\]
where $f_t(\dd m, \dd x)\in \pb(\R\times\R^2)$ is the law of $(\overline{m},\overline{X}_t)$. It can be shown that $w_t$ is a weak solution of the vorticity Equation \cref{eq:vorticity}. Moreover, the propagation of chaos result implies that the random measure 
\[ W^N_t := \frac{1}{N}\sum_{i=1}^N m^N_t \delta_{X^i_t} \in \mathcal{M}(\R^d),\]
converges weakly towards $w_t$ (in the space of measures). 

Compared to the classical setting, the main difficulty is that the Biot and Savart kernel \cref{eq:biotsavart} is singular in $(0,0)$. To deal with the problem, some of the historical stepping stones include the following. In \cite{marchioro_hydrodynamics_1982}, Marchioro and Pulvirenti use a regularised kernel and prove a moderate interaction result by coupling the trajectories. The result is improved by M\'el\'eard in \cite{meleard_trajectorial_2000, meleard_monte-carlo_2001} who proves the pathwise propagation of chaos with more general initial data. A different approach, without regularisation, is due to Osada in \cite{osada_propagation_1986}. Similarly to \cite{osada_propagation_1985}, it is based on the analytical study of the generator of the particle system. The propagation of chaos result holds under the assumptions of a large viscosity and a bounded initial data. Two recent works have improved these results using two different approaches. 
\begin{itemize}
    \item In \cite[Theorem 2.12]{fournier_propagation_2014}, Fournier, Hauray and Mischler use a martingale compactness method (similar to the one presented in Section \cref{sec:martingalecompactness}). To cope with the singularity of the Biot and Savart kernel, new entropy estimates are derived which, among other things, imply that any two particles do not stay too close to each other, see \cite[Lemma 3.3]{fournier_propagation_2014}. We also refer to the introduction of the article which contains a much more complete review of the existing works on the subject, also including deterministic models. 
    \item In \cite{jabin_quantitative_2018}, Jabin and Wang have analytically derived entropy bounds which imply the propagation of chaos in TV norm for the system with constant circulations $m^N_i=1$. Compared to the previous approach, the result is quantitative. The strategy is reviewed in Section \cref{sec:jabin}. 
\end{itemize}

To conclude, we discuss some extensions of these ideas to other important singular kernels derived from a Coulomb potential and already mentioned in Section \cref{sec:mckeantowardssingular}, namely when $K(x) = \xi_d x/|x|^d$ in dimension $d$, for a constant $\xi_d\in\R$. In dimension $d=2$, the attractive case $\xi_d>0$ is called the Keller-Segel kernel. Following the probabilistic methods of \cite{fournier_propagation_2014}, recent works on the corresponding particle model include \cite{fournier_stochastic_2017, godinho_propagation_2015, liu_propagation_2016, liu_propagation_2019}, see also the references therein. More analytical methods include \cite{haskovec_convergence_2011, berman_propagation_2019, bresch_mean-field_2019} and the references therein. 

Another natural extension concerns kinetic systems on the product space $\R^d\times\R^d$ for which the particle system is defined by the Newton equations with random noise: 
\[\dd X^i_t = V^i_t\dd t,\quad \dd V^i_t = K\star \mu_{\mathcal{X}^N_t}(X^i_t)\dd t + \sigma \dd B^i_t,\]
where $K(x) = \xi_d x/|x|^d$ on $\R^d$. The limit equation obtained as $N\to+\infty$ is called the Vlasov-Poisson-Fokker-Planck equation on $\R^d\times\R^d$: 
\[\partial_t f_t(x,v) + v\cdot \nabla_x f_t + (K\star \rho_t)(x)\cdot \nabla_v f_t(x,v) = \frac{\sigma^2}{2}\Delta_v f_t(x,v),\]
where $\rho_t(\dd x) = \int_{v\in\R^d} f_t(\dd x,\dd v)$. The propagation of chaos result via entropy bounds is proved in \cite{jabin_mean_2016}. More recently, the Vlasov-Poisson-Fokker-Planck is derived in \cite{carrillo_propagation_2019} from a regularised particle system with cutoff. Following this work, a better cutoff size is obtained in \cite{huang_mean-field_2020}. See also the references therein for a more detailed account of earlier works on the subject.

\subsubsection{The Landau equation}\label[II]{sec:landau}

The Landau operator is obtained as a \emph{grazing collision limit} \cite{toscani_grazing_1998, villani_review_2002} of the Boltzmann operator defined on the right-hand side of \cref{eq:Boltzmannphysics_summary} when the angular cross-section $\Sigma = \Sigma^\varepsilon(\theta)$ depends on a parameter $\varepsilon\to0$ such that 
\[\Sigma^\varepsilon(\theta)\underset{\varepsilon\to0}{\longrightarrow}0,\]
uniformly on any interval $[\theta_0,\pi]$, $\theta_0>0$ and 
\[\int_0^\pi \sin^2(\theta/2)\Sigma^{\varepsilon}(\theta)\dd\theta\underset{\varepsilon\to0}{\longrightarrow} \Lambda \in (0,+\infty).\]
In the spatially homogeneous case, it leads to the Landau equation: 
\begin{equation}\label[IIeq]{eq:landau}
    \partial_t f_t(v) = \nabla_v\cdot \int_{\R^3} a(v-v_*) \big[f_t(v_*)\nabla_vf_t(v)-f_t(v)\nabla_vf_t(v_*)\big]\dd v_*,
\end{equation}
where for $u\in\R^3$, we define the matrix $a(u) = \widetilde{\Phi}(|u|)\mathsf{P}(u)$, where $\widetilde{\Phi}$ is explicit in terms of the velocity cross section $\Phi$ and the constant $\Lambda$, and $\mathsf{P}(u) := I_3 - \frac{u\otimes u}{|u|^2}$ is the orthogonal projection matrix on $u^\perp$. As usual, we refer to the classical reviews \cite{villani_review_2002} and \cite{bellomo_macroscopic_2004} for a more complete overview of the Landau equation in kinetic theory. 

The formal derivation of \cref{eq:landau} from the Boltzmann equation has been made rigorous in a probabilistic framework in \cite{guerin_convergence_2003}. Similarly to the Boltzmann equation, the Landau equation is shown to be associated to a nonlinear martingale problem. Then using the strong pathwise martingale compactness method (see Section \cref{sec:provingcompactness_summary}), it is obtained as the limit of a Boltzmann particle system when both parameters $N$ and $\varepsilon$ converge, $N\to+\infty$ and $\varepsilon\to0$, see \cite[Theorem 4.1]{guerin_convergence_2003}. This procedure also gives a Monte Carlo algorithm for the approximation of the Landau equation, using a Bird simulation algorithm which will be discussed below. 

The Landau equation can also be obtained as the $N\to+\infty$ limit of a system of $N$ diffusion processes. In \cite{fontbona_measurability_2009}, the authors consider the particle system driven by $N^2$ independent Brownian motions $(B^{ij}_t)^{}_{i,j}$ and defined by the system of SDEs: 
\begin{equation}\label[IIeq]{eq:particlelandau}\dd X^i_t = \frac{1}{N}\sum_{j=1}^N b(X^i_t-X^j_t)\dd t + \frac{1}{\sqrt{N}}\sum_{j=1}^N \sigma(X^i_t-X^j_t)\dd B^{ij}_t,\end{equation}
for $i\in\{1,\ldots,N\}$, where $b:\R^3\to\R^3$ and $\sigma:\R^3\to\mathcal{M}_3(\R)$ are defined by 
\[b(u) = \nabla\cdot a(u), \quad \sigma(u)\sigma(u)^\mathrm{T} = a(u),\] 
for the matrix $a$ in \cref{eq:landau}. The nonlinear limit SDE is not a classical McKean-Vlasov system as it is driven by a space-time white noise instead of a standard Brownian motion. To prove the propagation of chaos, Fontbona, Gu\'erin and M\'el\'eard have developed a dedicated optimal coupling method which leads to a quantitative convergence estimate. A non quantitative result was obtained before in \cite{meleard_systemes_1988} using martingale methods and martingale measures (see Remark \cref{rem:martingalemeasure}). 

A more standard McKean-Vlasov system of the form \cref{eq:mckeanvlasov_summary} is given in \cite{fournier_particle_2009} with still $b=\nabla\cdot a$ and this time, for $x\in\R^3$ and $\mu\in\pb(\R^3)$, $\sigma(x,\mu)$ is the unique square root of the matrix: 
\[a\star \mu (x) = \int_{\R^3} a(x-y)\mu(\dd y).\] 
When the matrix $a$ is sufficiently smooth (roughly when $\widetilde{\Phi}(|u|) = |u|^2$), the propagation of chaos result thus follows from a standard synchronous coupling method, using some ad hoc preliminary estimates. A numerical scheme and its convergence analysis is also presented. When the velocity cross section is such that $\widetilde{\Phi}(|u|) = |u|^{2+\gamma}$ with $\gamma\in(-2,0)$ (this is the case of the moderately soft potentials), then $a$ (and thus $\sigma$) is not smooth. The propagation of chaos is proved in \cite{fournier_propagation_2016} using a new optimal coupling method for diffusion processes. The coupling is based on the observation that in order to couple the solutions of the two diffusion SDEs: 
\[\dd X^1_t = \sigma_1(t) \dd B^1_t,\quad \dd X^2_t = \sigma_2(t)\dd B^2_t,\]
for two different diffusion matrices $\sigma_1$ and $\sigma_2$, then an \emph{optimal} choice is:
\[B^2_t = \int_0^t U(a_1(s),a_2(s)) \dd B^1_s,\] 
where $a_k(s) = \sigma_k(s)\sigma_k(s)^{\mathrm{T}}$ for $k\in\{1,2\}$ and
\[U(a_1,a_2) := a_2^{-1/2}a_1^{-1/2}(a_1^{1/2}a_2a_1^{1/2})^{1/2}.\]
This coupling is a ``dynamical'' version of the optimal coupling between the normal laws $\mathcal{N}(0,a_1)$ and $\mathcal{N}(0,a_2)$ \cite{givens_class_1984}. Since the matrix $U(a_1,a_2)$ is orthogonal, it follows that $U(a_1,a_2) B^2_t$ is a standard Brownian motion. A non quantitative result using martingale methods is also shown in \cite{fournier_propagation_2016}. 

Yet another particle system has been proposed in \cite{carrapatoso_propagation_2015}. The SDE system \cref{eq:particlelandau} is the same as in \cite{fontbona_measurability_2009} but the Brownian motions are not independent, they satisfy $B^{ij}_t=-B^{ji}_t$. Contrary to the previous one, this particle system preserves the momentum and energy. The propagation of chaos result is proved using the pointwise empirical approach described in Section \cref{sec:abstractmischlermouhot_summary}. The same particle system is also studied in \cite{fournier_kac-like_2017} where the authors use the optimal coupling method of \cite{fournier_propagation_2016} for the case $\gamma \in [0,1]$. 

\subsubsection{DSMC for the Boltzmann equation}\label[II]{sec:dsmc}

As already discussed many times in this review, the propagation of chaos towards the Boltzmann equation of rarefied gas dynamics \cref{eq:Boltzmannphysics_summary} is a long standing problem which becomes extremely difficult in the unbounded cases \cref{eq:hardsphereB_summary} and \cref{eq:truemaxwellmolecules_summary} described in Section \cref{sec:boltzmann_summary}. In the easiest case of the Maxwell molecules with cutoff, the propagation of chaos follows for instance from Kac's Theorem \cref{thm:kac} but in fact, from any of the methods described in Section \cref{sec:boltzmannreview}. In the unbounded cases, most results are known only in the spatially homogeneous setting. For the hard-sphere cross section and for the (true) Maxwell molecules, the first complete and rigorous results are due to Sznitman \cite{sznitman_equations_1984} using martingale methods (see Theorem \cref{thm:sznitmanboltzmann}) and Murata \cite{murata_propagation_1977} using a coupling approach (in dimension two). For the true Maxwell molecules, in a series of papers \cite{desvillettes_probabilistic_1999}, \cite{fournier_markov_2001,fournier_monte-carlo_2001,fournier_monte-carlo_2001-1}, \cite{fournier_stochastic_2002}, the authors combined Tanaka's probabilitic representation of the Boltzmann equation with the martingale method of Sznitman and obtained existence results and particle approximation results respectively in dimension one, two and three. The strategy is based on a cutoff approximation with a vanishing cutoff parameter when $N\to+\infty$. Lately, the analytical approach developed in \cite{mischler_kacs_2013} has lead to quantitative results in both the hard sphere and true Maxwell molecules cases. The latest results on the subject are summarised in Section \cref{sec:kacprogram} and we refer the interested reader to the introduction of \cite{mischler_kacs_2013} for a more detailed review of known results. Probabilistic coupling methods have also been recently been developed to treat the case of Maxwell molecules and hard and moderately soft potentials, see \cite{fournier_rate_2016, xu_uniqueness_2018} and the references therein.  

The numerical treatment of the Boltzmann equation is also an old question, with many techniques available. The difficulty comes from the approximation of the collision integral, on the right-hand side of \cref{eq:Boltzmannphysics_summary}. In dimension three, this is an integral over a (3+2)-dimensional space which makes any naive deterministic quadrature method practically inefficient. We will not discuss how efficient deterministic methods could be implemented (it is an active research area, see for instance \cite{mouhot_fast_2006, dimarco_numerical_2014, pareschi_stability_2021} and the references therein) and we will focus on a brief overview of stochastic particle methods which can be seen as a natural application of the propagation of chaos property. 

\subsubsection*{An exact simulation algorithm.}

Since the Boltzmann equation is obtained as the limit of a system of particles interacting according to \cref{eq:boltzmanngenerator_summary}, a natural idea is to simulate this particle system on any time interval $[0,T]$ and to take (a possibly smoothened version of) its empirical measure as an approximation of the solution of the Boltzmann equation. It is important to note that it only makes sense to simulate cutoff mollified models so for physical cases of interest (hard spheres or Maxwell molecules), one also needs to introduce a cutoff approximation of the cross section. An advantage of this method is that the particle system can be simulated \emph{exactly} so the only errors comes from the $N$-particle discretization (that is, the convergence rate in the propagation of chaos) and the cutoff approximation. This method is called the \emph{Direct Simulation Monte Carlo} (DSMC) method and has been developed in the 60's by Bird \cite{bird_direct_1970}. We will discuss below Bird's algorithm but before that, we give the algorithmic form of Proposition \cref{prop:acceptreject} and Example \cref{example:semiparametric}: the Algorithm~\cref{algo:exact} below simulates a particle system in the semi-parametric cutoff case defined in Definition \cref{def:Boltzmannparammodel_summary}. With the notations of the definition, the interaction rate $\lambda$ is assumed to be bounded by a constant $\Lambda$ \cref{eq:uniformboundlambda_summary} and the semi-parametric post-collisional distribution $q$ (defined by \cref{eq:semiparametric_summary}) is bounded up to a factor $M$ by a distribution $q_0$ \cref{eq:semiparambound_summary}. 

\begin{center}
\begin{minipage}{.9\linewidth}
\begin{algorithm}[H]
 Set $t=0$ \;
 Draw the initial states $Z^1_0,\ldots, Z^N_0$ \;
 \While{$t\leq T$}{
  Draw $\tau$ from an exponential law with parameter $\Lambda M(N-1)/2$ \;
  Update each particle in $(Z^1_t,\ldots, Z^N_t)$ on $[t,t+\tau]$ according to $L^{(1)}$ \;
  Draw $(i,j)\in \{1,\ldots,N\}^2$ uniformly among the $N(N-1)/2$ pairs \;
  Draw $\theta\sim q_0(\theta)\nu(\dd\theta)$ \;
  Draw $\eta\in[0,1]$ uniformly \;
  \If{$\eta \leq  \frac{\lambda(Z^i_{t+\tau},Z^j_{t+\tau})q(Z^i_{t+\tau},Z^j_{t+\tau},\theta)}{\Lambda M q_0(\theta)}$}{
   $Z^i_{t+\tau} \gets \psi_1(Z^i_{t+\tau},Z^j_{t+\tau},\theta)$ \;
   $Z^j_{t+\tau} \gets \psi_2(Z^i_{t+\tau},Z^j_{t+\tau},\theta)$ \;
   }
  $t \gets t+\tau$  \;
 }
 \caption{Exact simulation}\label[II]{algo:exact}
\end{algorithm}
\end{minipage}
\end{center}

This algorithm and some variants can be found in \cite{graham_stochastic_1997}, \cite{fournier_monte-carlo_2001-1, fournier_monte-carlo_2001}, \cite{meleard_stochastic_1998} or in \cite{guerin_convergence_2003} for an application to the Landau equation. Note that if a process with generator $L^{(1)}$ can be simulated exactly (for instance if it is the generator of a transport operator) or in the spatially homogeneous case, then it is not necessary to discretize time and the output of Algorithm \cref{algo:exact} is exact: the generator of the particle system is \cref{eq:boltzmanngenerator_summary}. However, it is necessary to simulate a Poisson process with a parameter which is $\mathcal{O}(N)$; the accumulation of jumps on small time intervals may become difficult to handle when $N$ is very large. 

In the kinetic non spatially homogeneous case, the state space is $E=\R^d\times\R^d$ with the following assumptions:
\begin{itemize}
    \item the operator $L^{(1)}$ acts only on the space variable and includes the boundary conditions;
    \item the interactions are purely local: for $x,v,x_*,v_*\in\R^d$,
    \[\lambda((x,v),(x_*,v_*)) \equiv \lambda(v,v_*)\delta_{x,x_*} ;\]
    \item the post-collisional distribution depends only on the velocity variable: for $x,v,x_*,v_*\in\R^d$,
    \[q((x,v),(x_*,v_*), \theta) \equiv q(v,v_*,\theta).\]
\end{itemize}
As explained in Example \cref{example:mollifiedmodels}, the only way to treat this very singular case (due to the local interaction) is to consider a \emph{mollified} model. For instance, M\'el\'eard \cite{meleard_monte-carlo_2001} considers a bounded spatial domain which is divided into a finite number of cells of equal volume $\delta^d$, $\delta>0$. In the case of a torus $\mathbb{T}^d$ a possible choice is simply to consider a uniform spatial grid. Then the following mollified collision rate is considered:
\[\lambda^\delta ((x,v),(x_*,v_*)) \equiv \lambda(v,v_*) I_\delta(x,x_*),\]
where $I_\delta$ is the sum over the cells $G$ : 
\[I^\delta(x,y) := \frac{1}{\delta^d}\sum_{G} \1_{x,y \in G}.\]
It physically means that two particles are allowed to interact only when they are in the same cell. The scaling ensures that the purely local Boltzmann equation is (formally) recovered when $\delta\to0$ in the limit mollified Boltzmann equation. The rigorous proof of the propagation of chaos property for this model when $\delta\to0$ and $N\to+\infty$ can be found in \cite{meleard_stochastic_1998}. Since the simulation is exact, the propagation of chaos is also a convergence proof of Algorithm \cref{algo:exact}. Since the algorithm is by nature sequential in time (the collisions are treated sequentially one by one), a drawback of this method is that most of the collisions will be fictitious: the if-loop will almost never be entered into. This comes from the fact that the accept-reject scheme is as efficient as the bound on $\lambda$ and $q$ are small. 

\subsubsection*{The Bird algorithm.}

In the 60's, Bird \cite{bird_direct_1970} introduced a simulation algorithm of the Boltzmann equation of rarefied gas dynamics \cref{eq:Boltzmannphysics_summary} which can be understood as a time-discrete version of Algorithm \cref{algo:exact} with parallelized collisions over the cells. First, the time interval $[0,T]$ is discretized uniformly with a time step $\Delta t$ and the goal is to construct a time discrete approximation of the particles at the times $t_k = k\Delta t$ for $k\in\{0,\ldots,K\}$, $K\in\N$. A short heuristic description of the algorithm is the following. 
\begin{enumerate}
    \item The flow of $L^{(1)}$ and the boundary conditions are treated separately from the collision process. At each time step $t_k$, the positions are updated first and the positions at time $t_{k+1}$ are used to update the velocities from $t_k$ to $t_{k+1}$. 
    \item At each time step, each cell is treated independently: formally, it is equivalent to solve the spatially homogeneous problem in each cell during the time step $\Delta t$. 
    \item Instead of computing an exact simulation based on a Poisson process, a time counter is attached to each cell. Collision events are proposed and each time a collision is accepted, the time counter is incremented by a fixed time which is computed from the theoretical average time between two collisions. If $N_G$ denotes the number of particles in the cell $G$, then the parameter of the Poisson process which gives the (inverse of the) average time between two collisions in $G$ is bounded by:
    \[\frac{\Lambda}{N}\sum_{x^i,x^j\in G} I_\delta(x^i,x^j) = \frac{N_G(N_G-1)}{2}\frac{\Lambda}{N}\frac{1}{\delta^d}.\]
    Since the collision probability depends on the current state of the particles (pairs of particles do not collide with the same probability), the previous bound is used in an accept-reject scheme and for the computation of the time counter. Note that this method does not necessitate to compute the jump probabilities which is an expensive $\mathcal{O}(N_G^2)$ operation. Note also that it is possible to re-compute better bounds $\Lambda$ and $M$ at each iteration: a global bound is not necessary and the product $\Lambda M$ can be replaced by a bound over the $N(N-1)/2$ quantities $\sup_\theta \lambda(V^i_{k},V^j_{k})q(V^i_{k},V^j_{k},\theta)/q_0(\theta)$.
\end{enumerate}

\begin{center}
\begin{minipage}{.9\linewidth}
\begin{algorithm}[H]
  Draw the initial states $Z^1_0,\ldots, Z^N_0$ \;
  \For{$k=0$ \KwTo $K$}{
  Update each position $X^i_{t_k}$ according to $L^{(1)}$ until $t_{k+1}$ \;
  Set $V^i_k = V^i_{t_k}$ for $i\in\{1,\ldots,N\}$ \;
  Decompose the domain into disjoint equal cells of volume $\delta^d$ \; 
  \For{each cell $G$}{
    Set $t_c = t_k$ \;
    \While{$t_c\leq t_{k+1}$}{
        Set $N_G$ the number of particles in the cell $G$ \; 
        Draw uniformly two particles $V^i_{k}$ and $V^j_{k}$ in the cell $G$ \;
        Draw $\theta\sim q_0(\theta)\nu(\dd\theta)$ \;
        Draw $\eta\in[0,1]$ uniformly \;
        \If{$\eta \leq \frac{\lambda(V^i_{k},V^j_{k})q(V^i_{k},V^j_{k},\theta)}{\Lambda M q_0(\theta)}$}{
            Set $\Delta t_{ij} = {\left(\frac{N_G(N_G-1)}{2}\frac{\lambda(V^i_{k},V^j_{k})}{N}\frac{1}{\delta^d}\right)}^{-1}$ \;
            $V^i_{k} \gets \psi_1(V^i_{k},V^j_{k},\theta)$ \;
            $V^j_{k} \gets \psi_2(V^i_{k},V^j_{k},\theta)$ \;
            $t_c \gets t_c+\Delta t_{ij}$ \;
        }
    }
  }
  Set $V^i_{t_{k+1}} = V^i_{k}$ for $i\in\{1,\ldots,N\}$ \;
 }
 \caption{Bird algorithm}\label[II]{algo:bird}
\end{algorithm}
\end{minipage}
\end{center}

The convergence proof of the Bird algorithm is due to Wagner \cite[Theorem 4.1]{wagner_convergence_1992} using (non quantitative) martingale methods. The main result of Wagner is a propagation of chaos result via the empirical measure: Wagner proves that if the empirical measure of the initial state converges then this also holds true for the empirical measure of the output of the Bird algorithm at any later time (note that the algorithm actually defines a time continuous Markov process). The (heuristic) relationship between the limit of the Bird algorithm and the Boltzmann equation is explained in \cite[Section 5]{wagner_convergence_1992}. Algorithm \cref{algo:bird} is referred as the ``modified Bird algorithm with fictitious collisions''.  

\begin{remark}
    This method simulates the Boltzmann equation in \emph{weak form} (since it is based on the simulation of the post-collisional distribution). For the main application case \cref{eq:Boltzmannphysics_summary}, which is written in strong form, there is nothing else to do thanks to the invariance of the collision kernel by the pre- and post-collisional changes of variables, see Example \cref{example:spatiallyhomogeneousboltzmann}. 
\end{remark}

After Bird, Nanbu \cite{nanbu_direct_1980} proposed an algorithm which is roughly speaking a time discretization of the mean-field jump model described at the beginning of Section \cref{sec:meanfieldjump_summary}. At each time step, each particle updates its velocity by choosing a ``collision partner'' which does not update its state during this collision. As before, the collision is accepted or rejected with a probability which depends on the collision rate. The relationship with the Boltzmann equation is shown in Example \cref{example:nanbuboltzmann}. A drawback is that in the physical case of the Boltzmann equation of rarefied gas dynamics \cref{eq:Boltzmannphysics_summary} the algorithm does not preserve the energy and momentum. Another version was thus proposed by Babovsky \cite{babovsky_simulation_1986}: at each time step, the $N$-particle system is randomly uniformly separated into two groups of equal size from which we obtain $N/2$ randomly uniformly sampled collision pairs. Similarly to the Bird algorithm, if a collision is accepted, the two particles update their states. The main difference with Bird algorithm is that each particle can collide at most once per time step. This has a strong influence on the time accuracy. The convergence analysis of the Nanbu-Babovsky algorithm can be found in \cite{babovsky_convergence_1989}. A detailed review and comparison of the Bird and Nanbu-Babovsky algorithms can be found in \cite[Chapter 10]{cercignani_mathematical_1994} as well as several variants. We also refer the interested reader to the lecture notes~\cite{pareschi_introduction_2001}.

\subsection{Models of self-organization}\label[II]{sec:selforganization}

So far, we have been quite vague about what the particles represent. In this section, we present more concrete modelling problems which further motivate the study of particle systems. In the following examples, particles will be used to model large animal societies (Section \cref{sec:flocking}), neuronal networks (Section \cref{sec:neurons}) and socio-economic agents (Section \cref{sec:socioeconomicmodels}). Similarly to Statistical Physics models, the common feature of all these systems is the spontaneous emergence of a large scale complex global dynamics out of the simple and seemingly unorganized motion of many indistinguishable particles. The detailed study of such behaviour is not the primary interest of this review and the following will focus on the first step of the analysis which is the derivation of PDE models which can serve as a theoretical basis to explain \emph{self-organized} phenomena. In order to illustrate the potential complexity of this approach even for seemingly simple models, the next Section \cref{sec:kuramoto} is devoted to a brief overview of recent results on the famous Kuramoto model.  

There is a vast and growing literature on self-organization and collective dynamics models. Further much more detailed examples can be found in the books and review articles \cite{bellomo_active_2017, bellomo_active_2019, naldi_mathematical_2010, muntean_collective_2014, degond_mathematical_2018, albi_vehicular_2019, vicsek_collective_2012}. 

\subsubsection{Phase transitions and long-time behaviour: the example of the Kuramoto model}\label[II]{sec:kuramoto}

The Kuramoto model is the most classical model for synchronization phenomena between populations of \emph{oscillators}, which may be used to model a clapping crowd, a population of fireflies or a system of neurons to cite a few examples. Despite its formal simplicity, the Kuramoto model exhibits a complex long-time behaviour which has motivated a vast literature, see for instance the reviews \cite{acebron_kuramoto_2005, goncalves_large_2015} or the articles \cite{bertini_dynamical_2009, bertini_synchronization_2014} and the references therein. This section is focused on two recent works \cite{bertini_synchronization_2014, delgadino_diffusive-mean_2021} which prove, among other things, that the propagation of chaos does not always hold uniformly in time for the Kuramoto model and some of its variants. The main reason is a phase transition phenomenon. Both works actually prove some kind of large deviation results. Earlier results in this direction can be found in \cite{dai_pra_mckean-vlasov_1996, dawson_critical_1983, dawson_large_1987}.

Let $N$ \emph{oscillators} be defined by $N$ angles $\theta^i_t\in\R$ (defined modulo $2\pi$ so that they can actually be seen as elements of the circle) which satisfy the following McKean-Vlasov SDE:
\[\dd\theta^i_t = \xi_i\dd t -\frac{K}{N}\sum_{j=1}^N \sin(\theta^i_t-\theta^j_t)\dd t  + \dd B^i_t,\]
where $K\in\R$ is a real parameter of the model and $(\xi_i)_{i\in\{1,\ldots,N\}}$ are $N$ i.i.d. random variables which model the natural frequency of the oscillators (also called the \emph{disorder}). It is often assumed that at least the expectation $\E\xi_i$ is finite, in which case, up to a time translation, it reduces to the case where the natural frequencies have zero mean, see \cite{acebron_kuramoto_2005}. When a realization of the natural frequencies is chosen beforehand, then the model is said to be of \emph{quenched} type. At least when $\xi_i=0$ for all $i$, the propagation of chaos on any finite time interval follows immediately from McKean's Theorem~\cref{thm:mckean}. A natural question is therefore the long-time behaviour of the system and the uniform in time propagation of chaos. The limit Fokker-Planck equation can be shown to admit the following family of stationary solutions: 
\begin{equation}\label[IIeq]{eq:vonmisescircle}M_{\kappa,\theta_0}(\theta)\propto \exp(\kappa\cos(\theta-\theta_0)),\end{equation}
where $\kappa\geq0$ solves the compatibility equation $\kappa = 2KI_1(\kappa)/I_0(\kappa)$ and $I_0$ and $I_1$ are the modified Bessel functions of order 0 and 1. The parameter $\theta_0\in\R$ can be taken arbitrarily (by rotational invariance). The probability density function \cref{eq:vonmisescircle} is called the \emph{von Mises distribution with concentration parameter $\kappa$ and center $\theta_0$}. The trivial solution $\kappa=0$ is always a solution of the compatibility equation, it corresponds to the trivial disorder equilibrium where all the oscillators are asymptotically uniformly distributed over the circle. If $K>1$ then there exists also a unique nontrivial solution $\kappa>0$ of the compatibility equation and the associated family of stationary becomes \emph{asymptotically stable} meaning that, up to a negligible (in a certain sense, see below) set of initial conditions, there exists a $\theta_0\in\R$ which depends only on the initial condition and such that the solution of the Fokker-Planck equation converges towards the von Mises distribution \cref{eq:vonmisescircle} associated to this $\theta_0$. This phenomenon is called a \emph{phase transition} and a complete description of the long-time dynamics of the solution of the Fokker-Planck equation can be found in \cite{giacomin_global_2012, degond_phase_2015}. 

Consequently, if the propagation of chaos holds uniformly in time then the empirical measure $\mu_{\mathcal{X}^N_t}$ necessarily converges towards an element of the family \cref{eq:vonmisescircle} as $N,t\to+\infty$. This is not always the case as shown by the large deviation principle proved in \cite[Theorem~1.1]{bertini_synchronization_2014}. More precisely, let $K>1$ and let $\kappa>0$ be the unique nontrivial solution of the compatibility equation. Fix also a constant $T>0$. Assume that $f^N_0$ is $f_0$-chaotic for a $f_0$ is such that $\int_{\mathbb{S}^1} \exp(i\theta)f_0(\dd\theta) \ne 0$ (otherwise the solution of the Fokker-Planck equation does not converge towards an element of \cref{eq:vonmisescircle}). Then Bertini et al. show that there exist $\theta_0\in\R$ which depends only on the initial condition $f_0$ and a sequence of processes $(W^{N,T}_t)_{t\in[0,T]}$ which converges weakly to a standard Brownian motion such that for all $\varepsilon>0$:
\[\lim_{N\to+\infty} \mathbb{P}{\left(\sup_{\tau\in[C(K)/N,T]}\big\|\mu_{\mathcal{X}^N_{N\tau}}-M_{\kappa,\theta_0+D(K)W^{N,T}_\tau}\big\|_{H^{-1}}\leq\varepsilon\right)}=1,\]
where $C(K),D(K)>0$ depend only on $K$, the initial condition and $\varepsilon$. As a consequence, the propagation of chaos is \emph{not} uniform in time and breaks down at times proportional to $N$.

Another way to study the long-time behaviour of particle systems is to consider an appropriate scaling limit. For the Kuramoto model and more generally for McKean-Vlasov gradient systems, the natural scaling is the \emph{diffusive} scaling defined by \[f^{\varepsilon,N}_t(\mathbf{x}^N) := \varepsilon^{Nd}f^N_{t/\varepsilon^2}(\varepsilon\mathbf{x}^N) = \mathrm{Law}(\varepsilon X^1_{t/\varepsilon^2},\ldots,\varepsilon X^N_{t/\varepsilon^2}),\]
where $\varepsilon>0$ is the scaling parameter. In the case of the Kuramoto model, this is the law of a highly oscillating system with a frequency of order $\varepsilon^{-1}$ and $K=\mathcal{O}(\varepsilon^{-1})$. The authors of \cite{delgadino_diffusive-mean_2021} study a class of McKean-Vlasov gradient systems on the torus which generalizes the Kuramoto model. Using a gradient flow framework (see Section \cref{sec:gradientflows}), one of the main results of the article is an explicit counter example which proves that for some chaotic initial conditions, the two limits $N\to+\infty$ and $\varepsilon\to0$ do not commute above the phase transition. Consequently, the propagation of chaos cannot hold uniformly in time. The links between this result and log-Sobolev inequalities is explored in \cite{delgadino_phase_2021}, see also Remark \cref{rem:phasetransitionslogsobuniformpoc}.

\subsubsection{Swarming models}\label[II]{sec:flocking}

Over the last twenty years, there has been a growing interest in both the Mathematics and Physics communities for theorizing the underlying principles of large animal societies. Among the most common examples of such systems, flocks of birds, fish schools, large herds of mammals or ant colonies exhibit a collective coherent complex behaviour without any obvious exterior organizing principle such as a leader. Other examples can be found in the microscopic world (for instance colonies of bacteria or spermatozoa) or in human societies (for instance crowds phenomena or traffic flows). In all these systems, each individual can be roughly described as a kinetic particle $(X^i_t, V^i_t)$ and the underlying principles which model the global motion of the system should obey the Newton's laws (plus noise) $\dd X^i_t = V^i_t\dd t$ and $\dd V^i_t = F(\mathcal{X}^N_t)\dd t$, where $F$ is a force or a sum of forces. This section is devoted to the description of some examples of elementary mechanisms commonly used in \emph{swarming models}. Most of them are based on the assumption that particles have a sensing region and interact with the other particles which belong to this region. The easiest way to model this is to take an \emph{observation kernel} $K:\R_+\to\R_+$ which vanishes at infinity, for instance $K(r)=\1_{[0,R]}(r)$ for a fixed interaction radius $R$, and to consider that the sensing region of a particle at position $X^i_t$ depends on the map $x\mapsto K(|x-X^i_t|)$. Then one has to define which kind of behaviour a particle will adopt: for instance it can try to avoid the other particles in its sensing region or on the contrary to move closer to the center of mass of its neighbours. Alternatively, a particle can simply try to align its velocity with the velocities of the other particles in order to create a coherent motion, this is called a \emph{flocking model}. A gallery of models can be found for instance in the reviews \cite{vicsek_collective_2012} or \cite{albi_vehicular_2019}. Note that unlike classical physical systems, the particles are able to produce their own energy for self-propulsion, so there are no a priori conservation laws (apart from mass conservation). In the Physics literature, such particle system is called \emph{active matter}. 

The next objective is to consider large systems and thus to derive (rigorously) the $N\to+\infty$ limit. When propagation of chaos holds, this reduces the problem to the analysis of a single kinetic PDE. Following the principles of statistical physics, one can also try to compute the \emph{hydrodynamic} limit of the solution of the kinetic PDE to study the system on larger time and space scales. This naturally raises the problem of uniform in time propagation of chaos but for the examples below, we will focus on modelling aspects and we will not address this question; we refer the interested reader to the quoted references and to Section \cref{sec:kuramoto} for an example which demonstrates that the question can become very delicate. 

\subsubsection*{Attraction-Repulsion.} One of the first deterministic mathematical swarming models, due to D'Orsogna et al. \cite{dorsogna_self-propelled_2006}, is based on the combination of self-propulsion and an attraction-repulsion force. With the mean-field scaling introduced in \cite{carrillo_double_2009}, the model reads: 
\begin{align*}
    \frac{\dd V^i_t}{\dd t} &= (\alpha-\beta|V^i_t|^2)V^i_t - \frac{1}{N}\nabla_{x^i}\sum_{j\ne i} U\big(|X^i_t-X^j_t|\big),
\end{align*}
where $U(r)=-C_a\e^{-r/\ell_a}+C_r\e^{-r/\ell_r}$ is the Morse potential. The nonnegative constants $\alpha$, $\beta$, $C_a$, $\ell_a$, $C_r$, $\ell_r$ are respectively the propulsion coefficient, the friction coefficient, the strength of alignment, the typical alignment length, the strength of the repulsion and the typical repulsion length. Due to the propulsion and friction forces, each particle tends to adopt the fixed cruising speed $\sqrt{\alpha/\beta}$. Although entirely deterministic, the propagation of chaos is covered by \cite{bolley_stochastic_2011} and the limit PDE reads:
\[\partial_t f_t(x,v) + v\cdot\nabla_x f_t = -\nabla_v\cdot((\alpha-\beta|v|^2)vf_t) + (\nabla_x U\star \rho[f_t])\cdot\nabla_v f_t,\]
where $\rho[f_t](\dd x) = \int_{\R^d}f_t(\dd x,\dd v)$. The analysis of the limit kinetic PDE and its hydrodynamic limit in \cite{carrillo_double_2009} gives a rigorous theoretical explanation for the emergence of complex patterns such as rotating mills which were observed in numerical simulations only in \cite{dorsogna_self-propelled_2006}. 

\subsubsection*{Flocking.} The alignment mechanism introduced by Cucker and Smale \cite{cucker_mathematics_2007} reads: 
\begin{align*}
    \frac{\dd V^i_t}{\dd t} &=  \frac{1}{N}\sum_{j\ne i } K\big(|X^j_t-X^i_t|\big)(V^j_t-V^i_t), 
\end{align*}
where $K$ is an observation kernel which is typically taken equal to $K(r) = (1+|r|^2)^{-\gamma/2}$, $\gamma>0$. The main result is that if the observation kernel is large enough in the sense that $\int_0^{+\infty}K(r)\dd r = +\infty$, then the particle system satisfies for all $i,j\in\{1,\ldots,N\}$,
\[|V^i_t-V_\infty|\leq C_1\e^{-\lambda t},\quad |X^i_t-X^j_t|\leq C_2,\]
for some constants $C_1,C_2,\lambda>0$ and for an asymptotic velocity $V_\infty\in\R^d$. Note that since the momentum is preserved, $V_\infty = \frac{1}{N}\sum_{i=1}^N V^i_0$. This property is called \emph{flocking}. There is an extensive literature on the deterministic Cucker-Smale model that we will not discuss here, see the reviews \cite{naldi_particle_2010, muntean_derivation_2014,  albi_vehicular_2019}. 

On the other hand, there are various ways to add a stochastic component to the Cucker-Smale model. Maybe the most obvious way in this context, is to consider the McKean-Vlasov model introduced in \cite{ha_emergence_2009}: 
\[\dd V^i_t = \frac{1}{N}\sum_{j\ne i } K\big(|X^j_t-X^i_t|\big)(V^j_t-V^i_t) + \sigma\dd B^i_t,\]
for $N$ independent Brownian motions $(B^i_t)^{}_t$. In this case and despite the fact the drift is not globally Lipschitz and bounded, the propagation of chaos is proved in~\cite{bolley_stochastic_2011} using the synchronous coupling method (see also Section \cref{sec:mckeantowardssingular}) or in \cite{pedeches_asymptotic_2018} using martingale arguments (see also Section \cref{sec:martingalecompactness}). The limit Fokker-Planck equation reads:
\begin{equation}\label[IIeq]{eq:kineticfp}\partial_t f_t(x,v) + v\cdot\nabla_x f_t = -\nabla_v\cdot(\xi[f_t]f_t)+\frac{\sigma^2}{2}\Delta_v f_t,\end{equation}
with
\[\xi[f_t](x,v):=\int_{\R^d\times\R^d} K(|x'-x|)(v'-v)f_t(\dd x',\dd v').\]
More refined models can also be considered with a non constant diffusion matrix, with boundary conditions \cite{choi_propagation_2018} or when the observation kernel is anisotropic. In this last case, one can for instance consider an observation kernel $K(|X^j_t-X^i_t|)\equiv K_{V^i_t}(|X^j_t-X^i_t|)$ which also depends on the velocity of the $i$-th particle: this includes the biologically relevant case where the observation kernel is the indicator function of a cone of vision centered around the velocity of the particle. In this case, propagation of chaos is proved in \cite{choi_collective_2019} using a synchronous coupling argument. 

In \cite{ahn_stochastic_2010}, Ahn and Ha considered the Cucker-Smale with a random environmental noise:  
\[\dd V^i_t = \frac{1}{N}\sum_{j\ne i } K\big(|X^j_t-X^i_t|\big)(V^j_t-V^i_t) + \sigma\big(Z^i_t,\mu_{\mathcal{Z}^N_t}\big)\dd B_t,\]
where $Z^i_t = (X^i_t,V^i_t)$ and $(B_t)_t$ is a Brownian motion which is the same for all the particles (also called \emph{common noise}) and $\sigma$ is a possibly non constant diffusion matrix. In this case, the propagation of chaos does not hold in the usual sense. For general McKean-Vlasov systems of this form, given a realization of the common noise, a \emph{conditional} propagation of chaos property can be shown \cite{coghi_propagation_2016} by revisiting the classical arguments of Dobrushin \cite{dobrushin_vlasov_1979} in the deterministic case. However the limit law $f_t$ is not deterministic and satisfies a \emph{stochastic} PDE which depends on the common noise (roughly speaking, it is the PDE \cref{eq:kineticfp} where the Laplacian is replaced by a transport term involving the Brownian motion). For the Cucker-Smale system, this type of result can be found in~\cite{choi_cucker-smale_2019}. 

There exist many other Cucker-Smale models where the stochasticity is incorporated through a diffusive behaviour. For further examples, we refer the interested reader to the review \cite{cattiaux_stochastic_2018} and the references therein. Lately, \cite{friesen_stochastic_2020} proposed a stochastic Cucker-Smale model based on a Nanbu interaction mechanism (Example \cref{example:nanbuboltzmann}). The propagation of chaos for this model is proved using martingale arguments. 

\subsubsection*{Flocking with geometrical constraints.} 

In the 90's, Vicsek et al. \cite{vicsek_novel_1995} introduced a time discrete ``flocking algorithm'' using the minimal assumption that the particle move at a fixed constant speed. The Vicsek model has quickly become one of the most prominent models in the active matter literature. Several works have numerically exhibited the emergence of complex patterns at the particle level; see for instance \cite{chate_collective_2008} where the emergence of high-density band-like structures on a compact spatial domain is studied. From a mathematical point of view, since the speed of the particles is fixed, the velocity of each particle is defined by its \emph{orientation} which is an element of the unit sphere $\mathbb{S}^{d-1}$. The motion can thus be interpreted as a constrained dynamical system on a manifold. Following this idea, Degond and Motsch \cite{degond_continuum_2008} gave a mean-field time-continuous interpretation of the Vicsek model defined by a system of Stratonovich SDEs: 
\begin{align*}
    \dd X^i_t &= c_0 V^i_t \dd t \\
    \dd V^i_t &= \nu\big(|J^i_t|\big)\mathsf{P}(V^i_t)\Omega^i_t \dd t + \sqrt{2\sigma\big(|J^i_t|\big)}\mathsf{P}(V^i_t)\circ \dd B^i_t,
\end{align*}
where $c_0>0$ is the speed of the particles, $\nu,\sigma>0$ are respectively the intensity of the alignment and the strength of the diffusion and 
\[\Omega^i_t := \frac{J^i_t}{|J^i_t|}\in \mathbb{S}^{d-1},\quad J^i_t = \frac{1}{N}\sum_{j=1}^N K\big(|X^j_t-X^i_t|\big)V^j_t \in \R^d.\] 
Given $v\in\R^d$, the matrix $\mathsf{P}(v):=I_d-\frac{v\otimes v}{|v|^2}$ is the projection on the plane orthogonal to $v$. The SDE is written in the Stratonovich sense (indicated by the symbol $\circ$), so that for all $i$ and all $t\geq0$, $V^i_t\in \mathbb{S}^{d-1}$ provided that $|V^i_0|=1$. In this model, the alignment force exerted on particle $i$ belongs to the tangent hyperplane of the orientation $V^i_t$ and is directed towards the local average orientation $\Omega^i_t$. The strength of this force may depend on the norm of $J^i_t$ which plays the role of a (local) \emph{order parameter}: when the system is in a disordered state with all the orientations uniformly scattered on the sphere, then $|J^i_t|$ tends to zero as $N \rightarrow + \infty$. In the opposite case of a flocking state, $|J^i_t|$ concentrates around a fixed point of the sphere, with a concentration parameter which depends on the observation kernel. 

The propagation of chaos property is proved in \cite{bolley_mean-field_2012} in the case $\nu(|J|)=|J|$. The authors use the synchronous coupling method for a particle system which is defined a priori in the whole space $E=\R^d\times\R^d$. Consequently, a regularisation argument is needed for the projection matrix (which is singular at the origin) to prove local well-posedness. It is then shown that this system stays constrained on the manifold $E=\R^d\times \mathbb{S}^{d-1}$ (i.e. $|V^i_t|=1$ for all $t$) provided that it holds true at initial time. Since the singularity is never visited, the regularisation can therefore be removed, the particle system coincides with the original one and is well-defined globally in time. The result is extended in \cite{briant_cauchy_2021} in particular in the more singular case $\nu(|J|)=1$. The limit Fokker-Planck equation reads:
\begin{equation}\label[IIeq]{eq:kineticfpvicsek}\partial_t f_t(x,v) + c_0 v\cdot\nabla_x f_t = \nu\big(|J[f_t]|\big)\nabla_v\cdot(\mathsf{P}(v)\Omega[f_t]f_t) + \sigma\big(|J[f_t]|\big)\Delta_v f_t,\end{equation}
where 
\[J[f_t](x) = \int_{\R^d} K(|x'-x|)v' f_t(\dd x',\dd v')\in\R^d,\quad \Omega[f_t] = \frac{J[f_t]}{|J[f_t]|}\in\mathbb{S}^{d-1},\]
and $\Delta_v$, $\nabla_v\cdot$ denote respectively the Laplace-Beltrami and the divergence operators on the sphere $\mathbb{S}^{d-1}$. 

An analogous mean-field jump particle system is introduced in \cite{dimarco_self-alignment_2016} and the corresponding propagation of chaos result which leads to a BGK equation is proved in \cite{diez_propagation_2020}. Keeping the key assumption of the fixed speed, a Boltzmann interaction mechanism is proposed in \cite{bertin_boltzmann_2006, bertin_hydrodynamic_2009} and the propagation of chaos for various Boltzmann models is studied in \cite{carlen_kinetic_2013-1, carlen_kinetic_2013}. 

The behaviour of the spatially-homogeneous version of the kinetic Fokker-Planck PDE \cref{eq:kineticfpvicsek} is well-understood: well-posedness results and long-time convergence results are proved in \cite{figalli_global_2018, kang_dynamics_2016, briant_cauchy_2021} in the case $\nu(|J|)=1$ and phase transition phenomena are explored in depth in \cite{degond_phase_2015} in particular in the case $\nu(|J|)=|J|$. The stationary solutions of the spatially-homogeneous PDE belong to the family of von Mises distributions on the sphere $\mathbb{S}^{d-1}$ thus generalizing the framework of the Kuramato model to higher dimensions (the Kuramato model is equivalent to the one dimensional spatially-homogeneous Vicsek model). Finally, the hydrodynamic limit is derived in \cite{degond_continuum_2008, dimarco_self-alignment_2016}. However, the analysis of the spatially-inhomogeneous case remains mostly open. To the best of our knowledge, and despite some numerical evidence, a complete theory able to explain the phenomena reported at the particle level in \cite{chate_collective_2008} is still lacking. For Boltzmann models, very few is known in the mathematics literature even at the kinetic level \cite{carlen_boltzmann_2015}. 

The framework of the Vicsek model can also be used to model alignment mechanisms on other manifolds than the sphere. For instance, in dimension 3, the particles may be defined by their full \emph{body-orientation} which is a rotation matrix in $SO_3(\R)$, see the lecture notes \cite{giacomin_alignment_2019} for an extension of the Vicsek model to this case. In the liquid crystal literature, a different alignment mechanism called \emph{nematic} is used, which roughly speaking, corresponds to replacing the sphere $\mathbb{S}^{d-1}$ by the projective space $\mathbb{S}^{d-1}/\pm\mathrm{Id}$, see for instance \cite{degond_nematic_2020} and the references therein. 

\subsubsection*{Topological interactions.} There is experimental and numerical evidence \cite{ballerini_interaction_2008} to support the idea that in order to maintain cohesion in a bird flock, the interactions between the individuals are rather based on their \emph{rank} than on their relative distance. It means that given a particle $i$, the influence of a particle $j$ on $i$ at time $t$ depends on the rank $R[\mu_{\mathcal{X}^N_t}](X^i_t,X^j_t)\in\{1,\ldots,N\}$ of particle $j$ defined such that particle~$j$ is the $R[\mu_{\mathcal{X}^N_t}](X^i_t,X^j_t)$-th nearest neighbour of $i$: 
\[R[\mu_{\mathcal{X}^N_t}](X^i_t,X^j_t) := \#\big\{k\in\{1,\ldots,N\},\,\,|X^i_t-X^k_t|<|X^i_t-X^j_t|\big\}.\]
In a mean-field framework, it is more natural to use the normalised rank defined by $r[\mu_{\mathcal{X}^N_t}](X^i_t,X^j_t)=R[\mu_{\mathcal{X}^N_t}](X^i_t,X^j_t)/N$ where given $x,y,z\in\R^d$ and $\mu\in\pb(\R^d)$,
\begin{equation}\label[IIeq]{eq:topointeraction}r[\mu](x,y) :=\big\langle \mu,\psi(x,y,\cdot)\big\rangle,\quad \psi(x,y,z):=\1_{[0,1)}{\left(\frac{|x-z|}{|x-y|}\right)}.\end{equation}
All the models previously described can be alternatively defined using \emph{topological interactions} by replacing the metric observation kernel $K(|X^i_t-X^j_t|)$ by the rank-based observation kernel $K\big(r[\mu_{\mathcal{X}^N_t}](X^i_t,X^j_t)\big)$, where in this case $K:[0,1]\to\R_+$ is a smooth given function. This change has two consequences: first the interaction is no longer symmetric (this is not a real difficulty) and secondly, this adds a new source of nonlinearity but since it is of mean-field type (i.e. it only depends on the empirical measure), the limit can be easily derived, at least formally. Note however that since the function $\psi$ in \cref{eq:topointeraction} is not Lipschitz, an ad hoc argument is needed, for instance a regularisation procedure (see Section \cref{sec:mckeantowardssingular}). For the (deterministic) Cucker-Smale model, this is investigated in \cite{haskovec_flocking_2013}. For Boltzmann (Nanbu) interactions with a collision rate which depends on $K\big(r[\mu_{\mathcal{X}^N_t}](X^i_t,X^j_t)\big)$, several models are discussed in \cite{blanchet_topological_2016, blanchet_kinetic_2017} and a rigorous propagation of chaos result is proved in \cite{degond_propagation_2019}. 

\subsubsection{Neuron models}\label[II]{sec:neurons}

The modelling of (biological) neuronal networks has a long story that we do not intend to extensively review here. We will only give a glimpse on the subject by quoting some recent models relevant with our subject.

\subsubsection*{Mean-field jump models.} A neuron is mainly described by its membrane potential in $\R_+$ and maybe also by some other variables which depend on the considered model. In an abstract setting $X^i_t\in \R^d$ will denote the state of neuron $i\in\{1,\ldots,N\}$ at time $t$. The value of the membrane potential is typically linked to the jump rate of random events called \emph{spikes}. When a neuron spikes, its membrane potential is automatically reset at a default value and this spiking event increases the membrane potentials of the other (neighbouring) neurons. In a mean-field setting, this small potential increase is proportional to $1/N$. This small toy model is exactly a mean-field jump model with simultaneous jumps considered in Example \cref{example:simultaneousjumps}. Such model was considered first in \cite{de_masi_hydrodynamic_2015} and then in \cite{fournier_toy_2016}. The propagation of chaos can be proved using compactness or coupling methods. More recently, the question is also addressed in various very general cases which include diffusion models in \cite{andreis_mckeanvlasov_2018}. 

\subsubsection*{Diffusion models.} Another popular class of neuron models is based on McKean-Vlasov diffusion processes. In the abstract setting described in \cite{touboul_propagation_2014}, the neurons are clustered into $P(N)$ populations. Each population of neurons $\alpha$ has $N_\alpha$ neurons and $N = \sum_{\alpha=1}^{P(N)} N_\alpha$. Each population $\alpha$ is located at a position $r_\alpha\in\Gamma$ where $\Gamma$ is a nice space modelling the cerebral cortex. The spike of a neuron at location $r_\gamma$ produces a time continuous current which affects the other neurons at location $r_\alpha$ with a delay $\tau(r_\alpha,r_\gamma)\geq0$. The state of the neuron $X^i_t$  belonging to population $\alpha$ is thus governed by the SDE: 
\[\dd X^i_t = F(t,r_\alpha,X^i_t)\dd t + \frac{1}{P(N)}\sum_{\gamma=1}^{P(N)} \frac{1}{N_\gamma}\sum_{p(j)=\gamma}  b\big(r_\alpha,r_\gamma,X^i_t,X^j_{t-\tau(r_\alpha,r_\gamma)}\big)\dd t + \sigma(r_\alpha)\dd B^i_t,\]
where $b:\Gamma\times \Gamma\times \R^d \times \R^d \to \R^d $ is the current function, $p(j)\in\{1,\ldots,P(N)\}$ is the population index of particle $j$ and the functions $F:\R_+\times\Gamma\times\R^d\to\R^d$, $\sigma:\Gamma\to\mathcal{M}_d(\R)$ denote the intrinsic deterministic dynamics and the external noise exerted on the neuron. The limit Fokker-Planck equation is not of one of the types previously studied: it involves a time delay and an intricate spatial dependence which both raise well-posedness issues. On top of that, for classical neuron models such as the FitzHugh-Nagumo model, the parameters are not globally Lipschitz. The adaptation to this complex framework of the classical synchronous coupling method of Sznitman can be found in \cite{bossy_clarification_2015, touboul_propagation_2014}. 

\subsubsection*{Point processes models.} Finally, forgetting the details of the membrane potential, the neuronal activity can also be modelled by $N$ counting processes (i.e. non-decreasing integer-valued jump processes) with a jump parameter which depends on the number of past and neighbouring jumps. These processes are called \emph{(interacting) Hawkes processes} or \emph{self-exciting counting processes}. The state of the neuron $X^i_t\in\N$ is simply defined as its number of spikes up to time $t$. The mean-field analysis of such models has been initiated in \cite{delattre_hawkes_2016}. Shortly later, Chevallier \cite{chevallier_mean-field_2017} introduced a class of \emph{age dependent Hawkes processes} for which the jumping rate of the neuron $X^i_t$ depends on the elapsed time since the last spike, called the \emph{age} and denoted by: 
\[S^i_{t^-}:= t - \sup\{T^i\in X^i,\,\, T^i< t\},\]
where we write $T^i\in X^i$ when $T^i\in\R_+$ is a jump time of the counting process $(X^i_t)^{}_t$. Moreover each spiking event affects the jump rate of the other neurons. In summary, the jump rate of neuron $i$ is defined by:
\[\lambda^i_t := \psi{\left(S^i_{t^-}, \frac{1}{N}\sum_{j=1}^N \Big(\int_0^{t^-} H_{ij}(t-\tau)X^j(\dd \tau) + F_{ij}(t)\Big)\right)},\]
where $H_{ij}, F_{ij} : \R_+\to\R$ are random interaction functions, $\psi:\R_+\times\R\to\R_+$ is called the \emph{intensity function}, and $X^j$ denotes the random measure (or point process) associated to the process $(X^j_t)^{}_t$. In this expression, the communication function $H_{ij}$ models how the spike of a neuron $j$ at time $\tau$ affects the spike rate of neuron $i$ at time $t$. In our usual setting, it means that $X^i_t$ satisfies the SDE: 
\[X^i_t = \int_0^t\int_0^{+\infty}\1_{u\leq \psi\big(S^i_{s^-},\frac{1}{N}\sum_{j=1}^N(\int_0^{s^-} H_{ij}(s-\tau)X^j(\dd\tau) + F_{ij}(s))\big)}\mathcal{N}^i(\dd s, \dd u),\]
where the $\mathcal{N}^i$ are $N$ independent Poisson random measures on $\R_+\times\R_+$ with intensity $\dd s\otimes\dd u$. Using a synchronous coupling argument, it is shown in \cite{chevallier_mean-field_2017} that the limit $N\to+\infty$ exists and the distribution of the age $S^i_{t^-}$ of each neuron at time $t$ converges towards the solution of the PDE: 
\[\left\{
\begin{array}{l}
\displaystyle{\partial_t n(s,t) + \partial_s n(s,t) + \psi(s, m(t)+F_0(t))n(s,t) = 0,}\\
\displaystyle{n(0,t) = \int_0^{+\infty} \psi(s, m(t)+F_0(t))n(s,t)\dd s,}\\ 
\displaystyle{m(t) = \int_0^t h(t-\tau)n(0,\tau)\dd\tau},
\end{array}
\right.\]
where $F_0$ and $h$ denote the expectations of the functions $F_{ij}$ and $H_{ij}$. The solution $n(s,t)$ is the distribution of neurons with age $s$ at time $t$. This PDE was studied before by Pakdaman, Perthame and Salort \cite{pakdaman_adaptation_2014}. On this subject, see for instance~\cite{canizo_asymptotic_2018} and the references therein.

\subsubsection{Socio-economic models}\label[II]{sec:socioeconomicmodels}

In this section, the particles model interacting socio-economic agents (human beings) with all the variety of possible interactions that one can imagine: to give a flavour of some recent modelling trends, we present a selection of models for opinion dynamics, wealth distribution or rating score in games. More on the subject can be found in the book \cite{naldi_mathematical_2010}. The only modelling assumption is that an interaction involves only two agents so that all the models presented are Boltzmann models. Interactions which involve more than two but still a finite fixed number of agents could also be relevant in some situations but we will not discuss this point \cite{toscani_kinetic_2020}. The following parametric Boltzmann models are defined using the notations of Section \cref{sec:boltzmann_summary} with a generator of the form \cref{eq:Boltzmannnonsym_summary}. 

\subsubsection*{Opinion dynamics.} In the opinion formation considered model in \cite{toscani_kinetic_2006}, an opinion is a real number in $[-1,1]$, $L^{(1)}=0$ and an interaction between two agents with opinions $(z_1,z_2)$ leads to the post-collisional opinions:
\begin{align*}
\tilde{\psi}_1(z_1,z_2,\eta_1,\eta_2) &= z_1 -  \gamma P(|z_1|)(z_1-z_2) + \eta_1D(|z_1|), \\ 
\tilde{\psi}_2(z_1,z_2,\eta_1,\eta_2) &= z_2 -  \gamma P(|z_2|)(z_2-z_1) + \eta_2D(|z_2|),
\end{align*}
where $\gamma\geq0$ and the functions $P$ and $D$ model respectively the intrinsic tendency to the consensus and the diffusion. Typically, for extreme opinions, one expects $P$ and $D$ to be small. The parameters $(\eta_1,\eta_2)$ are independent zero mean random variables with a fixed variance $\sigma^2$. A similar model posed on the whole real line $\R$ and written in Nanbu form is studied in \cite{degond_continuum_2017} with $P(|z|) = 1$ and $D(|z|)=1$. The collision rate may depend on the individual opinions of the agents or of the difference between their opinions (typically, two agents with far-away opinions are less likely to interact). In \cite{degond_continuum_2017} collision rates which depend on a mean-field quantity are also considered. In both works, the authors study the long-time dynamics and the equilibrium distributions of the model. An important assumption is the grazing collision scaling $\gamma\to0$, $\sigma^2/\gamma\to \kappa$ for a fixed $\kappa$. This choice turns the Boltzmann equation into a more amenable (Landau) Fokker-Planck equation (see \cite{toscani_grazing_1998, villani_review_2002} and Section \cref{sec:landau}). Phase transitions phenomena for this equation are investigated in \cite{degond_continuum_2017} as well as non-spatially homogeneous versions of this model. 

\subsubsection*{Wealth distribution.} A model of wealth distribution model inspired from \cite{matthes_steady_2008} can be found in \cite{cortez_quantitative_2016}. The authors assume $E=\R$ with $L^{(1)}=0$, $\lambda=1$, $\tilde{\Theta}=\R^4$ and 
\[\tilde{\psi}_1\big(z_1,z_2,(L,R,\tilde{L},\tilde{R})\big) = Lz_1+Rz_2,\]
and
\[\tilde{\psi}_2\big(z_1,z_2,(L,R,\tilde{L},\tilde{R})\big) = \tilde{L}z_2+\tilde{R}z_1.\]
In this model, the state of a particle represents the wealth of an individual and the parameters $(L,R,\tilde{L},\tilde{R})$ specify how a trade between two individuals affect their wealth. This model generalises a famous model due to Kac \cite{kac_foundations_1956}. It is assumed that the parameters $(L,R,\tilde{L},\tilde{R})$ are distributed so that $\E[L+R]=\E[\tilde{L}+\tilde{R}]=1$. This model is called \emph{conservative}, which means that during a trade, the wealth of each agent is conserved in average. Several other examples of conservative and non conservative models are presented in \cite{matthes_steady_2008}. The rigorous propagation of chaos property is proved in \cite{cortez_quantitative_2016} using a coupling method (Section \cref{sec:couplingBoltzmann}). 

Lately, the authors of \cite{during_continuum_2020} introduced the non-conservative model 
\begin{align*}
\tilde{\psi}_1(z_1,z_2,R) &= R(z_1+z_2), \\ 
\tilde{\psi}_2(z_1,z_2,R) &= (1-R)(z_1+z_2),
\end{align*}
where $R$ is a parameter drawn from the uniform distribution on $[0,1]$. The model is originally written in a discrete time and discrete space setting (meaning that the wealths $z_1,z_2$ belong to $\N$). The continuum and mean-field limits are investigated using a martingale approach. Then, the limit Boltzmann equation is shown to admit several families of equilibria depending on the initial wealth distribution.  

\subsubsection*{Elo rating system.} In this example, the particles are players in a one-versus-one game, for instance during a chess competition or during a sport or e-sport event. Each player is characterised by its intrinsic strength $\rho$ (which is fixed) and a rating $r$. The goal of the Elo rating system is to update the ratings of the players at each game so that they match the intrinsic strengths of the players. Following the Elo system, a simple model for a game is a Boltzmann collision model of type \cref{eq:boltzmannpsitilde_summary} with $L^{(1)} =0$, which updates the ratings of two players $z_1=(r_1,\rho_1)$ and $z_2=(r_2,\rho_2)$ as follows: 
\begin{align*}
\tilde{\lambda}(r_1,r_2) &= \lambda w(|r_1-r_2|),\\
\tilde{\psi}_1(r_1,r_2,\theta) &= r_1 + \gamma\big(S(\rho_1,\rho_2,\theta)-b(r_1-r_2)\big), \\ 
\tilde{\psi}_2(r_1,r_2,\theta) &= r_2 - \gamma\big(S(\rho_1,\rho_2,\theta)-b(r_1-r_2)\big),
\end{align*}
where $\lambda,\gamma>0$ are given parameters, $S(\rho_1,\rho_2,\theta) \in\{-1,1\}$ is the \emph{score} of the game (1 means a win) and $b:\R\to(-1,1)$ is an odd increasing function which predicts the score of the game given the difference of ratings. The parameter $\theta\sim\nu(\dd\theta)$ is assumed to be such that 
\[\E_\nu\big[S(\rho_1,\rho_2,\theta)\big] = b(\rho_1-\rho_2),\] 
which means that the probability of a win for the player 1 is equal to
\[\mathbb{P}_\nu(S(\rho_1,\rho_2,\theta)=1) = \frac{1}{2}+\frac{1}{2}b(\rho_1-\rho_2).\]
The collision rate $w$ depends only on the absolute difference between the ratings (typically a game involves players with similar rating scores). The Boltzmann equation \cref{eq:boltzmannpsitilde_summary} reads in weak form: 
\begin{align*}
    &\frac{\dd}{\dd t}\iint_{\R^2}\varphi(r,\rho) f_t(\dd r,\dd \rho)\\
    &= \frac{\lambda}{2}\int_{\R^2\times\R^d\times\Theta} w(|r_1-r_2|)\Big\{\varphi\big(r_1 + \gamma\big(S(r_1,r_2,\theta)-b(r_1-r_2)\big),\rho_1\big)\\
    &\phantom{ \frac{\lambda}{2}\int_{\R^2\times\R^d\times\Theta} w(|r_1-r_2|)\Big\{}\quad+\varphi\big(r_2 - \gamma\big(S(r_1,r_2,\theta)-b(r_1-r_2)\big),\rho_2\big)\\
    &\phantom{ \frac{\lambda}{2}\int_{\R^2\times\R^d\times\Theta} w(|r_1-r_2|)\Big\{}\quad-\varphi(r_1,\rho_1)-\varphi(r_2,\rho_2)\Big\} f_t(\dd r_1,\dd \rho_1)f_t(\dd r_2,\dd \rho_2)\nu(\dd\theta).
\end{align*}
In the grazing collision limit $\gamma\to0$ and $\gamma\lambda\to\kappa$ for a fixed $\kappa>0$, a first order Taylor expansion gives, at least formally, the following equation in strong form: 
\begin{equation}\label[IIeq]{eq:pdeelo}\partial_t f_t(r,\rho) + \partial_r \big(a[f_t]f_t\big) = 0,\end{equation}
where
\[a[f_t](r,\rho) := \kappa\partial_r\int_{\R^2\times\R^2} w(|r_*-r|)\big(b(\rho_*-\rho)-b(r_*-r)\big)f_t(\dd r_*,\dd\rho_*).\]
This is the equation derived in \cite{jabin_continuous_2015} from a time-discrete model. A more elaborated model is proposed in \cite{during_boltzmann_2019} to incorporate a \emph{learning procedure} which increases the intrinsic strength of the players at each game. The long-time behaviour of the grazing collision limit Fokker-Planck equation is then investigated theoretically and numerically. In particular, the solution of \cref{eq:pdeelo} is shown to concentrate on the diagonal $\{\rho=r\}$ as expected. 

\subsection{Applications in data sciences and optimization}\label[II]{sec:datasciences}

Nowadays, the development of data sciences has pushed the development of ever more efficient algorithms. Typical tasks the are discussed below include sampling and filtering (Section \cref{sec:mcmc}), optimization (Section \cref{sec:optimization}) and the training of neural networks (Section \cref{sec:neuralnetworks}). All these situations are challenging, in particular due to the curse of dimensionality, to the high computational cost of naive methods or to the difficulty of finding a satisfactory theoretical framework to prove the convergence of the algorithms. To cope with these problems, various \emph{metaheuristic} methods based on the simulation of systems of particles have been developed. The models in Section \cref{sec:selforganization} illustrate how simple interaction mechanisms can lead to a complex behaviour. In this section, we explore some ideas to design good interaction mechanisms to be used to solve difficult numerical problems. The motivation is twofold: on the one hand, particle systems are easy to simulate and on the other hand, the mean-field theory gives a natural theoretical foundation for the convergence proof of the methods. 

\subsubsection{Some problems related to Monte Carlo integration}\label[II]{sec:mcmc}

Let $\pi$ be an unknown probability density function on a state space $E$ called the \emph{target distribution}. In Bayesian statistics, $\pi$ is typically a \emph{posterior} distribution which gives the distribution of the parameters of a model given the observations. To get an estimate of these parameters, one needs to compute various observables of the form $\langle \pi, \varphi\rangle$ for a test function $\varphi\in C_b(E)$. In general it is not possible to compute directly such an integral because the value of $\pi$ at each point can be computed only up to a multiplicative normalising constant or because the dimension of the state space is too high to use standard quadrature methods. The Monte Carlo paradigm is based on the law of large numbers: if $X^1,\ldots, X^N$ are $N$ independent $\pi$-distributed samples, then an asymptotic estimate of the observable is $\langle \mu_{\mathcal{X}^N}, \varphi\rangle$. However, constructing \emph{good} samples is not easy: in this section, we present a selection of known methods to achieve this goal and illustrate them with some applications. The underlying idea is to look at the samples as particles which are chaotic or, in a dynamical framework, which propagates chaos towards the target distribution: this thus provide many samples which becomes asymptotically i.i.d. and $\pi$-distributed. A classical reference on Monte Carlo methods is \cite{robert_monte_2004}.

\subsubsection*{Scaling limits of the Metropolis-Hastings algorithm.}

In a series of famous articles \cite{metropolis_monte_1949, metropolis_equation_1953, hastings_monte_1970} Metropolis, Hastings et al. have introduced an algorithm to construct a Markov chain which is ergodic with stationary distribution $\pi$. It aims to sample approximately $\pi$-distributed random variables, for a probability measure $\pi$ known up to a multiplicative factor (for instance a Gibbs measure with density $Z^{-1} \e^{-V}$ with respect to some non-negative measure, where the potential $V$ is known but~$Z$ can be very expensive to compute), which can be hard to sample from.
This renowned algorithm has become a building block for many more advanced methods. In its most basic form, it produces a single ergodic time-discrete Markov chain $(X_k)_{k\in\N}$ such that $\mathrm{Law}(X_k) \to \pi$ when $k\to+\infty$ and $X_{k},X_{k'}$ are asymptotically independent when $|k-k'|\to+\infty$. Although very efficient in simple cases, the convergence of the Metropolis-Hastings algorithm is often slow, in particular when~$\pi$ is multimodal. This is due to the sequential nature of the algorithm: typically, the desired $\pi$-distributed samples are extracted from the states of only one Markov chain at different times, well spaced in time and after an initial burn-in phase. 

Among the many extensions and improvements of the Metropolis-Hastings algorithm, the recent article \cite{clarte_collective_2021} studies a more efficient parallelised version of the Metropolis-Hasting algorithm which is directly inspired by the theory of propagation of chaos. The starting point is a map $E\times \pb(E)\to \pb(E), (x,\mu)\mapsto \Theta_\mu(\dd y | x)$ called the \emph{proposal distribution}. Let $\alpha_\mu$ and $h$ be the functions defined for all $\mu\in\pb(E)$, $x,y\in E$ and $u\in\R_+$ by
\[\alpha_\mu(x,y) := \frac{\Theta_\mu(y|x)\pi(x)}{\Theta_\mu(x|y)\pi(y)},\quad h(u) := \min(1,u).\]
Despite the dependence on $\mu$, $\alpha_\mu$ has to be easy enough to compute for this method to be numerically useful (nice typical examples include $\Theta_\mu(y|x) = K \star \mu (y)$ for some suitable kernel $K$). This step requires the knowledge of $\pi$ up to a multiplicative constant: this is particularly well-suited for Gibbs measures of the kind $\pi(\dd x) = Z^{-1} \e^{-V(x)} \mu(\dd x)$ with a potential $V$ and a non-negative measure $\mu$ which make the normalization constant (often called the partition function) $Z$ costly to compute. For numerical reasons, it will be often more amenable to work with $\log \alpha_\mu(x,y)$ instead of $\alpha_\mu(x,y)$ directly (in the previous example, this reduces the computation of the quotient $\pi(x)/ \pi(y)$ to the difference $V(x) - V(y)$)). The algorithm in \cite{clarte_collective_2021} constructs the Markov chain $\mathcal{X}^N_k= (X^1_k,\ldots,X^N_k)$ on $E^N$ such that each component $i\in \{1,\ldots,N\}$ is updated at step $k\in\N$ according to the transition kernel:
\[X^i_{k+1} \sim K_{\mu_{\mathcal{X}^N_k}}(X^i_k, \dd y),\]
where for $x\in E$  and $\mu\in\pb(E)$, the transition kernel is of the form 
\[K_{\mu}(x,\dd y) := 
\underbrace{\vphantom{\int_{z\in E}}h(\alpha_{\mu}(x,y))\Theta_{\mu}(\dd y|x)}_{\text{accept}}
+
\underbrace{\Big[1-\int_{z\in E}h(\alpha_{\mu}(x,z))\Theta_{\mu}(\dd z|x)\Big]\delta_x(\dd y)}_{\text{reject}}.\]

From an algorithmic point of view, at each iteration $k$ and for each particle $i$, a proposal $Y^i_k \sim \Theta_{\mu_{\mathcal{X}^N_k}}(\dd y | X^i_k)$ is sampled first; then the state of particle $i$ at the next iteration is set to $X^i_{k+1}=Y^i_k$ with probability $h(\alpha_{\mu_{\mathcal{X}^N_k}}(X^i_k,Y^i_k))$ (accept) and to $X^i_{k+1}=X^i_k$ otherwise (reject). 

The classical Metropolis-Hasting algorithm corresponds to the case where $\Theta_\mu$ does not depend on the measure argument $\mu$, in which case the previous construction simply gives $N$ independent Markov chains. When the proposal distribution depends on the empirical measure of the system, then this algorithm defines an interacting mean-field jump particle system in discrete time. Note that in this case $\pi^{\otimes N}$ is generally not a stationary distribution of the particle system. To get back to the traditional continuous time framework, it is possible to simply attach to each particle an independent Poisson process which triggers the jumps or a global Poisson process which triggers the simultaneous jumps of the $N$ particles. The result is a particle system of the form described in Section \cref{sec:meanfieldjump_summary}. Under appropriate Lipschitz regularity assumptions on $\Theta_\mu$ which are detailed in \cite{clarte_collective_2021}, then, when $N\to+\infty$, the propagation of chaos property holds towards the solution of the integro-differential equation: 
\begin{equation}\label[IIeq]{eq:pdecmc}
\partial_t f_t(x) = \int_{E} \pi(x)\Theta_{f_t}(y|x)h(\alpha_{f_t}(x,y)){\left(\frac{f_t(y)}{\pi(y)}-\frac{f_t(x)}{\pi(x)}\right)}\dd y.\end{equation}
This result is proved in \cite{clarte_collective_2021} using the optimal coupling argument described in Section \cref{sec:optimalcoupling}. Note that when $\Theta$ does not depend on its measure argument, then Equation \cref{eq:pdecmc} is nothing more than the forward Kolmogorov equation associated to the time continuous version of the Metropolis-Hasting Markov chain. In both the interacting and non interacting cases, it can readily be seen that $\pi$ is a stationary solution of \cref{eq:pdecmc}. The propagation of chaos also ensures the asymptotic independence of the particles as expected. 

This mean-field interpretation of the Metropolis-Hasting algorithm has two main advantages: first the exponential convergence of the solution of \cref{eq:pdecmc} towards $\pi$ with an explicit convergence rate can be deduced from a purely analytical study of Equation \cref{eq:pdecmc}. In \cite[Section 5]{clarte_collective_2021}, such result follows from the entropy-dissipation structure of the equation: for a given level of precision, the time-convergence estimate gives a value of $T$ such that $f_T$ is a good enough approximation of $\pi$, and the mean-field system of \cite{clarte_collective_2021} gives an approximation of $f_t$ (and a fortiori of $f_T$, hence of $\pi$), which is uniform for $t$ in $[0,T]$. Secondly, this analysis gives some rationale for the choice of the proposal distribution, which is critical in all Metropolis-Hasting based methods. In \cite{clarte_collective_2021}, the best convergence rate is obtained for $\Theta_\mu(\dd y|x) = K\star \mu(y)\dd y$ for a normalised symmetric observation kernel $K$ (typically gaussian) which approximates the Dirac distribution $\delta_0$. In this case, at the particle level, the proposal distribution is a random perturbation with law $K$ of the state of another uniformly sampled particle. Other choices of proposal distributions which produce good results in practice can be found in \cite[Section 3]{clarte_collective_2021}. 

In this example, the mean-field limit reduces the analysis of a complex particle system to the analysis of a (hopefully) simpler PDE. Another example of such idea can be found in \cite{jourdain_optimal_2015} where an algorithm similar to the one in \cite{clarte_collective_2021} is studied. The main difference is that at each time step, the proposals are accepted or rejected globally for the $N$ particles and not individually. In other words, the algorithm is a simple Metropolis-Hasting algorithm on a product space $E^N$ with a tensorized target distribution $\pi^{\otimes N}$. When $E=\R$, under a proper diffusive time-rescaling, each component of the chain (i.e. each particle) satisfies the propagation of chaos property when $N\to+\infty$. The limit nonlinear process is a McKean-Vlasov diffusion process whose law satisfies the Fokker-Planck equation 
\[\partial_t f_t = \partial_x \Big\{G(f_t)V'f_t + \Gamma(f_t)\partial_x f_t\Big\},\]
where $G$ and $\Gamma$ are explicit functions of $f_t$ and $V, V', V''$ and $V$ is a Gibbs potential such that $\pi = \e^{-V}$. The long-time convergence analysis of the solution can be found in \cite{jourdain_optimal_2014} using also the entropy-dissipation properties of the equation. 

\subsubsection*{Ensemble Kalman Sampling.} A common inverse problem is to estimate a parameter $x\in\R^d$ from a noisy observation $y\in\R^k$ which is given by: 
\[y = \mathcal{G}(x) + \eta,\]
where the \emph{model} is given by a function $\mathcal{G}:\R^d\to\R^k$ and a covariance matrix $\Gamma$ such that the noise $\eta$ is a Gaussian random variable $\mathcal{N}(0,\Gamma)$. In the present setting, the model parameters $\mathcal{G}$ and $\Gamma$ are known, usually given by the underlying physical setting, but in some cases one can also optimise them (for instance trying several noise amplitudes), in order to fit at best some experimental data. In the Bayesian framework, the parameter $x$ that we want to estimate is assumed to be a priori distributed according to a distribution $\pi_0$ on $\R^d$ (called the prior distribution). Then, after an observation $y$, our knowledge of $x$ is updated: the posterior distribution of $x$ knowing the observation $y$ is computed using Bayes' formula. For a given model $\mathcal{G}$ and $\Gamma$, the posterior distribution of $x$ is equal to 
\[\pi(x) \propto \exp(-\Phi(x))\pi_0(x),\]
where the likelihood function $\e^{-\Phi}$ is the Gibbs potential of the loss function: 
\[\Phi(x) := \frac{1}{2} |y-\mathcal{G}(x)|^2_\Gamma = \langle y-\mathcal{G}(x), \Gamma^{-1/2}(y-\mathcal{G}(u))\rangle.\]
To keep things simple, we will assume that $\pi_0$ is a centered Gaussian law with covariance matrix $\Gamma_0$. The target posterior distribution is thus (up to a normalisation constant): 
\begin{equation}\label[IIeq]{eq:targeteki}
\pi(x) \propto \exp(-\Phi_R(x)),\quad \Phi_R(x) := \Phi(x) + \frac{1}{2}|x|_{\Gamma_0}^2.\end{equation}
In order to reconstruct $x$, one can either draw samples from the posterior distribution $\pi$ or compute the points which maximise $\pi$, this method being known as the Maximum A Posteriori (MAP). The recent Ensemble Kalman Inversion (EKI) methods propose various metaheuristic diffusion interacting particle schemes to solve these sampling and optimization problems. Unlike the Metropolis-Hasting algorithm, these methods exploit the specific form of the target distribution. 

When the target distribution is in Gibbs form \cref{eq:targeteki}, a simple diffusion process with stationary distribution $\pi$, called the \emph{Langevin dynamics}, is given by the SDE: 
\[\dd X_t = -\nabla\Phi_R(X_t) \dd t + \sqrt{2}\dd B_t.\]
The law $f_t$ of $X_t$ solves the Fokker-Planck equation: 
\[\partial_t f_t = \nabla\cdot(f_t\nabla\Phi_R) + \Delta f_t.\]
Similarly to the Metropolis-Hasting case, it is possible to simply simulate a Langevin dynamics and use its ergodic properties to get samples from $\pi$. Note, however, that on a computer it is not possible to construct a time-continuous process and in practice the method thus relies on a discretization scheme which introduces a bias in the stationary distribution. For this reason, rather than being used as a direct sampling method, the discretized Langevin dynamics is more often plugged into the proposal distribution of a Metropolis-Hastings algorithm in order to correct this bias (it is then called the \emph{Metropolis Adjusted Langevin Algorithm}). Moreover, the Langevin dynamics requires to evaluate the gradient of $\Phi_R$ which can be impossible or very costly. In the present case the gradient of the potential reads: 
\begin{equation}\label[IIeq]{eq:nablaphir}
\nabla \Phi_R(x) = \nabla \mathcal{G}(x)\Gamma^{-1}(\mathcal{G}(x)-y) + \Gamma_0^{-1}u.\end{equation}

In \cite{garbuno-inigo_interacting_2020}, the authors introduce the following modified Fokker-Planck equation in order to speed up the convergence of the Langevin dynamics: 
\begin{equation}\label[IIeq]{eq:ekifp}\partial_t f_t = \nabla\cdot(f_t \mathrm{Cov}[f_t]\nabla\Phi_R) + \Tr\big(\mathrm{Cov}[f_t]\nabla^2 f_t\big),\end{equation}
where for $\mu\in\pb(E)$, $\mathrm{Cov}[\mu]$ is the covariance matrix: 
\[\mathrm{Cov}[\mu] := \int_{\R^d} \big(x-m[\mu]\big)\otimes \big(x-m[\mu]\big) \mu(\dd x),\quad m[\mu] := \int_{\R^d} x \mu(\dd x).\]
The nonlinear Fokker-Planck equation \cref{eq:ekifp} is the formal mean-field limit of the following McKean-Vlasov interacting particle system (called \emph{ensemble} in this context): 
\begin{equation}\label[IIeq]{eq:eks}\dd X^i_t = -\mathrm{Cov}\big[\mu_{\mathcal{X}^N_t}\big]\nabla\Phi_R(X^i_t)\dd t + \sqrt{2\mathrm{Cov}\big[\mu_{\mathcal{X}^N_t}\big]}\dd B^i_t,\end{equation}
for $i\in\{1,\ldots,N\}$ and where the $B^i_t$ are $N$ independent Brownian motions. In~\cite{garbuno-inigo_interacting_2020}, the system \cref{eq:eks} is called the \emph{Ensemble Kalman Sampler} (EKS) and the long-time behaviour of \cref{eq:ekifp} is studied using a gradient-flow approach. To obtain a derivative-free algorithm, the authors also use the following approximation, for $\mu\in\pb(E)$,
\begin{equation}\label[IIeq]{eq:ekilinearapproxG}\mathrm{Cov}[\mu]\nabla\mathcal{G}(x) \simeq \mathrm{Cov}[\mu,\mathcal{G}]:= \int_{\R^d} \big(x-m[\mu]\big)\otimes \big(\mathcal{G}(x)-\langle \mu,\mathcal{G}\rangle\big) \mu(\dd x).\end{equation}
Using \cref{eq:nablaphir} and the approximation \cref{eq:ekilinearapproxG} the EKS \cref{eq:eks} thus becomes derivative-free: 
\begin{equation}\label[IIeq]{eq:eksapprox}\dd X^i_t = -\mathrm{Cov}\big[\mu_{\mathcal{X}^N_t},\mathcal{G}\big]\Gamma^{-1}\big(\mathcal{G}(X^i_t)-y\big)\dd t -\mathrm{Cov}\big[\mu_{\mathcal{X}^N_t},\mathcal{G}\big]\Gamma_0^{-1}X^i_t\dd t +\sqrt{2\mathrm{Cov}\big[\mu_{\mathcal{X}^N_t}\big]}\dd B^i_t.\end{equation}
Unfortunately, the approximation \cref{eq:ekilinearapproxG} is exact only when $\mathcal{G}$ is linear and in general, the derivative-free EKS \cref{eq:eksapprox} does not converge towards the correct target distribution. In the linear case the propagation of chaos for the system \cref{eq:eks} is shown in \cite{ding_ensemble_2021}. Since the covariance matrix is a quadratic quantity, the Lipschitz assumptions of McKean's theorem do not hold. One of the methods described in Section \cref{sec:mckeantowardssingular} might be used; however the authors of \cite{ding_ensemble_2021} introduce a new bootstrapping method. The starting point is the classical synchronous coupling of Sznitman. Then, Ding and Li prove the following properties. 
\begin{enumerate}
    \item If $f_0$ has bounded moments of order $p\geq2$, then the nonlinear system and the particle system also have bounded moments of order $p\geq2$ on any finite time interval, see \cite[Lemma 5.2]{ding_ensemble_2021} and \cite[Proposition 5.4]{ding_ensemble_2021}. 
    \item Let $Y^i_t = X^i_t-\overline{X}^i_t$. The crucial property \cite[Lemma 5.4]{ding_ensemble_2021} states that if there exists $0\leq\alpha<1$ such that 
    \begin{equation}\label[IIeq]{eq:ekibootstraphyp}\E\big|Y^i_t\big|^2\leq CN^{-\alpha},\end{equation}
    then for any $\varepsilon>0$, 
    \begin{equation}\label[IIeq]{eq:ekibootstrap}
    \E\Big|Y^i_t-\frac{1}{N}\sum_{j=1}^N Y^j_t\Big|^2\leq C N^{-1/2-\alpha/2+\varepsilon}.\end{equation}
    \item Under the hypothesis \cref{eq:ekibootstraphyp} and using \cref{eq:ekibootstrap}, it is possible to prove \cite[Lemma 5.5]{ding_ensemble_2021}:
    \begin{equation}\label[IIeq]{eq:ekibootstrapend}\E\big|Y^i_t\big|^2\leq CN^{-1/2-\alpha/2+\varepsilon}.\end{equation}
    The proof is based on It\=o's formula and an explicit control of the quantity  
    \[\E\Big\|\mathrm{Cov}\big[\mu_{\mathcal{X}^N_t}\big]-\mathrm{Cov}\big[\mu_{\overline{\mathcal{X}}^N_t}\big]\Big\|^2 
    \]
    by $N^{-1}$ and $\E\Big|Y^i_t-\frac{1}{N}\sum_{j=1}^N Y^j_t\Big|^2$ (see \cite[Lemma B.2]{ding_ensemble_2021}). 
    \item From \cref{eq:ekibootstrapend} and \cref{eq:ekibootstraphyp}, by a bootstrapping argument starting from $\alpha=0$, it follows that \cref{eq:ekibootstrapend} holds with $\alpha=1-2\varepsilon$, which gives the optimal convergence rate up to $\varepsilon$.
\end{enumerate}
The proof crucially uses the linearity of $\mathcal{G}$. In \cite{ding_ensemble_2021-1}, the weakly nonlinear case where $\mathcal{G}(x) = Ax + g(x)$ for a small $g$ is investigated as well as the corresponding time-discrete algorithm. The present method as well as various other EKI methods are investigated numerically in \cite{reich_fokker-planck_2021}. A new methodology for nonlinear settings can be found in \cite{pavliotis_derivative-free_2022}. 

\subsubsection*{Filtering problems.} The two previous examples focus on a static target. Filtering can be understood as a ``dynamic sampling'' problem. An example of filtering problem which extends some of the notions that we have discussed is the famous \emph{Kalman filter}. The goal is to estimate a time-evolving signal $X_t\in\R^d$ which evolves according to the following SDE
\[\dd X_t = F(X_t)\dd t + \Sigma^{1/2}_1\dd B^1_t,\]
with known parameters $F:\R^d\to\R^d$, $\Sigma_1\in\mathcal{M}_d(\R)$ and $B^1_t$ a Brownian motion. The signal is not measured directly and it is only observed through the noisy linear transformation $Y_t\in\R^k$ defined by: 
\[\dd Y_t = GX_t + \Sigma^{1/2}_2\dd B^2_t,\]
with a known linear map $G:\R^d\to\R^k$ seen as a matrix, $\Sigma_2\in\mathcal{M}_k(\R)$ and $B^2_t$ an independent Brownian motion. The goal is to compute the conditional distribution $\pi_t$ of $X_t$ for any $t\geq0$ knowing the observed path $Y_{[0,t]}$, i.e. for any test function $\varphi\in C_b(\R^d)$, the goal is to compute: 
\[\langle \pi_t, \varphi\rangle := \E\big[\varphi(X_t)|\mathscr{F}_t\big],\]
where $\mathscr{F}_t = \sigma(Y_s, s\leq t)$. In the linear case $F(X_t) \equiv FX_t$ with $F:\R^d\to\R^d$ a linear map seen as matrix, the Bayes theorem implies that $\pi_t$ is a Gaussian law with mean $\widehat{X}_t$ and covariance matrix $P_t$ which satisfy the equations: 
\begin{align*}
\dd \widehat{X}_t &= F\widehat{X}_t\dd t + P_tG^\mathrm{T} \Sigma_2^{-1}(\dd Y_t - G\widehat{X}_t\dd t),\\
\frac{\dd}{\dd t}P_t &= FP_t + P_tF^\mathrm{T} - P_t G^\mathrm{T}\Sigma_2^{-1}G P_t + \Sigma_1.
\end{align*}
The equation on $P_t$ is a matrix-valued Riccati equation. The equation on $\widehat{X}_t$ is called the \emph{Kalman-Bucy} filter. Unfortunately, the solutions of these equations cannot be computed easily in general so an approximation method is needed. The key observation is their link with the conditional nonlinear McKean-Vlasov diffusion defined by:
\[\dd \overline{X}_t = F\overline{X}_t \dd t + \Sigma_1^{1/2}\dd W_t + \mathrm{Cov}[f_t]C^\mathrm{T}\Sigma_2^{-1}\big(\dd Y_t - G\overline{X}_t\dd t - \Sigma_2^{1/2}\dd V_t\big),\]
where $W_t, V_t$ are independent Brownian motions and $f_t = \mathrm{Law}(\overline{X}_t|\mathscr{F}_t)$. Then it can be shown that 
\[\widehat{X}_t = \E[\overline{X}_t|\mathscr{F}_t],\quad P_t = \mathrm{Cov}[f_t].\]
This readily suggests that the solutions of the Kalman-Bucy filter and the Riccati equation can be approximated by an interacting particle system, in this context called a \emph{particle filter}. The propagation of chaos thus appears as the crucial theoretical foundation of the method. The lack of Lipschitz regularity and the fact that the law is only defined conditionally to the random process $Y_t$ make things quite difficult and the result is not already covered by a theorem in the present review. Rigorous results are proved by Del Moral, Kurtzmann and Tugaut in \cite{del_moral_stability_2018} in the linear case and in \cite{del_moral_stability_2017} in the nonlinear case. The methodology of the proofs is non standard and the complexity of the model prevents us to give a faithful presentation here.   

The time continuous Kalman-Bucy filter that has been presented is one example but maybe not the most representative example of filtering problem. In practice, there are only time discrete processes, because they are part of a numerical simulation or because the signal is observed only at discrete times. A more traditional abstract filtering problem in discrete time, also called a \emph{state-space model}, is given by the two Markov chains with transition kernels: 
\[X_{k+1} \sim K(\dd x | X_{k}), \quad Y_{k+1} \sim g(\dd y | X_k),\]
The \emph{hidden Markov chain} $(X_k)_{k\in\N}$ is observed only through the \emph{observation process} $(Y_k)_{k\in\N}$ which is defined conditionally on $(X_k)_k$. The goal is to compute the conditional distribution $\pi_{k|k}$ for all $k\in\N$, defined for all $\varphi\in C_b(\R^d)$ by:
\[\langle \pi_{k|k}, \varphi\rangle = \E[\varphi(X_k)|Y_{0:k}],\]
where $Y_{0:k} = (Y_0,\ldots, Y_k)$. Bayes theorem gives the recursion formula: 
\[\pi_{k|k}(\dd x) \propto g(Y_k|x)\pi_{k|k-1}(\dd x),\quad\pi_{k|k-1}(\dd x) = \int_{\R^d} \pi_{k-1|k-1}(\dd z)K(\dd x|z).\]
In general it is not possible to obtain the expression of $\pi_{k|k}$ in closed form. For this reason, the class of \emph{Sequential Monte Carlo} (SMC) methods, also known as \emph{particle filters} aim at approximating it by an empirical measure $\mu_{\mathcal{X}^N_k}$ where $\mathcal{X}^N_k$ is understood as a time evolving particle system. Most often, the SMC methods rather rely on a weighted empirical measure, where the weights of the particles are obtained using an \emph{importance sampling} method. The convergence of the approximating empirical measure or of the importance weights is naturally related to propagation of chaos. The connection between the two domains is due to Del Moral \cite{del_moral_measure-valued_1998} at the end of the 90's. Since then, SMC methods have become increasingly popular with real-world applications in engineering, signal processing and more recently in machine learning to cite a few. For further details, we refer the interested reader to the short surveys \cite{kantas_overview_2009, crisan_survey_2002} for a practical introduction to the subject and to the larger monographs \cite{del_moral_feynman-kac_2004, del_moral_mean_2013} and \cite{doucet_sequential_2001} for the theoretical foundations, in particular the links with mean-field theory. 

\subsubsection{Agent Based Optimization}\label[II]{sec:optimization}

In its most abstract form, an optimization problem consists in finding the point $x_\star\in E \subset \R^d$, assumed to be unique, which minimizes a given function $G:E\to\R_+$. The problem is notoriously difficult in high dimensional spaces or when $G$ has many local minima. In the 90's, Kennedy and Eberhart \cite{kennedy_particle_1995}
introduced a class of optimization algorithms based on a swarm of interacting agents. The \emph{Particle Swarm Optimization} (PSO) methods are inspired by biological concepts: each agent (or particle) follows a set of simple rules which is a mix between an individual exploration behaviour of the state space and a collective exploitation of the swarm knowledge in order to efficiently find and converge to the global minimum of $G$. From an algorithmic point of view, the algorithm is appealing by its (relative) simplicity and its versatility as it does not requires expensive computations like the gradient of $G$. In the last decades, many variants and practical implementations of the original PSO algorithm have been proposed and a full inventory of these \emph{Swarm Intelligence} (SI) methods would go beyond the present review. Although these algorithms have proved their efficiency for notoriously difficult problems, their main drawback is their lack of theoretical mathematical foundations. Most of the SI methods are based on metaheuristic principles which can hardly be turned into rigorous convergence results, in particular when the number of agents involved becomes large. Lately, there has been a growing interest for the convergence analysis of SI methods using the tools developed in the kinetic theory community for mean-field particle systems in Physics or Biology. At this point in the present review, it becomes blatantly clear that a rigorous mean-field interpretation of SI methods could be of primary interest as it reduces the difficult analysis of a many particle system into the analysis of a single PDE for which many tools are already available to study its long-time convergence properties. 

Following these ideas, a very simple though quite efficient method has recently been introduced by Pinneau et al. \cite{pinnau_consensus-based_2017}. This method called \emph{Consensus Based Optimization} (CBO) is based on the following McKean-Vlasov particle system:
\begin{equation}\label[IIeq]{eq:cbopinneau}
    \dd X^i_t = -\lambda\Big(X^i_t - v\big[\mu_{\mathcal{X}^N_t}\big]\Big)H^\varepsilon\Big(G(X^i_t)- G\big(v\big[\mu_{\mathcal{X}^N_t}\big]\big)\Big)\dd t+ \sqrt{2}\sigma\big|X^i_t - v\big[\mu_{\mathcal{X}^N_t}\big]\big|\dd B^i_t,
\end{equation}
where $\lambda>0$, $\sigma\geq0$, $H^\varepsilon$ is a smoothened version of the Heaviside function $H(u)=\1_{u\geq0}$ and given $\mu\in\pb(\R^d)$, 
\[v[\mu]:= \frac{1}{\langle \mu, \omega^\alpha\rangle}\int_{\R^d} x\omega^\alpha(x)\mu(\dd x),\quad \omega^\alpha(x) := \exp(-\alpha G(x)), \quad\alpha>0.\] 
The quantity $v[\mu_{\mathcal{X}^N_t}]$ is a weighted average of the positions of the particles. Particles which are located near a minimum of $G$ have a larger weight. The drift term is thus an \emph{exploitation} term: it is a standard gradient relaxation (for a quadratic potential) towards the current weighted average position of the swarm. The diffusion term is an \emph{exploration} term which becomes as large as the particle is far from the current weighted average. To better understand the particular choice of the weight $\omega^\alpha$, recall the Laplace principle: it states that if a probability measure $f$ is absolutely continuous with respect to the Lebesgue measure and if $x_\star$ belongs to the support of $f$, then
\[\lim_{\alpha\to+\infty} {\left(-\frac{1}{\alpha}\log\langle f,\omega^\alpha\rangle\right)} = G(x_\star).\]
When applied to the mean-field limit solution of the Fokker-Planck equation:
\begin{equation}\label[IIeq]{eq:cbopde}\partial_t f_t(x) = -\lambda\nabla\cdot\Big(\big(x-v[f_t]\big)H^\varepsilon\big(G(x)-G(v[f_t])\big)f_t\Big) + \sigma^2\Delta\big(|x-v[f_t]|^2 f_t\big),\end{equation}
this result indicates that the Gibbs-like measure $\omega^\alpha f_t/\langle f_t,\omega^\alpha\rangle$ is close to $\delta_{x_\star}$ and the weighted average of the particles is thus expected to satisfy
\[v\big[\mu_{\mathcal{X}^N_t}\big]\underset{N\to+\infty}{\longrightarrow}v[f_t]\simeq x_\star.\]
Using this heuristics, the deterministic term in \cref{eq:cbopinneau} drives the dynamics of the particle system towards the current consensus point $v\big[\mu_{\mathcal{X}^N_t}\big]$ which is always close to $x_\star$ and which keeps concentrating as particles get closer to it since the noise amplitude is lower when the particles are close to $x_{\star}$. This informal reasoning at the particle level can be made rigorous for the limit equation \cref{eq:cbopde} whose solution is expected to be a good approximation of the one-particle distribution as $N\to+\infty$. In particular, it is possible to prove that a consensus is attained in the sense that $f_t\to\delta_{x_\star}$ as $t\to+\infty$.


The analytical study of the PDE \cref{eq:cbopde} and in particular the proof that a consensus is attained can be found in \cite{carrillo_analytical_2018}. However, the rigorous propagation of chaos result, which would be necessary to conclude that the particle system converges towards $\delta_{x_\star}$, remains open in the general case. A rigorous result is available in \cite{fornasier_consensus-based_2020} in the constrained case where $G$ is minimized over a compact submanifold of $\R^d$. The proof follows the classical Sznitman coupling approach. A crucial ingredient \cite[Lemma 3.1]{fornasier_consensus-based_2020} is the bound: 
\[\E \Big|v\big[\mu_{\overline{\mathcal{X}}^N_t}\big]-v[f_t]\Big|^2\leq CN^{-1},\]
where the system $\overline{\mathcal{X}}^N_t$ is i.i.d. with law $f_t$. Note that this bound is actually a large deviation estimate. 

Further developments on the CBO method can be found in \cite{carrillo_consensus-based_2021} where a modification of the diffusion coefficient is introduced in order to obtain dimension free convergence results. A review and a comparison of recent SI methods, including the CBO method and the original PSO algorithm, can be found in \cite{totzeck_trends_2021} and a numerical comparison can be found in the short note \cite{totzeck_numerical_2018}. We also quote the recent article~\cite{grassi_particle_2021} which gives a more unifying framework for the mean-field interpretation of PSO and CBO methods. In particular, a time-continuous mean-field interpretation of the original PSO algorithm is introduced which, unlike \cref{eq:cbopinneau}, is based on a \emph{kinetic} McKean-Vlasov diffusion system: 
\begin{align*}
    \dd X^i_t &= V^i_t \dd t,\\
    \dd V^i_t &= -\gamma V^i_t\dd t + \lambda_1(Y^i_t-X^i_t)\dd t + \lambda_2(Y^\mathrm{min}_t - X^i_t)\dd t \\
    &\qquad + \sigma_1 \diag(Y^i_t - X^i_t)\dd B^{1,i}_t + \sigma_2\diag(Y^\mathrm{min}_t - X^i_t)\dd B^{2,i}_t,\\ 
    \dd Y^i_t &= \nu(X^i_t - Y^i_t)\1_{G(X^i_t)\leq G(Y^i_t)}\dd t,\\
    Y^\mathrm{min}_t &= \argmin \big\{G(Y^1_t),\ldots, G(Y^N_t)\big\}, 
\end{align*}
where $(X^i_t,V^i_t)$ is the couple position-velocity, $Y^i_t$ is the best position of particle $i$ and $Y^\mathrm{min}_t$ is the best position of the whole system. The evolution of the velocity is a combination of a (technical) friction force, two drift forces towards the best positions $Y^i_t$ and $Y^\mathrm{min}_t$ and two noise terms with a norm which depends on the distance to the best positions. 

\subsubsection{Overparametrized Neural Networks}\label[II]{sec:neuralnetworks}

Training neural networks can be understood as an optimization task. Should the commonly used algorithms converge to the good optimum is in many cases still an open question. Recent independent works \cite{mei_mean_2018, rotskoff_trainability_2019, sirignano_mean_2020, chizat_global_2018} have shown that the training process of neural networks possesses a natural mean-field interpretation which gives new insights towards a rigorous theoretical justification to this convergence problem. 

For $k\in\N$, let $(X_k,Y_k)\in\R^p\times\R$ be a sequence of i.i.d. $\pi$-distributed random variables called the \emph{training data set}, where $\pi\in\pb(\R^p\times\R)$ is an unknown distribution. The random variable $X_k$ is an object (e.g. an image) and $Y_k$ is its label. A (single hidden layer) neural network composed of $N$ neurons is characterised by $N$ parameters $\boldsymbol{\theta}^N = (\theta^1,\ldots,\theta^N)\in (\R^d)^N$. The training task of the neural network consists in finding the parameters which minimize the risk functional: 
\[R^N(\boldsymbol{\theta}^N) := \E_{X,Y\sim\pi}\Big[\ell\big(Y,\widehat{y}(X,\boldsymbol{\theta}^N)\big)\Big],\] 
where for a data $x\in\R^p$, the predicted label $\widehat{y}$ is of the form: 
\[\widehat{y}(x,\boldsymbol{\theta}^N):= \big\langle \mu_{\boldsymbol{\theta}^N},\sigma(x,\cdot)\big\rangle.\]
The function $\sigma:\R^p\times\R^d\to\R$ is a given function called the \emph{activation function}. The \emph{loss function} $\ell:\R\times\R\to\R_+$ is taken equal to $\ell(y,\widehat{y}):=|y-\widehat{y}|^2$. Note that the risk functional depends only on the empirical measure so it can actually be rewritten $R^N(\boldsymbol{\theta}^N) = R(\mu_{\boldsymbol{\theta}^N})$, where the risk functional $R$ is defined on the whole set $\pb(\R^d)$ by 
\[\forall \mu\in \pb(\R^d),\quad R(\mu) := \iint_{\R^p\times\R} \ell\big(y,\langle \mu,\sigma(x,\cdot)\rangle\big) \pi(\dd x,\dd y).\]
Since the distribution $\pi$ is unknown, the parameters of the neural network are updated sequentially each time a new $\pi$-distributed data pair object-label is given. The most common updating rule is the (noisy) Stochastic Gradient Descent (SGD), which updates each parameter $i\in\{1,\ldots,N\}$ at iteration $k$ by following the gradient of the risk functional:
\begin{equation}\label[IIeq]{eq:sgd}\theta^i_{k+1} = \theta_k^i + 2s_k\big(Y_k - \widehat{y}(X_k,\boldsymbol{\theta}^N_k)\big)\nabla_{\theta}\sigma(X_k,\theta^i_k) +\sqrt{\frac{2s_k}{\beta}}W^i_k.\end{equation}
where $s_k\in\R_+$ is a step size, $\beta\in(0,+\infty]$ and $W^i_k$ are independent standard Gaussian random variables. In the noisy case $\beta<+\infty$, it is customary to add a confinement potential to the risk functional in order to ensure good convergence properties. We do not add it here to keep the presentation as light as possible. The whole point is to interpret \cref{eq:sgd} as the time discretization of a McKean-Vlasov particle system, where the particles are the parameters of the neural network $\theta^i_k$. Since the $(X_k,Y_k)$ are assumed to be i.i.d., the CLT suggests the approximation: 
\begin{align}
    -2\big(Y_k - \widehat{y}(X_k,\boldsymbol{\theta}^N_k)\big)\nabla_{\theta}\sigma(X_k,\theta^i_k) &= N\nabla_{\theta^i} \ell\big(Y_k,\widehat{y}(X_k,\boldsymbol{\theta}^N_k)\big)\nonumber \\
    &\simeq N\nabla_{\theta^i} R^N(\boldsymbol{\theta}^N_k) + \Sigma^{1/2}\big(\theta^i_k,\mu_{\boldsymbol{\theta}^N_k}\big)\widetilde{W}^i_k,\label[IIeq]{eq:sgdclt}
\end{align}
where $\widetilde{W}^i_k$ is a standard $d$-dimensional Gaussian random variable and the covariance matrix is defined by: 
\begin{align*}\Sigma\big(\theta^i_k,\mu_{\boldsymbol{\theta}^N_k}\big) &:=  N^2\E_{X,Y\sim\pi}\Big[\nabla_{\theta^i}\ell\big(Y,\widehat{y}(X,\boldsymbol{\theta}^N_k)\big)\nabla_{\theta^i}\ell\big(Y,\widehat{y}(X,\boldsymbol{\theta}^N_k)\big)^\mathrm{T}\Big]\\
&= \E_{X,Y\sim\pi}{\left[\big|\partial_{\widehat{y}}\ell\big(Y,\widehat{y}(X,\boldsymbol{\theta}^N_k)\big)\big|^2\nabla_{\theta}\sigma(X,\theta^i_k)\nabla_{\theta}\sigma(X,\theta^i_k)^\mathrm{T}\right]}.
\end{align*}
Since $R^N$ is actually a function of the empirical measure, the SGD dynamics \cref{eq:sgd} can be rewritten with our usual notations: 
\begin{equation}\label[IIeq]{eq:sgdmckeandiscrete}
    \theta^i_{k+1} = \theta_k^i + s_k b\big(\theta^i,\mu_{\boldsymbol{\theta}^N_k}\big) + \sqrt{s_k}\sigma_k\big(\theta^i_k,\boldsymbol{\theta}^N_k\big) G^i_k,
\end{equation}
where $G^i_k$ is a standard Gaussian random variable, 
\[\sigma_k\big(\theta^i_k,\boldsymbol{\theta}^N_k\big) := {\left(s_k\Sigma\big(\theta^i_k,\mu_{\boldsymbol{\theta}^N_k}\big) + \frac{2}{\beta }I_d\right)}^{1/2},\] 
and 
\[b\big(\theta^i,\mu_{\boldsymbol{\theta}^N_k}\big) = -N\nabla_{\theta^i} R^N(\boldsymbol{\theta}^N_k) = -\E_{X,Y\sim\pi}{\left[\partial_{\widehat{y}}\ell\big(Y,\widehat{y}(X,\boldsymbol{\theta}^N_k)\big)\nabla_{\theta}\sigma(X,\theta^i_k)\right]}.\]
Finally, taking a time-step $s_k = \varepsilon\xi(\varepsilon k)$ for $\xi$ a smooth function and $\varepsilon>0$ small, the Equation \cref{eq:sgdmckeandiscrete} becomes the standard Euler-Maruyama discretization of the (time inhomogeneous) McKean-Vlasov particle system:
\begin{equation}\label[IIeq]{eq:mckeanparticlesgd}\dd \theta^i_t = \xi(t)b\big(\theta^i,\mu_{\boldsymbol{\theta}^N_k}\big)\dd t + \sqrt{\frac{2\xi(t)}{\beta}}\dd B^i_t\Big.\end{equation}
The main difference with \cref{eq:mckeanvlasov_summary} is the time dependent coefficient $\xi(t)$ but it does not affect the argument of most of the techniques investigated in Section \cref{sec:mckeanreview}. In particular, the propagation of chaos results implies that in the limit $N\to+\infty$ and $\varepsilon\to0$ the distribution $f_t$ of the neurons satisfies the Fokker-Planck equation:
\begin{equation}\label[IIeq]{eq:sgdfp}\partial_t f_t(\theta) = -\xi(t)\nabla_\theta \cdot (b(\theta,f_t)f_t) + \xi(t)\Delta_\theta f_t.\end{equation}
This informal derivation is made rigorous in the following works. 
\begin{enumerate}
    \item In \cite{mei_mean_2018}, the authors prove the simultaneous double limit $N\to+\infty$ and $\varepsilon\to0$ from the rescaled empirical measure $\mu_{\boldsymbol{\theta}^N_{\lfloor t/\varepsilon\rfloor}}$ of the discrete SGD \cref{eq:sgd} to the time-continuous solution of \cref{eq:sgdfp}, without directly using the approximating time-continuous particle system \cref{eq:mckeanparticlesgd}. The key estimate \cite[Lemma 7.2 and Lemma 7.6]{mei_mean_2018} is a concentration inequality which controls the discrepancy between the rescaled SGD and a synchronously coupled system of nonlinear McKean-Vlasov diffusion processes. The Azuma-Hoffding inequality gives a quantitative bound for the analogous of the approximation \cref{eq:sgdclt} in this case. In this time-discrete framework, the synchronous coupling is obtained by taking the Gaussian random variables in \cref{eq:sgd} equal to the integral of the Brownian motion of the coupled McKean-Vlasov diffusion on each time step. The parameter $\varepsilon\equiv \varepsilon_N$ in the time step is linked to $N$: it can be taken equal to any inverse power $\varepsilon_N= N^{-\gamma}$, $\gamma>0$. A very similar coupling approach is used in \cite{de_bortoli_quantitative_2020} with the difference that the authors prove the propagation of chaos for the time-continuous particle system \cref{eq:mckeanparticlesgd} only. In the regime where the next order approximation in \cref{eq:sgdclt} is kept, the final diffusion matrix depends on $\Sigma$. Both works are based on the global Lipschitz and boundedness assumptions of McKean's Theorem \cref{thm:mckean}.
    \item In \cite{sirignano_mean_2020}, the authors use a compactness argument with ad hoc estimates to prove the convergence of the rescaled empirical measure of the SGD, without using the time continuous approximation \cref{eq:mckeanparticlesgd}. The proof is non quantitative and is written in the case $\beta=+\infty$ but it can accommodate more singular cases, without global Lipschitz assumptions but with the assumption of bounded moments for $\pi$ and the initial distribution. 
    \item In \cite{chizat_global_2018}, the authors solve a more general problem: using the fact that the functional $R^N(\boldsymbol{\theta}^N)\equiv R(\mu_{\boldsymbol{\theta}^N})$ defines a gradient flow on $(\R^d)^N$, they prove that as $N\to+\infty$ the empirical measure of this gradient flow converges towards the Wasserstein gradient flow defined by the risk functional $R$ on $\pb_2(\R^d)$. The proof is quite similar in spirit to what has been presented in Section \cref{sec:gradientflows} (i.e. a compactness argument for curves using Ascoli's theorem) but it is relatively simpler in this case because the framework is entirely deterministic (in particular, the empirical measure is a deterministic object).
\end{enumerate}
Of course, proving the propagation of chaos is only a first step (and in a sense the easiest one) towards the rigorous analysis of the optimization problem outlined above. As illustrated many times in this section, the goal is now to exploit the long-time convergence properties of the limit Fokker-Planck equation \cref{eq:sgdfp}. When $\ell(y,\widehat{y})=|y-\widehat{y}|^2$ a key observation is that this equation has a gradient flow structure. Using the fact that in this case: 
\[R(\mu) = R_0 + 2\int_{\R^d} V(\theta)\mu(\dd\theta) + \iint_{\R^d\times\R^d} W(\theta,\theta') \mu(\dd\theta)\mu(\dd\theta'),\]
where $R_0 = \E_{X,Y\sim\pi}[Y^2]$ and defining the potentials
\[V(\theta) := -\E_{X,Y\sim\pi}[Y\sigma(X,\theta)],\quad W(\theta,\theta') := \E_{X,Y\sim\pi}[\sigma(X,\theta)\sigma(X,\theta')],\]
then for $\theta\in\R^d$ and $\mu\in\pb(\R^d)$, the drift function is equal to
\[b(\theta,\mu) = -\nabla_\theta \frac{\delta R(\mu)}{\delta \mu}(\theta),\]
so that \cref{eq:sgdfp} is an evolutionary PDE in the sense of Definition \cref{def:evolutionarypde} and thus a gradient flow. This gradient flow structure is exploited in \cite{chizat_global_2018} and \cite{mei_mean_2018} to prove the long-time convergence of the SGD \cref{eq:sgd} and of the solution of \cref{eq:sgdfp} towards a global minimizer of $R$.

\subsection{Beyond propagation of chaos}\label[II]{sec:beyondpoc}

In this last section, we give a glimpse on some results which extend or complete the question of propagation of chaos. We discuss two natural directions: the fluctuation theory when the propagation of chaos property holds (Section \cref{sec:fluctuations}) and another other type of many-particle limit when the propagation of chaos does \emph{not} hold (Section \cref{sec:flemingviot}).

\subsubsection{Fluctuations}\label[II]{sec:fluctuations}

Propagation of chaos can be interpreted as a kind of law of large numbers where the empirical process $\mu_{\mathcal{X}^N_t}$ converges towards the deterministic limit $f_t$. The next stage is to consider the asymptotic behaviour when $N\to+\infty$ of the fluctuation process
\begin{equation}\label[IIeq]{eq:fluctuationprocess}\eta^N_t := \sqrt{N}{\left(\mu_{\mathcal{X}^N_t}-f_t\right)},
\end{equation}
thus giving a form of Central Limit Theorem. The first problem is to identify a suitable space to which $\eta^N_t$ and its (potential) limit belong. From its definition, $\eta^N_t$ belongs to the space of signed measures. It may not be the case for the limit and as we shall see, the ``good'' point of view is to look at $\eta^N_t$ as an element of a space of distributions. In this subsection, we denote by $\mathscr{H}'$ this space, defined as the dual of a space $\mathscr{H}$ of test functions. Then, the second problem is to identify and characterise the limit as a process in $\mathscr{H}'$. A choice is to study the limit of the finite dimensional distributions 
\begin{equation}\label[IIeq]{eq:fddfluctuations}{\left(\langle \eta^N_{t_1}, \varphi^1\rangle, \langle \eta^N_{t_2}, \varphi^2\rangle,\ldots, \langle \eta^N_{t_k}, \varphi^k\rangle\right)}\in \R^k, \end{equation}
for $\varphi^1,\ldots,\varphi^k\in \mathscr{H}$ and $(t_1,\ldots, t_k)\in \R_+^k$. If the limit exists, the finite dimensional distributions characterise a process $(\eta_t)_t\in C([0,T],\mathscr{H}')$. Another approach is the study of the asymptotic behaviour of the pathwise characteristic function 
\begin{equation}\label[IIeq]{eq:characteristicfunctionfluctuations}\E{\left[\e^{i\langle \eta^N_{[0,T]},\Phi\rangle}\right]},\end{equation}
where $\Phi$ is a test function on a suitable subset of the path space. The final step is to find the SDE (in $\mathscr{H}'$) which governs the evolution of the limit process $\eta_t$. The expected behaviour is a kind of infinite dimensional Ornstein-Uhlenbeck process. 

In the following, we briefly review the main results for the two classes of models studied before, the Boltzmann models and the mean-field McKean models. 

\subsubsection*{Boltzmann models.} The study of fluctuations for Boltzmann models has been initiated by Kac and McKean for McKean's 2-speed caricature of a Maxwellian gas~\cite{cohen_probabilistic_1973,mckean_fluctuations_1975}. The case of the three dimensional hard-sphere gas is also discussed in~\cite{mckean_fluctuations_1975}, within the framework of~\cite{grunbaum_propagation_1971}. The one-dimensional Kac model (Example~\cref{example:kacmodels}) is studied by Tanaka  \cite{tanaka_fluctuation_1982} in the equilibrium case and by Uchiyama \cite{uchiyama_fluctuations_1983} in the non equilibrium case. This last work is based on the following chain of arguments. 
\begin{enumerate}[(1)]
    \item Using the generator of the particle system, identify formally the limit generator of the real-valued process $h(\langle \eta^N_t,\varphi\rangle)$, where $h\in C^\infty_c(\R)$ and $\varphi\in\mathscr{S}(\R)$ belongs to the Schwartz space of functions rapidly decaying at infinity and $\eta^N_t\in \mathscr{S}'(\R)$ is seen as a tempered distribution. 
    \item Show that the sequence of laws of the processes $(\eta^N_t)_N$ is tight in the space $\pb(D([0,T],\mathscr{S}_\delta'(\R)))$ where $\mathscr{S}_\delta'(\R)\subset\mathscr{S}'(\R)$ is a subset of the space of tempered distributions. Check that any limit point concentrates on $\pb(C([0,T],\mathscr{S}_\delta'(\R)))$. 
    \item Identify any limit point as the solution of a martingale problem using the expression derived in the first step. 
\end{enumerate}
This method is then applied to a more realistic three-dimensional (cutoff) model in \cite{uchiyama_fluctuation_1983} (see also \cite{tanaka_probabilistic_1983}). The method of Uchiyama is extended to more general Boltzmann models in \cite{ferland_compactness_1992}. The limit of the characteristic functions \cref{eq:characteristicfunctionfluctuations} is studied for a Boltzmann model with simultaneous jumps in a countable state space in \cite{uchiyama_fluctuations_1988}.

\subsubsection*{Mean-field models.} The fluctuations of the simple one-dimensional McKean-Vlasov diffusion $b(x,\mu)=-\lambda x$, $\lambda>0$ and $\sigma=I_d$ are studied in \cite{tanaka_central_1981}. The starting point is the proof that for the pathwise version of \cref{eq:characteristicfunctionfluctuations} has an explicit limit: 
\begin{equation}\label[IIeq]{eq:characfunctionpathwisefluctu}
\lim_{N\to+\infty} \E{\left[\e^{i\langle \eta^N_{[0,T]},\xi\rangle}\right]}=\e^{-Q(\xi)/2},\end{equation}
for an explicit functional $Q(\xi)$, where $\xi$ belongs to a subspace of the space of test functions on $C([0,T],\R)$ built using the finite-dimensional polynomial functions: 
\[\xi(\omega) = \varphi_k(\omega(t_1),\ldots,\omega(t_k)),\]
where $\varphi_k$ is a polynomial. Then a SDE which governs the ``gaussian random field'' with characteristic function \cref{eq:characfunctionpathwisefluctu} is obtained in an appropriate space of distributions. The general linear case with $b(x,\mu)=\tilde{b}\star\mu(x)$ is investigated in \cite{tanaka_limit_1982}. Tanaka uses a method originally due to Braun and Hepp in a deterministic case which consists in studying the (pathwise) ``fluctuation field'' : 
\[\big\langle \eta^{N}_{[0,T]},\xi\big\rangle = \sqrt{N}{\left(\frac{1}{N}\sum_{i=1}^N \xi{\left(X^{i,N}_{[0,T]}\right)}-\big\langle f_{[0,T]}, \xi\big\rangle\right)},\]
where $f_{[0,T]}\in \pb(C([0,T],\R^d))$ and $\xi$ is a smooth function on the path space for a specific notion of differentiability. The idea is to write $X^{i,N}_{[0,T]}$ as the flow of a SDE which depends on $\mu_{\mathcal{X}^N_{[0,T]}}$. Then, under smoothness assumptions, a theorem due to Braun and Hepp which is generalised in \cite{tanaka_limit_1982} implies the convergence of the finite dimensional distributions \cref{eq:fddfluctuations} and/or of the (pathwise) characteristic function. A large deviation principle with an explicit rate function $I$ is also obtained.

The differentiability assumptions of \cite{tanaka_limit_1982} are weakened by Sznitman in \cite{sznitman_nonlinear_1984} using a Girsanov transform argument. The result is valid in $\R^d$ and in a bounded domain with reflecting boundary conditions. The method of Sznitman is employed in \cite{shiga_central_1985} for a mean-field jump process.

Following the ideas and results of Sznitman, Hitsuda and Mitoma \cite{hitsuda_tightness_1986} prove the tightness of the fluctuation process in a space of distributions (using a trajectorial representation and a synchronous coupling argument) and derive a SDE for the limit. The model is studied in dimension one only. The result is improved in \cite{fernandez_hilbertian_1997} where the authors identify a minimal (in a certain sense) space of distributions for the fluctuation process (a weighted negative Sobolev space). This approach is then carried out for a moderate interaction model in \cite{jourdain_propagation_1998} and for a very general jump-diffusion model in \cite{meleard_convergence_1998}. A detailed presentation can be found in M\'el\'eard's course~\cite[Section 5]{meleard_asymptotic_1996}.

\subsubsection{Measure-valued limits: an example}\label[II]{sec:flemingviot}

As explained many times in this review, the propagation of chaos property is equivalent to the convergence of the empirical process towards a deterministic limit. It means that the law of the limit is a Dirac delta. In some cases, propagation of chaos does not hold but the empirical process still has a limit when $N\to+\infty$. This limit is thus a (random) measure-valued process with a law which is not a Dirac delta. A classical reference on measure-valued processes is Dawson's course \cite{dawson_measure-valued_1993}. 

To give a flavour of the subject, let us give a semi-informal derivation of the most important measure-valued process, the famous \emph{Fleming-Viot process}, starting from the toy example of Section \cref{sec:toyexampleempiricalprocess}. We recall here its construction and highlight the differences which lead from the propagation of chaos to a measure valued limit. A similar presentation can be found in Dawson's course \cite{dawson_measure-valued_1993}. 
\begin{itemize}
    \item We assume that $E$ is compact, say $E=\mathbb{T}^d$ the torus in dimension $d$ and that the motion is a pure jump process, without deterministic drift (for simplicity).
    \item Instead of a constant jump rate $\lambda\equiv 1$, we speed up the process and take a jump rate $\lambda_N=N$ which depends on the number of particles. Compared to our usual setting in the Boltzmann case, it means that each \emph{pair} of particles update its state in average $\mathcal{O}(1)$ times during one unit of time. To prove the propagation of chaos, we assumed that each \emph{particle} updates its state in average $\mathcal{O}(1)$ times during one unit of time.
    \item The jump is still sampled from a linear jump transition measure: for $\mu\in\pb(E)$ and $x\in E$, $P_\mu(x,\dd y) = K_N\star\mu(\dd y)$ where $K_N:E\to E$ is a symmetric kernel. We assume this time that $K_N$ is a smooth mollifier when $N \to +\infty$, in the sense that there is $\sigma > 0$ such that for all $x \in E$ and all smooth $\phi$ on $E$, we have the Taylor expansion: 
    \begin{equation}\label[IIeq]{eq:taylorKN}\int_{E} \varphi(y)K_N(y-x)\dd y = \varphi(x) + \frac{\sigma}{N}\Delta\varphi(x) + \mathcal{O}{\left(\frac{1}{N^2}\right)},\end{equation}
\end{itemize}

With these modifications, the empirical process $(\mu_{\mathcal{X}^N_t})_t$ is a measure-valued Markov process with generator
\begin{equation}\label[IIeq]{eq:generatorFVN}\widehat{\mathcal{L}}_N\Phi(\mu) = N^2\iint_{E\times E} {\left\{\Phi{\left(\mu-\frac{1}{N}\delta_x+\frac{1}{N}\delta_y \right)}-\Phi(\mu)\right\}}(K_N\star\mu)(\dd y)\mu(\dd x),\end{equation}
where we assume that the test function $\Phi\in C_b(\pb(E))$ is a polynomial function 
\[\Phi(\mu) = \langle \mu^{\otimes k}, \varphi_k\rangle,\]
with $k\in\N$ and $\varphi_k\in C_b(E^k)$. We recall that since $E$ is compact, the set of polynomial functions on $\pb(E)$ is dense in $C_b(\pb(E))$. Note that a polynomial function can be extended to the space of signed measures. Following Dawson's course \cite{dawson_measure-valued_1993}, the first order derivative of the polynomial $\Phi$ (seen as a function on signed measures) at $\mu\in\pb(E)$ is defined as the function on $E$ :  
\begin{align}
\frac{\delta \Phi(\mu)}{\delta\mu} : x\in E\mapsto & \lim_{\varepsilon\to0} \frac{\Phi(\mu+\varepsilon\delta_x)-\Phi(\mu)}{\varepsilon}\nonumber \\
&\quad = \sum_{j=1}^k \int_{E^{k-1}}\varphi_k(x^1,\ldots,x^{j-1},x,x^{j+1},\ldots,x^k)\prod_{\ell\ne j}\mu\big(\dd x^\ell\big) \in\R,\label[IIeq]{eq:deltaphideltamu}\end{align}
and similarly, 
\begin{equation}\label[IIeq]{eq:deltaphideltamu2}
\frac{\delta \Phi^2(\mu)}{\delta^2\mu}:(x,y)\in E^2 \mapsto  \frac{\partial^2}{\partial\varepsilon_1\partial\varepsilon_2}\Phi(\mu+\varepsilon_1\delta_x+\varepsilon_2\delta y)\Big|_{\varepsilon_1=\varepsilon_2=0}\in\R.\end{equation}
Similarly to what we have presented in Section \cref{sec:limitgenerator}, the goal is to write an expansion of the generator \cref{eq:generatorFVN} as $N\to+\infty$. This time we work on the space of polynomials and we use the notion differentiability defined above. Reporting \cref{eq:taylorKN}, \cref{eq:deltaphideltamu} and \cref{eq:deltaphideltamu2} into \cref{eq:generatorFVN}, a direct computation gives the expansion as $N\to+\infty$ : 
\[\widehat{\mathcal{L}}_N\Phi(\mu) = \mathcal{L}_{\mathrm{FV}}\Phi(\mu) + R_N,\]
where $|R_N|=\mathcal{O}(1/N)$ and $\mathcal{L}_{\mathrm{FV}}$ is the Fleming-Viot generator defined by:
\begin{equation}\label[IIeq]{eq:generatorFV}\mathcal{L}_{\mathrm{FV}}\Phi(\mu) := \sigma\int_E\Delta{\left(\frac{\delta\Phi(\mu)}{\delta\mu}\right)}(x)\mu(\dd x) + \iint_{E\times E} \frac{\delta \Phi^2(\mu)}{\delta^2\mu}(x,y) Q_\mu(\dd x, \dd y),\end{equation}
where $Q_\mu(\dd x, \dd y) := \mu(\dd x)\otimes\delta_x(\dd y)-\mu(\dd x)\otimes\mu(\dd y)$. It can be proved that the Fleming-Viot generator \cref{eq:generatorFV} defines a $\pb(E)$-valued Markov process, called the Fleming-Viot process, which can also be characterised using the various points of view developed in the previous sections: the convergence of the $N$-particle semi-group, the infinite system of moment measures, the solution of a martingale problem. Everything is well detailed in Dawson's course \cite[Sections 2.5 to 2.9]{dawson_measure-valued_1993} in the slightly different situation where $K_N=\delta_0$ but the particles are subject to a Brownian noise between the jumps. The properties of the Fleming-Viot process are studied in the reference articles \cite{dawson_wandering_1982,donnelly_countable_1996}. 

In population dynamics, the space $E$ is the space of types (or alleles) and each jump is interpreted as the simultaneous death of an individual and the birth of a new individual with a type sampled uniformly among the population with a mutation given by $K_N$. The particle model is called the Moran model. The state space is often a discrete space. Historically, Fleming and Viot \cite{fleming_measure-valued_1979} derived the measure-valued limit using a suitable discretatization of the continuous state space and taking the limit in a martingale problem when both $N\to+\infty$ and the discretization step goes to zero. Alternatively to the Moran particle process, the Fleming-Viot process is also the measure-valued limit of the famous Wright-Fisher model. The main difference with the Moran process is that all the $N$ particles update their state at the same time. For an introduction to the limit $N\to+\infty$ in this case using martingale arguments, see \cite[Chapter 10, Section 4]{ethier_markov_1986} and the references therein. Finally, the lectures \cite{etheridge_introduction_2000} and \cite{etheridge_mathematical_2011} contain more recent references on the subject as well as many applications in mathematical biology. 

\section*{Acknowledgments} The authors wish to thank Pierre Degond for his precious advice and careful proofreading of this manuscript. The authors also thank Paul Thevenin for fruitful comments and discussions. The work of AD is supported by an EPSRC-Roth scholarship co-funded by the Engineering and Physical Sciences Research Council and the Department of Mathematics at Imperial College London. Finally, the authors also wish to thank the anonymous reviewers for their careful proofreading and useful comments and hints. 

\appendix

\section{Generator estimates against monomials}\label[II]{appendix:generatorestimates}

Generators estimates are required in particular in Section \cref{sec:limitsemigroup} and for compactness methods in Section \cref{sec:martingalecompactness} and Section \cref{sec:martingaleboltzmannreview}. For some polynomials $\mu \mapsto \langle\mu^{\otimes k},\varphi_k\rangle$, the purpose is to compare the generator $\mathcal{L}_N [ R_{\varphi_k} \circ \boldsymbol{\mu}_N ]$ of the empirical Markov process to the composition $\mathcal{L}_{\infty} {\left[ R_{ \varphi_k } \circ \boldsymbol{\mu}_N \right]}$, where the empirical map is $\boldsymbol{\mu}_N : \mathbf{x}^N \mapsto \mu_{\mathbf{x}^N}$. The generator $\mathcal{L}_{\infty}$ of the limit measure-valued process was defined in section \cref{sec:limitsemigroup}. This latter generator requires most of the time a specific formalism to be computed. We consider here the case of tensorized functions $\varphi_k = \varphi^1 \otimes\ldots\otimes\varphi^k$ : this relies on combinatorial and symmetry arguments, in a way which is reminiscent of \cite{jabin_mean_2016}. The first and most important example is $k = 2$ (see the compactness methods, where it is a key result). For mean-field generators of the form \cref{eq:Nparticlemeanfieldgenerator_summary}, the target generator against degree-$2$ monomials reads 
\[ \mathcal{L}_{\infty} {\left[ R_{ \varphi^1 \otimes \varphi^2 } \circ \boldsymbol{\mu}_N \right]} = R_{L_{\boldsymbol{\mu}_N} \varphi^1 \otimes \varphi^2} \circ \boldsymbol{\mu}_N + R_{\varphi^1 \otimes L_{\boldsymbol{\mu}_N} \varphi^2} \circ \boldsymbol{\mu}_N, \]
as in Section \cref{sec:limitsemigroup}.

\begin{lemma}[Quadratic estimates for mean-field generators] \label[II]{lemma:quadraticmeanfield} Let $\mathcal{L}_N$ be a mean-field generator of the form \cref{eq:Nparticlemeanfieldgenerator_summary}. Let $\varphi^1,\varphi^2\in \Dom(L_\mu)$ such that $\varphi^1\varphi^2\in  \Dom(L_\mu) $ for all $\mu\in\pb(E)$. Then it holds that
\[ \mathcal{L}_N {\left[ R_{\varphi^1 \otimes \varphi^2} \circ \boldsymbol{\mu}_N \right]} = R_{L_{ \boldsymbol{\mu}_N } \varphi^1 \otimes \varphi^2} \circ \boldsymbol{\mu}_N + R_{\varphi^1 \otimes L_{ \boldsymbol{\mu}_N } \varphi^2} \circ \boldsymbol{\mu}_N + \frac{1}{N} {\left\langle \boldsymbol{\mu}_N , \Gamma_{L_{\boldsymbol{\mu}_N}} {\left( \varphi^1 , \varphi^2 \right)} \right\rangle}. \]
\end{lemma}

\begin{proof}
Starting from $\mathcal{L}_N {\left[ R_{\varphi^1 \otimes \varphi^2} \circ \boldsymbol{\mu}_N \right]}$ at $\mathbf{x}^N=(x^1,\ldots,x^N)\in E^N$, let us compute 
\begin{align*}
&\mathcal{L}_N {\left[ R_{\varphi^1 \otimes \varphi^2} \circ \boldsymbol{\mu}_N \right]} {\left( \mathbf{x}^N \right)} = \sum_{i=1}^N L_{\mu_{\mathbf{x}^N}} \diamond_i {\left[ \mathbf{x}^N \mapsto {\left\langle \mu_{\mathbf{x}^N} , \varphi^1 \right\rangle} {\left\langle \mu_{\mathbf{x}^N} , \varphi^2 \right\rangle} \right]} {\left( \mathbf{x}^N \right)} \\
&= \sum_{i=1}^N {\left\{ \frac{1}{N^2} L_{\mu_{\mathbf{x}^N}} {\left[ \varphi^1 \varphi^2 \right]} {\left( x^i \right)} + \frac{1}{N^2} \sum_{\substack{j = 1 \\ j \neq i}}^N \varphi^2 {\left( x^j \right)} L_{\mu_{\mathbf{x}^N}} \varphi^1 {\left( x^i \right)} + \varphi^1 {\left( x^j \right)} L_{\mu_{\mathbf{x}^N}} \varphi^2 {\left( x^i \right)} \right\}} \\
&= \frac{1}{N^2} \sum_{i=1}^N {\left\{ L_{\mu_{\mathbf{x}^N}} {\left[ \varphi^1 \varphi^2 \right]} {\left ( x^i \right)} - \varphi^2 {\left( x^i \right)} L_{\mu_{\mathbf{x}^N}} \varphi^1  {\left( x^i \right)} - \varphi^1 {\left( x^i \right)} L_{\mu_{\mathbf{x}^N}} \varphi^2  {\left( x^i \right)} \right\}} \\
&\phantom{abcdef} + \frac{1}{N^2} \sum_{i,j=1}^N \varphi^1 {\left( x^j \right)} L_{\mu_{\mathbf{x}^N}} \varphi^2  {\left( x^i \right)} + \varphi^2 {\left( x^j \right)} L_{\mu_{\mathbf{x}^N}} \varphi^1  {\left( x^i \right)} \\
&= \frac{1}{N} {\left\langle \mu_{\mathbf{x}^N} , \Gamma_{L_{\mu_{\mathbf{x}^N}}} {\left( \varphi^1 , \varphi^2 \right)} \right\rangle} + R_{L_{\mu_{\mathbf{x}^N}} \varphi^1 \otimes \varphi^2} {\left( \mu_{\mathbf{x}^N} \right)} + R_{\varphi^1 \otimes L_{\mu_{\mathbf{x}^N}} \varphi^2} {\left( \mu_{\mathbf{x}^N} \right)},
\end{align*}
and the last term is exactly the desired expression.
\end{proof}

Once again, the carr\'e du champ controls the quadratic quantities. Let us try now to extend this estimate to any degree-$k$ monomial. A possible goal of this is to control the limit generator against polynomials, in order to approach its behaviour against any function by density. Unfortunately this fails here since the bound obtained still requires some growth comparison condition between $k$ and $N$ (see Section \cref{sec:limitgenerator}).

\begin{lemma}[Extension to large-degree monomials] \label[II]{lemma:largedegreemeanfield}
Fix $N \geq 2$. For every $j \geq 0$, let us define the operators $\Gamma^{(j+2)}_{L_\mu}: C_b {\left( E^{j+2} \right)} \rightarrow C_b {\left( E \right)}$ for $\mu\in\widehat{\pb}_N(E)$ by:
\begin{equation}\label[IIeq]{eq:jcarreduchamp}\Gamma^{(j+2)}_{L_{\boldsymbol{\mu}_N}} {\left( \varphi^1 , \ldots , \varphi^{j+2} \right)} = L_{\boldsymbol{\mu}_N} {\left[ \varphi^1 \ldots \varphi^{j+2} \right]} - \sum_{i=1}^{j+2} \varphi^i L_{\boldsymbol{\mu}_N} \prod_{ \substack{\ell = 1 \\ \ell \neq i} }^{j+2} \varphi^\ell, \end{equation}
where we implicitly assume that any product of test functions belong to the domain of the generator $L_\mu$ for all $\mu\in\widehat{\pb}_N(E)$. For $k\geq2$, let us assume that for any $0 \leq j \leq k-2$, there exists $C_j > 0$ such that for any $\{\ell_1,\ldots,\ell_{j+2}\}\subset\{1,\ldots,k\}$,
\begin{equation}\label[IIeq]{eq:boundjcarreduchamp} \sup_{\mu \in \widehat{\pb}_N(E)} \sup_{\|\varphi^1\|_{\infty}, \ldots,\|\varphi^k\|_{\infty} \leq 1} \Big\langle \mu^{\otimes k-1-j} , \Gamma^{(j+2)}_{L_{\mu}} {\left( \varphi^{\ell_1} , \ldots , \varphi^{\ell_{j+2}} \right)} \otimes \varphi^{\ell_{j+3}} \otimes \ldots \otimes \varphi^{\ell_k} \Big\rangle \leq C_j, \end{equation}
where $\{\ell_{j+3},\ldots,\ell_{k}\}=\{1,\ldots,k\}\setminus\{\ell_1,\ldots,\ell_{j+2}\}$. Then for $\mathcal{L}_N$ of the form \cref{eq:Nparticlemeanfieldgenerator_summary}, the following generator estimate holds with $\varphi_k = \varphi^1 \otimes \ldots \otimes \varphi^k$
\[ \mathcal{L}_N {\left[ R_{\varphi_k} \circ \boldsymbol{\mu}_N \right]} = \sum_{i=1}^k R_{\varphi^1 \otimes \ldots \otimes \varphi^{i-1} \otimes L_{\boldsymbol{\mu}_N} \varphi^i \otimes \varphi^{i+1} \otimes \ldots \otimes \varphi^k} \circ \boldsymbol{\mu}_N + \frac{1}{N} \sum_{j=0}^{k-2}\binom{k}{j+2}\frac{C_j}{N^j}.\]
In particular, if $C_j=\mathcal{O}(C^j)$ for a fixed $C>0$, then the remainder is controlled by $\mathcal{O}{\left(N^{-1}k^2{\left( 1 + \frac{C}{N} \right)}^k\right)}$.
\end{lemma}

Note that for $j=0$, the usual carr\'e du champ operator $\Gamma^{(2)}_{L_{\mu_{\mathbf{x}^N}}} = \Gamma_{L_{\mu_{\mathbf{x}^N}}}$ is recovered.

\begin{proof} Let us consider a tensorized $k$-particle test function $\varphi_k = \varphi^1 \otimes \ldots \otimes \varphi^k$ and $\mathbf{x}^N \in E^N$. The generator $\mathcal{L}_N$ is of the form \cref{eq:Nparticlemeanfieldgenerator_summary} so we have
\[ \mathcal{L}_N {\left[ R_{\varphi_k} \circ \boldsymbol{\mu}_N \right]} {\left( \mathbf{x}^N \right)} = \sum_{i=1}^N L_{\mu_{\mathbf{x}^N}} \diamond_i {\left[ \mathbf{x}^N \mapsto {\left\langle \mu_{\mathbf{x}^N}^{\otimes k} , \varphi^1 \otimes \ldots \otimes \varphi^k \right\rangle} \right]} {\left( \mathbf{x}^N \right)}. \]
We then use the linearity of $L_{\mu_{\mathbf{x}^N}}$ and the fact that it vanishes on constants. To compute the $L_{\mu_{\mathbf{x}^N}} \diamond_i {\left[ \cdot \right]}$ term, it is sufficient to develop the $\mu_{\mathbf{x}^N}^{\otimes k}$-sum and to discriminate on how many times $x^i$ appears. If there are $j$ occurrences, this leads to the sum 
\[   N^{-k} \sum_{\substack{ \{ \ell_1 , \ldots , \ell_{j} \} \\ \subset \{ 1 , \ldots , k \} }} \sum_{\substack{i_{\ell_{j+1}}, \ldots, i_{\ell_k} \\ i \notin \{ i_{\ell_{j+1}}, \ldots, i_{\ell_k} \}  }} L_{\mu_{\mathbf{x}^N}} {\left[ \varphi^{\ell_1} \ldots \varphi^{\ell_j} \right]} {\left( x^i \right)} \prod_{\ell \in \{ \ell_{j+1}, \ldots, \ell_k \} } \varphi^\ell {\left( x^{i_\ell} \right)},   \]
where we recall that for a given $\{\ell_1,\ldots,\ell_j\}\subset\{1,\ldots,k\}$, we write $\{\ell_{j+1},\ldots,\ell_k\}=\{1,\ldots,k\}\setminus \{\ell_1,\ldots,\ell_j\}$. The term $L_{\mu_{\mathbf{x}^N}} \diamond_i {\left[ \cdot \right]}$ is then obtained by summing over $1 \leq j \leq k$. Summing then over $i$ gives
\begin{equation}\label[IIeq]{eq:LNRphiksumSjk}\mathcal{L}_N {\left[ R_{\varphi_k} \circ \boldsymbol{\mu}_N \right]} {\left( \mathbf{x}^N \right)} = \sum_{j=1}^k S_j^k {\left(  \mu_{\mathbf{x}^N} \right)}, \end{equation}
using the shortcut
\begin{multline*}S_j^k {\left(  \mu_{\mathbf{x}^N} \right)} \\:= N^{-k} \sum_{\substack{ \{ \ell_1 , \ldots , \ell_{j} \} \\ \subset \{ 1 , \ldots , k \} }} \sum_{\substack{i_{\ell_{j}}, \ldots, i_{\ell_k} \\ i_{\ell_j} \notin \{ i_{\ell_{j+1}}, \ldots, i_{\ell_k} \}  }} L_{\mu_{\mathbf{x}^N}} {\left[ \varphi^{\ell_1} \ldots \varphi^{\ell_j} \right]} {\left( x^{i_{\ell_j}} \right)} \prod_{\ell \in \{ \ell_{j+1}, \ldots, \ell_k \} } \varphi^\ell {\left( x^{i_\ell} \right)}.\end{multline*}
Introduce now for $1 \leq j \leq k$
\begin{multline*}R_j^k {\left(  \mu_{\mathbf{x}^N} \right)} \\:= N^{-k} \sum_{\substack{ \{ \ell_1 , \ldots , \ell_{j} \} \\ \subset \{ 1 , \ldots , k \} }} \sum_{\substack{i_{\ell_{j}}, \ldots, i_{\ell_k} \\ i_{\ell_j} \in \{ i_{\ell_{j+1}}, \ldots, i_{\ell_k} \}  }} L_{\mu_{\mathbf{x}^N}} { \left[ \varphi^{\ell_1} \ldots \varphi^{\ell_j} \right] } {\left( x^{i_{\ell_j}} \right)} \prod_{\ell \in \{ \ell_{j+1}, \ldots, \ell_k \} } \varphi^\ell {\left( x^{i_\ell} \right)},\end{multline*}
so that
\begin{multline*}S_j^k {\left(  \mu_{\mathbf{x}^N} \right)} + R_j^k {\left(  \mu_{\mathbf{x}^N} \right)} \\= N^{-k} \sum_{\substack{ \{ \ell_1 , \ldots , \ell_{j} \} \\ \subset \{ 1 , \ldots , k \} }} \sum_{i_{\ell_{j}}, \ldots, i_{\ell_k}} L_{\mu_{\mathbf{x}^N}} {\left[ \varphi^{\ell_1} \ldots \varphi^{\ell_j} \right]} {\left( x^{i_{\ell_j}} \right)} \prod_{\ell \in \{ \ell_{j+1}, \ldots, \ell_k \} } \varphi^\ell {\left( x^{i_\ell} \right)}.\end{multline*}
Moreover $R_k^k {\left(  \mu_{\mathbf{x}^N} \right)} = 0$ and
\begin{align}
S_1^k {\left(  \mu_{\mathbf{x}^N} \right)} + R_1^k {\left(  \mu_{\mathbf{x}^N} \right)} &= N^{-k} \sum_{\ell_1=1}^k \sum_{ i_{\ell_1}, \ldots, i_{\ell_k} } L_{\mu} \varphi^{\ell_1} {\left( x^{i_{\ell_j}} \right)} \prod_{\ell \in \{ \ell_{2}, \ldots, \ell_k \} } \varphi^\ell {\left( x^{i_\ell} \right)} \nonumber\\
&= \sum_{i=1}^k \varphi^1 \otimes \ldots \otimes \varphi^{i-1} \otimes L_{\mu_{\mathbf{x}^N}} \varphi^i \otimes \varphi^{i+1} \otimes \ldots \otimes \varphi^k {\left(  \mu_{\mathbf{x}^N} \right)} \nonumber\\
&= \sum_{i=1}^k R_{\varphi^1 \otimes \ldots \otimes \varphi^{i-1} \otimes L_{\mu_{\mathbf{x}^N}} \varphi^i \otimes \varphi^{i+1} \otimes \ldots \otimes \varphi^k} {\left(  \mu_{\mathbf{x}^N} \right)}.\label[IIeq]{eq:Sk1Rk1}
\end{align}
An alternative way to write $R_j^k {\left(  \mu_{\mathbf{x}^N} \right)}$ is
\begin{multline*}R_j^k {\left(  \mu_{\mathbf{x}^N} \right)} \\= N^{-k} \sum_{\substack{ \{ \ell_1 , \ldots , \ell_{j+1} \} \\ \subset \{ 1 , \ldots , k \} }} \sum_{ i_{\ell_{j+1}}, \ldots, i_{\ell_k} } \sum_{m=1}^{j+1} {\left\{ \varphi^{\ell_m} L_{\mu} \prod_{ \substack{n = 1 \\ n \neq m} }^{j+1} \varphi^{\ell_n} \right\}} {\left( x^{i_{\ell_j}} \right)} \prod_{\ell \in \{ \ell_{j+2}, \ldots, \ell_k \} } \varphi^\ell {\left( x^{i_\ell} \right)}.\end{multline*}
Using the $j$-carr\'e du champ \cref{eq:jcarreduchamp}, we have the telescopic expression for $1\leq j< k$ :
\begin{align*}
&S^k_{j+1}+R^k_{j+1} - R^k_j \\
&= N^{-k} \sum_{\substack{ \{ \ell_1 , \ldots , \ell_{j+1} \} \\ \subset \{ 1 , \ldots , k \} }} \sum_{ i_{\ell_{j+1}}, \ldots, i_{\ell_k} } \Gamma^{(j+1)}_{L_{\mu_{\mathbf{x}^N}}} {\left( \varphi^1 , \ldots , \varphi^{j+1} \right)} {\left( x^{i_{\ell_{j+1}}} \right)} \prod_{\ell \in \{ \ell_{j+2}, \ldots, \ell_k \} } \varphi^l {\left( x^{i_\ell} \right)}    \\
&= N^{-j} \sum_{\substack{ \{ \ell_1 , \ldots , \ell_{j+1} \} \\ \subset \{ 1 , \ldots , k \} }} \Big\langle \mu_{\mathbf{x}^N}^{\otimes k-j} , \Gamma^{(j+1)}_{L_{\mu_{\mathbf{x}^N}}} {\left( \varphi^{\ell_1} , \ldots , \varphi^{\ell_{j+1}} \right)} \otimes \varphi^{\ell_{j+2}} \otimes \ldots \otimes \varphi^{\ell_k} \Big\rangle.
\end{align*} 
We then sum this expression over $1 \leq j \leq k-1$, we add $S^k_1+R^k_1$ and we use that $R_k^k {\left(  \mu_{\mathbf{x}^N} \right)}=0$. From \cref{eq:LNRphiksumSjk} and \cref{eq:Sk1Rk1}, we conclude that $\mathcal{L}_N {\left[ R_{\varphi_k} \circ \mu_N \right]}{\left( \mathbf{x}^N \right)}$ is equal, up to a remainder, to the expected generator 
\[ \sum_{i=1}^k R_{\varphi^1 \otimes \ldots \otimes \varphi^{i-1} \otimes L_{\mu_{\mathbf{x}^N}} \varphi^i \otimes \varphi^{i+1} \otimes \ldots \otimes \varphi^k} {\left(  \mu_{\mathbf{x}^N} \right)}, \]
the remainder being
\[ \sum_{j=0}^{k-2} N^{-1-j} \sum_{\substack{ \{ \ell_1 , \ldots , \ell_{j+2} \} \\ \subset \{ 1 , \ldots , k \} }} \Big\langle \mu_{\mathbf{x}^N}^{\otimes k-1-j} , \Gamma^{(j+2)}_{L_{\mu_{\mathbf{x}^N}}} {\left( \varphi^{\ell_1} , \ldots , \varphi^{\ell_{j+2}} \right)} \otimes \varphi^{\ell_{j+3}} \otimes \ldots \otimes \varphi^{\ell_k} \Big\rangle. \]
The final estimate then follows using the boundedness assumption \cref{eq:boundjcarreduchamp}, the number of combinations $\{ \ell_1 , \ldots , \ell_j \} \subset \{ 1 , \ldots , k \}$ and the binomial expansion.
\end{proof}

Consider now the situation of the Boltzmann models. The Boltzmann generator is described in Section \cref{sec:boltzmanngeneral}:
\[ \mathcal{L}_N \varphi_N = \frac{1}{N} \sum_{1 \leq i < j \leq N} L^{(2)} \diamond_{i,j} \varphi_N \]
where $L^{(2)}$ reads
\[ L^{(2)} \varphi_2 {\left( x^1 , x^2 \right)} = \lambda {\left( x^1 , x^2 \right)} \int_{E^2} {\left[ \varphi_2 {\left( x'^1 , x'^2 \right)} - \varphi_2 {\left( x^1 , x^2 \right)} \right]} \Gamma^{(2)} {\left( x^1 , x^2 , \dd x'^1 , \dd x'^2 \right)}, \]
where $\lambda$ and $\Gamma^{(2)}$ satisfy Assumption \cref{assum:L2_summary}. The symmetry properties imply a nice shape for the symmetrized version of $L^{(2)}$
\begin{align*}
 L^{(2)}_{\mathrm{sym}} \varphi_2 {\left( x^1 , x^2 \right)} &= \frac{L^{(2)}_{\mathscr{sym}} \varphi_2 {\left( x^1 , x^2 \right)} + L^{(2)}_{\mathscr{sym}} \varphi_2 {\left( x^2 , x^1 \right)}}{2} \\
 &= \frac{\lambda {\left( x^1 , x^2 \right)}}{2} \int_{E^2} \Big\{\varphi_2 {\left( x'^1 , x'^2 \right)} + \varphi_2 {\left( x'^2 , x'^1 \right)} \\
& \phantom{abcdefghijklm} - \varphi_2 {\left( x^1 , x^2 \right)} - \varphi_2 {\left( x^2 , x^1 \right)}\Big\} \Gamma^{(2)} {\left( x^1 , x^2 , \dd x'^1 , \dd x'^2 \right)}, 
\end{align*}
this implies $L^{(2)}_{\mathrm{sym}} {\left[ \varphi^1 \otimes \varphi^2 \right]} = L^{(2)}_{\mathrm{sym}} {\left[ \varphi^2 \otimes \varphi^1 \right]}$ for every  $\varphi^1 , \varphi^2 \in \mathcal{F}$. For the limit generator, this symmetry suggests to define $L_{\mu}$ as
\[ \forall \varphi \in \mathcal{F}, \forall x \in E,  \quad L_{\mu} \varphi ( x ) := {\left\langle \mu , L^{(2)}_{\mathrm{sym}} [ \varphi \otimes 1 ] ( x , \cdot )  \right\rangle} = {\left\langle \mu , L^{(2)}_{\mathrm{sym}} [ \varphi \otimes 1 ] ( \cdot , x )  \right\rangle}, \]
and equivalently $\varphi \otimes 1$ can be taken instead of $1 \otimes \varphi$ in the above definition. The needed estimate is now the following.

\begin{lemma}[Quadratic estimates for Boltzmann collisions] \label[II]{lemma:quadraticboltzmann}
The quadratic estimate for degree-$2$ monomials reads
\begin{align*}
\mathcal{L}_N {\left[ R_{\varphi^1 \otimes \varphi^2} \circ \boldsymbol{\mu}_N \right]} = R_{L_{ \boldsymbol{\mu}_N } \varphi^1 \otimes \varphi^2} \circ \boldsymbol{\mu}_N &+ R_{\varphi^1 \otimes L_{ \boldsymbol{\mu}_N } \varphi^2} \circ \boldsymbol{\mu}_N + \frac{1}{N} R_{L^{(2)}_{\mathrm{sym}} {\left[ \varphi^1 \otimes \varphi^2 \right]}} \circ \boldsymbol{\mu}_N \\
&+ \frac{1}{N} {\left\langle \boldsymbol{\mu}_N , \Gamma_{L_{\boldsymbol{\mu}_{\mathbf{x}^N}}} {\left( \varphi^1 , \varphi^2 \right)}  \right\rangle}.
\end{align*}
\end{lemma}

Note that compared to Lemma \cref{lemma:quadraticmeanfield}, an additional symmetrizing term appears.

\begin{proof} It is a direct computation. Let us start with
\[\mathcal{L}_N {\left[ R_{\varphi^1 \otimes \varphi^2} \circ \boldsymbol{\mu}_N \right]} {\left( \mathbf{x}^N \right)} = \frac{1}{N} \sum_{1 \leq i < j \leq N} L^{(2)} \diamond_{i,j} {\left[ \mathbf{x}^N \mapsto {\left\langle \mu_{\mathbf{x}^N} , \varphi^1 \right\rangle} {\left\langle \mu_{\mathbf{x}^N} , \varphi^2 \right\rangle} \right]} {\left( \mathbf{x}^N \right)}. \]
We then develop the expression inside the term $L^{(2)} \diamond_{i,j} {\left[ \cdot \right]}$. Since $L^{(2)} {\left[ 1 \otimes 1 \right]} = 0$, the only remaining terms are (up to a factor $N^{-2}$)
\begin{align*}
\varphi^1 \varphi^2 {\left( x^i \right)} &+ \varphi^1 \varphi^2 {\left( x^j \right)} + \varphi^1 {\left( x^i \right)} \varphi^2 {\left( x^j \right)} + \varphi^1 {\left( x^i \right)} \varphi^2 {\left( x^j \right)} \\
&+ {\left[ \varphi^1 {\left( x^i \right)} + \varphi^1 {\left( x^j \right)} \right]} \sum_{k \neq i,j} \varphi^2 {\left( x^k \right)} + {\left[ \varphi^2 {\left( x^i \right)} + \varphi^2 {\left( x^j \right)} \right]} \sum_{k \neq i,j} \varphi^1 {\left( x^k \right)}.
\end{align*} 
Applying $L^{(2)} \diamond_{i,j}$, the total expression is now (up to a factor $N^{-3}$) the sum over $1 \leq i < j \leq N$ of the terms 
\begin{align*} 
L^{(2)} &{\left[ \varphi^1 \varphi^2 \otimes 1 \right]} + L^{(2)} {\left[ 1 \otimes \varphi^1\varphi^2 \right]} + L^{(2)} {\left[ \varphi^1 \otimes \varphi^2 \right]} + L^{(2)} {\left[ \varphi^2 \otimes \varphi^1 \right]} \\
&+ {\left[ L^{(2)} {\left[ \varphi^1 \otimes 1 \right]} + L^{(2)} {\left[ 1 \otimes \varphi^1 \right]} \right]} \sum_{k \neq i,j} \varphi^2 {\left( x^k \right)} \\
&+ {\left[ L^{(2)} {\left[ \varphi^2 \otimes 1 \right]} + L^{(2)} {\left[ 1 \otimes \varphi^2 \right]} \right]} \sum_{k \neq i,j} \varphi^1 {\left( x^k \right)},
\end{align*}
where all the functions are evaluated at the point ${\left( x^i , x^j \right)}$. The property $\lambda ( x , x ) =0$ implies
\[ \sum_{1 \leq i < j \leq N} {\left[ L^{(2)} {\left[ \varphi^1 \otimes 1 \right]} + L^{(2)} {\left[ 1 \otimes \varphi^1 \right]} \right]} = \sum_{i,j = 1}^N L^{(2)}_{\mathrm{sym}} {\left[ \varphi^1 \otimes 1 \right]}. \]
Note also that up to a factor $N^{-3}$, $R_{L_{\mu_{\mathbf{x}^N}} \varphi^1 \otimes \varphi_2} {\left( \mu_{\mathbf{x}^N} \right)} + R_{\varphi^1 \otimes L_{\mu_{\mathbf{x}^N}} \varphi_2} {\left( \mu_{\mathbf{x}^N} \right)}$ equals the sum (evaluated at the point ${\left( x^i , x^j \right)}$) over $1 \leq i , j \leq N$ of
\begin{align*}
L^{(2)}_{\mathrm{sym}} &{\left[ \varphi^1 \otimes 1 \right]} \sum_{k \neq i,j} \varphi^2 {\left( x^k \right)} + L^{(2)}_{\mathrm{sym}} {\left[ \varphi^2 \otimes 1 \right]} \sum_{k \neq i,j} \varphi^2 {\left( x^k \right)} \\
&+ {\left[ \varphi^2 {\left( x^i \right)} + \varphi^2 {\left( x^j \right)} \right]} L^{(2)}_{\mathrm{sym}} {\left[ \varphi^1 \otimes 1 \right]} + {\left[ \varphi^1 {\left( x^i \right)} + \varphi^1 {\left( x^j \right)} \right]} L^{(2)}_{\mathrm{sym}} {\left[ \varphi^2 \otimes 1 \right]}.
\end{align*}
In the same way, ${\left\langle \mu_{\mathbf{x}^N} , \Gamma_{L_{\mu_{\mathbf{x}^N}}} {\left( \varphi^1 , \varphi^2 \right)}  \right\rangle}$ equals up to a factor $N^{-2}$ the sum (evaluated at the point ${\left( x^i , x^j \right)}$) over $1 \leq i , j \leq N$ of
\begin{multline*}L^{(2)}_{\mathrm{sym}} {\left[ \varphi^1 \varphi^2 \otimes 1 \right]} - {\left[ \varphi^2 {\left( x^i \right)} + \varphi^2 {\left( x^j \right)} \right]} L^{(2)}_{\mathrm{sym}} {\left[ \varphi^1 \otimes 1 \right]} \\- {\left[ \varphi^1 {\left( x^i \right)} + \varphi^1 {\left( x^j \right)} \right]} L^{(2)}_{\mathrm{sym}} {\left[ \varphi^2 \otimes 1 \right]}.\end{multline*}
At the end of the day, summing everything with the adequate power of $N$, one gets 
\begin{align*}
\mathcal{L}_N {\left[ R_{\varphi^1 \otimes \varphi^2} \circ \boldsymbol{\mu}_N \right]} = R_{L_{ \boldsymbol{\mu}_N } \varphi^1 \otimes \varphi^2} \circ \boldsymbol{\mu}_N &+ R_{\varphi^1 \otimes L_{ \boldsymbol{\mu}_N } \varphi^2} \circ \boldsymbol{\mu}_N + \frac{1}{N} R_{L^{(2)}_{\mathrm{sym}} {\left[ \varphi^1 \otimes \varphi^2 \right]}} \circ \boldsymbol{\mu}_N \\
&+ \frac{1}{N} {\left\langle \boldsymbol{\mu}_N , \Gamma_{L_{\mu_{\mathbf{x}^N}}} {\left( \varphi^1 , \varphi^2 \right)}  \right\rangle.}
\end{align*}
\end{proof}

\section{A combinatorial lemma}\label[II]{appendix:combinatoriallemma}

This combinatorial lemma is used in Section \cref{sec:martingalecompactness} to control the jumps of the limit process.

\begin{lemma} \label[II]{lem:Intersection}
Let $(\Omega,\mathscr{F},P)$ be a probability space, and consider two integers $2 \leq p \leq N$. Let $(A_i)_{1 \leq i \leq N}$ be a sequence of events in $\mathscr{F}$ such that $P(A_i) > 1/p$, and assume the existence of an integer $q \geq 1$ such that any intersection involving $(q+1)$ of the $A_i$ is $P$-negligible. Then 
\[ \frac{N}{p} < q.\]
\end{lemma}

As a corollary, from an infinite sequence $(A_n)_{n \geq 1}$ of events such that $P(A_n) > 1/p$ (for a given $p\geq2$), it is possible to build a non-negligible intersection involving an arbitrary large number of $A_n$.

\begin{proof} For $1 \leq j \leq q$, consider the set of $j$-intersections
\[ \mathcal{A}_j = {\left\{ \bigcap_{\ell = 1}^j A_{i_\ell}, \, i_1, \ldots, i_j \in \{1,\ldots,N\} \text{ pairwise distinct and } P{\left( \bigcap_{\ell = 1}^j A_{i_\ell} \right)} > 0 \right\}}. \]
From this, we construct a partition of $\bigcup_{1 \leq i \leq N} A_i$ which is composed of $j$ intersections. Let us first define the class of sets which are intersections of at most $j$ subsets. 
\[ \mathcal{R}_j := {\left\{ a \cap {\left(\bigcup_{k \geq j+1} \bigcup_{a' \in \mathcal{A}_k} a'\right)}^c, \,\, a \in \mathcal{A}_j \right\}}. \]
The intersections $a$ within $\bigcup_{1 \leq j \leq q} \mathcal{R}_j$ are pairwise disjoint, because the recovering of two $j$-intersections belongs at least to a $(j+1)$-intersection. Then, by definition of $q$, $\bigcup_{1 \leq i \leq N} A_i$ is $P$-a.s. covered by $\bigcup_{1 \leq j \leq q} \bigcup_{a \in \mathcal{R}_j} a$. As a consequence,
\begin{equation} \label[IIeq]{eq:totalMAss}
P{\left( {\bigcup_{1 \leq i \leq N} A_i} \right)} = \sum_{j=1}^q \sum_{a \in \mathcal{R}_j} P(a).
\end{equation}
For any $i$ and $a \in \mathcal{R}_j$, define now the contribution of $A_i$ to $a$ as $f_i(a) := P\left(a \cap A_i\right)$. Since
\[ \forall 1 \leq k \leq j, \quad P{\left( A_{i_k} \cap \bigcap_{l=1}^j A_{i_l}\right)} = P{\left( {\bigcap_{l=1}^j A_{i_l}} \right)} \] 
it is straightforward to check that
\begin{equation} \label[IIeq]{eq:sameContrib}
f_i(a) = 
\left\{
\begin{array}{ll}
0 \quad &\text{if } a \cap A_i = \emptyset   \\
P(a) \quad &\text{if } a \cap A_i \neq \emptyset 
\end{array}
\right.
\end{equation}
From the definition of $\mathcal{R}_j$, exactly $j$ of the $A_i$ positively contribute to an intersection $a \in \mathcal{R}_j$. Using this and \cref{eq:sameContrib}, it follows that
\[ \forall a \in \mathcal{R}_j, \quad \sum_{i=1}^N f_i(a) = j P(a). \] 
We sum this relation over $a \in \mathcal{R}_j$, then divide by $j$, and eventually sum over $1 \leq j \leq q$. Injecting this into \cref{eq:totalMAss} gives
\begin{equation} \label[IIeq]{eq:volumeBound}
 \sum_{i=1}^N \sum_{j=1}^q \frac{1}{j} \sum_{a \in \mathcal{R}_j} f_i(a) = P{\left( {\bigcup_{1 \leq i \leq N} A_i} \right)} \leq 1.
\end{equation} 
Since $A_k \subset\bigcup_{1 \leq i \leq N} A_i$ for every $1 \leq k \leq q$, the mass $P(A_i)$ shall be recovered as
\[ P(A_i) = \sum_{a \in \bigcup_{1 \leq j \leq q} \mathcal{R}_j} P(a \cap A_i) = \sum_{j=1}^q \sum_{a \in \mathcal{R}_j} f_i(a), \]
using the previous partition of $\bigcup_{1 \leq i \leq N} A_i$. Since $\frac{1}{j}$ remains bigger than $\frac{1}{q}$, \cref{eq:volumeBound} then leads 
\[ \frac{1}{q} \sum_{i=1}^N P(A_i) \leq 1. \]
The conclusion follows by the definition of $p$.
\end{proof} 

\section{Convergence in the Skorokhod space and tightness criteria}\label[II]{appendix:tightness} 

This section gathers the classical and less classical convergence results which are used in Section \cref{sec:mckeancompactnessreview} and Section \cref{sec:martingaleboltzmannreview} to prove propagation of chaos via compactness arguments. Probability remainders on the Skorokhod space, martingales, semimartingales and $D$-semimartingales can be found in Appendix \cref{appendix:skorokhod} and Appendix~\cref{appendix:martingales}. 

The following criterion due to Aldous \cite{aldous_stopping_1978} is the most classical result to prove the tightness of the laws of a sequence of c\`adl\`ag processes. 
\begin{theorem}[Aldous criterion]\label[II]{thm:aldous} For each $n\in\N$, let $(X^n_t)^{}_t$ be an adapted c\`adl\`ag process on the filtered probability space $(\Omega,\mathscr{F},(\mathscr{F}_t)_t,\mathbb{P})$. Assume that the sequence of processes satisfies the following conditions. 
\begin{enumerate}[(i)]
\item For all $N\in\N$ and for all $\varepsilon>0$ there exist $n_0\in\N$ and $K>0$ such that 
\[n\geq n_0 \Rightarrow \mathbb{P}{\left(\sup_{t\leq N} |X^n_t|>K\right)}\leq\varepsilon.\]
\item For all $N\in\N$ and for all $\varepsilon>0$ it holds that 
\begin{equation} \label[IIeq]{eq:AldousCondition}
\lim_{\theta\downarrow0}\limsup_{n} \sup_{S,T\in\mathscr{T}_N: S\leq T\leq S+\theta} \mathbb{P}(|X^n_T-X^n_S|\geq\varepsilon)=0,
\end{equation}
where $\mathscr{T}_N$ denotes the set of all $(\mathscr{F}_t)_t$-stopping times that are bounded by $N$. 
\end{enumerate}
Then the sequence of processes $(X^n_t)^{}_t$ is tight. 
\end{theorem}

\begin{proof} See \cite[Chapter VI, Section 4a]{jacod_limit_2003} or \cite[Chapter 3, Theorem 8.6]{ethier_markov_1986}.
\end{proof}

This criterion can be extended to a more general Polish space $(E,\rho)$ by replacing the first condition by the tightness of $(X^n_t)_{n\geq0}$ for each $t$ in a dense subset of $\R_+$. In the second condition \cref{eq:AldousCondition}, the norm $|X^n_T-X^n_S|$ has to be replaced by the distance $\rho(X^n_T,X^n_S)$.

The following theorem reduces the question of tightness in $\pb(D([0,T],\mathscr{E}))$ for an arbitrary space $\mathscr{E}$ to the simpler question of tightness in $\pb(D([0,T],\R))$. 

\begin{theorem}[Jakubowski] \label[II]{thm:jakubowski}
Let $\mathscr{E}$ be a completely regular topological space with metrisable compacts. Let $\mathcal{F}$ be a family of continuous functions on $\mathscr{E}$ which satisfies the following properties. 
\begin{enumerate}[(i)]
    \item $\mathcal{F}$ separates points in $\mathscr{E}$.
    \item $\mathcal{F}$ is closed under addition, i.e. if $\Phi_1,\Phi_2\in\mathcal{F}$, then $\Phi_1+\Phi_2\in\mathcal{F}$. 
\end{enumerate}
Let $T\in(0,+\infty]$. Let $(\mu_n)_{n\in\N}$ a family of probability measures in $\pb(D([0,T],\mathscr{E}))$. Then the family $(\mu_n)_{n\in\N}$ is tight if and only if the following properties hold. 
\begin{enumerate}[(i)]
    \item For all $\varepsilon>0$ and for all $t>0$ there exists a compact set $K_{t,\varepsilon}\subset\mathscr{E}$ such that for all $n\in\N$, 
    \[\mu_n(D([0,t],K_{t,\varepsilon}))>1-\varepsilon,\]
    and we can consider only $t=T$ when $T<+\infty$.
    \item The family $(\mu_n)_{n\in\N}$ is $\mathcal{F}$-weakly tight in the sense that for all $\Phi\in\mathcal{F}$, the family of probability measures $(\widetilde{\Phi}_\#\mu_n)_{n\in\N}$ is tight in $\pb(D([0,T],\R))$ where $\widetilde{\Phi}$ denotes the natural extension of $\Phi$ on $D([0,T],\mathscr{E})$ : 
    \[\widetilde{\Phi} : D([0,T],\mathscr{E}) \to D([0,T],\R),\quad \omega \mapsto \Phi\circ\omega.\]
\end{enumerate}
\end{theorem}

\begin{proof}
See \cite[Theorem 3.1 and Theorem 4.6]{jakubowski_skorokhod_1986}. 
\end{proof}

Theorem \cref{thm:jakubowski} is used in this review with $\mathscr{E}=\pb(E)$ and $\mathcal{F}$ the family of linear functions $\Phi(\mu)=\langle \varphi,\mu\rangle$ with $\varphi\in C_b(E)$. In this case, a similar result also appears in \cite[Theorem 2.1]{roellycoppoletta_criterion_1986}.

The following theorem gives a necessary and sufficient condition for the weak limit of a sequence of c\`adl\`ag processes to be almost surely continuous. 

\begin{theorem}[Continuity mapping in $D$] \label[II]{thm:continuitysara}
Given the Polish space $(E,\rho)$, let us define for $x$ in $D([0,T],E)$,
\[ J(x) := \int_0^{+ \infty} \e^{-t}{\left[1 \land \sup_{0 \leq s \leq t} \rho(x(s^-),x(s))\right]} \dd t. \]
Let $\big((X^n_t)^{}_t\big)_n$ be a sequence of adapted $E$-valued c\`adl\`ag processes which converges in law towards a c\`adl\`ag process $X$. Then $X$ is a.s. continuous if and only $J(X^n)$ converges in law towards $0$.
\end{theorem}

\begin{proof} See \cite[Chapter 3, Theorem 10.2]{ethier_markov_1986}.
\end{proof}

The basic tightness criterion for semimartingales is due to Rebolledo.

\begin{theorem}[Rebolledo criterion] Let $\big((X^n_t)_{t \geq 0}\big)_{n\geq0}$ be a sequence of c\`adl\`ag square integrable semimartingales. Let us write the decompositoin $X^n_t = A^n_t + M^n_t$, where $(M^n_t)_{t \geq 0}$ is a local square integrable martingale and $(A^n_t)_{t \geq 0}$ is an adapted finite variation paths process. If the two following conditions are fulfilled, then the sequences of processes $(M^n_t)_{t \geq 0}$, $([M^n]_t)_{t \geq 0}$ and $(X^n_t)^{}_t$ are tight.
\begin{enumerate}[(i)]
\item For every $t$ within a dense subset of $\R_+$, $(M^n_t)_{n \geq 0}$ and $(A^n_t)_{n \geq 0}$ are tight sequences.
\item Both processes $(\langle M^n\rangle)_{n \geq 0}$ and $(A^n)_{n \geq 0}$ satisfy condition \cref{eq:AldousCondition}. 
\end{enumerate}
\end{theorem}

\begin{proof}
See \cite[Theorem 2.3.2, Corollary 2.3.3]{joffe_weak_1986}.
\end{proof}

Finally, based on the Rebolledo criterion, a useful tightness criterion for $D$-semimartingales is proved in \cite{joffe_weak_1986}. This result is based on the following three assumptions. In the following $\big((\mathbf{X}^n_t)_{t \geq 0}\big)_{n\geq0}$ is a sequence of $\R^d$-valued $D$-semi-martingales and we use the notations of Definition \cref{def:semiMartD} and Lemma \cref{lemma:semiMartVariation}.

\begin{assumption} \label[II]{assum:TightnessH1}
There exist a constant $C>0$ and a sequence of positive adapted processes such that $(C^n_t)_{t\geq0}$ a.s. (recall that $\mathbf{b}$ and $\mathbf{a}$ are random processes):
\[ \forall t \geq 0, \forall \mathbf{x}^d \in \R^d, \quad \| \mathbf{b}(\mathbf{x}^d,t) \|^2 + \Tr \mathbf{a}(\mathbf{x}^d,t) \leq C^n_t \big(C + \| \mathbf{x}^d \|^2\big), \]
and for every $T>0$,
\begin{equation} \label[Ieq]{eq:bornegenetight}
\sup_{n \geq 0} \sup_{0 \leq t \leq T} \E [C^n_t] < + \infty, \quad \lim_{r \to + \infty} \sup_{n \geq 0} \mathbb{P}{\left( \sup_{0\leq t\leq T} C^n_t > r \right)} = 0.
\end{equation} 
\end{assumption}

\begin{assumption} \label[II]{assum:TightnessH2}
The initial sequence $(X^n_0)_{n\geq0}$ of random variables is such that
\[ \sup_{n \geq 0} \E \| \mathbf{X}^n_0 \|^2 < + \infty. \]
\end{assumption}

These two first assumptions are necessary to guaranty a $L^2$ Gronwall-like bound on $\mathbf{X}^n_t$ proved in \cite[Lemma 3.2.2]{joffe_weak_1986}. The next one is more technical but not difficult to check in practise.

\begin{assumption} \label[II]{assum:TightnessH3}
There exist a positive function $\alpha$ on $\R_+$ and a decreasing sequence of numbers $(\rho_n)_n$ such that $\lim_{t \to 0^+} \alpha(t) = 0$ and $\lim_{n \to +\infty} \rho_n = 0$, and for all $0<s<t$ and $n \geq 0$,
\[ A^n(t) - A^n(s) \leq \alpha(t-s) + \rho_n. \]
\end{assumption}

This assumption implies that the jumps of $A^n$ are smaller than $\rho_n$. Throughout the applications in the present article, we consider that $A(t) = t$ so this assumption is automatically fulfilled. 

\begin{theorem}[Joffe-Metivier criterion] \label[II]{thm:joffecriterion}
If Assumptions \cref{assum:TightnessH1}, \cref{assum:TightnessH2} and \cref{assum:TightnessH3} are verified, then the sequence $((\mathbf{X}^n_t)_{t \geq 0})_{n\geq0}$ of $D$-semimartingales is tight. If moreover convergence in law is assumed for the initial sequence in Assumption \cref{assum:TightnessH2}, then the canonical process is continuous in probability under the law of any limit point of the sequence.
\end{theorem}

\begin{proof}
See \cite[Proposition 3.2.3]{joffe_weak_1986} and \cite[Theorem 3.3.1]{joffe_weak_1986}. Under additional assumptions, this latter theorem also characterizes any limit point of the sequence as the solution of a martingale problem.
\end{proof}

\providecommand{\href}[2]{#2}
\providecommand{\arxiv}[1]{\href{http://arxiv.org/abs/#1}{arXiv:#1}}
\providecommand{\url}[1]{\texttt{#1}}
\providecommand{\urlprefix}{}


\medskip
Received xxxx 20xx; revised xxxx 20xx.
\medskip

\end{document}